\RequirePackage{fix-cm}
\pdfoutput=1

\documentclass[leqno, fleqn, centertags, 12pt]{article}

\usepackage{latexsym}
\usepackage{amsmath}
\usepackage{amssymb}
\usepackage{verbatim}

\usepackage[T1]{fontenc}

\usepackage{graphicx}

\usepackage{tikz-cd}

\usepackage[upright, widespace]{fourier}

\usepackage{fullpage}

\newskip\nineskipamount \nineskipamount=9pt plus 0pt minus 0pt
\newskip\zeroskipamount \zeroskipamount=0pt plus 0pt minus 0pt
\usepackage[noorphans,vskip=\nineskipamount]{quoting}

\makeatletter
\renewcommand{\@makefntext}[1]{\vspace*{0.5ex}\parindent=0em
\hspace*{-0.4em}
\hbox to 0.4em{\hss\@makefnmark}\hspace*{0.4em}{#1}
}
\makeatother

\newcounter{mysectionnumber}
\newcommand{\mysection}[2]{\setcounter{footnote}{0}
\setcounter{myparnum}{0}
\refstepcounter{mysectionnumber}
\vspace{27pt}{\Large {\themysectionnumber.} {#1}}\label{#2}\vspace*{15pt}}

\newcommand{\mynonumbersection}[1]{\vspace{27pt}{\Large {#1}}\vspace*{15pt}}

\newcommand{\myit}[1]{\textbf{\textit{#1}}\hspace{0.0em}}

\newcounter{myparnum}[mysectionnumber]
\renewcommand{\themyparnum}{\arabic{mysectionnumber}.\arabic{myparnum}}
\newcommand{\mypar}[2]{\refstepcounter{myparnum}{\vspace{\medskipamount}\textbf{{\themyparnum. #1}\label{#2}}\hspace{0.5em}}}

\newcounter{mylemmanum}[myparnum]
\newcommand{\myuppar}[1]{\vspace{\medskipamount}\textbf{#1}\hspace*{0.5em}}

\newcounter{myappendnumber}
\newcounter{myaparnum}[myappendnumber]
\newcommand{\myappend}[2]{\setcounter{footnote}{0}
\setcounter{myaparnum}{0}
\setcounter{myparnum}{0}
\refstepcounter{myappendnumber}
\vspace{21pt}{\Large A\dff.{\themyappendnumber.}\oss {#1}}\label{#2}\vspace*{15pt}}

\newcommand{\myapar}[2]{\refstepcounter{myaparnum}{\vspace{\medskipamount}\textbf{{\themyaparnum. #1}\label{#2}}\hspace{0.5em}}}

\renewcommand{\themyaparnum}{A\halfff\fff.\fff\themyappendnumber.\arabic{myaparnum}}

\newcounter{myapparnum}[mysectionnumber]
\newcommand{\proof}{\vspace{\medskipamount}{\textbf{{\emph{Proof}.}}\hspace*{1em}}}
\newcommand{\prooftitle}[1]{\vspace{\medskipamount}{\textbf{{\emph{#1}.}}\hspace*{1em}}}

\newcommand{\eproof}{ $\blacksquare$}

\newcommand{\dis}{\displaystyle}

\def\sss{\hspace{0.05em}\ }
\def\dss{\hspace{0.1em}\ }
\def\trs{\hspace{0.15em}\ }
\def\qss{\hspace{0.2em}\ }
\def\pss{\hspace{0.3em}\ }
\def\oss{\hspace{0.4em}\ }

\def\halfff{\hspace*{0.025em}}
\def\fff{\hspace*{0.05em}}
\def\dff{\hspace*{0.1em}}
\def\trf{\hspace*{0.15em}}
\def\qff{\hspace*{0.2em}}
\def\pff{\hspace*{0.3em}}
\def\off{\hspace*{0.4em}}

\newcommand{\nsp}{\hspace*{-0.1em}}
\newcommand{\nnsp}{\hspace*{-0.15em}}

\newcommand{\dnsp}{\hspace*{-0.2em}}

\renewcommand{\leq}{\leqslant}
\renewcommand{\geq}{\geqslant}
\newcommand{\id}{\mathop{\mbox{id}}\nolimits}

\newcommand{\zzz}{\mathbb{Z}}

\newcommand{\rrr}{\mathbb{R}}

\newcommand{\num}[1]{|\qff #1 \qff|}

\newcommand{\ttoo}{\hspace*{0.2em}\longrightarrow\hspace*{0.2em}}

\newcommand{\spin}{\operatorname{Spin}\dff}

\newcommand{\ttt}{\othertau\fff}
\newcommand{\aaa}{o\dff}

\begin{document}

\setlength{\baselineskip}{12pt plus 0pt minus 0pt}
\setlength{\parskip}{12pt plus 0pt minus 0pt}
\setlength{\abovedisplayskip}{12pt plus 0pt minus 0pt}
\setlength{\belowdisplayskip}{12pt plus 0pt minus 0pt}

\newskip\smallskipamount \smallskipamount=3pt plus 0pt minus 0pt
\newskip\medskipamount   \medskipamount  =6pt plus 0pt minus 0pt
\newskip\bigskipamount   \bigskipamount =12pt plus 0pt minus 0pt

\author{Nikolai\qss V.\qss Ivanov}
\title{Group\qss actions\qss on\qss complexes,\oss Kozsul\qss models,\oss presentations,\oss 
and\qss a\qss theorem\qss of\pss Coxeter}
\date{}

\footnotetext{\hspace*{-0.65em}\copyright\oss 
Nikolai\qss V.\qss Ivanov,\oss 2023.\trs }

\footnotetext{\hspace*{-0.65em}The picture of\dss a dodecahedron on\dss p.\qss 27\qss is\dss
adapted\dss by\trs F.\qss Letoutchaia\dss from\dss the image\qss
created\dss by\trs K.\qss Andr\'{e}\qss and distributed\dss at\qss 
https\halfff:/\!/\hspace*{-0.06em}commons.wikimedia.org/wiki/File:Dodecahedron.svg\qss under\dss
Creative Commons\dss Attribution-Share\dss Alike\qss 3.0\dss Unported\dss license,\oss
https\halfff:/\!/\hspace*{-0.06em}creativecommons.org/licenses/by-sa/3.0/deed.en.\oss}

\maketitle

\vspace{6ex}

{\renewcommand{\baselinestretch}{1.05}
\selectfont

\myit{\hspace*{0em}\large Contents}\vspace*{1.5ex} \\ 
\hbox to 0.8\textwidth{\myit{Introduction} \hfil 2}\hspace*{0.5em} \vspace*{1.5ex}\\
\hbox to 0.8\textwidth{\myit{\phantom{1}1.}\hspace*{0.5em} Generators and\sss relations for actions\sss transitive on\sss vertices \hfil 5}\hspace*{0.5em} \vspace*{0.25ex}\\
\hbox to 0.8\textwidth{\myit{\phantom{1}2.}\hspace*{0.5em} Kozsul\dss models of\dss actions\sss transitive on vertices \hfil 9}\hspace*{0.5em} \vspace*{0.25ex}\\
\hbox to 0.8\textwidth{\myit{\phantom{1}3.}\hspace*{0.5em} Implications\sss between edge and edge-loop relations \hfil 14}\hspace*{0.5em} \vspace*{0.25ex}\\
\hbox to 0.8\textwidth{\myit{\phantom{1}4.}\hspace*{0.5em} Presentations \hfil 18}\hspace*{0.5em} \vspace*{0.25ex}\\
\hbox to 0.8\textwidth{\myit{\phantom{1}5.}\hspace*{0.5em} Examples  \hfil 22}\hspace*{0.5em}\vspace*{1.5ex}\\
\hbox to 0.8\textwidth{\hspace*{4em} \emph{Fundamental\trs groups\dss of\pss CW-complexes}  \hfil 22}\hspace*{0.5em}\vspace*{0.25ex}\\
\hbox to 0.8\textwidth{\hspace*{4em} \emph{Symmetric\sss groups}  \hfil 23}\hspace*{0.5em}\vspace*{0.25ex}\\
\hbox to 0.8\textwidth{\hspace*{4em} \emph{Rotations of\trs a\sss regular dodecahedron}  \hfil 25}\hspace*{0.5em}\vspace*{0.25ex}\\
\hbox to 0.8\textwidth{\hspace*{4em} \emph{The binary\dss icosahedral\dss group}  \hfil 29}\hspace*{0.5em}\vspace*{1.25ex}\\
\hbox to 0.8\textwidth{\myit{\phantom{1}6.}\hspace*{0.5em} Coxeter's\qss 
implication \hfil 34}\hspace*{0.5em} \vspace*{0.25ex}\\
\hbox to 0.8\textwidth{\myit{\phantom{1}7.}\hspace*{0.5em} Actions\sss with several\sss orbits of\dss vertices \hfil 42}\hspace*{0.5em} \vspace*{1.5ex}\\
\myit{Appendices}\hspace*{0.5em}  \hspace*{0.5em} \vspace*{1.5ex}\\
\hbox to 0.8\textwidth{\myit{A.1.}\hspace*{0.5em} Coxeter's\dss proof\dss
of\trs his implication\hfil 51}\hspace*{0.5em} \vspace*{0.25ex}\\
\hbox to 0.8\textwidth{\myit{A.2.}\hspace*{0.5em} Coxeter's\qss implication\qss and\qss universal\qss 
central\qss extensions\hfil 52}\hspace*{0.5em}  \vspace*{0.25ex}\\
\hbox to 0.8\textwidth{\myit{A.3.}\hspace*{0.5em} Cayley\dss diagrams and\dss 
scaffoldings\hfil 56}\hspace*{0.5em} \vspace*{1.5ex}\\
\hbox to 0.8\textwidth{\myit{References}\hspace*{0.5em}\hfil 57}\hspace*{0.5em}  \vspace*{0.25ex}  

}

{\renewcommand{\baselinestretch}{1.05}
\selectfont

\newpage
\mynonumbersection{Introduction}

\myuppar{Actions and\sss presentations of\dss discrete groups.}
Let $G$ be a discrete group acting on a simply-connected\dss CW\dnsp-complex $Z$ 
preserving\sss the\dss CW-structure.\oss 
Such actions provide an efficient\sss tool\sss
for proving\sss that\sss the group $G$\sss is\dss finitely\sss presented
when\sss it\dss is.\oss
In more details,\oss if\dss the stabilizers of $0$\dnsp-cells are 
finitely\sss presented,\oss the stabilizers of $1$\dnsp-cells are
finitely\sss generated,\oss and\dss the number of\dss orbits of\dss
cells of\dss dimension $\leq\qff 2$\sss is\dss finite,\oss
then $G$\sss is\dss finitely\sss presented.\oss
We will\sss call\dss this\sss theorem\sss the\qss \emph{Brown\dss finiteness\dss theorem}.\oss
A\sss proof\dss of\dss this\sss theorem,\oss
based on\sss the\dss Bass--Serre\dss theory\qss \cite{s}\qss of\dss group acting on\sss trees,\oss
was given\sss in\dss 1983\dss by\dss K.S.\dss Brown\qss \cite{br}.\oss
Actually,\oss the full\qss \emph{Brown\dss theorem}\pss \cite{br}\qss 
is\dss more general\sss and\sss precise.\oss
If\dss we know presentations of\dss the stabilizers of $0$\dnsp-cells,\oss
generators of\dss the stabilizers of $1$\dnsp-cells,\oss
and\dss the gluing maps of $2$\dnsp-cells,\oss
then\dss Brown\dss theorem provides a presentation of\sss $G$\nnsp.\oss
No\sss finiteness assumptions are needed\sss in\sss general,\oss
but\sss under\sss the above assumptions we get\sss a finite presentation.\oss
A more elementary\sss proof\dss of\qss Brown's\dss theorem was given\sss by\trs
M.\dss Armstrong\qss \cite{a1},\qss \cite{a2}.\oss

In\sss fact,\qss at\sss least\sss in\sss the case of\dss actions\sss
transitive on $0$\dnsp-cells,\qss Brown\dss theorem\dss 
was known\sss before\sss his work.\oss
This special\sss case was implicitly\sss used\sss in\dss 1978\dss 
by\trs Hatcher\dss and\dss Thurston\qss \cite{ht}\qss
and\sss was explicitly\sss stated\dss 
by\dss Laudenbach\dss in an exposition\qss \cite{l}\qss of\dss their work.\oss
But\sss the paper\qss \cite{l}\qss is\dss focused on\sss the construction of\dss $Z$
and\sss the action of\sss $G$ when $G$\sss is\dss the mapping class group
of\dss a closed surface,\oss and\sss the proof\dss of\dss the general\dss theorem about\sss
actions and\sss presentations\dss is\dss hardly discussed.\oss
Apparently,\oss the papers\qss \cite{ht}\qss and\qss \cite{l}\qss were not\sss
known\sss to\dss Brown.\oss

\myuppar{Kozsul's\dss method.} 
A fairly\sss general\sss special\sss case of\trs Brown\dss finiteness\sss theorem
was proven\sss in\dss 1965\dss by\trs Kozsul\qss \cite{k}.\oss
Let\sss us assume\sss that\sss the action of\sss $G$\sss is\dss transitive on $0$\dnsp-cells.\oss
Suppose\sss that\sss the \dnsp$1$\dnsp-skeleton 
$X\off =\off Z_{\dff 1}$ of\sss $Z$\sss is\dss
a graph without\dss loops and\sss multiple edges
and\sss that\sss $X$\sss is\qss
\emph{locally\dss finite}\pss in\sss the sense\sss that\sss every\sss vertex\dss
is\dss an endpoint\sss of\dss only a\sss finite number of\dss edges.\oss
Suppose\sss further\sss that\sss the number of\dss orbits of $2$\dnsp-cells\dss
is\dss finite.\oss
If\dss also\sss the stabilizers of\dss vertices\qss 
are finitely\sss presented,\oss then $G$\sss is\dss finitely\sss presented.\oss
Kozsul\qss \cite{k}\qss uses a different\sss but\sss equivalent\dss language,\oss 
and\dss this claim\dss is\dss equivalent\sss to\qss \cite{k},\oss 
Theorem\qss 2\qss of\trs Chapter\qss 3.\footnotemark\oss
Note\sss that\sss under\sss these assumptions\sss the stabilizers of\dss
edges are automatically\sss finitely\sss presented.\oss

\footnotetext{Unfortunately,\oss in\qss \cite{k}\qss
the definition of\dss homotopies of\dss paths in\sss graphs contains a misprint.\oss
The expression\sss 
$(\trf a_{\dff 0}\dff,\qff \ldots\dff,\qff 
a_{\dff i}\dff,\qff a_{\dff i\dff +\dff 1}\dff,\qff
\ldots\dff,\qff a_{\dff n}\trf)$\sss
at\sss the\sss top of\trs p.\dss 37\dss
should\sss be replaced\dss by either\sss
$(\trf a_{\dff 0}\dff,\qff \ldots\dff,\qff 
a_{\dff i}\dff,\qff a_{\dff i\dff +\dff 2}\dff,\qff
\ldots\dff,\qff a_{\dff n}\trf)$\nnsp,\oss
as was suggested\sss in\qss \cite{i},\oss
or\halfff,\oss better\halfff,\oss by\sss
$(\trf a_{\dff 0}\dff,\qff \ldots\dff,\qff 
a_{\dff i}\dff,\qff a_{\dff i}\dff,\qff a_{\dff i\dff +\dff 1}\dff,\qff
\ldots\dff,\qff a_{\dff n}\trf)$\nnsp,\oss
giving a more general\dss result.\oss}

Like every other\sss proof\dss of\dss such\sss theorems,\oss
Kozsul\trs proof\dss begins by writing down a candidate\sss
to a presentation of\sss $G$ modulo\sss the stabilizer of\dss a vertex.\oss
This presentation defines a group $\mathbb{G}$\sss together with\sss
a canonical\dss homomorphism\sss
$\varphi\dff \colon\dff \mathbb{G}\qff \ttoo\qff G$\nnsp.\oss
In order\sss to prove\sss that\sss $\varphi$\sss is\dss an\sss isomorphism,\oss
Kozsul\dss uses $\mathbb{G}$ to construct\sss a new\sss graph
$\mathbb{X}$\sss together\sss with a morphism of\dss graphs\sss
$f\dff \colon\dff \mathbb{X}\qff \ttoo\qff X$\nnsp.\oss
Then\dss Kozsul\trs proves\sss that $f$ is\dss a covering of\dss graphs,\oss
and\sss the simply-connectivity assumption\sss implies\sss 
that $f$ is\dss actually an\sss isomorphism.\oss 
This\sss implies\sss that\sss $\varphi$\sss is\dss an\sss isomorphism,\oss
and\dss hence provides\sss the desired\dss presentation of\sss $G$ 
modulo\sss the stabilizer\halfff.\oss
We will\sss call\sss the graph $\mathbb{X}$\sss the\qss 
\emph{Kozsul\dss model}\pss of\sss $X$\nnsp.\oss
The construction of\dss the\dss Kozsul\dss model\sss $\mathbb{X}$\sss
is\dss similar\sss to\sss the construction of\dss the universal\sss covering\sss
$\widetilde{X}$\sss relative\sss to a graph of\dss groups\sss 
$(\trf G\fff,\qff Y\trf)$ in\sss the\dss Bass--Serre\dss theory\qss \cite{s},\pss
but\sss precedes\sss the\sss latter\sss by about\sss 3 years.\oss
Cf.\qss \cite{s},\oss Section\qss 5.3.

Kozsul's\dss candidate\sss to a presentation of\sss $G$ 
modulo\sss the stabilizer depends on several\sss choices.\oss
First,\oss one needs\sss to fix a vertex $v$ of\sss $X$\nnsp.\oss
Let\sss $G_{\dff v}$ be\sss the stabilizer of $v$\nnsp.\oss
Since\sss the action\dss is\dss transitive on vertices,\oss
fixing $v$\sss is\dss hardly significant.\oss
Let $E_{\dff v}$ be\sss the set\sss of\dss edges having $v$ as an endpoint.\oss
By\sss the\sss local\sss finiteness $E_{\dff v}$ is\dss finite.\oss
Next,\oss for every $e\qff \in\qff E_{\dff v}$ one needs\sss
to choose an element $s_{\dff e}\qff \in\qff G$\sss 
taking $v$\sss to\sss the other endpoint\sss of\sss $e$\nnsp.\oss
Kozsul's\sss generators are symbols $g_{\dff e}$\sss
mapped\dss by $\varphi$\sss to\sss the corresponding elements $s_{\dff e}$\nsp.\oss
Finally,\oss one needs\sss to choose a\sss finite set\sss of\dss generators of\dss $G_{\dff v}$\nsp.\oss 
The\sss local\dss finiteness of\sss $X$ 
and\dss the finiteness of\dss the number of\dss orbits of\sss $2$\dnsp-cells ensure\sss
that\sss there are only\sss finitely many generators and\sss relations.\oss

\myuppar{Looking\dss for\sss natural\dss presentations.}
Let\sss us\sss temporarily\sss 
ignore\sss the finiteness properties
and\sss instead\dss look\sss for presentations which are natural\dss in some sense.\oss
Then\sss there\dss is\dss no need\sss to assume\sss that\sss $X$\sss is\dss locally\sss finite
and\sss that\sss there\dss is\dss a finite number of\dss orbits of\sss $2$\dnsp-cells.\oss
As\sss the set\sss of\dss generators of\sss $G_{\dff v}$
one can simply\sss take\sss the whole group\sss $G_{\dff v}$\nsp.\oss
The constructions of\dss the group $\mathbb{G}$ and of\dss the\dss Kozsul\trs model\sss still\sss work,\oss
and\sss we get\sss a fairly canonical\dss presentation.\oss
The set\sss of\dss generators\dss is\dss essentially\sss $E_{\dff v}$\nnsp,\oss
while\sss the relations depend on\sss the choice of\dss elements 
$s_{\dff e}\dff,\qff e\qff \in\qff E_{\dff v}$\nsp.\oss
Such dependence seems\sss to be inevitable.\oss
The most\sss interesting\sss relations,\pss
which we call\dss the\qss \emph{loop\dss relations},\oss 
are determined\dss by\sss loops used\dss to glue\sss the $2$\dnsp-cells.\oss 
A key role\dss is\dss also played\dss by relations defined\dss by\sss loops
of\dss length $2$ obtained\dss by\sss following an edge $e\qff \in\qff E_{\dff v}$
and\dss then\sss following\sss the same edge backwards,\oss
which we call\dss the\qss \emph{edge-loop relations}.\footnotemark\oss

\footnotetext{See\trs Lemma\qss \ref{symmetry}\qss below and\dss the remark\sss following
its proof\halfff.\oss
Kozsul\qss \cite{k}\qss did\sss not\sss discussed\dss these relations,\oss
although\sss in\sss his approach they are necessary and automatically\sss included.\oss}

This constructions of\dss $\mathbb{G}$ and of\dss the\dss
Kozsul\dss model\sss $\mathbb{X}$\sss are discussed\sss in\dss Sections\qss
\ref{one}\qss and\qss \ref{models}\qss respectively.\oss
The main\sss results\dss is\trs Corollary\qss \ref{isomorphism-complexes}\qss
to\sss the effect\sss that\sss the maps\sss
$f\dff \colon\dff \mathbb{X}\qff \ttoo\qff X$
and\sss
$\varphi\dff \colon\dff \mathbb{G}\qff \ttoo\qff G$
are isomorphisms.\oss
Since $\mathbb{G}$\sss is\dss constructed\dss in\sss terms of\sss $G_{\dff v}$
with\sss the help of\dss explicit\sss generators and\sss relations,\oss
this gives us a presentation of\sss $G$ modulo $G_{\dff v}$\nsp.\oss
If\dss there\dss is\dss a natural\sss choice of\dss elements
$s_{\dff e}\dff,\qff e\qff \in\qff E_{\dff v}$\nsp,\oss
then\sss the corresponding\sss presentation\dss is\dss also natural.\oss

\myuppar{Back\sss to finite presentations.}
The presentation\sss from\sss the previous subsection\dss is\dss very\sss redundant.\oss
The relations are far\sss from\sss being\sss independent,\qss and,\qss moreover\halfff,\qss
can be used\dss to express some generators in\sss terms of\dss the others.\oss
In\dss Section\qss \ref{implications}\qss we prove\sss the main\sss implications 
between\sss these relations.\oss
The proofs are simple,\oss but,\oss 
to\sss the best\sss knowledge of\dss the author\halfff,\oss 
these implications were not\sss even stated explicitly\sss before.\oss
In\dss Section\qss \ref{presentations}\qss we show\sss that\sss
many\sss relations can\sss be\sss turned\sss into definitions of\dss some generators;\oss
after\sss this\sss these relations can be discarded,\oss
and\dss then\sss the remaining\sss relations should\sss be rewritten\sss 
in\sss terms of\dss these definitions.\oss
When\sss the group $G_{\dff v}$\sss is\dss finitely\sss presented
and\dss there are only a finite number of\dss orbits of\sss $1$\dnsp-cells and $2$\dnsp-cells,\oss
then\sss the resulting\sss presentation\dss is\dss finite.\oss
This procedure\dss is\dss not\sss canonical\dss because
it\sss requires choosing\sss which\sss relations\sss will\dss be\sss turned\sss into definitions.\oss
We call\sss the collection of\dss necessary\sss choices a\qss \emph{scaffolding},\oss
and\sss in\dss Appendix\qss \ref{cayley-diagrams}\qss we relate scaffoldings\sss
to\dss Cayley\dss diagrams.\oss
Often\sss there\dss is\dss a natural\sss scaffolding\halfff,\oss
and\sss the resulting\sss presentation\dss is\dss also natural.\oss
Surprisingly,\oss for well-known groups\sss this procedure often\sss 
results in\sss their\sss well-known\dss presentations.\oss
We illustrate\sss this in\dss Section\qss \ref{examples},\oss
where we show in details how\sss this procedure works for\sss
the fundamental\sss groups of\dss CW-complexes,\oss
symmetric groups,\oss the group of\dss rotations of\dss a dodecahedron,\oss
and\dss the binary\sss icosahedral\dss group.

\myuppar{The\dss Coxeter's\dss implication.}
Using\sss the action of\dss the binary\sss icosahedral\dss group\sss 
on\sss
the regular dodecahedron\sss results\sss in\sss its presentation
with\sss generators\sss $g\fff,\qff r$ subject\sss to\sss the relations\vspace{3pt}\vspace{-1.5pt}
\[
\quad
g^{\dff 2}
\off =\off 
r^{\dff 3}
\off =\off
(\dff r\fff g\trf)^{\dff 5}
\quad
\mbox{and}\quad
r^{\dff 6}\off =\dff\off 1
\qff.
\]

\vspace{-9pt}\vspace{-1.5pt}
The substitutions\sss
$z
\off =\off 
g^{\dff 2}
\off =\off 
r^{\dff 3}$\dnsp,\qss
$s
\off =\off 
r^{\dff -\dff 1}$\dnsp,\oss
and\dss
$t
\off =\off
r\fff g$\sss
turns\sss this presentation\sss into\vspace{3pt}\vspace{-1.5pt}
\[
\quad
s^{\dff 3}
\off =\off 
t^{\dff 5}
\off =\off 
(\dff s\dff t\trf)^{\dff 2}
\off =\off
z\quad\
\mbox{and}\quad\
z^{\dff 2}
\off =\off
1
\qff
\]

\vspace{-9pt}\vspace{-1.5pt}
with\sss the generators $s\fff,\qff t\fff,\qff z$\nnsp.\oss
This\dss is\dss
the classical\trs Coxeter\qss \cite{c1}\qss presentation\qss
({\fff}he denotes\sss the binary\sss icosahedral\dss group by\sss 
$\langle\trf 2\fff,\qff 3\fff,\qff 5 \trf\rangle$\nsp).\oss
Coxeter\qss \cite{c1}\qss also proved a quite remarkable fact\sss
that\sss the relations\dss
$s^{\dff 3}
\off =\off 
t^{\dff 5}
\off =\off 
(\dff s\dff t\trf)^{\dff 2}
\off =\off
z$\sss
automatically\sss imply\sss that\sss
$z^{\dff 2}\off =\off 1$\nnsp.\oss
The present\sss author was not\sss able\sss to resist\sss
the\sss temptation\sss to discuss\sss this implication\sss in details.\oss

Coxeter's\qss \cite{c1}\trs proves\sss this implication by an\sss  
ingenious unmotivated\sss calculation.\oss
We reproduce\sss this proof\dss in\dss Appendix\qss \ref{coxeter-proof}.\oss
In\dss Section\qss \ref{coxeter-implication}\qss we use\sss the methods of\trs
Sections\qss \ref{one}\dss --\dss \ref{examples}\qss in order\sss to give
a conceptual\dss proof\dss of\trs Coxeter's\dss implication.\oss
In particular\halfff,\oss in\sss this proof\dss the exponent\sss $2$\sss
in\sss the relation $z^{\dff 2}\off =\off 1$\sss appears as\sss 
the\dss Euler\dss characteristic of\dss the sphere $S^{\fff 2}$\dnsp.\oss
But\sss some mystery\sss still\sss remains.\oss
In\dss Appendix\qss \ref{central-extensions}\qss we present\sss another
proof\dss of\trs Coxeter's\dss implication.\oss
It\dss is\dss based on\sss the\sss theory of\dss the universal\sss central\sss extensions,\oss
and\sss in\sss this proof\dss the exponent\sss $2$ appears as\sss the order 
of\dss the fundamental\dss group\sss
$\pi_{\dff 1}\dff(\trf SO\dff(\dff 3\dff)\trf)$\nnsp.\oss

\myuppar{Actions with several\sss orbits of\dss the vertices.}
This\dss is\dss the\sss topic of\trs Section\qss \ref{several}.\oss
The proofs in\sss this case are very\sss similar\sss to\sss the proofs
for actions\sss transitive on vertices,\oss and we\sss limit\sss
ourselves by\sss the arguments where\sss the changes are not\sss completely obvious.\oss

\myuppar{Notations.}
In\sss that\sss follows $X$\sss is\qss \emph{a graph without\dss loops and\sss multiple edges},\oss
and $G$\sss is\dss a group acting on $X$\sss from\sss the\sss left.\oss
The above CW-complex $Z$ will\dss not\sss appear explicitly.\oss
For a vertex $v$ of\dss $X$ we denote by\sss $G_{\dff v}$\sss its stabilizer\halfff.\oss
For an edge $e$  
we understand\dss by\sss its\qss \emph{stabilizer}\dss $G_{\dff e}$
the intersection $G_{\dff e}\off =\off G_{\dff v}\dff \cap\trf G_{\dff w}$
of\dss the stabilizers of\dss endpoints $v\fff,\qff w$ of\sss $e$\nnsp.\oss

\newpage
\mysection{Generators\qss and\qss relations\qss for\qss actions\qss transitive\qss on\qss vertices}{one}

\myuppar{Generators.}
In\dss Section\qss \ref{one} -- \ref{examples}\qss the action of\dss $G$ on $X$\sss is\dss assumed\dss 
to be\sss transitive on\sss the set\sss of\dss vertices of\dss $X$\sss
and $X$ is\dss assumed\sss to be a $1$\dnsp-dimensional simplicial\sss complex,\oss 
i.e.\qss a graph without\dss loops and\sss multiple edges.\oss
In\dss Section\qss \ref{one} -- \ref{presentations}\qss
we will\sss assume\sss that\sss a vertex $v$ of\dss $X$\sss is\dss fixed.\oss
Let\sss 
$E\off =\off E_{\fff v}$\dss 
be\sss the set\sss of\dss edges of\dss $X$\sss having $v$ as an endpoint\halfff.\oss
It\dss is\dss suggestive\sss to\sss think\sss that\sss an edge \nsp$e\qff \in\qff E$\sss is\dss
directed\sss from $v$\sss to\sss its other endpoint,\oss
which we denote by \nsp$\ttt(\dff e\trf)$\dnsp. 

Let us choose for every\dss $e\qff \in\qff E$\dss 
an element $s_{\dff e}\qff \in\qff G$\dss such\dss that\sss 
$s_{\dff e}\dff(\dff v\trf)
\off =\off
\ttt(\dff e\trf)$\nnsp.\oss
Since $X$ has no\sss loops or multiple edges,\dss 
$e$\sss is\dss uniquely determined\dss by\sss $s_{\dff e}$\nsp.\oss 
For each\dss $e\qff \in\qff E$\dss 
let\sss $g_{\dff e}$ be a symbol\sss corresponding\dss to\sss the edge $e$\nnsp,\oss
thought\sss as a generator of\sss $G$ modulo $G_{\dff v}$\nsp.\oss 
Let\vspace{3pt}\vspace{-0.5pt}
\[
\quad
\mathcal{S}
\off =\off 
\{\qff s_{\dff e}\qff |\qff e\qff \in\qff E \qff\}
\quad\
\mbox{and}\quad\
\mathcal{G}\off =\off \{\qff g_{\dff e}\qff |\qff e\qff \in\qff E \qff\} 
\qff.
\]

\vspace{-9pt}\vspace{-0.5pt}
Let\dss $\mathcal{F}$\dss be\sss the free group on\dss the set\sss $\mathcal{G}$\dss 
(or\halfff,\oss what\dss is\dss essentially\dss the same,\oss on $E$),\oss 
and\dss let\dss $\mathcal{F}\fff \ast\dff G_{\dff v}$\dss be,\oss as usual,\oss 
the free product\sss of\dss 
$\mathcal{F}$ and $G_{\dff v}$\nsp.\oss 
There\dss is\dss an obvious homomorphism\vspace{3pt}\vspace{-0.5pt}
\[
\quad
\psi\dff \colon\dff 
\mathcal{F}\fff \ast\dff 
G_{\dff v}\qff \ttoo\qff G
\qff
\]

\vspace{-9pt}\vspace{-0.5pt}
equal\dss to\sss the inclusion\dss 
$G_{\dff v}\qff \ttoo\qff G$\dss on $G_{\dff v}$ and\dss 
taking\dss $g_{\dff e}$\dss to\dss $s_{\dff e}$\nsp.\oss 
While\sss the group\dss $\mathcal{F}\dff \ast\dff G_{\dff v}$\dss is\dss canonically determined\dss 
by\dss the choice of\trs the vertex $v$\nnsp,\oss the homomorphism $\psi$
depends on\dss the choice of\dss elements\dss $s_{\dff e}$\nsp.\oss 
The freedom in choosing elements\dss $s_{\dff e}$\dss such\dss that\dss 
$s_{\dff e}\dff(\dff v\trf)$\dss is\dss the other endpoint\sss of\dss $e$ will\dss 
play a role in\qss Section\qss 3,\oss
and,\oss especially,\oss in\sss applications.\oss

\myuppar{Relations.}
Our\sss first\sss goal\dss is\dss to prove\sss that\dss the homomorphism $\psi$\sss is\dss surjective,\oss 
and\dss to provide a set $R$ normally\sss generating\sss 
the kernel\sss $\ker\qff \psi$ of\sss $\psi$\dnsp.\oss 
This set\sss $R$\sss should\dss be\sss thought\sss as a 
set\sss of\dss relations of\trs the group\sss $G$\sss modulo\sss the subgroup\sss $G_{\dff v}$
with\sss generators $g_{\dff e}$\nsp, $e\qff \in\qff E$\nnsp.\oss 
If\dss a presentation of\dss $G_{\dff v}$\sss is\dss already known,\oss 
then a presentation of\dss $G$ can\dss be obtained\dss by\sss adding\dss to it\dss 
new generators\dss $g_{\dff e}$\nnsp, $e\qff \in\qff E$\nnsp,\oss 
and\dss the elements of\sss $R$ as additional relations.\oss 
Since we\sss think about\dss the elements of\sss $R$ as\sss the
relations of\dss the group $G$ modulo $G_{\dff v}$\nsp,\oss
we will\dss usually\dss write\sss them\dss in\dss the form of\dss
equalities,\oss with\dss $a\off =\off b$\dss
meaning\dss that\sss  
$a\dff \cdot\dff b^{\dff -\dff 1}
\qff \in\pff 
\ker\qff \psi$\dnsp.\oss

It\sss turns out\sss that\sss $R$\sss consists of\dss
two\sss types of\dss elements,\oss
which we will\sss call\dss the\qss \emph{edge relations}\qss
and\dss the\qss \emph{loop relations}.\oss
The edge relations depend only on\sss the choice of\dss elements $s_{\dff e}$\sss
and\sss the action of\dss $G_{\dff v}$ on $E$\nnsp,\oss and\dss hence
are\qss \emph{local}\pss in a definite sense.\oss
In contrast,\oss the\sss loop relations are\sss global\sss and depend
on a choice of\dss a set\sss $\mathcal{L}$\sss of\trs loops in\sss the graph $X$\nnsp.\oss
In\sss this section we will\sss describe\sss relations of\dss both\sss types,\oss
and\sss in\dss Section\qss \ref{models}\qss we will\sss prove\sss that\sss
they normally\sss generate\sss $\ker\qff \psi$\dnsp.\oss
In\dss Section\qss \ref{implications}\qss
and\qss \ref{presentations}\qss we will\sss 
explain how\dss one can simplify\sss the resulting\sss presentation.\oss

\mypar{Lemma.}{transitivity} 
\emph{If\pss $X$\dss is\dss connected,\oss 
then\dss $\psi\dff(\trf \mathcal{F}\trf)$\dss 
acts\sss transitively\sss on\dss the set\dss of\qss vertices of\pss $X$\nnsp.\oss}

\proof 
Let $w$ be a vertex of\dss $X$.\oss 
Since $X$ is connected,\oss 
we can connect $v$ with $w$ by
a sequence\qss 
$v\off =\off v_{\dff 0}\dff,\off 
v_{\dff 1}\dff,\off 
\ldots\dff,\off 
v_{\dff n}\off =\off w$\qss 
of\dss vertices $v_{\dff i}$ such\dss that\trs 
$v_{\dff i}\dff,\pff v_{\dff i\dff +\dff 1}$\dss 
are connected\dss by\sss an edge for\qss 
$0\qff \leq\qff i\qff \leq\qff n\qff -\qff 1$\nnsp.\oss
We will\dss prove\sss that 
$w\off =\pff g\dff(\dff v\trf)$\sss for some 
$g\qff \in\qff \psi\dff(\trf \mathcal{F}\trf)$\dss by an induction by $n$\nnsp.\oss 
For\dss $n\off =\off 1$\dss this holds by\sss the choice of\dss elements\dss $s_{\dff e}$\nsp.\oss
If\sss $g\dff(\dff v\trf)\off =\off v_{\dff n\dff -\dff 1}$  and
$g\qff \in\qff \psi\dff (\trf \mathcal{F}\trf)$\nnsp,\oss 
then 
$v\off =\off g^{\dff -\dff 1}\dff (\dff v_{\dff n\dff -\dff 1}\dff)$\dss 
is\dss connected\dss by\sss an edge $e$ with 
$g^{\dff -\dff 1}\dff (\dff v_{\dff n}\dff )$\nnsp.\oss
Therefore\dss\vspace{3pt}
\[
\quad
g^{\dff -\dff 1}\dff (\dff v_{\dff n}\dff )
\off =\off
s_{\dff e}\dff \cdot\dff  g^{\dff -\dff 1}\dff (\dff v_{\dff n\dff -\dff 1}\dff)
\]

\vspace{-12pt}\vspace{3pt} 
by\dss the choice of  $s_{\dff e}$\nsp.\oss 
Hence
$v_{\dff n}
\qff =\qff 
g\dff \cdot\dff  s_{\dff e}\dff \cdot\dff g^{\dff -\dff 1}\dff (\dff v_{\dff n\dff -\dff 1}\dff)
\qff =\qff 
g\dff \cdot\dff  s_{\dff e}\trf  (\dff v\trf)$\nsp.\oss
The\dss lemma\dss follows.\oss  \eproof

\mypar{Corollary\halfff.}{generators} 
\emph{If\qss  $X$\sss is\dss connected,\oss 
then\dss $G$\dss is\dss generated\dss by\trs 
$\psi\dff(\trf \mathcal{F}\trf)$\dss and\dss $G_{\dff v}$\nsp.\oss 
Moreover\halfff,\oss 
$G\off =\off \psi\dff(\trf \mathcal{F}\trf)\dff \cdot\dff G_{\dff v}$\nsp.\oss 
In\dss particular\halfff,\oss
the homomorphism\qss 
$\psi\dff \colon\dff 
\mathcal{F}\fff \ast\dff 
G_{\dff v}\qff \ttoo\qff G$\qss is\dss surjective.\oss}  \eproof

\proof
If\sss $g\qff \in\qff G$\nnsp,\oss
then\sss $g\trf(\dff v\trf)\off =\off h\trf(\dff v\trf)$\sss
for some\sss $h\qff \in\qff \psi\dff(\trf \mathcal{F}\trf)$\sss
by\sss the\sss lemma.\oss  \eproof

\myuppar{The edge relations.}
The\qss \emph{edge relations}\pss
correspond\dss to\sss the pairs\dss
$(\dff e\fff,\pff t\trf)$\nnsp,\oss
where\dss $e\qff \in\qff E$\dss and\dss $t\qff \in\qff G_{\dff v}$\nsp.\oss
Given such $e$ and $t$\nnsp,\oss 
the vertex\dss $t\dff \cdot\dff s_{\dff e}\dff (\dff v\trf)$\dss is\dss 
connected\dss with $v$ by\dss the edge\sss $t\trf(\dff e\dff)$\dss 
(since\sss $s_{\dff e}\dff (\dff v\trf)$\sss is\dss connected by\sss $e$ with $v$\nnsp,\oss 
and\dss $t\qff \in\qff G_{\dff v}$\dss fixes $v$\nsp).\oss
Therefore\vspace{1.5pt}
\[
\quad
t\dff \cdot\dff s_{\dff e}\dff (\dff v\trf)
\off =\off
s_{\dff d}\dff (\dff v\trf)
\qff,
\]

\vspace{-12pt}\vspace{1.5pt}
where\dss
$d\off =\off t\trf(\dff e\dff)\qff \in\qff E$\nnsp.\oss
It\dss follows\sss that\dss the element\trs\vspace{1.5pt}
\[
\quad
k\trf(\dff e\fff,\pff t\trf)
\off =\off
s_{\dff d}^{\dff -\dff 1}\dff \cdot\dff t\dff \cdot\dff s_{\dff e}
\]

\vspace{-12pt}\vspace{1.5pt}
fixes $v$ and\dss hence belongs\sss to $G_{\dff v}$\nsp.\oss
The\qss \emph{edge relation}\qss corresponding\dss to\sss the pair\dss
$(\dff e\fff,\pff t\trf)$\dss is\vspace{1.5pt}
\[
\quad
g_{\dff d}^{\dff -\dff 1}\dff \cdot\dff t\dff \cdot\dff g_{\dff e}
\off =\off
k\trf(\dff e\fff,\pff t\trf)
\off.
\]

\vspace{-12pt}\vspace{1.5pt}
The left-hand\sss side of\trs this relation\dss is\dss an element\sss of\trs 
$\mathcal{F}\fff \ast\dff G_{\dff v}$\dss and\dss the right-hand side\dss is\dss 
an element\sss of\qss $G_{\dff v}\qff \subset\pff \mathcal{F}\dff\ast\dff G_{\dff v}$\nsp.\oss
We will\sss denote\sss this relation\dss by\trs
$E\dff(\dff e\fff,\pff t\trf)$\nnsp.\oss
The relation\dss $E\dff(\dff e\fff,\pff t\trf)$\dss obviously\dss holds in\sss $G$\nnsp.\oss
Speaking\sss more formally\halfff,\oss 
this means\sss that\dss the element\vspace{1.5pt}
\[
\quad
E\dff(\dff e\fff,\pff t\trf)
\off =\off
\left(\dff
g_{\dff d}^{\dff -\dff 1}\dff \cdot\dff t\dff \cdot\dff g_{\dff e}
\dff\right)^{\dff -\dff 1}
\cdot\dff
k\trf(\dff e\fff,\pff t\trf)
\]

\vspace{-12pt}\vspace{1.5pt}
of\qss 
$\mathcal{F}\dff\ast\dff G_{\dff v}$\dss
belongs\sss to\sss the kernel\sss of\trs $\psi$\nnsp.\oss

\myuppar{The\sss loop relations.}
The\qss \emph{loop relations}\pss
correspond\dss to\sss simplicial\dss loops in\sss $X$\sss
starting and ending at\dss the vertex $v$\nnsp,\oss
i.e.\qss to\sss sequences\sss
$l
\off =\off
\{\trf v_{\dff i} \qff\}_{\trf 0\qff \leq\qff i\qff \leq\qff n}$
of vertices\dss $v_{\dff i}$\dss of\trs $X$\sss such\dss that\qss 
$v_{\trf 0}\off =\off v_{\dff n}\off =\off v$\qss and
$v_{\dff i}$ is connected with $v_{\dff i\dff +\dff 1}$ by\sss an edge for all\qss 
$0\qff \leq\qff i\qff \leq\qff n\qff -\qff 1$\nnsp.\oss 
The number $n$\sss is\dss called\dss the\qss \emph{length}\qss of\trs the loop $l$\nnsp.\oss
We claim\dss that\dss there\dss is\dss a unique sequence\dss 
$s_{\dff 1}\dff,\off s_{\dff 2}\dff,\off\ldots\dff,\off s_{\dff n}$\dss 
of\dss elements of\trs $\mathcal{S}$ such that\vspace{3pt}
\begin{equation}
\label{loop}
\hspace{0.5em}\begin{array}{l}
s_{\dff 1}\dff(\dff v_{\trf 0}\trf)
\off =\off 
v_{\dff 1}\qff,\vspace{9pt}\\
s_{\dff 1}\dff \cdot\dff s_{\dff 2}\trf(\dff v_{\trf 0}\trf)
\off =\off 
v_{\dff 2}\qff,\vspace{9pt}\\
s_{\dff 1}\dff \cdot\dff s_{\dff 2}\dff \cdot\dff s_{\dff 3}\dff(\dff v_{\trf 0}\trf)
\off =\off 
v_{\dff 3}\qff,\vspace{6pt}\\
\ldots\ldots\vspace{6pt}\\
s_{\dff 1}\dff \cdot\dff s_{\dff 2}\dff \cdot\dff
\ldots\dff \cdot\dff 
s_{\dff n}\trf(\dff v_{\trf 0}\trf)
\off =\off 
v_{\dff n}\qff.
\end{array}
\end{equation}

\vspace{-12pt}\vspace{3pt}
Indeed,\oss
by\dss the definition of\trs $\mathcal{S}$\sss there exists an element\trs 
$s_{\dff 1}\qff \in\qff \mathcal{S}$\dss such\dss that\sss\vspace{2.5pt}
\[
\quad
s_{\dff 1}\dff(\dff v_{\trf 0}\trf)
\off =\off
s_{\dff 1}\dff(\dff v\trf)
\off =\off 
v_{\dff 1} 
\]

\vspace{-12pt}\vspace{2.5pt}
Since $v_{\dff 1}$ is\dss connected with $v_{\dff 2}$
by an edge,\oss
the vertex\dss 
$v_{\trf 0}
\off =\off
s_{\dff 1}^{\qff -\dff 1}\dff(\dff v_{\trf 1}\trf)$\dss
is connected to\dss 
$s_{\dff 1}^{\qff -\dff 1}\dff(\dff v_{\trf 2}\trf)$\dss 
by an edge and\dss hence\sss there exists\dss
$s_{\dff 2}\qff \in\qff \mathcal{S}$\dss such\dss that\vspace{2.5pt}
\[
\quad
s_{\dff 2}\dff(\dff v_{\trf 0}\trf)
\off =\off
s_{\dff 2}\dff(\dff v\trf)
\off =\off 
s_{\dff 1}^{\qff -\dff 1}\dff(\dff v_{\trf 2}\trf)
\qff,
\]

\vspace{-12pt}\vspace{2.5pt}
i.e.\qss
$s_{\dff 1}\dff \cdot\dff s_{\dff 2}\trf(\dff v_{\trf 0}\trf)
\off =\off 
v_{\dff 2}$\nsp.\oss
Similarly\halfff,\oss since $v_{\trf 2}$ is\dss connected with $v_{\trf 3}$
by an edge,\pss\vspace{2.5pt}
\[
\quad
 v_{\trf 0}
\off =\off
(\dff s_{\dff 1}\dff \cdot\dff s_{\dff 2}\trf)^{\dff -\dff 1}\dff(\dff v_{\trf 2}\trf)
\]

\vspace{-12pt}\vspace{2.5pt}
is connected to\dss 
$(\dff s_{\dff 1}\dff \cdot\dff s_{\dff 2}\trf)^{\dff -\dff 1}\dff(\dff v_{\trf 3}\trf)$\dss
by an edge and\dss hence\sss there exists\dss
$s_{\dff 3}\qff \in\qff \mathcal{S}$\dss such\dss that\vspace{2.5pt}
\[
\quad
s_{\dff 3}\dff(\dff v_{\trf 0}\trf)
\off =\off
s_{\dff 3}\dff(\dff v\trf)
\off =\off 
(\dff s_{\dff 1}\dff \cdot\dff s_{\dff 2}\trf)^{\dff -\dff 1}\dff(\dff v_{\trf 3}\trf)
\qff,
\]

\vspace{-12pt}\vspace{2.5pt}
i.e.\qss
$s_{\dff 1}\dff \cdot\dff s_{\dff 2}\dff \cdot\dff s_{\dff 3}\trf(\dff v_{\trf 0}\trf)
\off =\off 
v_{\dff 3}$\nsp.\oss
By continuing\dss to argue in\dss this way\dss
we will\sss see\sss that\dss there exists a sequence\dss 
$s_{\dff 1}\dff,\off s_{\dff 2}\dff,\off\ldots\dff,\off s_{\dff n}$\dss
with\dss the required\dss properties.\oss
The same argument\sss shows\sss that\trs
$v_{\dff 0}\off =\off v$\dss is\dss connected\dss with\dss
$(\dff s_{\dff 1}\dff \cdot\dff s_{\dff 2}
\dff \cdot\dff 
\ldots
\dff \cdot\dff 
s_{\dff i\dff -\dff 1}\trf)^{\dff -\dff 1}\dff(\dff v_{\trf i}\trf)$\dss
by an edge for every\dss $i\qff \leq\qff n$\nnsp.\oss
In view of\qss (\ref{loop})\qss this implies\sss that\sss $s_{\dff i}$\dss
is\dss the unique element\sss of\sss $\mathcal{S}$\sss such\dss that\vspace{2.5pt}
\[
\quad
s_{\dff i}\dff(\dff v_{\trf 0}\trf)
\off =\off
s_{\dff i}\dff(\dff v\trf)
\off =\off 
(\dff s_{\dff 1}\dff \cdot\dff s_{\dff 2}
\dff \cdot\dff 
\ldots
\dff \cdot\dff 
s_{\dff i\dff -\dff 1}\trf)^{\dff -\dff 1}\dff(\dff v_{\trf i}\trf)
\qff.
\]

\vspace{-12pt}\vspace{2.5pt}
It\dss follows\dss that\dss the sequence\dss 
$s_{\dff 1}\dff,\off s_{\dff 2}\dff,\off\ldots\dff,\off s_{\dff n}$\qss
is\dss unique.\oss

The\sss last\sss equality\dss in\qss (\ref{loop})\qss means\sss that\dss
the product\dss
$s_{\dff 1}\dff \cdot\dff s_{\dff 2}
\dff \cdot\dff 
\ldots
\dff \cdot\dff 
s_{\dff n}$\dss
takes\dss
$v_{\dff 0}\off =\off v$\trs
to\dss
$v_{\dff n}\off =\off v$\nnsp.\oss
Therefore\dss
$s_{\dff 1}\dff \cdot\dff s_{\dff 2}
\dff \cdot\dff 
\ldots
\dff \cdot\dff 
s_{\dff n}
\qff \in\qff G_{\dff v}$\nsp.\qff\oss
Let\trs
$g_{\dff 1}\dff,\off g_{\dff 2}\dff,\off\ldots\dff,\off g_{\dff n}$\dss
be\sss the generators of\trs $\mathcal{F}$\dss corresponding\dss to\sss the elements\dss 
$s_{\dff 1}\dff,\off s_{\dff 2}\dff,\off\ldots\dff,\off s_{\dff n}$\dss of\sss $\mathcal{S}$\dss
(more formally\halfff,\oss if\qss
$s_{\dff i}\off =\off s_{\dff e_{\dff i}}$\nsp,\oss
then\qss
$g_{\dff i}\off =\off g_{\dff e_{\dff i}}$\nsp).\oss
The\qss \emph{loop relation}\qss corresponding\dss 
to\sss the\sss loop\sss $l$\dss is\vspace{1.5pt}\vspace{0.25pt}
\[
\quad
g_{\dff 1}\dff \cdot\dff g_{\dff 2}
\dff \cdot\dff 
\ldots
\dff \cdot\dff 
g_{\dff n}
\off =\off
s_{\dff 1}\dff \cdot\dff s_{\dff 2}
\dff \cdot\dff 
\ldots
\dff \cdot\dff 
s_{\dff n}
\qff.
\]

\vspace{-10.5pt}\vspace{0.25pt}
The left-hand\sss side of\trs this relation\dss is\dss an element\sss of\trs 
$\mathcal{F}\qff \subset\pff \mathcal{F}\dff\ast\dff G_{\dff v}$\dss and\dss the right-hand side\dss is\dss 
an element\sss of\qss $G_{\dff v}\qff \subset\pff \mathcal{F}\dff\ast\dff G_{\dff v}$\nsp.\oss
We will\sss denote\sss this relation\dss by\trs
$L\trf(\dff l\qff)$\nnsp.\oss
The relation\dss $L\trf(\dff l\qff)$\dss obviously\dss holds in\sss $G$\nnsp.\oss
Speaking more formally\halfff,\oss this means\sss that\dss the element\vspace{1.5pt}
\[
\quad
L\trf(\dff l\qff)
\off =\off
\left(\dff
g_{\dff 1}\dff \cdot\dff g_{\dff 2}
\dff \cdot\dff 
\ldots
\dff \cdot\dff 
g_{\dff n}
\dff\right)^{\dff -\dff 1} 
\cdot\dff
\left(\dff
s_{\dff 1}\dff \cdot\dff s_{\dff 2}
\dff \cdot\dff 
\ldots
\dff \cdot\dff 
s_{\dff n}
\dff\right)
\]

\vspace{-12pt}\vspace{1.5pt}
of\qss 
$\mathcal{F}\dff\ast\dff G_{\dff v}$\dss
belongs\sss to\sss the kernel\sss of\trs $\psi$\nnsp.\oss

\myuppar{The edge-loop relations.}
For $e\qff \in\qff E$\dss let\sss $l_{\dff e}$ be\sss the\sss loop of\trs
length $2$ starting at\sss $v$\nnsp,\oss
following $e$ and\dss immediately\dss returning\dss to $v$ along $e$\nnsp.\oss
More formally\halfff,\oss this\dss is\dss 
the\sss loop\dss $v\fff,\pff \ttt(\dff e\trf)\fff,\pff v$\nnsp.\oss
We call\sss such\dss loops\sss the\qss \emph{edge-loops}.\oss
The corresponding\dss relation\dss $L\trf(\dff l_{\dff e}\qff)$\dss
has\sss the form\vspace{3pt}
\[
\quad
g_{\dff e}\dff \cdot\dff g_{\dff a}
\off =\off
s_{\dff e}\dff \cdot\dff s_{\dff a}
\qff,
\]

\vspace{-12pt}\vspace{3pt}
where\sss 
$a\off =\off s_{\dff e}^{\dff -\dff 1}\dff(\dff e\trf)$\nnsp.\oss 
We will\sss call\sss these\sss loop
relations\sss the\qss \emph{edge-loop relations}.\oss
As we will\sss see,\oss their\sss role\dss is\dss more similar\sss
to\sss the role of\dss edge relation\sss than\sss to\sss the role of\dss
other\sss loop relations.\oss

\myuppar{Introducing\dss relations\dss in\dss $\mathcal{F}\dff\ast\dff G_{\dff v}$\nsp.}
\emph{From\sss now\sss on we will\sss assume\sss that\sss $X$\sss is\dss connected.\oss}
Let\sss $\mathcal{L}$\dss be a collection of\trs loops in\sss $X$\sss based at\sss $v$\nnsp.\oss
Let\sss us\sss consider\sss the quotient\sss 
$\mathbb{G}
\off =\off
\mathcal{F}\dff\ast\dff G_{\dff v}\dff/\dff K$\nnsp,\oss 
where $K$\sss is\dss normally\dss generated\dss 
by all\sss elements\dss $E\dff(\dff e\fff,\pff t\trf)$\dss 
and\sss elements\dss $L\dff(\dff l\qff)$ for\sss loops\sss $l\qff \in\qff \mathcal{L}$\dss
and\dss for all\sss edge-loops\sss 
$l\off =\off l_{\dff e}\dff,\off e\qff \in\qff E$\nnsp.\oss
Let\vspace{1.5pt}
\[
\quad
\rho\dff \colon\dff
\mathcal{F}\dff\ast\dff G_{\dff v}\qff \ttoo\qff \mathbb{G}
\]

\vspace{-12pt}\vspace{1.5pt}
be\sss the canonical\dss projection.\oss
Since our relations hold in $G$\nnsp,\oss there\dss is\dss
a unique homomorphism\dss 
$\varphi\dff \colon\dff \mathbb{G}\qff \ttoo\qff G$
such\dss that\sss $\psi\off =\off \varphi\dff \circ\dff \rho$\nnsp.\oss
Since\sss $\psi$\sss is\dss surjective by\dss Corollary\qss \ref{generators},\qss
$\varphi$\sss is\dss also surjective.\oss
Our goal\dss is\dss to prove\sss that\sss $\varphi$\sss is\dss an isomorphism\dss if\trs
the collection\sss $\mathcal{L}$\sss includes sufficiently\dss many\dss loops\qss 
(see\trs Theorem\qss \ref{simply-connected}\qss for\dss the precise meaning 
of\qss ``sufficiently\dss many''). 

For every\sss $e\qff \in\qff E$\sss let\sss 
$\mathfrak{g}_{\dff e}\off=\off \rho\trf(\trf g_{\dff e}\dff)\qff \in\qff \mathbb{G}$\nnsp.\oss 
Then 
$\varphi\dff(\trf \mathfrak{g}_{\dff e}\dff)
\off =\off 
\psi\dff(\trf g_{\dff e}\dff)
\off =\off 
s_{\dff e}$\nsp.\oss
Let
$\mathfrak{G}
\off =\off
\rho\trf(\trf \mathcal{G}\trf)$  
and\dss
$\mathbb{G}_{\dff v}
\off =\off 
\rho\dff(\trf G_{\dff v}\trf)\qff \subset\qff \mathbb{G}$\nnsp.\oss
Since\sss 
$\psi\off =\off \varphi\dff \circ\dff \rho$\sss is\dss equal\sss on\sss $G_{\dff v}$\dss 
to\sss the inclusion $G_{\dff v}\qff \ttoo\qff G$\nnsp,\oss 
the maps\dss 
$G_{\dff v}\qff \ttoo\qff \mathbb{G}_{\dff v}$\dss 
and\dss 
$\mathbb{G}_{\dff v}\qff \ttoo\qff G_{\dff v}$\dss 
induced by\sss $\rho$\sss and\sss $\varphi$\sss respectively\halfff,\oss 
are isomorphisms.\oss
Similarly\halfff,\oss the fact\sss that\sss
$\psi\dff(\trf g_{\dff e}\dff)
\off =\off 
s_{\dff e}$ 
implies\sss that\sss $\rho$\sss induces a\sss bijection\dss 
$\mathcal{G}\qff \ttoo\qff \mathfrak{G}$\dss 
and\dss the fact\sss that\sss
$\varphi\dff(\trf \mathfrak{g}_{\dff e}\dff)
\off =\off 
s_{\dff e}$
implies\sss that\sss $\varphi$\sss induces a\sss bijection\dss 
$\mathfrak{G}\qff \ttoo\qff \mathcal{S}$\nsp\dnsp.\oss

\mypar{Lemma.}{exchange} 
$\mathbb{G}_{\dff v}\dff \cdot\qff \mathfrak{G}
\off \subset\off 
\mathfrak{G}\dff \cdot\qff \mathbb{G}_{\dff v}$\nsp.\oss

\proof 
Let\sss $t\qff \in\qff G_{\dff v}$ and\sss $e\qff \in\qff E$\nnsp.\oss 
Then\dss the element\sss $E\dff(\dff e\fff,\pff t\trf)$\sss belongs\sss to $K$\nnsp,\oss 
and\dss hence\sss the corresponding\dss relation\sss holds in\sss $\mathbb{G}$\nnsp.\oss
More precisely\halfff,\oss\vspace{2pt}
\[
\quad
\rho\dff \left(\trf g_{\dff d}^{\dff -\dff 1}\dff \cdot\dff t\dff \cdot\dff g_{\dff e}\trf\right)
\off =\off
\rho\dff \left(\trf s_{\dff d}^{\dff -\dff 1}\dff \cdot\dff t\dff \cdot\dff s_{\dff e}\trf\right)
\off =\off
\rho\dff \bigl(\trf k\trf(\dff e\fff,\pff t\trf)\trf\bigr)
\off,
\]

\vspace{-9pt}
where\dss
$d\off =\off t\trf(\dff e\dff)$\nnsp.\oss
Let\trs
$\tau\off =\dff\off \rho\dff \left(\trf t\qff\right)$\dss
and\dss
$\kappa\off =\off \rho\dff \bigl(\trf k\trf(\dff e\fff,\pff t\trf)\trf\bigr)$\dss
and\dss rewrite\sss the\sss last\sss equality\sss as\vspace{3pt}
\[
\quad
\left(\trf \mathfrak{g}_{\dff d}\trf\right)^{\dff -\dff 1} \cdot\dff \tau\dff \cdot\dff \mathfrak{g}_{\dff e}
\off =\off
\kappa
\qff.
\]

\vspace{-12pt}\vspace{3pt}
Therefore\dss
$\tau\dff \cdot\dff \mathfrak{g}_{\dff e}
\off =\off
\mathfrak{g}_{\dff d}\dff \cdot\dff \kappa$\nnsp.\oss
It\sss follows\sss that\sss 
$\tau\dff \cdot\trf \mathfrak{G}
\off \subset\off
\mathfrak{G}\dff \cdot\qff \mathbb{G}_{\dff v}$\nsp.\oss
The\sss lemma\sss follows.\oss  \eproof

\myuppar{Remark.} 
The proof\dss of\qss Lemma\qss \ref{exchange}\qss used only\dss 
the fact\dss that\dss the element\trs 
$\left(\trf g_{\dff d}\trf\right)^{\dff -\dff 1} \cdot\dff t\dff \cdot\dff g_{\dff e}$\dss 
is\dss forced\dss by\dss the relation\dss $E\dff(\dff e\fff,\pff t\trf)$\dss 
to be equal\dss to an element\sss of\trs $G_{\dff v}$\nsp.\oss 
The fact\dss that\dss
it\dss is\dss forced\dss to be equal\dss to\sss the specific element\trs 
$s_{\dff d}^{\dff -\dff 1}\dff \cdot\dff t\dff \cdot\dff s_{\dff e}$\dss 
played\dss no role.\oss

\mysection{Kozsul\qss models\qss of\pss actions\qss transitive\qss on\qss vertices}{models}

\myuppar{Kozsul\trs models.}
The\qss \emph{Kozsul\dss model}\pss of\trs the graph\dss $X$\dss together\dss
with\dss the action of\trs $G$\dss on\dss $X$\dss
is\dss a\sss graph\dss $\mathbb{X}$\dss
together\sss with an action of\dss $\mathbb{G}$\dss on\dss $\mathbb{X}$\nnsp.\oss
The graph\dss $\mathbb{X}$\dss is\dss a\qss ``relaxed''\trs version of\dss $X$\dss
which accounts only\dss for\dss the relations used\dss to define\sss $\mathbb{G}$\nnsp.\oss
Naturally,\oss it\sss depends on\sss the choice of\dss elements $s_{\dff e}$\sss
and\dss the collection of\trs loops $\mathcal{L}$\nnsp.\oss
There\dss is\dss 
a canonical\sss $\mathbb{G}$\dnsp-equivariant\dss map\dss
$f\dff \colon\dff \mathbb{X}\qff \ttoo\qff X$\nnsp,\oss
where\dss $\mathbb{G}$\dss acts on\dss $X$\dss via\sss the homomorphism\dss $\varphi$\nnsp.\oss
We will\sss see\sss that\sss $f$\dss is\dss an isomorphisms of\dss graphs\dss if\trs and\dss only\trs if\trs
the system of\dss relations used\dss to define\dss $\mathbb{G}$\dss is\dss complete,\oss
i.e.\qss if\trs and\dss only\trs if\trs
$\varphi\dff \colon\dff
\mathbb{G}\qff \ttoo\qff G$\dss
is\dss an\sss isomorphism.\oss

\myuppar{The vertices of\dss $\mathbb{X}$\nnsp.}
For a\sss graph\dss $Z$\trs we will\sss denote\sss by\trs $Z_{\trf 0}$\dss 
be\sss the set\sss of\dss vertices of\trs $Z$\nnsp.\oss
We will\dss define\sss first\dss $\mathbb{X}_{\trf 0}$\dss
and\dss a canonical\dss map\dss
$f\dff \colon\dff
\mathbb{X}_{\trf 0}\qff \ttoo\qff X_{\trf 0}$\nsp.\oss
Then we will\dss define edges of\trs $\mathbb{X}$\nnsp,\oss
check\dss the correctness of\trs this definition,\oss
and check\dss that\dss $f$\dss maps edges\sss to edges.\oss

Since $G_{\dff v}$ is\dss the stabilizer of\sss $v$ and $G$ acts\sss 
transitively on\dss $X_{\trf 0}$\nsp,\oss we can identify\dss 
$G/G_{\dff v}$\dss with\dss $X_{\trf 0}$\nsp.\oss 
Taking\dss this as a clue,\oss let\trs 
$\mathbb{X}_{\trf 0}
\off =\off 
\mathbb{G}\fff/\fff \mathbb{G}_{\dff v}$\dss 
be\sss the set\sss of\dss vertices 
of\trs the future complex\qss (graph)\dss $\mathbb{X}$\nnsp.\oss
The group\dss $\mathbb{G}$\sss acts on\dss $\mathbb{X}_{\dff 0}$\dss 
in\dss the usual\dss manner\halfff.\oss
We will\sss denote\sss this action,\oss  
like\sss the action of\trs $G$\dss
on\dss $X_{\dff 0}$\nsp,\oss
by\dss $(\trf \gamma\fff,\pff y\trf)\qff \longmapsto\qff \gamma\trf(\dff y\trf)$\nnsp,\oss
where\dss $\gamma\qff \in\qff \mathbb{G}$\dss and\dss 
$y\qff \in\pff \mathbb{X}_{\dff 0}$\nsp.\oss 

Since $\varphi\dff(\trf \mathbb{G}_{\dff v}\dff)\pff =\off G_{\dff v}$\nsp,\oss 
the homomorphism 
$\varphi\dff \colon\dff \mathbb{G}\qff \ttoo\qff G$\dss 
induces a map 
$f\dff \colon\dff 
\mathbb{G}\fff/\fff \mathbb{G}_{\dff v}
\qff \ttoo\qff 
G/G_{\dff v}$\nnsp,\oss 
which we may\sss consider also as a map
$f\dff \colon\dff 
\mathbb{X}_{\dff 0}\qff \ttoo\qff X_{\trf 0}$\nnsp.\oss
The coset $\mathbb{G}_{\dff v}\qff \in\qff \mathbb{G}/\fff \mathbb{G}_{\dff v}$\nnsp,\oss
thought\sss as a vertex of\dss the future complex $\mathbb{X}$\nnsp,\oss
will\dss be denoted\dss by $v^{\fff *}$\nnsp.\oss
The vertex $v^{\fff *}$ will\dss play a role similar\dss to\sss that\sss of\trs the vertex $v$ of\trs $X$\nnsp.\oss
In\dss these terms\sss the map\sss $f$\sss 
is\dss given\sss by the formula\vspace{3pt}
\[
\quad 
f\dff(\trf \gamma\trf(\dff v^{\dff *}\trf)\trf)
\off =\off
\varphi\dff(\trf \gamma\trf)\trf(\dff v\trf)
\qff,
\] 

\vspace{-12pt}\vspace{3pt}
where\dss $\gamma\qff \in\qff \mathbb{G}$\nnsp.\pss 
In\dss particular\halfff,\oss this formula\sss leads\sss to a correctly defined map.\oss
We claim\sss that\sss $f$\dss is\dss \emph{$\mathbb{G}$\dnsp-equivariant}\pss 
with respect\dss to\sss the natural\sss action of\trs $\mathbb{G}$\sss on\dss 
$\mathbb{X}_{\dff 0}\off =\off \mathbb{G}\fff/\fff \mathbb{G}_{\dff v}$\dss 
and\dss the action of\dss $\mathbb{G}$ on $X_{\trf 0}$ via\dss $\varphi$\nnsp.\oss
Indeed,\oss every\sss 
$z\qff \in\qff \mathbb{X}_{\dff 0}$\sss has\sss the form\sss
$z\off =\off \beta\trf(\dff v^{\dff *}\trf)$\nnsp,\oss 
where\sss $\beta\qff \in\qff \mathbb{G}$\nnsp,\oss
and\vspace{4.5pt}
\[
\quad
f\dff(\trf \gamma\trf(\dff z\trf) \trf) 
\off =\off
f\dff(\trf \gamma\trf \beta\trf(\dff v^{\dff *}\trf) \trf)
\off =\off
\varphi\dff(\trf \gamma\trf \beta\dff)\trf(\dff v\trf)
\]

\vspace{-34.5pt}
\[
\quad
\phantom{f\dff(\trf \gamma\trf(\dff z\trf) \trf) 
\off =\off
f\dff(\trf \gamma\trf \beta\trf(\dff v^{\dff *}\trf) \trf)
\off }
=\off
\varphi\dff(\trf \gamma\trf)\qff \varphi\dff(\trf \beta\dff)\trf(\dff v\trf)
\off =\off
\varphi\dff(\trf \gamma \trf)\dff f\dff(\trf \beta\trf(\dff v^{\dff *}\trf)\trf)
\off =\off
\varphi\dff(\trf \gamma \trf)\dff(\trf f\dff(\trf z\trf)\trf)
\qff.
\]

\vspace{-12pt}\vspace{3pt}
The $\mathbb{G}$\dnsp-equivariance\sss follows.\oss
The surjectivity of\sss $\varphi$\sss implies\sss that\trs
$f\dff \colon\dff 
\mathbb{X}_{\trf 0}\qff \ttoo\qff X_{\trf 0}$\dss
is\dss surjective.

\myuppar{The edges of\dss $\mathbb{X}$\nnsp.}
In order\sss to define\sss the edges,\oss 
it\dss is\dss sufficient\dss to define for every\dss vertex 
$x$ a set\sss $N\trf(\dff x\trf)$ of\trs its\qss \emph{neighbors},\oss 
i.e.\qss of\dss vertices connected\dss with\sss $x$\sss by\sss an edge\qss
(and\dss hence $\neq\qff x$\nsp),\oss 
and\dss then check\dss that\dss the resulting\dss 
relation of\trs being\qss
\emph{neighbors}\qss is\dss symmetric\fff:\oss 
if\qss $y\qff \in\qff N\trf(\dff x\trf)$\nnsp,\oss
then\qss $x\qff \in\qff N\trf(\trf y\trf)$\nnsp.\oss
To begin\sss with,\oss we set\sss 
$N\trf(\dff v^{\fff *}\dff)
\off =\off
\mathfrak{G}\dff \cdot\fff v^{\fff *}$\dnsp.\oss
Note\sss that\sss\vspace{3pt}
\[
\quad
f\dff(\trf \mathfrak{g}_{\dff e}\dff \cdot\dff v^{\fff *}\dff)
\off =\off
\varphi\dff(\trf \mathfrak{g}_{\dff e}\trf)\dff (\dff v\trf)
\off =\off
s_{\dff e}\dff (\dff v\trf)
\]

\vspace{-12pt}\vspace{3pt}
and\dss hence $f$ induces a\sss bijection\sss
$N\trf(\dff v^{\fff *}\dff)
\qff \ttoo\qff
N\trf(\dff v\trf)$\nnsp.\oss
In\sss particular\halfff,\dss $v^{\fff *}\dff \not\in\pff N\trf(\dff v^{\fff *}\dff)$\nnsp.\oss
The definition of\sss $N\dff(\dff x\trf)$ for other vertices\sss 
$x\qff \in\qff \mathbb{X}_{\trf 0}$\sss is\dss dictated\dss by\sss
the need\dss to have a canonical\sss action of\dss $\mathbb{G}$ on\sss $\mathbb{X}$\nnsp.\oss
Namely,\oss every\dss vertex of\trs $\mathbb{X}$\sss has\sss the form
$\gamma\trf(\dff v^{\dff *}\trf)$\nnsp,\oss where\dss $\gamma\qff \in\qff \mathbb{G}$\nnsp,\qss
and we set\vspace{3pt}
\[
\quad
N\trf(\trf \gamma\trf(\dff v^{\dff *}\trf)\trf) 
\off =\dff\off
\gamma\dff \cdot\dff \mathfrak{G}\trf(\dff v^{\dff *}\trf)
\qff.
\]

\vspace{-12pt}\vspace{3pt}
If\dss this definition\dss is\dss correct,\oss
i.e.\qss does not\sss depend on\sss the choice of\sss $\gamma$
such\dss that\sss
$x\off =\off \gamma\trf(\trf v^{\fff *}\trf)$\nnsp,\oss
as we will\sss see in a moment,\oss
then\sss the relation 
$y\qff \in\qff N\trf(\trf x\trf)$
is\dss invariant\dss under\dss the action of\sss $\mathbb{G}$\nnsp.\oss

\vspace{-12pt}\vspace{3pt}\vspace{-1pt}
\mypar{Lemma.}{correctness-of-N}
\emph{The definition\sss of\pss $N\trf(\dff x\trf)$\dss is\dss correct\halfff,\oss
i.e.\qss $N\trf(\dff x\trf)$\dss does not\sss depend on\dss the 
choice of\qss $\gamma$\qss such\dss that\qss
$x\off =\off \gamma\trf(\trf v^{\fff *}\trf)$\nnsp.\oss}

\proof
If\trs
$\gamma_{\fff 1}\trf (\dff v^{\fff *}\trf)
\off =\off
\gamma_{\dff 2}\trf (\dff v^{\fff *}\trf)$\nnsp,\oss
then\dss
$\gamma_{\dff 2}^{\dff -\dff 1}\qff \gamma_{\fff 1}\trf (\dff v^{\fff *}\trf)
\off =\off
v^{\fff *}$\dss
and\dss therefore\dss
$\delta
\off =\off 
\gamma_{\dff 2}^{\dff -\dff 1}\dff \gamma_{\fff 1}
\qff \in\qff
\mathbb{G}_{\dff v}$\nsp.\oss
Let\qss $\mathfrak{g}\qff \in\qff \mathfrak{G}$\nnsp.\oss
Lemma\qss \ref{exchange}\qss implies\sss that\trs
$\delta\qff \mathfrak{g}\off =\off \mathfrak{f}\qff \varepsilon$\dss
for some\dss
$\beta\qff \in\qff \mathfrak{G}$\dss
and\dss
$\varepsilon\qff \in\qff \mathbb{G}_{\dff v}$\nsp.\oss
Therefore\vspace{4.5pt}
\[
\quad
\gamma_{\fff 1}\qff \mathfrak{g}\trf(\dff v^{\fff *}\trf)
\off =\off
\gamma_{\dff 2}\qff \delta\qff \mathfrak{g}\trf(\dff v^{\fff *}\trf)
\off =\off
\gamma_{\dff 2}\qff \mathfrak{f}\qff \varepsilon\trf(\dff v^{\fff *}\trf)
\off =\off
\gamma_{\dff 2}\qff \mathfrak{f}\trf(\dff v^{\fff *}\trf)
\qff,
\]

\vspace{-12pt}\vspace{4.5pt}
where we used\dss the fact\dss that\trs
$\varepsilon\qff \in\qff \mathbb{G}_{\dff v}$\dss  
and\dss hence\dss
$\varepsilon\trf(\dff v^{\fff *}\trf)\off =\off v^{\fff *}$\nnsp.\oss
If\dss follows\dss that\vspace{4.5pt}
\[
\quad
\gamma_{\dff 1}\dff \cdot\qff \mathfrak{G}\trf(\dff v^{\fff *}\trf)
\off \subset\off
\gamma_{\dff 2}\dff \cdot\qff \mathfrak{G}\trf(\dff v^{\fff *}\trf)
\qff.
\]

\vspace{-12pt}\vspace{4.5pt}
By\sss interchanging $\gamma_{\dff 1}$
and $\gamma_{\dff 2}$ we see\sss that\dss the opposite inclusion
also\sss true and\dss hence\sss these\sss two sets are equal.\oss
The\sss lemma\sss follows.\oss  \eproof

\mypar{Lemma.}{symmetry} 
\emph{Let\qss $x\fff,\pff y\qff \in\qff \mathbb{X}_{\dff 0}$\nsp.\oss 
If\pss $x\qff \in\qff N\trf(\dff y\trf)$\nnsp,\oss
then\qss $y\qff \in\qff N\trf(\trf x\trf)$\nnsp.\oss}

\proof
Let\trs us choose\dss $\gamma\qff \in\qff \mathbb{G}$\dss
such\dss that\sss
$x\off =\off \gamma\trf(\dff v^{\fff *}\trf)$\nnsp.\oss 
Since\sss $y\qff \in\qff N\trf(\dff x\trf)$\nnsp,\oss
there exists an edge\sss $e\qff \in\qff E$ such\dss that\trs
$y\off =\off \gamma\qff \cdot\qff \mathfrak{g}_{\dff e} \qff(\dff v^{\fff *}\trf)$\nnsp.\oss
Let\trs $l$\dss be\sss the edge-loop\dss
$v\fff,\pff s_{\dff e} \trf(\dff v\trf)\fff,\pff v$\nnsp.\oss
The corresponding\dss relation\dss $L\dff(\dff l\qff)$\dss 
has\sss the form\sss
$g_{\dff e}\dff \cdot\dff g_{\dff d}
\off =\off
s_{\dff e}\dff \cdot\dff s_{\dff d}$\nsp,\oss
where\dss $d\qff \in\qff E$\dss and\dss
$s_{\dff e}\dff \cdot\dff s_{\dff d}\qff \in\qff G_{\dff v}$\nsp.\oss
By\sss applying\dss $\rho$\dss to\sss this relation,\oss 
we see\sss that\dss 
$\mathfrak{g}_{\dff e}\dff \cdot\dff \mathfrak{g}_{\dff d}
\off =\off
\rho\trf(\dff s_{\dff e}\dff \cdot\dff s_{\dff d}\trf)
\pff \in\off
\mathbb{G}_{\dff v}$\nsp.\oss
Since\dss $\mathbb{G}_{\dff v}$\dss fixes\sss $v^{\dff *}$\sss 
under\dss the action of\trs $\mathbb{G}$\sss on\dss 
$\mathbb{X}_{\dff 0}\off =\off \mathbb{G}/\mathbb{G}_{\dff v}$\nnsp,\oss 
this implies\sss that\vspace{2.5pt}
\[
\quad
x
\off =\off 
\gamma\trf(\dff v^{\fff *}\trf)
\off =\off 
\gamma\dff \cdot\trf
\rho\trf\left(\dff s_{\dff e}\dff \cdot\dff s_{\dff d}\trf\right)\trf(\dff v^{\fff *}\trf)
\off =\off
\gamma\trf \cdot\trf \mathfrak{g}_{\dff e}\dff \cdot\dff \mathfrak{g}_{\dff d}\trf(\dff v^{\fff *}\trf).
\]

\vspace{-12pt}\vspace{2.5pt}
Since\dss
$y\off =\off \gamma\trf \cdot\trf \mathfrak{g}_{\dff e} \qff(\dff v^{\fff *}\trf)$\dss
and\qss
$\mathfrak{g}_{\dff d}\qff \in\qff \mathbb{G}_{\dff E}$\nsp,\oss
this means\sss that\trs
$x\qff \in\qff N\trf(\trf y\trf)$\nnsp.\oss
The lemma follows.\oss  \eproof

\myuppar{Remark.}
The above proof\dss depends only\sss on\dss the edge-loop relations,\oss
and\dss is\dss the only\dss place in\dss the construction of\trs $\mathbb{X}$\dss
and\sss in\dss the proofs of\dss its properties where\sss these relations are used.\oss

\myuppar{The\sss graph\sss $\mathbb{X}$ and\sss the map\sss
$f\dff \colon\dff \mathbb{X}\qff \ttoo\qff X$\nnsp.}
By\qss Lemmas\qss \ref{correctness-of-N}\qss and\qss \ref{symmetry}\qss
the relation\dss $y\qff \in\qff N\trf(\trf x\trf)$\dss
is\dss correctly\sss defined and symmetric.\oss
Therefore we can define\sss $\mathbb{X}$\sss as\sss the\sss graph\sss
having\sss $\mathbb{X}_{\trf 0}$\sss as\sss the set\sss of\dss vertices
and\dss with edges connecting\dss pairs of\dss vertices\dss
$x\fff,\pff y$\dss such\dss that\trs
$y\qff \in\qff N\trf(\trf x\trf)$\nnsp.\oss
By\dss the very\sss definition\dss the relation 
$y\qff \in\qff N\trf(\trf x\trf)$
is\dss invariant\dss under\dss the action of\sss $\mathbb{G}$\nnsp.\oss
Hence\sss $\mathbb{G}$\sss canonically\sss acts on\dss $\mathbb{X}$\dnsp.
Also,\oss by\dss the definition of\trs the neighbors
the map\dss
$f\dff \colon\dff 
\mathbb{X}_{\trf 0}\qff \ttoo\qff X_{\trf 0}$\dss 
maps\sss neighbors of\dss $v^{\fff *}$\dss in\dss $\mathbb{X}$\dss 
to\sss neighbors of\dss $v$\dss in\dss $X$\nnsp.\oss
In view of\trs $\mathbb{G}$\dnsp-equivariance of\dss $f$\dss
this implies\sss that\dss $f$\dss is\dss a\sss map of\dss graphs\dss
$\mathbb{X}\qff \ttoo\qff X$\nnsp,\oss
i.e.\qss maps edges\sss to edges.\oss
Since\dss $X$\dss is\dss a\sss graph\sss without\dss loops,\qss
$\mathbb{X}$\dss is\dss also a\sss graph\sss without\dss loops.\oss

\mypar{Lemma.}{local-isomorphism} 
\emph{The map\qss 
$f\dff \colon\dff 
\mathbb{X}\qff \ttoo\qff X$\qss is\dss a\sss local\dss isomorphism of\trs graphs,\oss 
i.e.\qss for every\dss vertex\dss
$z$\dss of\pss $\mathbb{X}$\dss the map\dss $f$\dss maps\sss the set\dss of\qss edges of\pss 
$\mathbb{X}$\dss having\dss $z$\sss as\sss an\sss endpoint\dss bijectively\sss 
onto\sss the set\sss of\qss edges of\pss $X$\dss having\dss
$f\dff (\dff z\trf)$ as\sss an\sss endpoint\halfff.\oss}

\proof
In\sss terms of\dss neighbors\sss this means\sss that\sss $f$ 
induces a\sss bijection\sss
$N\trf(\trf z\trf)
\qff \ttoo\qff
N\trf(\trf f\trf(\trf z\trf)\trf)$
for every\sss $z\qff \in\qff \mathbb{X}_{\trf 0}$\nsp.\oss
By\sss the definition of\sss 
the neighbors in $\mathbb{X}$\sss
this\dss is\dss the case for\sss $z\off =\off v^{\fff *}$\dnsp.\oss
Now\sss the\sss $\mathbb{G}$\dnsp-equivariance of\sss $f$\sss
implies\sss the general\sss case.\oss  \eproof\vspace{-0.5pt}

\mypar{Lemma.}{connectedness} 
\emph{The graph\sss $\mathbb{X}$\sss is\dss connected.\oss}\vspace{-0.5pt}

\proof
It\dss is\dss sufficient\dss to prove\sss that\sss for every\sss
$\gamma\qff \in\qff \mathbb{G}$\sss the vertex\sss 
$\gamma\trf(\dff v^{\fff *}\trf)$\dss 
is\dss connected\sss with\sss $v$\sss by\sss a path.\oss
By\dss the very\sss definition of\sss $\mathbb{G}$
every $\gamma\qff \in\qff \mathbb{G}$\sss
can\sss be presented
as a product\sss of\dss several\sss elements 
$\mathfrak{g}_{\dff e}$ and several\sss elements of\trs $G_{\dff v}$\nsp.\oss
We can\dss write\trs $\gamma$\dss in\dss the form\dss
$\gamma\off =\off g\qff \beta$\nnsp,\oss
where\dss $\beta$\dss is\dss a\sss shorter\sss product\qss 
({\fff}it\dss may\dss be even empty\fff)\qss
and either\dss $g\off =\off \mathfrak{g}_{\dff e}$\qss for some\dss
$e\qff \in\qff E$\nnsp,\oss
or\trs $g\qff \in\qff \mathbb{G}_{\dff v}$\nsp.\oss
Using an\sss induction\dss by\dss the\sss length of\trs the product\halfff,\oss
we can assume\sss that\sss $\beta\trf(\dff v^{\fff *}\trf)$\sss
is\dss connected\sss with $v^{\fff *}$ by\sss a path.\oss
It\dss follows\dss that\trs
$\gamma\trf(\dff v^{\fff *}\trf)
\off =\off
g\qff \beta\trf(\dff v^{\fff *}\trf)$\dss
is\dss connected\sss with\sss $g\trf(\dff v^{\fff *}\trf)$\sss by\sss a path.\oss
If\trs $g\off =\off \mathfrak{g}_{\dff e}$\nsp,\oss 
then\dss $g\trf(\dff v^{\fff *}\trf)$\dss is\dss a\sss neighbor of\dss $v^{\fff *}$\sss
and\dss hence\dss is\dss connected\sss with $v^{\fff *}$ by an edge.\oss
By\sss concatenating\dss the path connecting\dss $\gamma\trf(\dff v^{\fff *}\trf)$\dss
with\sss $g\trf(\dff v^{\fff *}\trf)$\dss and\sss this edge,\oss 
we\sss get\sss a path connecting\dss $\gamma\trf(\dff v^{\fff *}\trf)$\dss with\sss $v^{\fff *}$\nnsp.\oss 
If $g\qff \in\qff \mathbb{G}_{\dff v}$\nsp,\oss 
then\dss 
$g\trf(\dff v^{\fff *}\trf)
\off =\off 
v^{\fff *}$\nnsp,\oss and hence\sss
the original\dss path already\sss connects\dss 
$\gamma\trf(\dff v^{\fff *}\trf)$\dss with\dss $v^{\fff *}$\nnsp.\oss 
This completes\sss the proof\halfff.\oss  \eproof

\myuppar{The map\sss $f$\sss as a covering\halfff.}
The surjectivity of\dss
$f\dff \colon\dff 
\mathbb{X}_{\trf 0}\qff \ttoo\qff X_{\trf 0}$\sss
together\dss with\dss Lemma\qss \ref{local-isomorphism}\qss
means\sss that\trs
$f\dff \colon\dff 
\mathbb{X}\qff \ttoo\qff X$\dss
is\dss an analogue of\dss covering\sss maps for\sss graphs.\oss
In\sss particular\halfff,\oss the geometric realization\dss
$\num{f\fff}\dff \colon\dff 
\num{\mathbb{X}}\qff \ttoo\qff \num{X}$\dss
is\dss a\sss covering of\trs topological\sss spaces.\oss
Also,\oss Lemma\qss \ref{local-isomorphism}\qss immediately\sss
implies\sss the following\dss paths\dss lifting\dss property.\oss
Suppose\sss that\sss 
$v_{\dff 0}\dff,\off v_{\dff 1}\dff,\off\ldots\dff,\off v_{\dff n}$\dss 
is\dss a\sss path\sss
in\sss $X$ and\sss that\sss $u_{\qff 0}$\sss is\dss
a vertex of\dss $\mathbb{X}$ such\dss that\sss
$f\dff(\trf u_{\qff 0}\trf)\off =\off v_{\dff 0}$\nsp.\oss
Then\sss there\dss is\dss a\sss unique path\sss in $\mathbb{X}$
of\qss the\dss form\sss
$u_{\dff 0}\dff,\off u_{\dff 1}\dff,\off\ldots\dff,\off u_{\dff n}$
such\dss that\sss 
$f\dff(\trf u_{\dff i}\trf)\off =\off v_{\dff i}$\sss 
for every\sss $i$\nnsp.\oss

\mypar{Lemma.}{lifted-loops}
\emph{Let\sss $l$ be\sss a\sss path\dss
$v_{\dff 0}\dff,\off v_{\dff 1}\dff,\off\ldots\dff,\off v_{\dff n}$\dss
in\qss $X$\nnsp.\oss 
Suppose\sss that\sss $l$\sss is\qss a\dss loop based at\sss $v$\nnsp,\oss
i.e.\qss that\trs $v_{\dff 0}\off =\off v_{\dff n}\off =\off v$\dnsp.\oss
If\pss $l\qff \in\qff \mathcal{L}$\nnsp,\oss 
then\dss the\dss lift\dss of\trs $l$\sss 
starting\dss at\sss $v^{\dff *}$
ends\dss at\sss $v^{\dff *}$\nsp.\oss}

\proof 
In\dss fact\halfff,\oss the construction of\trs $\mathbb{X}$\dss
was\dss motivated\dss by\dss the desire\sss to have\sss this property\halfff,\oss
and\dss the proof\dss amounts\sss to verification\dss that\sss definitions
work as intended.\oss
Let\trs $l$\dss be\sss the loop\qss
$v_{\dff 0}\dff,\off v_{\dff 1}\dff,\off\ldots\dff,\off v_{\dff n}$\nsp.\oss 
The corresponding\dss loop relation\dss $L\dff(\dff l\qff)$\dss is\vspace{3pt}
\[
\quad
g_{\dff 1}
\dff \cdot\dff 
g_{\dff 2}
\dff \cdot\dff 
\ldots
\dff \cdot\dff 
g_{\dff n}
\off =\off
s_{\dff 1}\dff \cdot\dff s_{\dff 2}
\dff \cdot\dff 
\ldots
\dff \cdot\dff 
s_{\dff n}
\qff,
\]

\vspace{-12pt}\vspace{3pt}
where\sss the elements\dss $s_{\dff i}\qff \in\qff \mathcal{S}$\dss
are determined\dss by\qss (\ref{loop}),\oss
the right\dss hand side\dss is\dss considered as an element\sss of\trs $G_{\dff v}$\nsp,\oss 
and each\dss $g_{\dff i}$\dss is\dss the generator of\trs $\mathcal{F}$\dss corresponding\dss
to $s_{\dff i}$\nsp,\dss 
$0\qff \leq\qff i\qff \leq\qff n$\nnsp.\oss 
By applying\dss 
$\rho\dff \colon\dff 
\mathcal{F}\dff\ast\dff G_{\dff v}
\qff \ttoo\qff 
\mathbb{G}$\dss to\sss this relation we see\sss that\vspace{3pt}
\[
\quad
\mathfrak{g}_{\dff 1}
\dff \cdot\dff 
\mathfrak{g}_{\dff 2}
\dff \cdot\dff 
\ldots
\dff \cdot\dff 
\mathfrak{g}_{\dff n}
\off =\off
s_{\dff 1}\dff \cdot\dff s_{\dff 2}
\dff \cdot\dff 
\ldots
\dff \cdot\dff 
s_{\dff n}
\qff,
\]

\vspace{-12pt}\vspace{3pt}
where\dss 
$\mathfrak{g}_{\dff i}
\off =\off
\rho\dff(\trf g_{\dff i}\trf)
\qff \in\qff 
\mathfrak{G}$\dss 
and\dss the right\dss hand side\dss is\dss considered 
as an element\sss of\trs $\mathbb{G}_{\dff v}$\nsp.\oss 
Since\sss the subgroup\dss $\mathbb{G}_{\dff v}$\dss is\dss the stabilizer of\trs 
$v^{\dff *}$\dss in\dss $\mathbb{X}$\trs by\dss the construction of\trs $\mathbb{X}$\nnsp,\oss 
this implies that\vspace{3pt}
\begin{equation}
\label{fix}
\quad
\mathfrak{g}_{\dff 1}
\dff \cdot\dff 
\mathfrak{g}_{\dff 2}
\dff \cdot\dff 
\ldots
\dff \cdot\dff 
\mathfrak{g}_{\dff n}\trf
(\trf v^{\dff *}\trf)
\off =\off
v^{\dff *}
\qff.
\end{equation} 

\vspace{-12pt}\vspace{3pt}
We claim that the sequence $u_{\dff 0}\dff,\qff \ldots\dff,\qff u_{\dff n}$ of the vertices\vspace{3pt}
\begin{equation}
\label{loop1}
\hspace{0.5em}\begin{array}{l}
u_{\dff 0}
\off =\off
v^{\dff *}\dff,\vspace{12pt}\\
u_{\dff 1}
\off =\off
\mathfrak{g}_{\dff 1}\dff(\dff v^{\dff *}\trf)\qff,\vspace{12pt}\\
u_{\dff 2}
\off =\off
\mathfrak{g}_{\dff 1}\dff \cdot\dff \mathfrak{g}_{\dff 2}\trf(\dff v^{\dff *}\trf)\qff,\vspace{7.5pt}\\
\ldots\ldots\vspace{7.5pt}\\
u_{\dff n}
\off =\off
\mathfrak{g}_{\dff 1}\dff \cdot\dff \mathfrak{g}_{\dff 2}\dff \cdot\dff
\ldots\dff \cdot\dff 
\mathfrak{g}_{\dff n}\trf(\dff v^{\dff *}\trf)
\qff
\end{array}
\end{equation}

\vspace{-7.5pt}
is\dss a path\dss in\dss $\mathbb{X}$ starting\sss at\dss $v^{\dff *}$
and\dss lifting\dss $l$\nnsp.\qff\oss
Let\trs $0\qff \leq\qff i\qff \leq\qff n\qff -\qff 1$\nnsp.\oss 
Since $\mathfrak{g}_{\dff i\dff +\dff 1}\qff \in\qff \mathfrak{G}$\nnsp,\oss we have\vspace{4.5pt} 
\[
\quad
u_{\dff i\dff +\dff 1}
\off =\off 
\mathfrak{g}_{\dff 1}\dff \cdot\dff \mathfrak{g}_{\dff 2}\dff \cdot\dff
\ldots\dff \cdot\dff 
\mathfrak{g}_{\dff i\dff +\dff 1}\trf(\dff v^{\dff *}\trf)
\off \in\off  
\mathfrak{g}_{\dff 1}\dff \cdot\dff \mathfrak{g}_{\dff 2}\dff \cdot\dff
\ldots\dff \cdot\dff 
\mathfrak{g}_{\dff i}\dff \cdot\qff
\mathfrak{G}\trf(\trf v^{\dff *}\trf)
\qff,
\]

\vspace{-7.5pt}
and\dss hence\sss $u_{\dff i\dff +\dff 1}$\sss is\dss a 
neighbor of\qss 
$\mathfrak{g}_{\dff 1}\dff \cdot\dff \mathfrak{g}_{\dff 2}\dff \cdot\dff
\ldots\dff \cdot\dff 
\mathfrak{g}_{\dff i}\trf(\dff v^{\dff *}\trf)
\off =\off
u_{\dff i}$\nnsp.\oss 
Therefore,\pss $u$\dss is\dss indeed a path.\oss 
Since\dss 
$u_{\dff 0}
\off =\off
v^{\dff *}$\nnsp,\oss this path starts at\trs $v^{\dff *}$\nnsp.\oss 
Next,\vspace{4.5pt}
\[
\quad
f\dff(\dff u_{\dff i}\trf)
\off =\off
f\dff \left(\dff 
\mathfrak{g}_{\dff 1}\dff \cdot\dff \mathfrak{g}_{\dff 2}\dff \cdot\dff
\ldots\dff \cdot\dff 
\mathfrak{g}_{\dff i}
\trf(\dff v^{\dff *}\trf)
\dff\right) 
\off =\off
\varphi\dff \left(\dff 
\mathfrak{g}_{\dff 1}\dff \cdot\dff \mathfrak{g}_{\dff 2}\dff \cdot\dff
\ldots\dff \cdot\dff 
\mathfrak{g}_{\dff i}
\qff\right)
\trf(\dff v\trf)
\]

\vspace{-33pt}
\[
\quad
\phantom{f\dff(\dff u_{\dff i}\trf)
\off }
=\off
\varphi\dff \left(\dff 
\mathfrak{g}_{\dff 1}\dff\right)
\dff \cdot\qff 
\varphi\dff \left(\dff \mathfrak{g}_{\dff 2}\dff\right)
\dff \cdot\dff
\ldots\dff \cdot\qff 
\varphi\dff \left(\dff \mathfrak{g}_{\dff i}
\qff\right)
\trf(\dff v\trf) 
\off =\off
s_{\dff 1}\dff \cdot\dff s_{\dff 2}
\dff \cdot\dff 
\ldots
\dff \cdot\dff 
s_{\dff i} 
\trf(\dff v\trf)
\off =\off
v_{\dff i}
\qff,
\]

\vspace{-7.5pt}
and\dss hence\sss the path\sss $u$\dss is\dss indeed\sss a\dss lift\sss of\trs $l$\nnsp.\oss  
By \qss(\ref{fix})\qss the endpoint\dss $u_{\dff n}$\dss is\dss equal\dss to\sss $v^{\dff *}$\dnsp.\oss 
Hence $u$\dss is\dss a\dss lift\sss of\trs $l$\dss and\dss is\dss a closed\sss path.\oss
Lemma\qss \ref{local-isomorphism}\qss implies\sss that\dss there\dss is\dss only\sss one\dss
lift\halfff.\oss  \eproof

\mypar{Theorem.}{simply-connected}
\emph{Let\pss $\num{X}^{\dff +}$ be\sss the result\sss of\qss glueing\dss 
$2$\dnsp-cells\sss to\sss the geometric realization\dss $\num{X}$\dss of\trs $X$\dss 
along\dss the geometric realization of\trs all\trs loops\sss of\qss
the form\dss $g\dff(\trf l\qff)$\dss with\dss 
$g\qff \in\qff G$\dss and\dss 
$l\qff \in\qff \mathcal{L}$\nnsp.\oss 
If\pss $\num{X}^{\dff +}$\dss is\dss simply-connected,\oss 
then every\trs lift\trs of\trs a closed\dss path\dss in\dss $X$\dss starting\sss
at\sss $v$\sss to a path\dss in\qss $\mathbb{X}$\dss starting\sss at\dss $v^{\dff *}$\sss
is\dss a closed\dss path.\oss}

\proof 
For a path $q$ in $X$ or $\mathbb{X}$ we denote by 
$\num{q}$ the corresponding\dss geometric path in\dss $\num{X}$\dss 
or\dss $\num{\mathbb{X}}$\dss respectively\halfff.\oss 
Clearly\halfff,\pss $\num{q}$\dss is\dss a closed\dss path\dss 
if\trs and\dss only\trs if\trs $q$\dss is.\oss 
If\sss $r$\sss is\dss the path\sss
$u_{\dff 0}\dff,\off 
u_{\dff 1}\dff,\off\ldots\dff,\off 
u_{\dff n}$\nsp,\oss
we denote by\sss $r^{\dff -\dff 1}$\sss the path\sss
$u_{\dff n}\dff,\off u_{\dff n\dff -\dff 1}\dff,\off\ldots\dff,\off u_{\dff 0}$\nsp.\oss
The fact\sss that\sss $r^{\dff -\dff 1}$\sss is\dss indeed a path\sss implicitly depends on\dss
Lemma\qss \ref{symmetry},\oss and\dss hence on edge-loop relations.

Let\dss $p$\dss be a closed\dss path in\dss $X$\dss starting at\dss $v$\nnsp.\oss 
By\dss the elementary\dss theory of\trs CW-complexes\sss 
the assumption\sss that\sss $\num{X}^{\dff +}$ is\dss simply-connected\dss 
implies that\sss 
$\num{p}$\sss is\dss homotopic in $\num{X}$ 
to a product\sss of\dss paths of\trs the form\dss
$\num{r\qff g\dff(\trf l\qff)\qff r^{\dff -\dff 1}}$\nsp,\oss 
where $l\qff \in\qff \mathcal{L}$\nnsp,\qss 
$g\qff \in\qff G$\nnsp,\oss 
and $r$\sss is\dss a path connecting $v$ with\dss the starting\sss vertex 
$g\dff(\dff v\trf)$ of\sss $g\dff(\trf l\trf)$\nnsp.\oss
Therefore\sss it\dss is\dss sufficient\dss to prove\sss that\dss lift\sss of\dss
every\sss such\sss path
$r\qff g\dff(\trf l\qff)\qff r^{\dff -\dff 1}$\sss 
is\dss closed.\oss
This\dss is\dss the only argument\sss depending on algebraic\sss topology.

Let\sss $r^{\fff \sim}$\sss be\sss the\sss lift\sss of\sss $r$ starting at\sss $v^{\fff *}$\dnsp.\oss
The surjectivity of\dss
$f\dff \colon\dff 
\mathbb{X}_{\trf 0}\qff \ttoo\qff X_{\trf 0}$\sss 
implies\sss that\sss the endpoint\sss of\sss $r^{\fff \sim}$\sss
is\dss equal\dss to\sss $\beta\trf(\dff v^{\fff *}\trf)$
for some\sss $\beta\qff \in\qff \mathbb{G}$\nnsp.\oss
At\sss the same\sss time $f$\sss takes\sss this vertex\sss to
$g\dff(\dff v\trf)$\nnsp,\oss
and\dss hence
$g\dff(\dff v\trf)
\off =\off
f\dff(\trf \beta\dff(\trf v^{\dff *}\trf)\trf)
\off =\dff\off
\varphi\dff(\trf \beta\trf)\dff(\dff v\trf)$\nnsp.\oss
It\sss follows\sss that
$g^{\dff -\dff 1}\dff \varphi\dff(\trf \beta\trf)$\sss
fixes $v$\nnsp,\oss i.e.\dss 
$g^{\dff -\dff 1}\dff \varphi\dff(\trf \beta\trf)
\qff \in\qff
G_{\dff v}$\nsp.\oss
Since\sss $\varphi$\sss induces an isomorphism\dss 
$\mathbb{G}_{\dff v}\qff \ttoo\qff G_{\dff v}$\nsp,\oss
there exists\sss $h\qff \in\qff \mathbb{G}_{\dff v}$
such\dss that\sss
$g^{\dff -\dff 1}\dff \varphi\dff(\qff \beta\trf)
\off =\off
\varphi\dff(\trf h\trf)$
and\dss hence\sss
$g
\off =\off 
\varphi\dff(\qff \beta\dff \cdot\dff h^{\dff -\dff 1}\dff)$\nnsp.\oss 
Let\sss 
$\gamma\off =\off \beta\dff \cdot\dff h^{\dff -\dff 1}$\nnsp.\oss
Then\vspace{4.5pt}
\[
\quad
\varphi\dff(\trf \gamma\trf)\off =\off g
\quad
\mbox{and}\quad
\gamma\trf (\trf v^{\dff *}\trf)
\off =\off 
\beta\dff \cdot\dff h^{\dff -\dff 1}\dff (\trf v^{\dff *}\trf)
\off =\off
\beta\dff (\trf v^{\dff *}\trf)
\qff,
\]

\vspace{-12pt}\vspace{4.5pt}
because\dss $h\qff \in\qff \mathbb{G}_{\dff v}$\nsp.\oss
Therefore\sss the endpoint\sss of\dss the\sss lift\sss $r^{\sim}$\dss 
has\sss the form\dss
$\gamma\trf (\trf v^{\dff *}\trf)$\dss for some\dss $\gamma\qff \in\qff\mathbb{G}$\dss 
such\dss that\trs 
$\varphi\dff(\trf \gamma\trf)\off =\off g$\nnsp.\oss
Let\sss $l^{\dff \sim}$\sss be\sss the\sss lift\sss of\sss $l$ starting at\sss 
$v^{\fff *}$\dnsp\dnsp.\oss
By\trs Lemma\qss \ref{lifted-loops}\qss the path $l^{\dff \sim}$ ends also at\sss 
$v^{\fff *}$\dnsp\dnsp,\oss i.e.\qss is\dss a\sss loop.\oss
The equivariance of\sss $f$ implies\sss that\sss 
$\gamma\trf(\trf l^{\dff \sim}\trf)$\sss
is\dss the\sss lift\sss of\sss $g\dff(\trf l\qff)$
starting and ending at\sss $\gamma\trf(\dff v^{\fff *}\trf)$\nnsp.\oss
The\sss lift\sss of\trs 
$r\trf g\dff(\trf l\qff)\qff r^{\dff -\dff 1}$\dss 
starting\sss at\dss $v^{\dff *}$\dss can\sss be obtained\dss by\dss following\dss 
first\dss the\sss lift\dss $r^{\sim}$\dss of\dss $r$\dss 
to its endpoint\dss $\gamma\trf (\trf v^{\dff *}\trf)$\nnsp,\oss 
then\dss tracing\dss the closed\dss path\dss $\gamma\trf (\trf l^{\dff\sim}\trf)$\dss
which\dss lifts $g\dff(\dff l\trf)$\nnsp,\oss 
and\dss finally\dss retracing\dss $r^{\sim}$\dss back\halfff,\oss  
i.e.\qss following\dss $(\dff r^{\sim}\dff)^{\dff -\dff 1}$\nnsp.\oss 
Clearly\halfff,\oss this path returns\dss to\sss 
$v^{\dff *}$\dnsp,\oss
i.e.\qss the\sss lift\sss of\sss
$r\trf g\dff(\trf l\qff)\qff r^{\dff -\dff 1}$\sss
is\dss closed.\oss
The\sss lemma\sss follows.\oss  \eproof

\mypar{Corollary\halfff.}{isomorphism-complexes} 
\emph{Under\dss the assumptions of\trs the\dss theorem\dss the map\dss 
$f\dff \colon\dff \mathbb{X}\qff \ttoo\qff X$\dss is\dss an isomorphism of\dss graphs
and\dss the\sss homomorphism\dss
$\varphi\dff \colon\dff \mathbb{G}\qff \ttoo\qff G$\dss
is\dss an\dss isomorphism\sss of\qss groups.\oss}

\proof
The map $f$\sss is\dss surjective on vertices.\oss
By\trs Lemma\qss \ref{local-isomorphism}\qss it\dss is\dss a\sss local\dss isomorphism.\oss
Therefore,\oss if\dss it\dss is\dss not\sss an\sss isomorphism,\oss
then\sss it\dss is\dss not\sss injective on\sss vertices.\oss
Since $\mathbb{X}$\sss is\dss connected\dss by\trs Lemma\qss \ref{connectedness},\oss
in\sss this case\sss there exists a non-closed\sss path $p$\sss in $\mathbb{X}$
such\sss that\sss $f\dff(\trf p\trf)$\sss is\dss closed.\oss 
Replacing,\oss if\dss necessary,\dss $p$\sss by\sss the path $g\trf(\dff p\trf)$\sss
for some $g\qff \in\qff \mathbb{G}$\nnsp,\oss
we may assume\sss that\sss $p$ starts at\sss $v^{\fff *}$\dnsp\dnsp.\oss
But\sss then $p$ should\dss be closed\dss by\trs Theorem\qss \ref{simply-connected}.\oss
The contradiction proves\sss that $f$ is\dss an\sss isomorphism.\oss 
In\sss turn,\oss this implies\sss that\sss $\varphi$ induces a\sss
bijection\sss 
$\mathbb{G}/\mathbb{G}_{\dff v}
\qff \ttoo\qff
G/G_{\dff v}$\nsp.\oss
Since $\varphi$ induces an\sss isomorphism
$\mathbb{G}_{\dff v}\qff \ttoo\qff G_{\dff v}$\nsp,\oss
it\sss follows\sss that\sss $\varphi$\sss is\dss an\sss isomorphism.\oss  \eproof\vspace{-0.125pt}

\mysection{Implications\qss between\qss edge\qss
and\qss edge-loop\qss relations}{implications}

\myuppar{Another\sss form of\dss edge relations.}
The edge relation\dss $E\dff(\dff e\fff,\pff t\trf)$\dss can\dss be rewritten\sss in\sss
the form\vspace{1.5pt}
\[
\quad
t\dff \cdot\dff g_{\dff e}
\off =\off
g_{\dff t\dff(\dff e\trf)}\dff \cdot\dff 
k\trf(\dff e\fff,\pff t\trf)
\qff.
\]

\vspace{-12pt}\vspace{1.5pt}
Clearly\halfff,\oss the element\sss $k\trf(\dff e\fff,\pff t\trf)$\sss is\dss
uniquely\sss determined\dss by\dss the equality\halfff.\oss
Therefore\sss the edge relation\dss $E\dff(\dff e\fff,\pff t\trf)$\dss
holds in a quotient\sss group of\trs $\mathcal{F} \ast\dff G_{\dff v}$\dss
if\trs and\dss only\dss if\trs there exists an element\trs $h\qff \in\qff G_{\dff v}$\dss
such\dss that\sss
$t\dff \cdot\dff g_{\dff e}
\off =\off
g_{\dff t\dff(\dff e\trf)}\dff \cdot\dff 
h$\dss
in\dss this quotient\halfff.\oss

\mypar{Lemma.}{edge-product} 
\emph{The relations\qss $E\dff(\dff e\fff,\pff t_{\dff 1}\trf)$\dss
and\qss $E\dff(\dff t_{\dff 1}\dff(\dff e\trf)\fff,\pff t_{\dff 2}\trf)$\dss 
together\dss imply\qss
$E\dff(\dff e\fff,\pff t_{\dff 2}\dff t_{\dff 1}\trf)$\nnsp.\oss}\vspace{-3.5pt}

\proof
If\qss $E\dff(\dff e\fff,\pff t_{\dff 1}\trf)$\sss
and\dss $E\dff(\dff t_{\dff 1}\dff(\dff e\trf)\fff,\pff t_{\dff 2}\trf)$\sss hold,\oss
then\dss there exist\dss
$h_{\dff 1}\fff,\pff h_{\dff 2}\qff \in\qff G_{\dff v}$\dss
such\dss that\vspace{3pt}
\[
\quad
t_{\dff 1}\dff \cdot\dff g_{\dff e}
\off =\off
g_{\trf t_{\dff 1}\dff(\dff e\trf)}\dff \cdot\dff 
h_{\dff 1}
\hspace{1.5em}\mbox{and}\hspace{1.5em}
t_{\dff 2}\dff \cdot\dff g_{\trf t_{\dff 1}\dff(\dff e\trf)}
\off =\off
g_{\trf t_{\dff 2}\dff t_{\dff 1}\dff(\dff e\trf)}\dff \cdot\dff 
h_{\dff 2}
\qff.
\] 

\vspace{-12pt}\vspace{3pt}
It\dss follows\dss that\vspace{3pt}
\[
\quad
\left(\dff t_{\dff 2}\dff t_{\dff 1}\dff\right)\dff \cdot\dff g_{\dff e}
\off =\off
t_{\dff 2}\dff \cdot\dff 
g_{\trf t_{\dff 1}\dff(\dff e\trf)}\dff \cdot\dff 
h_{\dff 1}
\off =\off\
g_{\trf t_{\dff 2}\dff t_{\dff 1}\dff(\dff e\trf)}\dff \cdot\dff 
h_{\dff 2}\dff \cdot\dff h_{\dff 1}
\off =\off
g_{\trf t_{\dff 2}\dff t_{\dff 1}\dff(\dff e\trf)}\dff \cdot\dff 
\left(\dff h_{\dff 2}\dff h_{\dff 1}\dff\right)
\qff.
\]

\vspace{-12pt}\vspace{3pt}
Since\dss $h_{\dff 2}\dff h_{\dff 1}\qff \in\qff G_{\dff v}$\nsp,\oss
this implies\dss
$E\dff(\dff e\fff,\pff t_{\dff 2}\dff t_{\dff 1}\trf)$\nnsp.\oss  \eproof

\mypar{Lemma.}{edge-inverse}
\emph{The relations\qss $E\dff(\dff e\fff,\pff t^{\dff -\dff 1}\trf)$
and\pss $E\dff(\dff t^{\dff -\dff 1}\dff(\dff e\trf)\fff,\pff t\trf)$
are equivalent\halfff.\oss}

\proof
Let\trs $d\off =\off t^{\dff -\dff 1}\dff(\dff e\trf)$\nnsp.\oss
Then\dss $t\dff(\dff d\trf)\off =\off e$\nnsp.\oss
The relation\dss 
$E\dff(\dff t^{\dff -\dff 1}\dff(\dff e\trf)\fff,\pff t\trf)$\dss
means\dss that\vspace{3pt}
\[
\quad
t\dff \cdot\dff g_{\dff d}
\off =\off
g_{\dff t\dff(\dff d\trf)}\dff \cdot\dff 
h
\off =\off
g_{\dff e}\dff \cdot\dff 
h
\]

\vspace{-12pt}\vspace{3pt}
for some\dss $h\qff \in\qff G_{\dff v}$\nsp.\oss
Similarly\halfff,\oss the relation\dss
$E\dff(\dff e\fff,\pff t^{\dff -\dff 1}\trf)$\dss
means\dss that\vspace{3pt}
\[
\quad
t^{\dff -\dff 1}\dff \cdot\dff g_{\dff e}
\off =\off
g_{\dff t^{\dff -\dff 1}\dff(\dff e\trf)}\dff \cdot\dff 
k
\off =\off
g_{\dff d}\dff \cdot\dff 
k
\]

\vspace{-9pt}
for some\dss $k\qff \in\qff G_{\dff v}$\nsp.\oss
By\dss taking inverses we see\sss that\trs
$t^{\dff -\dff 1}\dff \cdot\dff g_{\dff e}
\off =\off
g_{\dff d}\dff \cdot\dff k$\qss is\dss
equivalent\dss to\vspace{3pt}
\[
\quad
g_{\dff e}^{\dff -\dff 1}\dff \cdot\dff t
\off =\off
k^{\dff -\dff 1}\dff \cdot\dff
g_{\dff d}^{\dff -\dff 1} 
\]

\vspace{-9pt}
and\dss hence\sss to\dss
$t\dff \cdot\dff g_{\dff d}
\off =\off
g_{\dff e}\dff \cdot\dff k^{\dff -\dff 1}$\dnsp.\oss
Since $G_{\dff v}$\sss is\dss invariant\dss under\dss taking inverses,\oss
it\dss follows\sss that\dss the relations\qss
$E\dff(\dff e\fff,\pff t^{\dff -\dff 1}\trf)$\dss
and\qss $E\dff(\dff t^{\dff -\dff 1}\dff(\dff e\trf)\fff,\pff t\trf)$\dss
are equivalent\halfff.\oss
This completes\sss the proof\halfff.\oss  \eproof

\mypar{Lemma.}{edge-triple-product} 
\emph{Suppose\sss that\trs $e\qff \in\qff E$ and\qss 
$t\fff,\pff r\fff,\pff u\qff \in\qff G_{\dff v}$\nsp.\oss 
Let\trs $d\off =\off t\dff(\dff e\trf)$\nnsp.\oss
If\qss $u\dff(\dff d\qff)\off =\off r\dff(\dff e\trf)$\nnsp,\oss 
then\dss the relations\qss 
$E\dff(\dff e\fff,\pff t\trf)$\nnsp,\pss 
$E\dff(\dff e\fff,\pff r\trf)$\nnsp,\oss 
and\pss 
$E\dff(\dff e\fff,\pff r^{\dff -\dff 1}\fff u\dff t\trf)$\dss 
together imply\qss $E\dff(\dff d\fff,\pff u\trf)$\nnsp.\oss}

\proof
The relations\dss
$E\dff(\dff e\fff,\pff t\trf)$\nnsp,\pss 
$E\dff(\dff e\fff,\pff r\trf)$\dss 
mean\dss that\vspace{4.5pt}
\[
\quad
t\dff \cdot\dff g_{\dff e}
\off =\off
g_{\dff t\dff(\dff e\trf)}\dff \cdot\dff 
h
\off =\off
g_{\dff d}\dff \cdot\dff 
h
\hspace{1.5em}\mbox{and}\hspace{1.5em}
r\dff \cdot\dff g_{\dff e}
\off =\off
g_{\dff r\dff(\dff e\trf)}\dff \cdot\dff 
h'
\] 

\vspace{-7.5pt}
for some\dss $h\fff,\pff h'\qff \in\qff G_{\dff v}$\nsp.\oss
Since\dss 
$u\dff t\trf(\dff e\trf)
\off =\off 
u\dff(\dff d\trf)
\off =\off 
r\trf(\dff e\trf)$\dss 
and\dss hence\sss 
$r^{\dff -\dff 1}\dff u\dff t\dff(\dff e\trf)
\off =\off 
e$\nnsp,\oss
the relation\dss
$E\dff(\dff e\fff,\pff r^{\dff -\dff 1}\fff u\dff t\trf)$\dss
means\sss that\sss 
$(\dff r^{\dff -\dff 1}\dff u\dff t\qff)\dff \cdot\dff g_{\dff e}
\off =\off
g_{\dff e}\dff \cdot\dff 
k$
for some\dss $k\qff \in\qff G_{\dff v}$\nsp.\oss 
It\dss follows\dss that\vspace{4.5pt}
\[
\quad
u\dff \cdot\dff
g_{\dff d}\dff \cdot\dff 
h
\off =\off
u\dff \cdot\dff
t\dff \cdot\dff 
g_{\dff e}
\off =\off
\left(\dff u\dff t\qff\right)\dff \cdot\dff
g_{\dff e}
\off =\off
r\qff
(\trf r^{\dff -\dff 1}\dff u\dff t\qff)\dff \cdot\dff g_{\dff e}
\]

\vspace{-34.5pt}
\[
\quad
\phantom{u\dff \cdot\dff
g_{\dff d}\dff \cdot\dff 
h
\off =\off
u\dff \cdot\dff
t\dff \cdot\dff 
g_{\dff e}
\off =\off
\left(\dff u\dff t\qff\right)\dff \cdot\dff
g_{\dff e}
\off }
=\off
r\dff \cdot\dff
g_{\dff e}\dff \cdot\dff 
k
\off =\off
g_{\dff r\dff(\dff e\trf)}\dff \cdot\qff 
h'\dff \cdot\dff
k
\off =\off
g_{\dff u\dff(\dff d\qff)}\dff \cdot\qff
(\trf h'\dff k\qff)
\]

\vspace{-12pt}\vspace{4.5pt}
and\dss hence\dss 
$u\dff \cdot\dff
g_{\dff d}
\off =\off
g_{\dff u\dff(\dff d\qff)}\dff \cdot\dff
(\trf h'\dff k\dff h^{\dff -\dff 1}\trf)$\nnsp.\oss
Since\dss
$h'\dff k\dff h^{\dff -\dff 1}\qff \in\qff G_{\dff v}$\nsp,\oss
this implies $E\dff(\dff d\fff,\pff u\trf)$\nnsp.\oss  \eproof

\mypar{Lemma.}{no-inversion-edges} 
\emph{Suppose\sss that\dss $e\qff \in\qff E$\nnsp,\dss
$h\qff \in\qff G_{\dff v}$\nsp,\pss and\sss $h\trf(\dff e\trf)\off =\off e$\nnsp.\oss
Let\sss
$a\off =\off s_{\dff e}^{\dff -\dff 1}\dff(\dff e \trf)$\sss
and\sss $s\off =\off s_{\dff e}$\nsp.\oss
Then\dss the relations\sss $E\dff(\dff e\fff,\pff h\trf)$
and\dss $L\trf(\trf l_{\dff e}\trf)$\sss together\sss imply\sss
$E\dff(\dff a\fff,\pff s^{\dff -\dff 1}\dff h\dff s \trf)$\nnsp.\oss}

\proof
Let\sss
$u\off =\off s^{\dff -\dff 1}\dff h\dff s$\nnsp.\oss
Then $u\dff(\dff a\trf)\off =\off a$\nnsp,\oss and\dss 
$E\dff(\dff e\fff,\pff h\trf)$\dss
and\dss
$E\dff(\dff a\fff,\pff u\trf)$\dss
have\sss the form\vspace{3pt}
\[
\quad
h\dff \cdot\dff g_{\dff e}
\off =\off 
g_{\dff e}\dff \cdot\dff
k\dff(\dff e\fff,\pff h\trf)
\quad\
\mbox{and}\quad\
u\dff \cdot\dff g_{\dff a}
\off =\off 
g_{\dff a}\dff \cdot\dff
k\dff(\dff a\fff,\pff u\trf)
\]

\vspace{-12pt}\vspace{3pt}
respectively\halfff.\oss
Recall\dss that\dss $L\trf(\trf l_{\dff e}\trf)$\dss
is\dss the relation\dss
$g_{\dff e}\dff \cdot\dff g_{\dff a}
\off =\off
s_{\dff e}\dff \cdot\dff s_{\dff a}$\nsp,\oss
or\halfff,\oss what\dss is\dss the same\vspace{3pt}
\begin{equation}
\label{e-a-relation}
\quad
g_{\dff a}
\dff \cdot\dff 
s_{\dff a}^{\dff -\dff 1}
\off =\off
g_{\dff e}^{\dff -\dff 1}
\dff \cdot\dff 
s_{\dff e}
\pff.
\end{equation}

\vspace{-12pt}\vspace{3pt}
By\dss the definition,\oss\vspace{1.5pt}
\[
\quad
k\trf(\dff a\fff,\pff u\trf)
\off =\off 
s_{\trf u\dff(\dff a\trf)}^{\dff -\dff 1}
\dff \cdot\dff 
u
\dff \cdot\dff 
s_{\dff a}
\off =\off 
s_{\dff a}^{\dff -\dff 1}
\dff \cdot\dff 
u
\dff \cdot\dff
s_{\dff a}
\off =\off
s_{\dff a}^{\dff -\dff 1}
\dff \cdot\dff 
s_{\dff e}^{\dff -\dff 1}
\dff \cdot\dff 
h
\dff \cdot\dff 
s_{\dff e}
\dff \cdot\dff 
s_{\dff a}
\]

\vspace{-12pt}\vspace{3pt}
and\dss hence\dss $E\dff(\dff a\fff,\pff u\trf)$\dss
is\dss equivalent\dss to\vspace{3pt}
\[
\quad
s_{\dff e}^{\dff -\dff 1}
\dff \cdot\dff 
h
\dff \cdot\dff 
s_{\dff e}
\dff \cdot\dff 
g_{\dff a}
\off =\off 
g_{\dff a}
\dff \cdot\dff
s_{\dff a}^{\dff -\dff 1}
\dff \cdot\dff
s_{\dff e}^{\dff -\dff 1}
\dff \cdot\dff 
h
\dff \cdot\dff 
s_{\dff e}
\dff \cdot\dff 
s_{\dff a}
\]

\vspace{-12pt}\vspace{3pt}
or\halfff,\oss what\dss is\dss the same,\oss to\vspace{3pt}
\[
\quad
s_{\dff e}^{\dff -\dff 1}
\dff \cdot\dff 
h
\dff \cdot\dff 
s_{\dff e}
\dff \cdot\dff 
g_{\dff a}
\dff \cdot\dff 
s_{\dff a}^{\dff -\dff 1}
\off =\off 
g_{\dff a}
\dff \cdot\dff
s_{\dff a}^{\dff -\dff 1}
\dff \cdot\dff
s_{\dff e}^{\dff -\dff 1}
\dff \cdot\dff 
h
\dff \cdot\dff 
s_{\dff e}
\pff.
\]

\vspace{-12pt}\vspace{3pt}
By\dss taking\sss into account\qss (\ref{e-a-relation})\qss
we can rewrite\sss this relation as\vspace{3pt}
\[
\quad
s_{\dff e}^{\dff -\dff 1}
\dff \cdot\dff 
h
\dff \cdot\dff 
s_{\dff e}
\dff \cdot\dff 
g_{\dff e}^{\dff -\dff 1}
\dff \cdot\dff 
s_{\dff e}
\off =\off 
g_{\dff e}^{\dff -\dff 1}
\dff \cdot\dff 
s_{\dff e}
\dff \cdot\dff
s_{\dff e}^{\dff -\dff 1}
\dff \cdot\dff 
h
\dff \cdot\dff 
s_{\dff e}
\pff,
\]

\vspace{-12pt}\vspace{3pt}
or\halfff,\oss
after\dss the obvious cancellations,\oss as\sss\vspace{2.25pt}
\[
\quad
{\dff e}^{\dff -\dff 1}
\dff \cdot\dff 
h
\dff \cdot\dff 
s_{\dff e}
\dff \cdot\dff 
g_{\dff e}^{\dff -\dff 1}
\off =\off 
g_{\dff e}^{\dff -\dff 1}
\dff \cdot\dff 
h 
\pff.
\]

\vspace{-12pt}\vspace{2.25pt}
Since\dss
$k\trf(\dff e\fff,\pff h\trf)
\off =\off 
s_{\dff e}^{\dff -\dff 1}
\dff \cdot\dff 
h
\dff \cdot\dff
s_{\dff e}$\nsp,\oss
we see\sss that\sss
$E\dff(\dff a\fff,\pff u\trf)$\sss
is\dss equivalent\dss to\sss 
$k\trf(\dff e\fff,\pff h\trf)
\dff \cdot\dff 
g_{\dff e}^{\dff -\dff 1}
\off =\off 
g_{\dff e}^{\dff -\dff 1}
\dff \cdot\dff 
h$
and\dss hence\sss to\dss
$E\dff(\dff e\fff,\pff h\trf)$\nnsp.\oss
This completes\sss the proof\halfff.\oss  \eproof

\mypar{Lemma.}{rotating-edge-loops}
\emph{Suppose\sss that\dss $e\qff \in\qff E$\dss
and\dss $h\qff \in\qff G_{\dff v}$\nsp.\oss
Let}\sss\vspace{4.5pt}
\[
\quad
c\off =\off h\trf(\dff e\trf)\dff,\quad
a\off =\off s_{\dff e}^{\dff -\dff 1}\dff(\dff e \trf)\dff,\quad
b\off =\off s_{\dff c}^{\dff -\dff 1}\dff(\dff c \trf)\dff,\quad
\mbox{\emph{and}}\quad
t
\off =\off 
s_{\dff c}^{\dff -\dff 1}\dff \cdot\dff h\dff \cdot\dff s_{\dff e}
\pff.
\]

\vspace{-12pt}\vspace{4.5pt}
\emph{Then $t\qff \in\qff G_{\dff v}\dff,\off b\off =\off t\trf(\dff a\trf)$ and\sss
the relations\sss
$E\dff(\dff e\fff,\pff h\trf)$\nnsp,\dss
$E\dff(\dff a\fff,\pff t\trf)$\nnsp,\dss
$L\trf(\trf l_{\dff e}\trf)$ 
together\sss imply\sss $L\trf(\trf l_{\dff c}\trf)$\nnsp.}

\proof
Clearly\halfff,\qss $t\qff \in\qff G_{\dff v}$ and\sss
$b\off =\off t\trf(\dff a\trf)$\nnsp.\oss
The relations\dss
$E\dff(\dff e\fff,\pff h\trf)$\dss
and\dss
$E\dff(\dff a\fff,\pff t\trf)$\dss
are\vspace{4.5pt}
\[
\quad
g_{\trf c}
\off =\off
h\dff \cdot\dff g_{\dff e}\dff \cdot\dff k\trf(\dff e\fff,\pff h\trf)^{\dff -\dff 1}
\quad
\mbox{and}\quad\
g_{\trf b}
\off =\off
t\dff \cdot\dff g_{\dff a}\dff \cdot\dff k\trf(\dff a\fff,\pff t\trf)^{\dff -\dff 1}
\]

\vspace{-7.5pt}
respectively\halfff.\oss
The equalities\sss 
$k\trf(\dff e\fff,\pff h\trf)
\off =\off 
s_{\trf c}^{\dff -\dff 1}\dff \cdot\dff h\dff \cdot\dff s_{\dff e}$\dss
and\dss
$k\trf(\dff a\fff,\pff t\trf)
\off =\off 
s_{\trf b}^{\dff -\dff 1}\dff \cdot\dff t\dff \cdot\dff s_{\dff a}$\dss 
imply\sss that\vspace{4.5pt}
\begin{equation}
\label{formula}
\quad
k\trf(\dff e\fff,\pff h\trf)^{\dff -\dff 1}
\qff \cdot\qff
t
\off =\off
s_{\dff e}^{\dff -\dff 1}\dff \cdot\dff h^{\dff -\dff 1}\dff \cdot\dff s_{\trf c}
\dff \cdot\dff 
s_{\dff c}^{\dff -\dff 1}\dff \cdot\dff h\dff \cdot\dff s_{\dff e}
\off =\off
1
\pff.
\end{equation}

\vspace{-7.5pt}
The relation\dss $L\trf(\trf l_{\dff e}\trf)$\sss is\dss
$g_{\dff e}\dff \cdot\dff g_{\dff a}
\off =\off
s_{\dff e}\dff \cdot\dff s_{\dff a}$\nsp.\oss
By\sss using\qss (\ref{formula})\qss and\dss then\dss $L\trf(\trf l_{\dff e}\trf)$\dss
with we see\sss that\vspace{3pt}
\[
\quad
g_{\dff c}\dff \cdot\dff g_{\dff b}
\off =\off
h\dff \cdot\dff g_{\dff e}\dff \cdot\dff k\trf(\dff e\fff,\pff h\trf)^{\dff -\dff 1}
\qff \cdot\qff
t\dff \cdot\dff g_{\dff a}\dff \cdot\dff k\trf(\dff a\fff,\pff t\trf)^{\dff -\dff 1}
\]

\vspace{-33pt}
\[
\quad
\phantom{g_{\dff c}\dff \cdot\dff g_{\dff d}
\off }
=\off
h\dff \cdot\dff g_{\dff e}\dff \cdot\dff g_{\dff a}\dff \cdot\dff k\trf(\dff a\fff,\pff t\trf)^{\dff -\dff 1}
\off =\off
h\dff \cdot\dff s_{\dff e}\dff \cdot\dff s_{\dff a}\dff \cdot\dff k\trf(\dff a\fff,\pff t\trf)^{\dff -\dff 1}
\]

\vspace{-33pt}
\[
\quad
\phantom{g_{\dff c}\dff \cdot\dff g_{\dff d}
\off }
=\off
h\dff \cdot\dff s_{\dff e}\dff \cdot\dff s_{\dff a}\dff \cdot\dff
s_{\dff a}^{\dff -\dff 1}\dff \cdot\dff t^{\dff -\dff 1}\dff \cdot\dff s_{\dff b}
\off =\off
h\dff \cdot\dff s_{\dff e}\dff \cdot\dff t^{\dff -\dff 1}\dff \cdot\dff s_{\dff b}
\]

\vspace{-33pt}
\[
\quad
\phantom{g_{\dff c}\dff \cdot\dff g_{\dff d}
\off }
=\off
h\dff \cdot\dff s_{\dff e}\dff \cdot\dff 
s_{\dff e}^{\dff -\dff 1}\dff \cdot\dff h^{\dff -\dff 1}\dff \cdot\dff s_{\dff c}
\dff \cdot\dff s_{\dff b}
\off =\off
s_{\dff c}
\dff \cdot\dff s_{\dff b}
\pff.
\]

\vspace{-9pt}
Therefore\dss
$g_{\dff c}\dff \cdot\dff g_{\dff b}
\off =\off
s_{\dff c}
\dff \cdot\dff s_{\dff b}$\nsp,\oss
i.e.\qss the relation\dss
$L\trf(\trf l_{\dff c}\trf)$\sss holds.\oss  \eproof

\mypar{Lemma.}{inverting-edge-loops}
\emph{Let\dss $e\qff \in\qff E$\sss and\dss
$a\off =\off s_{\dff e}^{\dff -\dff 1}\dff(\dff e \trf)$\nnsp.\oss
Let\sss $h\off =\off s_{\dff e}\dff \cdot\dff s_{\dff a}$\nsp.\oss
Then\sss $h\qff \in\qff G_{\dff v}$ and\dss
if\pss the relation\dss
$E\dff(\dff e\fff,\pff h^{\dff -\dff 1}\trf)$\dss
holds,\oss
then\dss the relations\qss $L\trf(\trf l_{\dff e}\trf)$\sss
and\pss $L\trf(\trf l_{\dff a}\trf)$\sss
are equivalent\halfff.\oss}

\proof
Clearly\halfff,\qss 
$s_{\dff a}\dff(\dff v\trf)
\off =\off
s_{\dff e}^{\dff -\dff 1}\dff(\dff v \trf)$
and\dss hence\dss
$h\dff(\dff v\trf)
\off =\off
s_{\dff e}\dff \cdot\dff s_{\dff a}\trf(\dff v\trf)
\off =\off
v$\nnsp,\oss
i.e.\qss $h\qff \in\qff G_{\dff v}$\nsp.\oss
This proves\sss the first\sss claim of\dss the\sss lemma.\oss
Let\dss
$d
\off =\off 
s_{\dff a}^{\dff -\dff 1}\dff(\dff a\trf)$\nnsp.\oss
Then\dss\vspace{4.5pt}
\[
\quad
d
\off =\off
s_{\dff a}^{\dff -\dff 1}\dff\left(\qff s_{\dff e}^{\dff -\dff 1}\dff(\dff e\qff)\trf\right)
\off =\off
s_{\dff a}^{\dff -\dff 1}\dff \cdot\dff s_{\dff e}^{\dff -\dff 1}\qff(\dff e\trf)
\off =\off
(\trf s_{\dff e}\dff \cdot\dff s_{\dff a}\trf)^{\dff -\dff 1}\qff(\dff e\trf)
\off =\off
h^{\dff -\dff 1}\qff(\dff e\trf)
\qff.
\]

\vspace{-7.5pt}
The relation\dss $E\dff(\dff e\fff,\pff h^{\dff -\dff 1}\trf)$\dss
implies\dss that\sss \vspace{4.5pt}
\[
\quad
h^{\dff -\dff 1}\dff \cdot\dff g_{\dff e}
\off =\off
g_{\dff h^{\dff -\dff 1}\qff(\dff e\trf)}\dff \cdot\qff
k\trf(\dff e\fff,\pff h^{\dff -\dff 1}\trf)
\off =\off
g_{\dff d}\dff \cdot\qff
k\trf(\dff e\fff,\pff h^{\dff -\dff 1}\trf)
\qff.
\]

\vspace{-7.5pt}
Since\dss
$k\trf(\dff e\fff,\pff h^{\dff -\dff 1}\trf)
\off =\off
s_{\trf d}^{\dff -\dff 1}\dff \cdot\dff h^{\dff -\dff 1}\dff \cdot\dff s_{\dff e}$\nsp,\oss
it\dss follows\dss that\sss 
$g_{\dff d}
\off =\off
h^{\dff -\dff 1}\dff \cdot\dff g_{\dff e}
\qff \cdot\qff
s_{\trf e}^{\dff -\dff 1}\dff \cdot\dff h\dff \cdot\dff s_{\dff d}$
and\dss hence\vspace{4.5pt}
\begin{equation}
\label{gs}
\quad
g_{\trf d}
\dff \cdot\dff 
s_{\dff d}^{\dff -\dff 1}
\off =\off
h^{\dff -\dff 1}
\dff \cdot\dff 
g_{\dff e}
\dff \cdot\dff 
s_{\dff e}^{\dff -\dff 1}\dff \cdot\dff h
\qff.
\end{equation}

\vspace{-7.5pt}
The relation\dss $L\trf(\trf l_{\dff a}\trf)$\sss is\dss the relation\dss
$g_{\dff a}\dff \cdot\dff g_{\dff d}
\off =\off
s_{\dff a}\dff \cdot\dff s_{\dff d}$\nsp,\oss
or\halfff,\oss
what\dss is\dss the same,\vspace{4.5pt}
\[
\quad
g_{\dff a}
\dff \cdot\dff 
g_{\dff d}
\dff \cdot\dff 
s_{\dff d}^{\dff -\dff 1}
\off =\off
s_{\dff a}
\pff.
\]

\vspace{-7.5pt}
In\dss view of\pss (\ref{gs})\qss it\dss is\dss equivalent\dss to\vspace{4.5pt}
\[
\quad
g_{\dff a}
\dff \cdot\dff 
h^{\dff -\dff 1}
\dff \cdot\dff 
g_{\dff e}
\dff \cdot\dff 
s_{\dff e}^{\dff -\dff 1}\dff \cdot\dff h
\off =\off
s_{\dff a}
\pff.
\]

\vspace{-7.5pt}
But\dss
$s_{\dff a}\off =\off s_{\dff e}^{\dff -\dff 1}\dff \cdot\dff h$\dss
and\dss hence\dss $L\trf(\trf l_{\dff a}\trf)$\sss
is\dss equivalent\dss to\sss 
$g_{\dff a}
\dff \cdot\dff 
h^{\dff -\dff 1}
\dff \cdot\dff 
g_{\dff e}
\off =\off
1$\nnsp,\oss
which\dss is,\oss in\dss turn,\oss is\dss equivalent\dss to\dss
$g_{\dff e}\dff \cdot\dff g_{\dff a}
\off =\off
h$\nnsp.\oss
Since\dss 
$h
\off =\off
s_{\dff e}\dff \cdot\dff s_{\dff a}$\nsp,\oss
the last\dss relation\dss is\dss nothing\sss else\sss but\trs
$L\trf(\trf l_{\dff e}\trf)$\nnsp.\oss
It\dss follows\dss that\trs
$L\trf(\trf l_{\dff e}\trf)$\dss
and\trs
$L\trf(\trf l_{\dff a}\trf)$\dss
are equivalent\halfff.\oss  \eproof

\mysection{Presentations}{presentations}

\myuppar{Turning\sss edge relations into definitions.}
Suppose\sss that\sss 
$e\qff \in\qff E$\nnsp,\dss $h\qff \in\qff G_{\dff v}$\nsp,\qss
and\dss let\dss
$c\off =\off h\dff(\dff e\trf)$\nnsp.\oss
The relation\dss $E\dff(\dff e\fff,\pff h\trf)$\dss
can\dss be rewritten as\vspace{3pt}
\begin{equation}
\label{edge-definition}
\quad
g_{\trf c}
\off =\off
h\dff \cdot\dff g_{\dff e}\dff \cdot\qff k\trf(\dff e\fff,\pff h\trf)^{\dff -\dff 1}
\pff
\end{equation}

\vspace{-9pt}
and\dss interpreted as a\qss \emph{definition}\qss of\trs the generator\dss 
$g_{\trf c}$\dss
in\dss terms of\trs the generator\dss $g_{\dff e}$\dss
and\sss elements\dss
$h\fff,\qff k\trf(\dff e\fff,\pff h\trf)\qff \in\qff G_{\dff v}$\nsp.\oss
Of\dss course,\oss different\sss choices of\trs
$h\qff \in\qff G_{\dff v}$\sss
such\dss that\sss
$c\off =\off h\dff(\dff e\trf)$\dss
lead\dss to different\dss definitions of\trs $g_{\dff c}$\nnsp,\oss
but\dss they\dss turn out\dss to be equivalent\dss if\trs the relations
$E\dff(\dff e\fff,\pff t\trf)$ with $t\qff \in\qff G_{\dff e}$
hold.\oss
Indeed,\oss if\trs
$h\fff,\pff r\qff \in\qff G_{\dff v}$ 
and\dss 
$c\off =\off h\dff(\dff e\trf)\off =\off r\dff(\dff e\trf)$\nnsp,\oss 
then 
$r^{\dff -\dff 1}\fff h\fff,\pff
h^{\dff -\dff 1}\fff r
\qff \in\qff G_{\dff e}$\nsp.\oss
By\sss applying\trs
Lemma\qss \ref{edge-triple-product}\qss to\dss $h\fff,\pff r$\dss
and\dss $u\off =\off 1$\dss and observing\dss that\dss 
$E\dff(\dff e\fff,\pff 1\trf)$\dss always holds,\oss
we see\sss that\trs
$E\dff(\dff e\fff,\pff r\trf)$\dss
implies\dss $E\dff(\dff e\fff,\pff h\trf)$\dss
if\pss 
$E\dff(\dff e\fff,\pff r^{\dff -\dff 1}\fff h\trf)$\dss
holds.\oss
Similarly\halfff,\pss
$E\dff(\dff e\fff,\pff h\trf)$\dss
implies\dss $E\dff(\dff e\fff,\pff r\trf)$\dss
if\pss 
$E\dff(\dff e\fff,\pff h^{\dff -\dff 1}\fff r\trf)$\dss
holds.\oss
This proves our claim.\oss
 
Let\dss us assume\sss that\dss the relations
$E\dff(\dff e\fff,\pff t\trf)$ with $t\qff \in\qff G_{\dff e}$
hold and\sss
define $g_{\dff c}$ by\dss the formula\qss (\ref{edge-definition})\qss 
for every\sss edge
$c\qff \in\qff E$ of\dss the form\dss
$c\off =\off h\dff(\dff e\trf)$\dss
with\dss $h\qff \in\qff G_{\dff v}$\nsp,\oss
i.e.\dss for every\sss edge $c$ in\dss the $G_{\dff v}$\dnsp-orbit\sss of\dss $e$\nnsp.\oss
Then\dss the relations
$E\dff(\dff e\fff,\pff h\trf)$ with $h\qff \in\qff G_{\dff v}$\sss
hold\dss by\dss the definition.\oss
Moreover\halfff,\oss
Lemma\qss \ref{edge-triple-product}\qss implies\sss that\dss
the relations
$E\dff(\dff c\fff,\pff h\trf)$
with\sss $h\qff \in\qff G_{\dff v}$
and\sss $c$\sss
in\dss the $G_{\dff v}$\dnsp-orbit\sss of\dss $e$\sss also hold.\oss
Similarly\halfff,\oss
Lemma\qss \ref{rotating-edge-loops}\qss implies\sss that\dss
relations\dss $L\trf(\trf l_{\dff c}\trf)$\sss
for such $c$\sss follow\dss from\dss 
$L\trf(\trf l_{\dff e}\trf)$\nnsp.\oss

Therefore,\oss for each orbit $O$ of\trs the action of\sss $G_{\dff v}$
on\sss $E$ one needs\sss 
the generator $g_{\dff e}$ only\dss for one representative\sss $e$\sss
of\sss $O$
and\sss
among\dss the relations $E\dff(\dff c\fff,\pff h\trf)$ and\dss
$L\trf(\trf l_{\dff c}\trf)$ with $c\qff \in\qff O$
and\dss
$h\qff \in\qff G_{\dff v}$\sss
one needs only\dss the relations $E\dff(\dff e\fff,\pff h\trf)$ 
with $h\qff \in\qff G_{\dff v}$\sss
and\dss the relation\dss $L\trf(\trf l_{\dff e}\trf)$\nnsp.\oss
Moreover\halfff,\oss
Lemmas\qss \ref{edge-product}\qss and\qss \ref{edge-inverse}\qss imply\dss
that\dss one needs only\dss the relations $E\dff(\dff e\fff,\pff h\trf)$ with $h$\sss
belonging\dss to a set\sss of\dss generators of\sss $G_{\dff v}$
and\dss the relation\dss $L\trf(\trf l_{\dff e}\trf)$\nnsp.\oss
Of\dss course,\oss the\sss loop relations involving discarded\dss generators
should\dss be rewritten\sss in\dss terms of\dss remaining ones.\oss

\myuppar{Inversions.}
Recall\trs that\sss every edge $e\qff \in\qff E$\sss has $v$ as an endpoint\sss
and\sss that\sss we denote by $\ttt(\dff e\trf)$
the other endpoint\sss of $e$\nnsp.\oss 
An element $g\qff \in\qff G$\sss is\dss said\sss to be an\qss
\emph{inversion}\qss of\dss $e\qff \in\qff E$\sss if\dss
$g\dff(\dff e\trf)\off =\off e$ and $g$ interchanges\sss
the endpoints of\sss $e$\nnsp,\oss
i.e.\dss $g\dff(\dff v\trf)\off =\off \ttt(\dff e\trf)$
and $g\dff(\trf \ttt(\dff e\trf)\trf)\off =\off v$\nnsp.\oss
Clearly\halfff,\oss if\dss $g$\sss is\dss an\sss inversion of\dss $e$\sss
and $h\qff \in\qff G_{\dff v}$\nsp,\oss
then $h\dff g\dff h^{\dff -\dff 1}$\dss is\dss an\sss inversion of\sss $h\dff(\dff e\trf)$\nnsp.\oss
Hence,\oss if\dss $e$\sss admits an\sss inversion,\oss
then every\sss edge in\dss the $G_{\dff v}$\dnsp-orbit\dss 
of\sss $e$\dss admits an\sss inversion.

\mypar{Lemma.}{orbits-inversions}
\emph{Suppose\sss that\dss an edge\dss $e\qff \in\qff E$\dss does not\sss
admits an inversion.\oss
If\trs $s\dff(\dff v\trf)\off =\off \ttt(\dff e\trf)$\nnsp,\oss
then\dss the edge\dss $s^{\dff -\dff 1}\dff(\dff e\trf)$\dss
does not\dss belongs\sss to\sss the $G_{\dff v}$\dnsp-orbit\sss of\dss $e$\nnsp.\oss}

\proof
Let\dss $a\off =\off s^{\dff -\dff 1}\dff(\dff e\trf)$\nnsp.\oss
Since\dss $s^{\dff -\dff 1}\dff(\trf \ttt(\dff e\trf)\trf)\off =\off v$\nnsp,\oss
the edge $a$\sss belongs\sss to\dss $E$\nnsp.\oss
If\dss $a$\sss belongs\sss to\sss the $G_{\dff v}$\dnsp-orbit\sss of\dss $e$\nnsp,\oss
then $h\dff(\dff a\trf)\off =\off e$\sss for some
$h\qff \in\qff G_{\dff v}$\nsp.\oss
Clearly,\dss
$h\dff s^{\dff -\dff 1}\dff(\trf \ttt(\dff e\trf)\trf)\off =\off v$\nnsp.\oss
Also,\pss $s^{\dff -\dff 1}\dff(\dff v\trf)\off =\off \ttt(\dff a\trf)$\dss
and\dss hence\dss
$h\dff s^{\dff -\dff 1}\dff(\dff v\trf)\off =\off \ttt(\dff e\trf)$\nnsp.\oss
It\dss follows\dss that\dss $h\dff s^{\dff -\dff 1}$\sss is\dss
an\sss inversion of\dss $e$\nnsp,\oss contrary\dss to\sss the assumption.\oss
This proves\sss the lemma.\oss   \eproof

\mypar{Lemma.}{correctness-involution}
\emph{Suppose\sss that\trs 
$d\off =\off u\trf(\dff e\trf)$\nnsp,\oss 
where\dss $e\qff \in\qff E$\dss 
and\dss $u\qff \in\qff G_{\dff v}$\nsp.\oss
If\qss $s\dff(\dff v\trf)\off =\off \ttt(\dff e\trf)$\dss
and\dss $r\trf(\dff v\trf)\off =\off \ttt(\dff d\trf)$\nnsp,\oss
then\dss $s^{\dff -\dff 1}\dff(\dff e\trf)$\dss 
and\dss $r^{\dff -\dff 1}\dff(\dff d\trf)$\dss 
belong\dss to\sss the same\dss $G_{\dff v}$\dnsp-orbit\halfff.\oss}

\proof
If\trs $d\off =\off e$\nnsp,\oss
then\dss $r\trf(\dff v\trf)\off =\off \ttt(\dff e\trf)$\dss
and\dss 
$s^{\dff -\dff 1}\dff r\trf(\dff v\trf)\off =\off v$\nnsp.\oss
Therefore\dss $s^{\dff -\dff 1}\dff r\qff \in\qff G_{\dff v}$\nsp.\oss
Since\vspace{1.5pt}
\[
\quad
s^{\dff -\dff 1}\dff r\qff
\left(\dff r^{\dff -\dff 1}\dff(\dff e\trf)\dff\right)
\off =\off
s^{\dff -\dff 1}\dff(\dff e\trf)
\qff,
\]

\vspace{-12pt}\vspace{1.5pt}
in\dss this case\dss
$s^{\dff -\dff 1}\dff(\dff e\trf)$\dss 
and\dss 
$r^{\dff -\dff 1}\dff(\dff e\trf)
\off =\off
r^{\dff -\dff 1}\dff(\dff d\trf)$\dss 
belong\dss to\sss the same\dss $G_{\dff v}$\dnsp-orbit\halfff.\oss
This proves\sss the\sss lemma\sss in\dss the case\dss
$d\off =\off e$\nnsp.\oss
In\dss the general\sss case,\oss
let\dss
$t\off =\off u^{\dff -\dff 1}\dff r$\nnsp.\oss
Then\dss\vspace{1.5pt}
\[
\quad
t^{\dff -\dff 1}\dff(\dff e\trf)
\off =\off
r^{\dff -\dff 1}\dff u\trf(\dff e\trf)
\off =\off 
r^{\dff -\dff 1}\dff(\dff d\trf)
\qff.
\]

\vspace{-12pt}\vspace{1.5pt}
Clearly\halfff,\oss 
$t\trf(\dff v\trf)
\off =\off
u^{\dff -\dff 1}\dff r\dff(\dff v\trf)
\off =\off
u^{\dff -\dff 1}\dff(\trf \ttt(\dff d\trf)\trf)
\off =\off
\ttt(\dff e\trf)$\qss
and\dss hence\sss the already\dss proved special\sss case 
implies\sss that\trs
$s^{\dff -\dff 1}\dff(\dff e\trf)$\dss 
and\dss 
$t^{\dff -\dff 1}\dff(\dff e\trf)
\off =\off 
r^{\dff -\dff 1}\dff(\dff d\trf)$\dss 
belong\dss to\sss the same\dss $G_{\dff v}$\dnsp-orbit\halfff.\oss  \eproof

\myuppar{An involution on\dss the set\sss of\dss orbits.}
Let\dss us\sss define an\sss involution\dss $\iota$\dss on\dss the set\sss
of\trs $G_{\dff v}$\dnsp-orbits on\dss $E$\dss as follows.\oss
Let\dss $\varepsilon$\dss be\sss the orbit\sss of\dss $e$\nnsp.\oss
If\dss $e$\dss admits an\sss inversion,\oss then\sss we set\trs
$\iota\dff(\dff \varepsilon\dff)\off =\off \varepsilon$\nnsp.\oss
If\sss $\varepsilon$ does\sss not\sss admits an\sss inversion,\oss
then\sss we\sss take as\dss $\iota\dff(\dff \varepsilon\dff)$\dss
the orbit\sss of\trs $s^{\dff -\dff 1}\dff(\dff e\trf)$\nnsp,\oss
where $s\qff \in\qff G$\sss is\dss such\dss that\sss
$s\dff(\dff v\trf)\off =\off \ttt(\dff e\trf)$\nnsp.\oss
By\qss Lemma\qss \ref{correctness-involution}\qss this definition\dss
is\dss correct\halfff.\oss
Let\dss us\dss check\dss that\dss $\iota$\dss is\dss an\sss involution.\oss
Suppose\sss that\sss $s\dff(\dff v\trf)\off =\off \ttt(\dff e\trf)$
and\dss let\sss
$a\off =\off s^{\dff -\dff 1}\dff(\dff e\trf)$\nnsp.\oss
Then\dss $s^{\dff -\dff 1}\dff(\dff v\trf)\off =\off \ttt(\dff a\trf)$\dss
and\vspace{1.5pt}
\[
\quad
\left(\trf s^{\dff -\dff 1} \trf\right)^{-\dff 1}\dff (\dff a\trf)
\off =\off
s\dff \left(\trf s^{\dff -\dff 1}\dff(\dff e\trf)\trf\right)
\off =\off
e
\qff.
\]

\vspace{-10.5pt}
It\dss follows\dss that\trs
$\iota\dff(\dff \iota\dff(\dff \varepsilon\dff)\dff)$\dss
is\dss equal\dss to\sss the orbit\sss of\dss $e$\nnsp,\oss
i.e.\qss to\dss $\varepsilon$\nnsp.\oss
Hence\dss $\iota\dff \circ\trf \iota\off =\off \id$\nnsp,\oss
i.e.\qss $\iota$\dss is\dss an\sss involution.\oss
Lemma\qss \ref{orbits-inversions}\qss implies\sss that\trs
$\iota\dff(\dff \varepsilon\dff)\off =\off \varepsilon$\trs
if\trs and\dss only\trs if\dss $e$\sss admits an\sss inversion.\oss

\myuppar{Turning\sss some edge-loop relations into definitions.}
Let\sss $e\qff \in\qff E$\sss and\dss 
$a\off =\off s_{\dff e}^{\dff -\dff 1}\dff(\dff e \trf)$\nnsp.\oss 
Suppose\sss that\sss
$e$ does not\sss admit\sss an inversion.\oss 
Lemma\qss \ref{orbits-inversions}\qss
implies\sss that\sss $a$ and $e$ belong\dss to different\sss $H$\dnsp-orbits.\oss
The relation\dss $L\trf(\trf l_{\dff e}\trf)$\sss 
can\sss be rewritten as\vspace{1.5pt}
\[
\quad
g_{\dff a}
\off =\off
g_{\dff e}^{\dff -\dff 1}\dff \cdot\dff
s_{\dff e}\dff \cdot\dff s_{\dff a}
\qff,
\]

\vspace{-10.5pt}
where 
$s_{\dff e}\dff \cdot\dff s_{\dff a}\qff \in\qff G_{\dff v}$\nsp.\oss
If\trs the edge $e$\sss is\dss 
the representative of\dss an orbit $O$\nnsp,\oss
then one can\sss take $a$ as\sss the representative of\trs
the orbit\sss $\iota\trf(\dff O\dff)\off \neq\off O$
and\sss interpret\sss $L\trf(\trf l_{\dff e}\trf)$
as\sss a\qss \emph{definition}\qss of\sss 
$g_{\trf a}$
in\dss terms of\sss $g_{\dff e}$
and 
$s_{\dff e}\dff \cdot\dff s_{\dff a}\qff \in\qff G_{\dff v}$\nsp.\oss
If\trs the representative of\sss $\iota\trf(\dff O\dff)$\sss is\dss
some other edge
$b\qff \in\qff E$\nnsp,\oss
then one can combine\sss the definition of\sss $g_{\dff b}$
in\dss terms of\sss $g_{\dff a}$ with\sss 
the definition of\sss $g_{\dff a}$
in\dss terms of\sss $g_{\dff e}$
and\sss get\sss a definition of\sss $g_{\dff b}$
in\dss terms of\sss $g_{\dff e}$\nsp.\oss
But\dss it\dss is\dss more natural\dss to replace $b$ by $a$
as\sss the representative.\oss
With\dss this\sss interpretation\dss  
$L\trf(\trf l_{\dff e}\trf)$\sss
holds\sss by\dss the definition.\oss
Lemma\qss \ref{inverting-edge-loops}\qss implies\sss that\trs
$L\trf(\trf l_{\dff a}\trf)$\sss also holds.\oss
Moreover\halfff,\oss Lemma\qss \ref{no-inversion-edges}\qss
implies\sss that\dss the edge relations
$E\dff(\dff a\fff,\pff r\trf)$
with\sss $r\qff \in\qff G_{\dff v}$\sss
follow\dss from\dss the edge relations
$E\dff(\dff e\fff,\pff h\trf)$
with\sss $h\qff \in\qff G_{\dff v}$\nsp.\oss
Therefore we can discard\dss the generator $g_{\dff a}$
and\dss the edge relations\dss 
$E\dff(\dff a\fff,\pff r\trf)$
with\sss $r\qff \in\qff G_{\dff v}$\nsp.\oss

\myuppar{Edges admitting\sss an\sss inversion.}
As above,\oss let\sss
$e\qff \in\qff E$\sss
and\dss let\sss
$a\off =\off s_{\dff e}^{\dff -\dff 1}\dff(\dff e \trf)$\nnsp.\oss 
Suppose\sss now\dss that
$e$ admits an inversion.\oss
Then $a$ and $e$ belong\dss to\sss the same
$G_{\dff v}$\dnsp-orbit\halfff.\oss
If\sss $e$ is\dss used\sss as a representative,\oss 
then $g_{\dff a}$ is\dss already\dss defined\sss in\dss terms 
of\sss $g_{\dff e}$\nsp.\oss
Moreover\halfff,\oss if\sss $e$ admits\sss an inversion,\pss 
one can choose as $s_{\dff e}$ such an\sss inversion,\pss
and\sss with such a choice $a\off =\off e$\nnsp.\oss
So,\oss when $e$ admits an\sss inversion,\qss
$L\trf(\trf l_{\dff e}\trf)$ cannot\sss be\sss
turned\sss into a definition and\sss
needs\sss to be rewritten\dss in\sss terms of\dss $g_{\dff e}$\nsp.\oss
If\dss $s_{\dff e}$\sss is\dss an\dss inversion
of\dss $e$\nnsp,\oss
then\sss $L\trf(\trf l_{\dff e}\trf)$\sss 
is\dss simply\dss the relation\dss \vspace{1.5pt}
\[
\quad
g_{\dff e}\dff \cdot\dff g_{\dff e}
\off =\off
s_{\dff e}\dff \cdot\dff s_{\dff e} 
\pff.
\]

\vspace{-10.5pt}
In\dss general,\pss
$a
\off =\off 
h\dff(\dff e\trf)$\dss
for some\dss $h\qff \in\qff G_{\dff v}$
and\dss 
$L\trf(\trf l_{\dff e}\trf)$\dss 
is\dss equivalent\dss to\dss
$g_{\dff e}\dff \cdot\qff 
h
\dff \cdot\dff 
g_{\dff e}
\off =\off
s_{\dff e}\dff \cdot\qff 
h
\dff \cdot\dff 
s_{\dff e}$\nsp.\oss
We\sss will\sss not\sss use\sss this fact\sss and\dss leave\sss 
its verification\sss to\sss the reader\halfff.\oss

\myuppar{Scaffoldings.}
A\pss \emph{scaffolding}\pss for\sss the action of\sss $G$ on $X$ consists\sss of\trs
a\sss set\dss $E_{\dff 0}$\dss of\trs representatives of\trs $G_{\dff v}$\dnsp-orbits\sss in $E$\nnsp,\oss
a\dss family\sss of\dss sets\sss $\mathcal{T}_{\fff e}$\sss  
of\trs representatives of\dss cosets in\dss 
$G_{\dff v}/G_{\dff e}$\dss for\dss $e\qff \in\qff E_{\dff 0}$\nsp,\oss
and\sss a\sss family\dss $s_{\dff e}$\nsp,\dss $e\qff \in\qff E$\dss of\dss elements of\trs $G$\dss
as in\dss Section\qss \ref{one},\oss
i.e.\qss such\dss that\trs 
$s_{\dff e}\dff (\dff v\trf)
\off =\off 
\ttt(\dff e\trf)$\dss
for every\dss $e\qff \in\qff E$\nnsp.\oss
A\sss scaffolding\dss is\dss said\dss to be\qss \emph{regular}\pss if\trs the following\sss
three conditions hold.\oss\vspace{-6pt}
\begin{itemize}

\item[({\fff}i{\fff})]\quad\
If\dss $e\qff \in\qff E_{\dff 0}$\sss 
admits an\sss inversion,\oss then\sss $s_{\dff e}$\sss is\dss an\sss inversion of\dss $e$\nnsp.\oss

\item[({\fff}ii{\fff})]\quad\
If\dss $e\qff \in\qff E_{\dff 0}$\sss does not\sss 
admits an\sss inversion,\oss then\dss
$a
\off =\off 
s_{\dff e}^{\dff -\dff 1}\dff(\dff e\trf)
\pff \in\pff
E_{\dff 0}$\dss
and\dss $s_{\dff a}\off =\off s_{\dff e}^{\dff -\dff 1}$\nnsp.

\item[({\fff}iii{\fff})]\quad\
If\qss $e\qff \in\qff E_{\dff 0}$\nsp,\pss $u\qff \in\qff \mathcal{T}_{\dff e}$\nsp,\oss
and\dss $d\off =\off u\dff(\dff e\trf)$\nnsp,\oss
then\dss 
$s_{\dff d}
\off =\off
u\dff s_{\dff e}\dff u^{\dff -\dff 1}$\nnsp.\oss

\end{itemize}

\vspace{-6pt}
Since $E_{\dff 0}$\sss is\dss a set\sss of\dss representatives of\sss $G_{\dff v}$\dnsp-orbits,\oss
the involution\sss $\iota$\sss induces an\sss involution\sss $\iota_{\dff 0}$ on\sss $E_{\dff 0}$\nsp.\oss
Let\sss $E_{\dff 1}$\sss be a set\sss of\dss representatives of\dss the orbits of\dss
the involution\sss $\iota_{\dff 0}$\nsp.\oss

The property\qss (iii)\qss is\dss intended\sss for simplifying\dss 
the process of\dss rewriting\sss the\sss loop relations
in\sss terms of\dss generators $g_{\dff e}$\sss with $e\qff \in\qff E_{\dff 0}$\nsp.\oss
If\qss 
$s_{\dff d}
\off =\off
u\dff s_{\dff e}\dff u^{\dff -\dff 1}$\nnsp,\oss
then\vspace{3pt}
\[
\quad
k\trf(\dff e\fff,\pff u\trf)
\off =\off
s_{\dff d}^{\dff -\dff 1}\dff u\qff s_{\dff e}
\off =\off
u\qff s_{\dff e}^{\dff -\dff 1}\dff u\dff ^{\dff -\dff 1}\dff u\qff s_{\dff e}
\off =\off
u\qff s_{\dff e}^{\dff -\dff 1}\dff s_{\dff e}
\off =\off
u
\qff.
\]

\vspace{-9pt}
Hence $E\dff(\dff e\fff,\pff u\trf)$\sss 
is\dss equivalent\dss to\sss
$g_{\dff d}^{\dff -\dff 1}\dff \cdot\dff u\dff \cdot\dff g_{\dff e}
\off =\off
u$\sss
and one can define $g_{\dff d}$ as
$u\dff \cdot\dff g_{\dff e}\dff \cdot\dff u^{\dff -\dff 1}$\nsp\dnsp.\oss

\mypar{Lemma.}{coherent-frames}
\emph{There exist\sss regular\dss scaffoldings.\oss}

\proof
Let\dss us\dss begin\dss by\sss
choosing\sss an arbitrary\sss representative of\dss every\sss $G_{\dff v}$\dnsp-orbit\sss
consisting of\dss edges admitting\sss an\sss inversion.\oss
For every\sss such\sss a\sss representative\sss $e$\dss let\sss $s_{\dff e}$\sss
be an\sss inversion of\dss $e$\nnsp.\oss
We will\sss consider\dss the other orbits in\sss pairs of\trs the form\dss
$\varepsilon\fff,\pff \iota\dff(\dff \varepsilon\dff)$\nnsp.\oss
Let\sss us choose\sss from each such\sss pair one orbit\halfff,\pss
say\halfff,\qss $\varepsilon$\nnsp.\oss 
Let\dss us\sss choose\sss an arbitrary\sss representative\sss $e$\sss
of\trs this orbit\sss
and an arbitrary\sss
element\sss $s_{\dff e}\qff \in\qff G$\sss
such\sss that\sss
$s_{\dff e}\dff(\dff v\trf)\off =\off \ttt(\dff e\trf)$\nsp.\oss 
Then\sss we\sss take\sss the edge\sss
$a\off =\off s_{\dff e}^{\dff -\dff 1}\dff(\dff e\trf)$\sss
as\sss the representative of\trs the orbit\dss
$\iota\dff(\dff \varepsilon\dff)$\nnsp,\oss
and\dss set\trs
$s_{\dff a}\off =\off s_{\dff e}^{\dff -\dff 1}$\nnsp.\oss
Since\dss 
$s_{\dff a}\dff(\dff v\trf)
\off =\off
s_{\dff e}^{\dff -\dff 1}\dff(\dff v\trf)
\off =\off 
\ttt(\dff a\trf)$\nnsp,\oss
the element\dss $s_{\dff a}$\sss satisfies\sss 
the condition\sss from\dss Section\qss \ref{one}.\oss
Also,\pss \vspace{1.5pt}
\[
\quad
s_{\dff a}^{\dff -\dff 1}\dff(\dff a\trf)
\off =\off
s_{\dff e}\dff(\trf s_{\dff e}^{\dff -\dff 1}\dff(\dff e\trf)\trf)
\off =\off
e
\]

\vspace{-12pt}\vspace{1.5pt}
and\dss $s_{\dff a}^{\dff -\dff 1}\off =\off s_{\dff e}$\nsp.\oss
Clearly\halfff,\oss the set\trs $E_{\dff 0}$\sss of\dss selected\sss
representatives of\trs $G_{\dff v}$\dnsp-orbits\sss and\dss the
elements\dss $s_{\dff e}$\nsp,\dss $e\qff \in\qff E_{\dff 0}$\dss
satisfy\dss the conditions\qss ({\fff}i{\fff})\qss and\qss ({\fff}ii{\fff}),\oss 
as also\sss the condition\qss ({\fff}iii{\fff})\qss 
for\dss $u\off =\off 1$\nnsp.\oss

Let\dss us\dss define\sss $s_{\dff d}$\sss for\qss
$d\qff \not\in\qff E_{\dff 0}$\nnsp.\oss
Since every\qss 
$d\qff \in\qff E$\qss 
can\sss be uniquely\dss written in\dss the form\dss 
$d\off =\off u\dff(\dff e\trf)$\dss 
with\dss $e\qff \in\qff E_{\dff 0}$\dss and\dss $u\qff \in\qff \mathcal{T}_{\dff e}$\nsp,\oss 
we can simply\sss set\qss\vspace{1.5pt}
\[
\quad
s_{\dff d}
\off =\off
u\dff s_{\dff e}\dff u^{\dff -\dff 1}
\qff.
\]

\vspace{-12pt}\vspace{1.5pt}
We need\dss to verify\dss that\dss this\dss is\dss a\sss legitimate choice
of\sss $s_{\dff d}$\nsp,\oss 
i.e.\qss that\sss 
$s_{\dff d}\dff(\dff v\trf)\off =\off \ttt(\dff d\trf)$\nnsp.\oss 
In order\dss to prove\sss this,\oss note\sss that\trs
$u\dff(\dff v\trf)\off =\off v$\dss
because\dss
$u\qff \in\qff \mathcal{T}_{\dff e}\qff \subset\qff H$\nnsp,\oss
and\dss that\trs
$s_{\dff e}\dff(\dff v\trf)\off =\off \ttt(\dff e\trf)$\dss
by\dss the choice of\sss $s_{\dff e}$\nsp.\oss
Also,\pss 
$d\off =\off u\dff(\dff e\trf)$\dss implies\sss that\trs 
$u\dff(\trf \ttt(\dff e\trf)\trf)\off =\off \ttt(\dff d\trf)$\nnsp.\oss
Hence\vspace{1.5pt}
\[
\quad
s_{\dff d}\dff(\dff v\trf)
\off =\off 
u\dff s_{\dff e}\dff u\dff ^{\dff -\dff 1}\dff(\dff v\trf)
\off =\off 
u\dff s_e\dff(\dff v\trf)
\off =\off 
u\trf(\trf \ttt(\dff e\trf)\trf)
\off =\off
\ttt(\dff d\trf)
\qff. 
\]

\vspace{-12pt}\vspace{1.5pt}
Therefore\dss 
$s_{\dff d}
\off =\off
u\dff s_{\dff e}\dff u^{\dff -\dff 1}$\dss is\dss indeed\sss 
a\dss legitimate choice\sss of\dss $s_{\dff d}$\nsp.\oss  \eproof

\mypar{Theorem.}{second-simplification}
\emph{Suppose\sss that\dss we are working\trs with a regular scaffolding\sss and\dss
the assumptions of\qss Theorem\qss \ref{simply-connected}\qss hold\dss
for a collection of\pss loops\sss $\mathcal{L}$\nnsp.\oss
For every\dss $e\qff \in\qff E_{\dff 1}$\dss let\trs $\mathcal{H}_{\dff e}$\dss
be a set\sss of\qss generators of\pss $G_{\dff e}$\nsp.\oss
Then\dss the group\dss $G$\dss can\dss be obtained\dss from\dss $G_{\dff v}$\dss
by\sss adding a\sss generator\dss $g_{\dff e}$\dss for every\dss
$e\qff \in\qff E_{\dff 1}$\dss
and\dss the following\dss relations.\oss}\vspace{-9pt}
\begin{itemize}

\item[$\phantom{E}(\trf E\trf)$]\quad\
\emph{The edge relations\pss 
$E\dff(\dff e\fff,\pff t\trf)$\qss for\qss $e\qff \in\qff E_{\dff 1}$\dss
and\dss $t\qff \in\qff \mathcal{H}_{\dff e}$\nsp.\oss} 

\item[$(\trf EL\trf)$]\quad\
\emph{The\sss edge-loop relations\pss $L\trf(\trf l_{\dff e}\trf)$\qss 
for edges\dss $e\qff \in\qff E_{\dff 1}$ admitting an inversion.\oss}

\item[$\phantom{E}(\trf L\trf)$]\quad\
\emph{The\qss loop relations\pss
$L\dff(\trf l\qff)$\qss for\qss $l\qff \in\qff \mathcal{L}$\nnsp,\oss
rewritten\sss in\dss
terms of\qss $g_{\dff e}$\dss with\dss
$e\qff \in\qff E_{\dff 1}$\nsp.\oss}

\end{itemize}

\vspace{-6pt}
\proof
Let\sss us consider\sss first\sss the\sss statement\sss resulting\sss from\sss 
replacing\sss $E_{\dff 1}$ by\sss $E_{\dff 0}$ in\sss the\sss theorem.\oss
This statement\sss follows\sss from\dss Corollary\qss \ref{isomorphism-complexes}\qss
together\sss with\sss the procedure of\dss turning edge relations into definitions
explained at\sss the beginning of\dss this section.\oss
Next,\oss let\sss us\sss turn\sss the edge-loop relations $L\trf(\trf l_{\dff e}\trf)$\sss
for edges $e\qff \in\qff E_{\dff 1}$ not\sss admitting an\sss inversion\sss
into definitions.\oss
Then\sss the edge-loop relation $L\trf(\trf l_{\dff e}\trf)$\sss for such an edge $e$ holds
by\sss the definition,\oss and,\oss as we saw,\oss 
implies\sss that\sss 
the relations $L\trf(\trf l_{\dff a}\trf)$ and\sss $E\dff(\dff a\fff,\pff r\trf)$\sss
for\sss 
$a\off =\off \iota_{\dff 0}\dff(\dff e\trf)$
and arbitrary\sss $r\qff \in\qff G_{\dff v}$\sss
also hold.\oss
Therefore after\sss this\dss it\dss is\dss sufficient\sss to keep\sss
the edge-loop relations $L\trf(\trf l_{\dff e}\trf)$\sss
only\sss for edges $e\qff \in\qff E_{\dff 1}$ admitting an\sss inversion.\oss
The\sss theorem\sss follows.\oss  \eproof

\newpage
\mysection{Examples}{examples}

\vspace{1.75pt}
\myuppar{Fundamental\dss groups of\pss CW-complexes.}
The\sss first\sss example relates\sss the above\sss theory\sss 
with\dss the standard\sss way\sss of\dss finding\sss presentations of\trs
the fundamental\dss groups of\qss CW-complexes.\oss
Let\dss $K$\dss be a connected\dss CW-complex.\oss
There\dss is\dss a\sss well\dss known standard\sss procedure replacing\dss $K$\dss
by\sss a homotopy\sss equivalent\trs CW-complex\sss $K'$ such\dss that\sss
$K'$ has only\sss one $0$\dnsp-cell\sss and\dss 
the\sss loops defined\dss by $1$\dnsp-cells are
not\sss contractible in\dss $K'$\nnsp.\oss 
We will\sss assume\sss that\sss the\trs CW-complex\sss $K$ already\dss has\sss
these properties.\oss\vspace{2pt}

Let\dss $x\qff \in\qff K$\dss be\sss the only\sss $0$\dnsp-cell\sss of\trs $K$\nnsp,\pss
and\trs let\dss $G\off =\off \pi_{\trf 1}\dff (\trf K\fff,\qff x\trf)$\nnsp.\oss
The fundamental\dss group\dss $G$\dss acts on\dss the universal\sss cover\dss
$K^{\fff \sim}$\sss of\trs $K$\nnsp.\oss
Let\dss $X$\dss be\sss the $1$\dnsp-skeleton of\trs $K$\nnsp.\oss
Since\sss the loops defined\dss by\dss the $1$\dnsp-cells of\trs $K$\dss
are not\sss contractible,\oss every\sss $1$\dnsp-cell\sss of\trs $X$\dss
connects\sss two different $0$\dnsp-cells.\oss
This allows\sss to consider\dss $X$\dss as a $1$\dnsp-dimensional\sss
simplicial\sss complex.\oss
The group\dss $G$\dss acts on\dss this complex,\oss
and one can use\sss the\sss theory\sss of\qss 
Sections\qss \ref{one}\dss --\dss \ref{presentations}\qss 
to construct\sss a\sss presentation of\trs $G$\nnsp.\vspace{2pt}

Since\dss $K$\dss has only\sss one $0$\dnsp-cell,\oss the action of\trs $G$\dss
is\dss transitive on\dss the set\sss of\dss vertices of\trs $X$\nnsp.\oss
Let\sss $v\qff \in\qff  X$\dss be a vertex\sss in\dss the preimage\sss of\sss $x$\dss
end\trs $E$\dss be\sss the set\sss of\dss edges having\sss $v$\sss as one of\trs
the endpoints.\oss
The group\dss $G$\dss acts freely\sss on\dss $K^{\fff \sim}$\sss and\dss hence on\dss $X$\nnsp.\oss
Therefore\sss the stabilizer\dss $G_{\dff v}$\dss of\dss $v$\dss
and\dss the stabilizers\dss $G_{\dff e}$\dss of\dss edges\dss $e\qff \in\qff E$\dss
are\dss trivial.\oss
In\sss particular\halfff,\oss each set\dss $\mathcal{T}_{\fff e}$\dss
of\dss representatives of\dss cosets\dss $G_{\dff v}/G_{\dff e}$\dss 
is\dss equal\dss to\dss $\{\trf 1\qff\}$\dss  
and\dss there are no edge relations.\oss
Since\sss the action of\trs $G$\dss on\dss $K^{\fff \sim}$\sss is\dss free,\oss
no edge\dss $e\qff \in\qff E$\dss admits an\sss inversion and\dss
for every\dss $e\qff \in\qff E$\dss 
there\dss is\dss a\qss \emph{unique}\qss element\sss $s_{\dff e}$\sss 
such\dss that\dss $s_{\dff e}\dff(\dff v\trf)$\dss 
is\dss the endpoint\sss of\dss $e$\sss other\dss than\sss $v$\nnsp.\vspace{2pt}

The edges\dss $e\qff \in\qff E$\dss correspond\dss to\sss the\dss lifts of\sss
$1$\dnsp-cells of\dss $K$\dss containing\sss $v$\sss in\dss their\sss boundary\halfff.\oss
One can\dss trace a $1$\dnsp-cell\dss in\dss two different\sss directions,\oss
and\dss this\sss leads\sss to\sss two such\dss lifts.\oss
Tracing\sss a $1$\dnsp-cell\sss in one of\trs the directions defines
an element\sss of\trs the fundamental\dss group\dss 
$G$\nnsp,\oss
and\dss by\dss the\sss theory\sss of\dss covering spaces\sss this element\dss
is\dss nothing\sss else but $s_{\dff e}$\nsp,\oss
where\sss $e$\sss is\dss the corresponding\sss edge.\oss
Since\sss the elements of\trs 
$G$\dss 
defined\dss by\dss tracing\sss
a $1$\dnsp-cell\dss in\dss two directions are\sss the inverses of\dss each other\halfff,\pss
$s_{\dff a}\off =\off s_{\dff e}^{\dff -\dff 1}$\dss
for\sss the edges\dss $a\fff,\pff e$\dss corresponding\dss to a $1$\dnsp-cell.\oss
It\dss follows\dss that\trs
$E_{\dff 0}\off =\off E$\dss together\dss with\dss the families\dss
$\mathcal{T}_{\fff e}\off =\off \{\trf 1\qff\}$\dss
and\dss
$s_{\dff e}$\nsp,\dss $e\qff \in\qff E$\dss is\dss a regular\sss scaffolding.\oss
We will\dss work\sss with\dss this\sss scaffolding.\oss\vspace{2pt}

Every\dss $G$\dnsp-orbit\sss of\trs $2$\dnsp-cells\sss of\trs $K^{\fff \sim}$\sss
contains a\sss $2$\dnsp-cell\sss such\dss that\sss $v$\sss belongs\sss to its boundary\halfff.\oss
By\sss starting\sss at\sss $v$\sss and\dss following\dss the boundary\sss of\dss
such a sell\dss we get\sss a\sss loop in\dss $X$\nnsp.\oss
Let\dss $\mathcal{L}$\dss be a set\sss of\dss such\dss loops,\oss
one for each\dss $G$\dnsp-orbit\halfff.\oss
Since\sss the result\sss of\dss glueing\dss the $2$\dnsp-cells of\trs $K^{\fff \sim}$\sss
to\dss $X$\dss is\dss simply-connected,\pss
$\mathcal{L}$\dss satisfies\sss the assumptions of\qss Theorem\qss \ref{simply-connected}.\oss
Therefore one can use\trs Theorem\qss \ref{second-simplification}\qss
together\sss with\dss its\sss complement\dss to construct\sss a\sss presentation of\trs $G$\nnsp.\oss
The resulting\dss presentation\sss has one\sss generator\dss for every\dss $1$\dnsp-cell\sss
of\trs $K$\dss and one relation\sss for every\dss $2$\dnsp-cell\sss of\trs $K$\nnsp.\oss
A\sss trivial\sss verification shows\sss that\dss this\dss
presentation\dss is\dss nothing else\sss
but\dss the standard\sss presentation of\dss
$G\off =\off \pi_{\trf 1}\dff (\trf K\fff,\qff x\trf)$\dss
defined\dss by\dss the\trs CW-structure of\trs $K$\nnsp.\oss

\myuppar{Symmetric\sss groups.}
Let\dss $\Sigma_{\dff n}$\sss be\sss the group of\dss permutations of\trs
the set\trs
$I_{\dff n}
\off =\off
\{\qff 1\fff,\pff 2\fff,\pff \ldots\fff,\pff n\qff\}$\nnsp.\oss
As usual,\oss we will\sss denote\sss by\dss $(\trf i\fff,\qff j\trf)$\dss
the\sss transposition of\trs two different\sss elements\dss
$i\fff,\qff j\qff \in\pff I_{\dff n}$\nsp,\oss
i.e.\qss the map\dss
$I_{\dff n}\qff \ttoo\qff I_{\dff n}$\dss interchanging\sss $i$\sss and\sss $j$\sss
and\dss fixing\sss all\sss other elements of\trs $I_{\dff n}$\nsp.\oss
The standard\sss presentation of\trs $\Sigma_{\dff n}$\sss
has\dss generators\dss
$\sigma_{\dff 1}\dff,\off \sigma_{\dff 2}\dff,\off \ldots\dff,\off \sigma_{\dff n\dff -\dff 1}$\dss
corresponding\halfff,\oss respectively\halfff,\oss to\sss the\sss transpostions\dss
$(\trf 1\fff,\qff 2\trf)\fff,\off (\trf 2\fff,\qff 3\trf)\dff,\off \ldots\dff,\off
(\trf n\qff -\qff 1\fff,\qff n\trf)$\nnsp.\oss
The relations are\sss the following\halfff.\oss\vspace{3pt}
\[
\quad
\sigma_{\dff i}^{\dff 2}
\off =\off
1
\hspace{1.0em}\mbox{for every}\hspace{0.8em}
i\qff \leq\qff n\qff -\qff 1
\qff,
\]

\vspace{-36pt}
\[
\quad
\sigma_{\dff i}\qff \sigma_{\dff i\dff +\dff 1}\qff \sigma_{\dff i}
\off =\off
\sigma_{\dff i\dff +\dff 1}\qff \sigma_{\dff i}\qff \sigma_{\dff i\dff +\dff 1}
\hspace{1.0em}\mbox{for every}\hspace{0.8em}
i\qff \leq\qff n\qff -\qff 2
\qff,
\hspace{1.0em}\mbox{and}
\]

\vspace{-36pt}
\[
\quad
\sigma_{\dff i}\qff \sigma_{j}
\off =\off
\sigma_{j}\qff \sigma_{\dff i}
\hspace{1.0em}\mbox{for}\hspace{1.0em}
\num{i\qff -\qff j}\qff \geq\qff 2
\qff.
\]

\vspace{-9pt}
We will\sss show\dss how\dss the\sss theory\sss of\qss
Sections\qss \ref{one}\dss --\dss \ref{presentations}\qss
naturally\dss leads\sss to\sss this presentation.\oss

Using an\sss induction\sss by\sss $n$\nnsp,\oss
we can assume\sss that\dss the group\dss $\Sigma_{\dff n\dff -\dff 1}$\dss
admits such a presentation.\oss
Let\dss us\sss take as\dss $X$\dss the $1$\dnsp-skeleton of\trs the
$(\fff n\dff -\dff 1\fff)$\dnsp-dimensional\sss simplex\sss with\sss the vertices\dss
$1\fff,\pff 2\fff,\pff \ldots\fff,\pff n$\nnsp,\oss
and\dss let\dss us\sss take as $v$\sss the vertex $1$\nnsp.\oss
The\sss group\dss $G\off =\off \Sigma_{\dff n}$\dss naturally\sss acts on\dss $X$\nnsp,\oss
and\dss the stabilizer\dss $G_{\dff 1}$\dss of\dss $v\off =\off 1$\dss is\dss the group\dss
of\dss permutations of\trs
the set\trs
$\{\qff 2\fff,\pff 3\fff,\pff \ldots\fff,\pff n\qff\}$\nnsp.\oss
The\sss latter\dss group\dss is\dss canonically\dss isomorphic\sss to\dss $\Sigma_{\dff n\dff -\dff 1}$\dss
and\dss hence has a presentation\sss with\sss generators\dss
$\sigma_{\dff 2}\dff,\off \sigma_{\dff 3}\dff,\off \ldots\dff,\off \sigma_{\dff n\dff -\dff 1}$\dss
corresponding\dss the\sss same\sss transpositions and subject\dss to\sss
the same relations as in\dss the standard\sss presentation of\trs $\Sigma_{\dff n}$\nsp.\oss

The edges of\trs $X$\dss having $v$ as an endpoint\sss are\sss
the edges\dss
$e_{\dff 2}\dff,\off e_{\dff 3}\dff,\off \ldots\dff,\off e_{\dff n}$\dss
connecting\dss $1$\sss with\dss
$2\fff,\pff 3\fff,\pff \ldots\fff,\pff n$\dss respectively\halfff.\oss
Let\dss $E$\dss be\sss the set\sss of\trs these edges.\oss
The group\sss $G_{\dff 1}$\sss acts\sss transitively\sss on\dss $E$\dss
and\dss hence we can\dss take\dss 
$E_{\dff 0}\off =\off \{\trf e_{\dff 2}\trf\}$\sss 
is\dss a set\sss of\dss representatives of\trs
$G_{\dff 1}$\dnsp-orbits.\oss
Let\dss $e\off =\off e_{\dff 2}$\dss and\dss let\dss
$G_{\dff 1\fff 2}$\dss be\sss the stabilizer of\dss $e$\sss in\dss $G_{\dff 1}$\nsp.\oss
Clearly\halfff,\pss $G_{\dff 1\fff 2}$\dss is\dss the group\dss
of\dss permutations of\trs
the set\trs
$\{\qff 3\fff,\pff \ldots\fff,\pff n\qff\}$\nnsp.\oss
The\sss latter\dss group\dss is\dss canonically\dss isomorphic\sss to\dss 
$\Sigma_{\dff n\dff -\dff 2}$\dss
and\dss is\dss generated\dss by\dss the set\trs
$\mathcal{H}_{\dff 1\fff 2}
\off =\off
\{\qff \sigma_{\dff 3}\dff,\off \ldots\dff,\off \sigma_{\dff n\dff -\dff 1} \qff\}$\nnsp,\oss
where\dss
$\sigma_{\dff i}\off =\off (\trf i\fff,\qff i\qff +\qff 1\trf)$\nnsp.\oss
The\sss transpostions\dss
$(\trf 2\fff,\qff i\trf)$\dss
for\dss $3\qff \leq\qff i\qff \leq\qff n$\dss
together with\sss $1$\sss
form\sss a\sss set\dss $\mathcal{T}_{\dff e}$\dss 
of\trs representatives of\dss cosets\dss in\sss $G_{\dff 1}/G_{\dff 1\fff 2}$\nsp.\oss

For\dss $2\qff \leq\qff i\qff \leq\qff n\qff -\qff 1$\dss
let\dss $s_{\dff i}\off =\off (\trf 1\fff,\qff i\trf)$\nnsp.\oss
Then\dss 
$s_{\dff i}\dff (\dff v\trf)
\off =\off
s_{\dff i}\dff(\dff 1\trf)
\off =\off
i$\nnsp,\oss
i.e.\qss $s_{\dff i}\dff (\dff v\trf)$\sss 
is\dss the endpoint\sss of\dss $e_{\dff i}$\dss different\dss from\dss $v$\dss
Hence we may\sss set\dss $s_{\dff e_{\dff i}}\off =\off s_{\dff i}$\nsp.\oss 
Clearly\halfff,\pss\vspace{3pt}
\[
\quad
(\trf 1\fff,\qff i\trf)
\off =\off
(\trf 2\fff,\qff i\trf)\dff \cdot\dff
(\trf 1\fff,\qff 2\trf)\dff \cdot\dff
(\trf 2\fff,\qff i\trf)^{\dff -\dff 1}\qss
\]

\vspace{-12pt}
and\dss hence\vspace{1pt}
\begin{equation}
\label{s-transpositions}
\quad
s_{\dff i}
\off =\off
(\trf 2\fff,\qff i\trf)\dff \cdot\dff
s_{\dff 2}\dff \cdot\dff
(\trf 2\fff,\qff i\trf)^{\dff -\dff 1}
\qff.
\end{equation}

\vspace{-9pt}
Clearly\halfff,\pss $s_{\dff 2}\off =\off (\trf 1\fff,\qff 2\trf)$\dss
is\dss an\sss inversion of\trs the edge\dss $e\off =\off e_{\dff 2}$\nnsp.\oss
Together\dss with\qss (\ref{s-transpositions})\qss this implies\sss that\dss
the set\trs
$E_{\dff 0}\off =\off \{\trf e_{\dff 2}\trf\}$\nnsp,\oss
the set\trs
$\mathcal{T}_{\dff e}
\off =\off
\{\pff
(\trf 2\fff,\qff i\trf)
\qff \mid\pff
3\qff \leq\qff i\qff \leq\qff n
\pff\}$\dss
of\dss representatives of\dss cosets\sss in\dss 
$G_{\dff 1}/G_{\dff 1\fff 2}$\nsp,\oss
and\dss the\sss family\sss of\dss elements\dss
$s_{\dff e_{\dff i}}\off =\off s_{\dff i}$\dss form\sss a\sss 
regular\sss scaffolding.\oss

The $2$\dnsp-skeleton of\trs the
$(\fff n\dff -\dff 1\fff)$\dnsp-dimensional\sss simplex\dss is\dss
simply-con\-nect\-ed,\oss and\dss $G$\dss acts
on\dss the set\sss of\trs $2$\dnsp-simplices\sss transitively\halfff.\oss
Therefore\sss the assumptions of\qss Theorem\qss \ref{simply-connected}\qss hold\dss
for\dss the set\trs
$\mathcal{L}$\dss
consisting of\trs the single\sss loop\dss
$1\fff,\pff 2\fff,\pff 3\fff,\pff 1$\nnsp.\oss
Now\sss we are ready\dss to
apply\qss Theorem\qss \ref{second-simplification}.\oss

By\qss Theorem\qss \ref{second-simplification}\qss the group $G$
can\dss be obtained\dss from $G_{\dff 1}$ by\sss adding\sss
one generator\sss $g_{\trf 2}\off =\off g_{\dff e_{\dff 2}}$\nsp,\oss which\sss we will\dss
also denote\sss by $\sigma_{\dff 1}$\nsp,\oss
the edge relations\pss 
$E\dff(\dff e_{\trf 2}\fff,\pff \sigma_{\dff i}\trf)$\qss for\dss 
$i\qff \geq\qff 3$\nnsp,\oss
the edge-loop relation\dss $L\trf(\trf l_{\dff e_{\trf 2}}\trf)$\nnsp,\oss
and\dss the\sss loop relation corresponding\dss to\sss the\sss loop\dss
$1\fff,\pff 2\fff,\pff 3\fff,\pff 1$\nnsp.\oss
The following\dss three\sss lemmas show\dss that\dss these relations
are nothing else\sss but\dss the standard\dss relations of\trs $\Sigma_{\dff n}$\dss
involving\dss $\sigma_{\dff 1}$\nsp.\oss
Since\sss the other standard\dss relations of\trs $\Sigma_{\dff n}$\dss
are\sss the relations of\trs $G_{\dff 1}$\nsp,\oss
this will\sss complete\sss the induction step and\dss prove\sss
that\dss the standard\dss presentation described above\dss is\dss
indeed a\sss presentation of\trs the\sss group\dss $\Sigma_{\dff n}$\nsp.\oss

\mypar{Lemma.}{symmetric-edge}
\emph{The edge relation\dss
$E\dff(\dff e_{\trf 2}\fff,\pff \sigma_{\dff i}\trf)$\dss
with\dss
$i\qff \geq\qff 3$\dss is\dss equivalent\dss to\dss
$\sigma_{\dff 1}\dff \cdot\dff \sigma_{\dff i}
\off =\off
\sigma_{\dff i}\dff \cdot\dff \sigma_{\dff 1}$\nsp.\oss}

\proof
The edge relation\dss 
$E\dff(\dff e_{\dff 2}\fff,\pff \sigma_{\dff i}\trf)$\dss
is\dss the relation\vspace{3pt}
\[
\quad
g_{\dff k}^{\dff -\dff 1}\dff \cdot\dff \sigma_{\dff i}\dff \cdot\dff g_{\trf 2}
\off =\off
s_{\dff k}^{\dff -\dff 1}\dff \cdot\dff \sigma_{\dff i}\dff \cdot\dff s_{\trf 2}
\pff,
\]

\vspace{-9pt}
where $k\off =\off \sigma_{\dff i}\trf(\trf 2\trf)$\nnsp.\oss
If\qss $i\qff \geq\qff 3$\nnsp,\oss
then\dss
$\sigma_{\dff i}\trf(\trf 2\trf)\off =\off 2$\dss and\dss hence\dss
$E\dff(\dff e_{\trf 2}\fff,\pff \sigma_{\dff i}\trf)$\dss 
actually\dss is\vspace{3pt}
\[
\quad
g_{\trf 2}^{\dff -\dff 1}\dff \cdot\dff \sigma_{\dff i}\dff \cdot\dff g_{\trf 2}
\off =\off
s_{\trf 2}^{\dff -\dff 1}\dff \cdot\dff \sigma_{\dff i}\dff \cdot\dff s_{\trf 2}
\pff.
\]

\vspace{-9pt}
Also,\oss if\qss $i\qff \geq\qff 3$\nnsp,\oss
then\dss
$s_{\trf 2}^{\dff -\dff 1}\dff \cdot\dff \sigma_{\dff i}\dff \cdot\dff s_{\trf 2}
\off =\off
\sigma_{\dff i}$\dss
and\dss hence\dss
$E\dff(\dff e_{\trf 2}\fff,\pff \sigma_{\dff i}\trf)$\dss
means\sss that\vspace{3pt}
\[
\quad
g_{\trf 2}^{\dff -\dff 1}\dff \cdot\dff \sigma_{\dff i}\dff \cdot\dff g_{\trf 2}
\off =\off
\sigma_{\dff i}
\pff.
\]

\vspace{-9pt}
This\dss is\dss equivalent\dss to\qss
$g_{\trf 2}\dff \cdot\qff \sigma_{\dff i}
\off =\off
\sigma_{\dff i}\dff \cdot\qff g_{\trf 2}$\nsp,\oss
and,\oss in\dss terms of\dss
$\sigma_{\dff 1}$\nsp,\qff\oss 
to\dss 
$\sigma_{\dff 1}\dff \cdot\qff \sigma_{\dff i}
\off =\off
\sigma_{\dff i}\dff \cdot\qff \sigma_{\dff 1}$\nsp.\oss  \eproof

\mypar{Lemma.}{symmetric-edge-loop}
\emph{The edge-loop relation\qss
$L\trf(\trf l_{\dff e_{\trf 2}}\trf)$\dss
is\dss equivalent\dss to\dss
$\sigma_{\dff 1}^{\dff 2}
\off =\off
1$\nnsp.\oss}

\proof
Since\sss $s_{\trf 2}$\dss is\dss an\sss inversion of\dss $e_{\trf 2}$\dss
and\dss $s_{\trf 2}\dff \cdot\dff s_{\trf 2}\off =\off 1$\nnsp,\oss
the edge-loop relation\dss $L\trf(\trf l_{\dff e_{\trf 2}}\trf)$\dss has\sss the form\dss
$g_{\trf 2}\dff \cdot\dff g_{\trf 2}\off =\off 1$\nnsp,\oss
and,\oss in\dss terms of\dss
$\sigma_{\dff 1}$\nsp,\qff\oss 
the form\dss 
$\sigma_{\dff 1}\dff \cdot\qff \sigma_{\dff 1}
\off =\off
1$\nnsp.\oss  \eproof

\mypar{Lemma.}{symmetric-loop}
\emph{If\qss the relation\qss 
$L\trf(\trf l_{\dff e_{\trf 2}}\trf)$\dss
holds,\oss
then\dss the\dss loop relation\dss
corresponding\dss to\sss the\sss loop\dss
$1\fff,\pff 2\fff,\pff 3\fff,\pff 1$\dss
is\dss equivalent\dss to\qss 
$\sigma_{\trf 2}\dff \cdot\qff
\sigma_{\dff 1}\dff \cdot\qff
\sigma_{\trf 2}
\off =\off
\sigma_{\dff 1}\dff \cdot\qff
\sigma_{\trf 2}\dff \cdot\qff 
\sigma_{\dff 1}$\nsp.\oss}

\proof
As we will\sss see,\oss
in\dss its\sss original\dss form\dss this relations\sss
involves not\sss only\dss $\sigma_{\dff 1}\off =\off g_{\trf 2}$\nsp,\oss
but\sss also\sss the generator\dss 
$g_{\dff 3}\off =\off g_{\dff e_{\dff 3}}$\nsp.\oss
In view of\qss (\ref{s-transpositions})\qss the element\dss $g_{\qff 3}$\dss
is\dss defined as follows\qss\vspace{3pt}\vspace{-0.5pt}
\[
\quad
g_{\qff 3}
\off =\off
(\trf 2\fff,\qff 3\trf)
\dff \cdot\dff
g_{\trf 2}\dff \cdot\dff
(\trf 2\fff,\qff 3\trf)^{\dff -\dff 1}
\off =\off
(\trf 2\fff,\qff 3\trf)
\dff \cdot\dff
g_{\trf 2}\dff \cdot\dff
(\trf 2\fff,\qff 3\trf)
\pff
\]

\vspace{-9pt}\vspace{-0.5pt}
(see remarks preceding\trs
Lemma\qss \ref{coherent-frames}).\oss
Since\dss $\sigma_{\trf 2}\off =\off (\trf 2\fff,\qff 3\trf)$\nnsp,\oss
this means\sss that\trs\vspace{3pt}\vspace{-0.5pt}
\begin{equation}
\label{g3}
\quad
g_{\qff 3}
\off =\off
\sigma_{\trf 2}\dff \cdot\dff
g_{\trf 2}\dff \cdot\dff
\sigma_{\trf 2}
\pff.
\end{equation}

\vspace{-9pt}\vspace{-0.5pt}
Let\dss us\dss turn\dss to\sss the our\dss loop relation.\oss
An easy\dss verification shows\sss that\vspace{3pt}\vspace{-0.5pt}
\[
\hspace{0.5em}\begin{array}{l}
s_{\trf 2}\trf(\dff 1\trf)
\off =\off 
2\qff,\vspace{12pt}\\
s_{\trf 2}\dff \cdot\dff s_{\trf 3}\trf(\dff 1\trf)
\off =\off 
3\qff,\vspace{12pt}\\
s_{\trf 2}\dff \cdot\dff s_{\trf 3}\dff \cdot\dff s_{\trf 2}\trf(\dff 1\trf)
\off =\off 
1\qff,
\end{array}
\]

\vspace{-9pt}\vspace{-0.5pt}
and\dss hence\sss the\sss loop relation\sss corresponding\dss to\sss
the\dss loop\dss
$1\fff,\pff 2\fff,\pff 3\fff,\pff 1$\qss is\qss\vspace{3pt}\vspace{-0.5pt}
\[
\quad
g_{\trf 2}
\dff \cdot\dff 
g_{\trf 3}
\dff \cdot\dff 
g_{\trf 2}
\off =\off
s_{\trf 2}
\dff \cdot\dff 
s_{\trf 3}
\dff \cdot\dff 
s_{\trf 2}
\pff
\]

\vspace{-9pt}\vspace{-0.5pt}
(see\dss Section\qss \ref{one}).\oss
Another\sss easy\dss verification shows\sss that\vspace{3pt}\vspace{-0.5pt}
\[
\quad
s_{\trf 2}
\dff \cdot\dff 
s_{\trf 3}
\dff \cdot\dff 
s_{\trf 2}
\off =\off
(\trf 1\fff,\qff 2\trf)\qff
(\trf 1\fff,\qff 3\trf)\qff
(\trf 1\fff,\qff 2\trf)
\off =\dff\off
(\trf 2\fff,\qff 3\trf)
\off =\off
\sigma_{\trf 2}
\]

\vspace{-9pt}\vspace{-0.5pt}
and\dss hence our\sss relation\sss means\sss that\vspace{3pt}\vspace{-0.5pt}
\[
\quad
g_{\trf 2}
\dff \cdot\dff 
g_{\trf 3}
\dff \cdot\dff 
g_{\trf 2}
\off =\off
\sigma_{\trf 2}
\pff.
\]

\vspace{-9pt}\vspace{-0.5pt}
Using\qss (\ref{g3})\qss and\dss the notation\dss $\sigma_{\dff 1}$\dss for\dss $g_{\trf 2}$\dss
turns\sss this relation\sss into\vspace{3pt}\vspace{-0.5pt}
\[
\quad
\sigma_{\dff 1}\dff \cdot\dff 
(\qff \sigma_{\trf 2}\dff \cdot\dff
\sigma_{\dff 1}\dff \cdot\dff
\sigma_{\trf 2}\pff)\dff \cdot\qff 
\sigma_{\dff 1}
\off =\off
\sigma_{\trf 2}
\pff.
\]

\vspace{-9pt}\vspace{-0.5pt}
If\pss $L\trf(\trf l_{\dff e_{\trf 2}}\trf)$\dss holds,\oss
then\dss $\sigma_{\dff 1}^{\dff 2}\off =\off 1$\qss
by\qss Lemma\qss \ref{symmetric-edge-loop}\qss
and\dss hence\sss the\sss last\dss relation\dss is\dss equivalent\dss to\vspace{3pt}\vspace{-0.5pt}
\[
\quad
\sigma_{\trf 2}\dff \cdot\qff
\sigma_{\dff 1}\dff \cdot\qff
\sigma_{\trf 2}
\off =\off
\sigma_{\dff 1}\dff \cdot\qff
\sigma_{\trf 2}\dff \cdot\qff 
\sigma_{\dff 1}
\pff.
\]

\vspace{-9pt}\vspace{-0.5pt}
This completes\sss the proof\halfff.\oss  \eproof\vspace{0.25pt}

\myuppar{Rotations of\dss a\sss regular dodecahedron.}
Let\dss $\mathcal{D}$\dss be\sss the group of\dss 
orientation-preserving\sss symmetries of\dss a regular dodecahedron\dss $D$\nnsp.\oss
We will\sss consider\dss the action of\trs $\mathcal{D}$\dss 
on\dss the graph\dss $X$\dss defined\dss by\dss
the vertices and edges of\trs $D$\sss in an obvious way\halfff.\oss
In\dss fact\halfff,\pss $\mathcal{D}$\dss can be defined\dss 
in a purely\sss combinatorial\dss manner\halfff.\oss
The faces of\trs $D$\dss lead\dss to\sss $12$\sss cycles of\trs the\sss length\sss $5$\sss in\dss $X$\nnsp.\oss
These cycles can be oriented\sss in such a\sss way\dss
that\dss the\sss orientations of\dss every\dss two cycles having a common edge induce opposite
orientations of\trs this edge.\oss
There are\sss two ways\sss to choose such orientations.\oss
Let\dss is\dss fix one of\trs them.\oss
Then\dss $\mathcal{D}$\dss can\sss be defined as\sss the group 
of\dss automorphisms of\trs the graph\dss $X$\dss
preserving\dss the collection of\trs these\sss $12$\sss oriented cycles.\oss

The group\dss $\mathcal{D}$\dss acts\sss transitively\sss 
on\dss the set\sss of\dss vertices of\trs $X$\nnsp.\oss
Let\dss us\dss fix a vertex\dss $v$\dss of\trs $X$\nnsp,\oss
and\dss let\trs 
$E
\off =\off
\{\trf
e_{\dff 1}\dff,\pff e_{\trf 2}\dff,\off e_{\trf 3}
\trf\}$\dss
be\sss the set\sss of\dss edges\sss having $v_{\dff 0}$ as an endpoint\halfff.\oss
Let\trs $w_{\dff 1}\dff,\pff w_{\trf 2}\dff,\off w_{\trf 3}$\dss
be\sss the other endpoints of\trs the edges\trs 
$e_{\dff 1}\dff,\pff e_{\trf 2}\dff,\off e_{\trf 3}$\dss
respectively\halfff.\oss
The stabilizer\sss $G_{\dff v}$\sss of\dss $v$\sss is\dss a\sss cyclic
group of\dss order\sss $3$\sss cyclically\dss permuting\dss the edges\dss
$e_{\dff 1}\dff,\pff e_{\trf 2}\dff,\off e_{\trf 3}$\nsp.\oss
Moreover\halfff,\oss\vspace{3pt}
\begin{equation}
\label{action-h}
\quad
h\trf(\dff e_{\dff 1}\trf)\off =\off e_{\trf 3}\qff,\hspace{1.0em}
h\trf(\dff e_{\trf 3}\trf)\off =\off e_{\trf 2}\qff,
\hspace{1.0em}\mbox{and}\hspace{1.0em}
h\trf(\dff e_{\trf 2}\trf)\off =\off e_{\trf 1}
\end{equation}

\vspace{-9pt}
for a unique\sss $h\qff \in\qff G_{\dff v}$\nsp,\oss
and such a element\dss $h$\dss generates\sss $G_{\dff v}$\nsp.\oss
The group $G_{\dff v}$ acts\sss transitively\sss on\sss $E$
and\dss hence we can\sss take\sss 
$E_{\dff 0}\off =\off \{\trf e_{\dff 1}\qff\}$\sss 
as a set\sss of\dss representatives of\trs
$G_{\dff v}$\dnsp-orbits.\oss
The stabilizer\sss $G_{\dff e}$\sss of\trs the edge\sss
$e\off =\off e_{\dff 1}$\dss in\sss $G_{\dff v}$\sss is\dss trivial\sss and\dss hence
we can\sss take\sss the group\dss $G_{\dff v}$\dss itself\dss as\sss the set\sss
$\mathcal{T}_{\dff e}$\sss of\dss representatives 
of\dss cosets of\sss $G_{\dff e}\off =\off 1$\nnsp.\oss
For each\dss $i\off =\off 1\fff,\pff 2\fff,\pff 3$\dss there\dss is\dss
a unique element\dss $s_{\dff i}\qff \in\qff \mathcal{D}$\trs 
leaving\dss the edge\sss $e_{\dff i}$\sss
invariant\sss and\sss interchanging\dss its endpoints.\oss
Clearly\halfff,\pss 
$s_{\dff i}\dff (\dff v\trf)
\off =\off
s_{\dff i}\dff(\dff v_{\dff 0}\trf)
\off =\off
w_{\dff i}$\dss
and\dss we may\sss set\dss $s_{\dff e_{\dff i}}\off =\off s_{\dff i}$\nsp.\oss
The uniqueness of\trs $s_{\dff i}$\dss implies\sss that\vspace{3pt}
\begin{equation}
\label{s-dodecahedron}
\quad
s_{\trf 2}\off =\off h^{\dff 2}\dff \cdot\dff s_{\dff 1}\dff \cdot\dff h^{\dff -\dff 2}
\hspace{1.2em}\mbox{and}\hspace{1.2em}
s_{\trf 3}\off =\off h\dff \cdot\dff s_{\dff 1}\dff \cdot\dff h^{\dff -\dff 1}
\pff.
\end{equation}

\vspace{-9pt}
It\dss follows\dss that\dss the set\trs 
$E_{\dff 0}\off =\off \{\trf e_{\dff 1}\qff\}$\nnsp,\oss
the set\trs $G_{\dff v}\off =\off \{\qff 1\fff,\pff h\fff,\pff h^{\dff 2}\qff\}$\dss
taken as\dss $\mathcal{T}_{\dff e}$\nnsp,\oss
and\dss the\sss family\sss of\dss elements\dss
$s_{\dff e_{\dff i}}\off =\off s_{\dff i}$\dss 
form\sss a\sss regular\dss scaffolding.\oss
In addition,\oss the rotation\dss $s_{\dff i}\dff \cdot\dff s_{\dff i}$\dss 
is\dss equal\dss to\sss the identity\sss on\dss the edge\sss $e_{\dff i}$\sss and\dss is\dss
orientation-preserving.\oss
Therefore\dss $s_{\dff i}\dff \cdot\dff s_{\dff i}\off =\off 1$\dss
and\dss hence\dss $s_{\dff i}^{\dff -\dff 1}\off =\off s_{\dff i}$\dss
for\dss $i\off =\off 1\fff,\pff 2\fff,\pff 3$\nnsp.\oss

The surface of\trs the dodecahedron\dss $D$\dss is\dss homeomorphic\sss to\sss
the $2$\dnsp-sphere and\dss hence\dss is\dss simply-connected.\oss
The\sss group $\mathcal{D}$ acts on\dss the set\sss of\dss faces of\trs $D$\dss
transitively\halfff,\oss and\dss hence\sss 
the assumptions of\trs Theorem\qss \ref{simply-connected}\qss hold\dss
for\trs
$\mathcal{L}\off =\off \{\trf l\qff\}$\nnsp,\oss
where\sss $l$\trs is\dss a\sss loop
defined\dss by\sss a\sss face of\trs $D$\dss
having\sss $v$\sss as\sss its vertex.\oss 
There are\sss three such\dss faces,\oss and\dss for exactly\sss one of\trs them\dss 
$l$\trs has\sss the form\dss\vspace{3pt}
\[
\quad
v\dff,\off 
w_{\dff 1}\dff,\off
a\dff,\off
b\dff,\off
w_{\dff 2}\dff,\off
v
\]

\vspace{-9pt}
for some vertices\dss $a\fff,\pff b$\trs of\trs $X$\nnsp.\oss
We will\dss assume\sss that\trs $l$\trs is\dss this\dss loop.\oss

Now\sss we are ready\dss to
apply\qss Theorem\qss \ref{second-simplification}.\oss
By\qss Theorem\qss \ref{second-simplification}\qss the group $\mathcal{D}$
can\dss be obtained\dss from $G_{\dff v}$ by\sss adding\sss
one generator\sss 
$g_{\trf 1}
\off =\off 
g_{\dff e_{\dff 1}}$\nsp,\oss 
the edge-loop relation\dss $L\trf(\trf l_{\dff e_{\dff 1}}\trf)$\nnsp,\oss
and\dss the\sss loop relation corresponding\dss to\sss the\sss loop\dss
$l$\nnsp.\oss
There are no edge relations because\dss $G_{\dff e}\off =\off 1$\nnsp.\oss
The edge-loop relation\dss $L\trf(\trf l_{\dff e_{\dff 1}}\trf)$\dss
has\sss the form\dss
$g_{\trf 1}\dff \cdot\dff g_{\trf 1}\off =\off s_{\dff 1}\dff \cdot\dff s_{\dff 1}$\nnsp.\oss
Since\dss $s_{\dff 1}\dff \cdot\dff s_{\dff 1}\off =\off 1$\nnsp,\oss
the relation\dss $L\trf(\trf l_{\dff e_{\dff 1}}\trf)$\dss
is\dss equivalent\dss to\dss 
$g_{\trf 1}\dff \cdot\dff g_{\trf 1}\off =\off 1$\nnsp,\oss
or\halfff,\oss what\trs is\trs the\sss same,\pss $g_{\trf 1}^{\dff 2}\off =\off 1$\nnsp.\oss

The main\dss task\dss is\dss to compute\sss the\sss loop relation\dss $L\dff(\trf l\trf)$\nnsp.\oss
In\dss its\sss original\dss form\dss
the\sss loop relation\dss $L\dff(\trf l\trf)$\dss 
involves not\sss only\dss $g_{\trf 1}$\nsp,\oss
but\sss also\sss the\sss generators\dss
$g_{\trf 2}\off =\off g_{\dff e_{\trf 2}}$\dss
and\qss
$g_{\trf 3}\off =\off g_{\dff e_{\trf 3}}$\nsp.\oss
In\dss view of\pss (\ref{s-dodecahedron})\pss
and\dss remarks preceding\trs Lemma\qss \ref{coherent-frames}\qss they\sss are\sss
redefined\sss as\sss the elements\vspace{3pt}
\begin{equation}
\label{g-dodecahedron}
\quad
g_{\trf 2}\off =\off h^{\dff 2}\dff \cdot\dff g_{\dff 1}\dff \cdot\dff h^{\dff -\dff 2}
\hspace{1.2em}\mbox{and}\hspace{1.2em}
g_{\trf 3}\off =\off h\dff \cdot\dff g_{\dff 1}\dff \cdot\dff h^{\dff -\dff 1}
\pff.
\end{equation}

\begin{figure}[h!]
\hspace*{7em}
\includegraphics[width=0.64\textwidth]{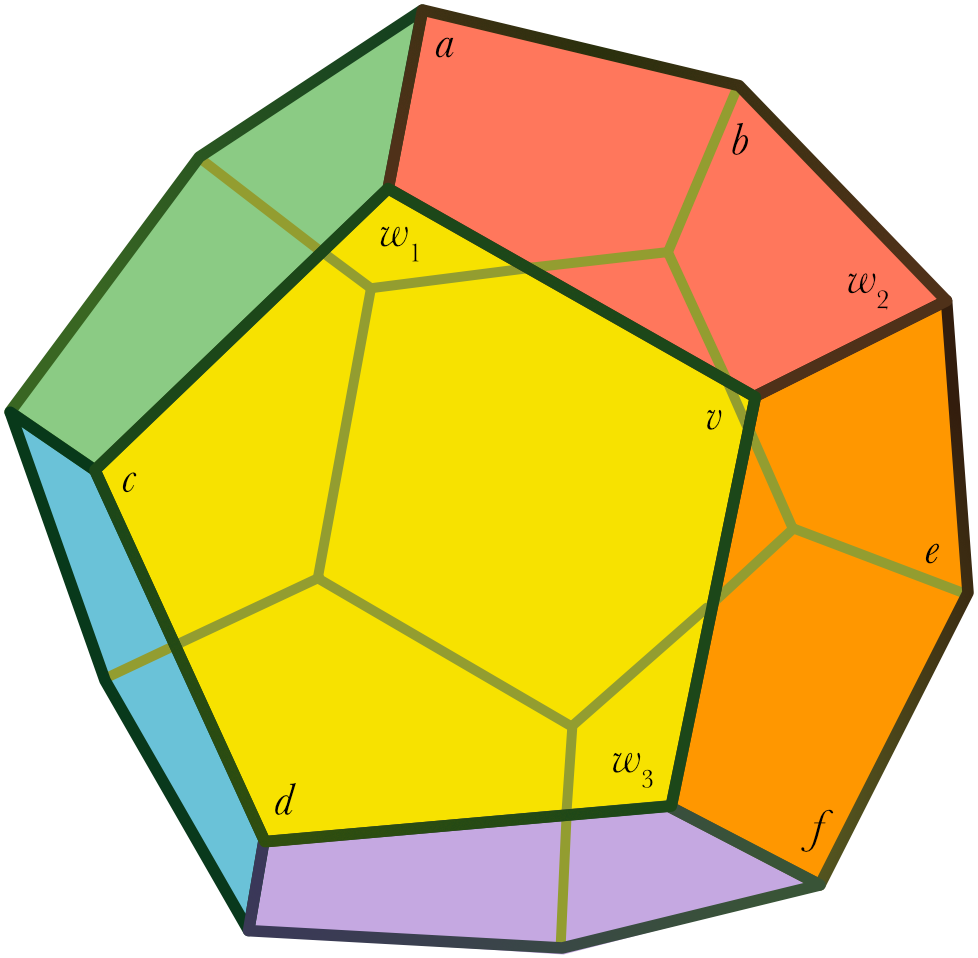}
\end{figure}\vspace{24pt}

Every\dss $s_{\dff i}$\dss flips\sss the edge\dss $v\dff w_{\dff i}$\dss
and\dss interchanges\sss two pentagons\sss adjacent\dss to\sss this edge.\oss
We will\sss assume\sss that\dss the vertices are marked as on\dss
the above picture.\oss 
Then\sss this picture allows\sss to easily\sss
determine\sss the action on\dss the marked vertices.\oss 
For example,\pss $s_{\trf 1}\trf(\dff a\trf)\off =\off w_{\dff 3}$\dss
and\dss $s_{\trf 1}\trf(\dff b\trf)\off =\off d$\nnsp.\oss
We need\dss to find\dss
$r_{\trf 1}\dff,\off
r_{\trf 2}\dff,\off
r_{\trf 3}\dff,\off
r_{\trf 4}\dff,\off
r_{\trf 5}
\off \in\off
\{\trf s_{\dff 1}\dff,\off s_{\dff 2}\dff,\off s_{\dff 3} \qff\}$\dss
such\dss that\vspace{3pt}
\[
\quad
r_{\trf 1}\trf(\dff v\trf)
\off =\off 
w_{\dff 1}
\qff,\quad\
r_{\trf 1}\dff \cdot\dff r_{\trf 2}\trf(\dff v\trf)
\off =\off 
a\qff,\quad\
r_{\trf 1}\dff \cdot\dff r_{\trf 2}\dff \cdot\dff r_{\trf3}\trf(\dff v\trf)
\off =\off 
b\qff,\quad\
\]

\vspace{-36pt}
\[
\quad
r_{\trf 1}\dff \cdot\dff r_{\trf 2}\dff \cdot\dff r_{\trf 3}\dff \cdot\dff r_{\trf 4}\trf(\dff v\trf)
\off =\off 
w_{\dff 2}\qff,\quad\
\mbox{and}\quad\
r_{\trf 1}\dff \cdot\dff r_{\trf 2}\dff \cdot\dff r_{\trf 3}\dff \cdot\dff 
r_{\trf 4}\dff \cdot\dff r_{\trf 5}\qff(\dff v\trf)
\off =\off 
v\qff.
\]

\vspace{-9pt}
Such elements $r_{\dff i}$ are uniquely\sss determined\qss
(see\dss Section\qss \ref{one},\oss where\sss they\dss were denoted\dss by\sss $s_{\dff i}$\nsp),\oss
and\dss in\dss our calculations
we will\dss not\dss mention\dss them explicitly\halfff.\oss
We\sss will\sss simply\dss will\dss write down\sss the above equalities in\dss terms of\trs
the elements\sss $s_{\dff i}$\nsp.\oss
The fact\dss that\dss $s_{\dff i}^{\dff -\dff 1}\off =\off s_{\dff i}$\dss 
allows\sss to shorten\dss the formulas.\oss
To begin\dss with,\pss 
$s_{\trf 1}\trf(\dff v\trf)
\off =\off 
w_{\dff 1}$\nsp.\oss
Clearly\halfff,\pss 
$s_{\dff 1}\trf(\dff a\trf)
\off =\off 
w_{\dff 3}$\dss
and\dss hence\vspace{3pt}
\[
\quad
s_{\trf 1}\dff \cdot\dff 
s_{\trf 3}\trf(\dff v\trf)
\off =\off 
a
\qff.
\]

\vspace{-9pt}
Next\halfff,\pss
$s_{\trf 3}\dff \cdot\dff 
s_{\trf 1}\trf(\dff b\trf)
\off =\off 
s_{\trf 3}\trf(\dff d\trf)
\off =\off 
w_{\dff 2}$\dss
and\dss hence\vspace{3pt}
\[
\quad
s_{\trf 1}\dff \cdot\dff 
s_{\trf 3}\dff \cdot\dff 
s_{\trf 2}\trf(\dff v\trf)
\off =\off 
b
\qff.
\]

\vspace{-9pt}
Next\halfff,\pss
$s_{\trf 2}\dff \cdot\dff s_{\trf 3}\dff \cdot\dff s_{\trf 1}\trf(\dff w_{\dff 2}\trf)
\off =\off
s_{\trf 2}\dff \cdot\dff s_{\trf 3}\trf(\dff c\trf)
\off =\off
s_{\trf 2}\trf(\dff e\trf)
\off =\off 
w_{\dff 1}$\qss
and\dss hence\vspace{3pt}
\[
\quad
s_{\trf 1}\dff \cdot\dff 
s_{\trf 3}\dff \cdot\dff 
s_{\trf 2}\dff \cdot\dff 
s_{\trf 1}\trf(\dff v\trf)
\off =\off 
w_{\dff 2}
\qff.
\]

\vspace{-9pt}
Finally\halfff,\pss
$s_{\trf 1}\dff \cdot\dff 
s_{\trf 2}\dff \cdot\dff 
s_{\trf 3}\dff \cdot\dff 
s_{\trf 1}\trf(\dff v\trf)
\off =\off
s_{\trf 1}\dff \cdot\dff 
s_{\trf 2}\dff \cdot\dff 
s_{\trf 3}\trf(\dff w_{\dff 1}\trf)
\off =\off
s_{\trf 1}\dff \cdot\dff 
s_{\trf 2}\trf(\qff f\trf)
\off =\off 
s_{\trf 1}\trf(\dff a\trf)
\off =\off 
w_{\dff 3}$\qss
and\dss hence\vspace{3pt}
\[
\quad
s_{\trf 1}\dff \cdot\dff 
s_{\trf 3}\dff \cdot\dff 
s_{\trf 2}\dff \cdot\dff 
s_{\trf 1}\dff \cdot\dff 
s_{\trf 3}\trf(\dff v\trf)
\off =\off 
v
\qff.
\]

\vspace{-9pt}
Therefore our\dss loop relation\dss is\vspace{3pt}
\[
\quad
g_{\trf 1}\dff \cdot\dff 
g_{\trf 3}\dff \cdot\dff 
g_{\trf 2}\dff \cdot\dff 
g_{\trf 1}\dff \cdot\dff 
g_{\trf 3}
\off =\off 
s_{\trf 1}\dff \cdot\dff 
s_{\trf 3}\dff \cdot\dff 
s_{\trf 2}\dff \cdot\dff 
s_{\trf 1}\dff \cdot\dff 
s_{\trf 3}
\qff.
\]

\vspace{-9pt}
We need\dss to compute\dss
$s_{\trf 1}\dff \cdot\dff 
s_{\trf 3}\dff \cdot\dff 
s_{\trf 2}\dff \cdot\dff 
s_{\trf 1}\dff \cdot\dff 
s_{\trf 3}$\nsp.\oss
Since\sss this element\dss fixes $v$ and\dss hence belongs\sss to\dss $G_{\dff v}$\nsp,\oss
it\dss is\dss determined\dss by\dss its action on any\sss of\trs
the vertices\sss $w_{\dff i}$\nsp.\oss
For example,\oss\vspace{3pt}
\[
\quad
s_{\trf 1}\dff \cdot\dff 
s_{\trf 3}\dff \cdot\dff 
s_{\trf 2}\dff \cdot\dff 
s_{\trf 1}\dff \cdot\dff
s_{\trf 3}\trf(\dff w_{\dff 3}\trf)
\off =\off
s_{\trf 1}\dff \cdot\dff 
s_{\trf 3}\dff \cdot\dff 
s_{\trf 2}\dff \cdot\dff 
s_{\trf 1}\trf(\dff v\trf)
\]

\vspace{-36pt}
\[
\quad
\phantom{s_{\trf 1}\dff \cdot\dff 
s_{\trf 3}\dff \cdot\dff 
s_{\trf 2}\dff \cdot\dff 
s_{\trf 1}\dff \cdot\dff
s_{\trf 3}\trf(\dff w_{\dff 3}\trf)
\off }
=\off
s_{\trf 1}\dff \cdot\dff 
s_{\trf 3}\dff \cdot\dff 
s_{\trf 2}\trf(\dff w_{\dff 1}\trf)
\]

\vspace{-36pt}
\[
\quad
\phantom{s_{\trf 1}\dff \cdot\dff 
s_{\trf 3}\dff \cdot\dff 
s_{\trf 2}\dff \cdot\dff 
s_{\trf 1}\dff \cdot\dff
s_{\trf 3}\trf(\dff w_{\dff 3}\trf)
\off }
=\off
s_{\trf 1}\dff \cdot\dff 
s_{\trf 3}\trf(\dff e\trf)
\off =\off
s_{\trf 1}\trf(\dff c\trf)
\off =\off 
w_{\dff 2}
\]

\vspace{-9pt}
and\dss hence\dss
$s_{\trf 1}\dff \cdot\dff 
s_{\trf 3}\dff \cdot\dff 
s_{\trf 2}\dff \cdot\dff 
s_{\trf 1}\dff \cdot\dff
s_{\trf 3}
\off =\off
h$\nnsp.\qff\oss
Now\dss we can\dss rewrite our\sss relation as\vspace{3pt}
\[
\quad
g_{\trf 1}\dff \cdot\dff 
g_{\trf 3}\dff \cdot\dff 
g_{\trf 2}\dff \cdot\dff 
g_{\trf 1}\dff \cdot\dff 
g_{\trf 3}
\off =\off 
h
\qff.
\]

\vspace{-9pt}
By\dss using\qss (\ref{g-dodecahedron})\qss
we can eliminate\dss $g_{\trf 2}$\dss
and\dss $g_{\trf 3}$\dss from\dss this relation and\dss get\vspace{4.5pt}
\[
\quad
g_{\trf 1}\dff \cdot\dff 
\left(\trf h\dff \cdot\dff g_{\dff 1}\dff \cdot\dff h^{\dff -\dff 1} \trf\right)\dff \cdot\dff 
\left(\trf h^{\dff 2}\dff \cdot\dff g_{\dff 1}\dff \cdot\dff h^{\dff -\dff 2} \trf\right)\dff \cdot\dff 
g_{\trf 1}\dff \cdot\dff 
\left(\trf h\dff \cdot\dff g_{\dff 1}\dff \cdot\dff h^{\dff -\dff 1} \trf\right)
\off =\off 
h
\qff.
\]

\vspace{-7.5pt}
Let\trs $g\off =\off g_{\trf 1}$\nsp.\oss
Since\dss $h^{\dff -\dff 2}\off =\off h$\nnsp,\oss
our relation\dss is\dss equivalent\dss to\vspace{3pt}
\[
\quad
g\dff \cdot\dff 
h\dff \cdot\dff 
g\dff \cdot\dff 
h\dff \cdot\dff 
g\dff \cdot\dff 
h\dff \cdot\dff 
g\dff \cdot\dff 
h\dff \cdot\dff 
g\dff \cdot\dff 
h
\off =\off 
1
\qff,
\]

\vspace{-9pt}
or\halfff,\oss what\dss is\dss the same,\pss to\dss
$\left(\dff g\dff h\trf\right)^{\fff 5}\off =\off 1$\nnsp.\oss
This\dss is\dss the\sss last\dss relation.\oss

\myuppar{The presentation of\dss $\mathcal{D}$\dnsp.}
The resulting\dss presentation\dss
has\sss two\sss generators\dss $g\dff,\pff h$\dss and\dss relations\vspace{1.375pt}
\[
\quad
g^{\dff 2}\off =\off 1\qff,\quad\
h^{\dff 3}\off =\off 1\qff,\quad\
\mbox{and}\quad\
\left(\dff g\dff h\trf\right)^{\fff 5}\off =\off 1
\qff.
\]

\vspace{-10.625pt}
This\dss is\dss the presentation of\trs $G$\dss
discovered\dss by\trs Hamilton\dss in\dss 1856,\oss
who studied\dss it\dss in\dss terms of\trs the regular\dss icosahedron.\oss
See\qss \cite{cm},\oss Section\qss 6.4.\oss
Our\dss methods apply\dss to\sss the\sss icosahedron equally\dss well,\oss
but\sss an example with\dss non-triangular\dss loops seems\sss to be\sss instructive.\oss
The first\dss two relations result\dss from\dss two symmetries of\trs the dodecahedron\dss $D$\nnsp.\oss
Namely\halfff,\oss the rotation\dss by\dss the angle\sss $\pi$\sss 
about\dss the\sss line passing\dss through\dss the center of\trs $D$\dss 
and\dss the midpoint\sss of\dss $e_{\dff 1}$\sss is\dss a symmetry\sss of\trs $D$\nnsp,\oss
as aslo\sss the rotation\dss by\dss the angle\sss $2\dff \pi/3$\sss 
about\dss the\sss line passing\dss through\dss the center of\trs $D$\dss
and\dss the point\sss $v$\nnsp.\oss 
They\dss have\sss the order\sss $2$\sss and\dss $3$\sss respectively\halfff.\oss
The last\dss relation\sss has a similar\sss meaning.\oss
Let\dss us\dss compute\sss the action of\trs $g\dff h$\dss
on\dss the face\dss $v\dff w_{\dff 1}\dff a\trf b\dff w_{\dff 3}$:\oss\vspace{3pt}\vspace{-2.5pt}
\[
\quad
g\dff h\trf(\dff v\trf)\off =\off g\trf(\dff v\trf)\off =\off w_{\dff 1}
\qff,
\]

\vspace{-36pt}
\[
\quad
g\dff h\trf(\dff w_{\dff 1}\trf)\off =\off g\trf(\dff w_{\dff 3}\trf)\off =\off a
\qff,
\]

\vspace{-36pt}
\begin{equation}
\label{rotating-a-face}
\quad
g\dff h\trf(\dff a\trf)\off =\off g\trf(\dff d\trf)\off =\off b
\qff,
\end{equation}

\vspace{-36pt}
\[
\quad
g\dff h\trf(\dff b\trf)\off =\off g\trf(\dff c\trf)\off =\off w_{\trf 2}
\qff,
\]

\vspace{-36pt}
\[
\quad
g\dff h\trf(\dff w_{\trf 2}\trf)\off =\off g\trf(\dff w_{\dff 1}\trf)\off =\off v
\qff,
\]

\vspace{-9pt}
It\dss follows\dss that\dss $g\dff h$\dss leaves\sss the face\dss
$v\dff w_{\dff 1}\dff a\trf b\dff w_{\dff 3}$\dss 
invariant\sss and,\pss moreover\halfff,\oss 
is\dss the clockwise\qss
({\fff}when one looks\sss at\dss
the dodecahedron\dss from\sss the outside)\qss 
rotation of\trs this face by\dss the angle\dss $2\dff \pi/5$\dss 
about\dss the line passing\dss through\dss the center of\trs $D$\dss
and\dss the center of\trs the pentagon\dss
$v\dff w_{\dff 1}\dff a\trf b\dff w_{\dff 3}$\nsp.\oss
The last\dss relation\sss reflects\sss the fact\dss that\dss this\dss is\dss
a symmetry\sss of\trs $D$\dss of\dss order\dss $5$\nnsp.\oss

\myuppar{The\sss group $\spin(\dff 1\dff)$\nnsp.}
Recall\dss that\dss the fundamental\dss group of\trs the group\dss $SO\dff(\dff 3\dff)$\dss
is\dss isomorphic\sss to\dss $\zzz/2$\dss and\dss hence\sss the universal\sss covering space of\trs
$SO\dff(\dff 3\dff)$\dss is\dss a\sss two-sheeted covering\sss space.\oss
It\dss is\dss a\sss topological\dss group in a canonical\dss way\halfff.\oss
We will\sss denote it\dss by\sss $\spin(\dff 1\dff)$\nnsp.\oss
There\dss is\dss a canonical\dss homomorphism\dss
$\spin(\dff 1\dff)\qff \ttoo\qff SO\dff(\dff 3\dff)$\dss
with\dss the kernel\sss contained\sss in\dss the center of\trs $\spin(\dff 1\dff)$\dss
and\dss isomorphic\sss to\dss $\zzz/2$\nnsp.\oss
Let\dss $c$\dss be\sss the non-trivial\sss element\sss of\trs this\sss kernel.\oss

Let\dss us\dss consider\dss the group of\dss rotations of\trs $\rrr^{\dff 3}$\dss
about\sss a\sss fixed axis.\oss
It\dss is\dss isomorphic\sss to\dss $SO\dff(\dff 2\dff)$\dss and\dss the\sss loop
of\dss rotations\sss by\dss the angles\dss $2\dff \pi\fff t$\nnsp,\pss
$0\qff \leq\qff t\qff \leq\qff 1$\nnsp,\oss represents\sss the generator of\trs
the fundamental\dss group of\trs $SO\dff(\dff 3\dff)$\nnsp.\oss
It\dss follows\sss that\dss the preimage of\trs this group of\dss rotations
in\dss $\spin(\dff 1\dff)$ is\dss also\sss isomorphic\sss to $SO\dff(\dff 2\dff)$\nnsp,\oss
and\dss that\dss the map\dss
$SO\dff(\dff 2\dff)\qff \ttoo\qff SO\dff(\dff 2\dff)$\dss
induced\dss by\dss
$\spin(\dff 1\dff)\qff \ttoo\qff SO\dff(\dff 3\dff)$\dss
is\dss the map\dss $r\off \longmapsto\off r^{\dff 2}$\dnsp.\oss
This suggests\sss to\sss think about\dss the elements of\trs this preimage
as rotations about\dss the same axis\sss by\sss angles from\dss $0$\dss
to\dss $4\dff \pi$\nnsp,\oss
with\dss the rotation\sss by\dss the angle\dss $2\dff \pi$\dss being\dss
the element\sss $c$\dss generating\dss the kernel.\oss

\myuppar{The binary\dss icosahedral\dss group.}
We define\sss it\dss in\dss terms of\dss a\sss regular dodecahedron\dss $D$\dss
with\dss the center at\trs $0\qff \in\qff \rrr^{\dff 3}$\dnsp.\oss
The\qss \emph{binary\dss icosahedral\dss group}\qss 
is\dss the preimage\dss $\mathcal{D}^{\fff \sim}$\dss of\trs the\sss group\sss $\mathcal{D}$\sss
of\dss orientation-preserving\sss symmetries of\trs $D$\nnsp.\oss
It\dss acts on\dss $D$\dss by\dss the canonical\dss homomorphism\dss
$\mathcal{D}^{\fff \sim}\qff \ttoo\qff \mathcal{D}$\dss and\dss will\sss use\sss this action\dss to
find a presentation of\trs $\mathcal{D}^{\fff \sim}$\dnsp.\oss
The\sss task\dss is\dss very\sss similar\dss to\sss the case of\trs the group\dss $\mathcal{D}$\dss
itself\halfff,\oss and\sss we will\sss use\sss the same notations\sss to\sss
the extent\dss possible.\oss

Let\dss $h$\dss be\sss the rotation\dss by\dss the angle\sss $2\dff \pi/3$\sss about\dss
the axis passing\dss through\dss $0$\dss and\dss the vertex\sss $v$\nnsp.\oss 
Then\dss the equalities\qss (\ref{action-h})\qss still\dss hold,\oss
but\dss $h^{\dff 3}$\dss is\dss not\dss the identity\halfff,\oss
but\dss the rotation\dss by\dss the angle\dss $2\dff \pi$\nnsp,\oss
i.e.\qss $h^{\dff 3}\off =\off c$\nnsp.\oss
It\dss follows\sss that\dss the stabilizer\dss $G_{\dff v}$\dss of\trs the vertex\sss $v$\sss
in\dss $\mathcal{D}^{\fff \sim}$\dss is\dss the cyclic group of\dss order\sss $6$\sss
generated\dss by\dss $h$\nnsp.\oss
As before,\qss $G_{\dff v}$ acts\sss transitively\sss on\dss $E$\dss
and\dss hence we can\dss take\dss 
$E_{\dff 0}\off =\off \{\trf e_{\dff 1}\qff\}$\sss 
is\dss a set\sss of\dss representatives of\trs
$G_{\dff v}$\dnsp-orbits.\oss
The stabilizer\dss $G_{\dff e}$\dss of\trs the edge\dss
$e\off =\off e_{\dff 1}$\dss in\sss $G_{\dff v}$\sss
is\dss equal\dss to\sss the kernel\sss of\trs 
$\spin(\dff 1\dff)\qff \ttoo\qff SO\dff(\dff 3\dff)$\nnsp,\oss
i.e.\qss is\dss the subgroup\dss $\{\qff 1\fff,\pff h^{\dff 3}\qff\}$\dss
of\trs $G_{\dff v}$\nsp.\oss
Therefore we can\dss take\dss
$\{\qff 1\fff,\pff h\fff,\pff h^{\dff 2}\qff\}$\dss
as\sss the set\dss
$\mathcal{T}_{\dff e}$\dss of\dss representatives of\dss cosets is\trs $G_{\dff v}/G_{\dff e}$\dss
and\dss the one-element\sss set\trs $\{\qff h^{\dff 2}\qff\}\off =\off \{\trf c\trf\}$\dss
as\sss the set\dss $\mathcal{H}_{\dff e}$\dss of\dss generators of\trs $G_{\dff e}$\nsp.\oss 

For\dss $i\off =\off 1\fff,\pff 2\fff,\pff 3$\trs let\dss $s_{\dff i}\qff \in\qff \spin(\dff 1\dff)$\dss
be\sss the rotation\dss by\dss the angle\sss $\pi$\sss about\dss
the axis passing\dss through\dss $0$\dss and\dss the midpoint\sss of\trs the edge\dss $e_{\dff i}$\nsp.\oss
Then\sss $s_{\dff i}$\sss is\sss the\dss unique\sss rotation\dss by\dss the angle\sss $\pi$\sss
which\dss interchanges\sss the endpoints of\dss $e_{\dff i}$\dss
({\fff}the other such element\dss is\dss a\sss rotation\dss by\dss $3\dff \pi$\nsp)\qss
and\dss we may\sss set\dss $s_{\dff e_{\dff i}}\off =\off s_{\dff i}$\nsp.\oss
The equalities\qss (\ref{s-dodecahedron})\qss still\dss hold\dss
because an element\sss of\trs $\spin(\dff 1\dff)$\dss conjugate\sss to a rotation\dss
by\sss an angle $\theta$ is\dss also a rotation\dss by\dss the same angle $\theta$\nnsp.\oss
It\dss follows\dss that\dss the set\trs 
$E_{\dff 0}\off =\off \{\trf e_{\dff 1}\qff\}$\nnsp,\oss
the set\trs 
$\mathcal{T}_{\dff e}\off =\off \{\qff 1\fff,\pff h\fff,\pff h^{\dff 2}\qff\}$\nnsp,\oss
and\dss the\sss family\sss of\dss elements\dss
$s_{\dff e_{\dff i}}\off =\off s_{\dff i}$\dss 
form\sss a\sss regular scaffolding.\oss
In addition,\oss the rotation\dss $s_{\dff i}\dff \cdot\dff s_{\dff i}$\dss 
is\dss equal\dss a rotation\sss by\dss the angle\dss $2\dff \pi$\dss
and\dss hence\dss
$s_{\dff i}\dff \cdot\dff s_{\dff i}\off =\off c$\dss
and\dss $s_{\dff i}^{\dff -\dff 1}\off =\off s_{\dff i}\dff \cdot\dff c$\dss
for\dss $i\off =\off 1\fff,\pff 2\fff,\pff 3$\nnsp.\oss

By\qss Theorem\qss \ref{second-simplification}\qss the group\sss $\mathcal{D}^{\fff \sim}$\sss
can\dss be obtained\dss from\dss the group\sss $G_{\dff v}$\sss by\sss adding\sss
one generator\sss 
$g_{\trf 1}
\off =\off 
g_{\dff e_{\dff 1}}$\nsp,\oss 
the edge relation\dss $E\dff(\dff e_{\dff 1}\fff,\pff c\trf)$\nnsp,\oss
the edge-loop relation\dss $L\trf(\trf l_{\dff e_{\dff 1}}\trf)$\nnsp,\oss
and\dss the\sss loop relation\dss $L\dff(\trf l\trf)$\dss corresponding\dss to\sss the\sss same\sss loop\dss
$l$\dss as before.\oss

Let\dss us\dss consider\dss the relation\dss $E\dff(\dff e_{\dff 1}\fff,\pff c\trf)$\dss first\halfff.\oss
The element\trs $c\qff \in\qff G^{\dff \sim}\qff \subset\qff \spin(\dff 1\dff)$\dss 
belongs\sss to\sss the kernel\sss of\trs
$\spin(\dff 1\dff)\qff \ttoo\qff SO\dff(\dff 3\dff)$\dss
and\dss hence\dss $c\dff(\dff e_{\dff 1}\dff)\off =\off e_{\dff 1}$\nsp.\oss
Therefore\sss the relation\dss $E\dff(\dff e_{\dff 1}\fff,\pff c\trf)$\dss is\vspace{3pt}
\[
\quad
g_{\dff 1}^{\dff -\dff 1}\dff \cdot\dff c\dff \cdot\dff g_{\dff 1}
\off =\off
s_{\dff 1}^{\dff -\dff 1}\dff \cdot\dff c\dff \cdot\dff s_{\dff 1}
\pff.
\]

\vspace{-9pt}
But\trs 
$s_{\dff 1}^{\dff -\dff 1}\dff \cdot\dff c\dff \cdot\dff s_{\dff 1}
\off =\off
s_{\dff 1}\dff \cdot\dff c\dff \cdot\dff c\dff \cdot\dff s_{\dff 1}
\off =\off
s_{\dff 1}\dff \cdot\dff s_{\dff 1}
\off =\off
c$\dss
and\dss hence\dss $E\dff(\dff e_{\dff 1}\fff,\pff c\trf)$\dss is equivalent\dss to\vspace{3pt}
\[
\quad
g_{\dff 1}^{\dff -\dff 1}\dff \cdot\dff c\dff \cdot\dff g_{\dff 1}
\off =\off
c
\pff,
\]

\vspace{-9pt}
or\halfff,\oss what\dss is\dss the same,\oss to\dss
$c\dff \cdot\dff g_{\dff 1}\off =\off g_{\dff 1}\dff \cdot\dff c$\nnsp.\oss

The edge-loop relation\dss $L\trf(\trf l_{\dff e_{\dff 1}}\trf)$\dss
is\dss
$g_{\trf 1}\dff \cdot\dff g_{\trf 1}\off =\off s_{\dff 1}\dff \cdot\dff s_{\dff 1}$\nnsp.\oss
Since\dss $s_{\dff 1}\dff \cdot\dff s_{\dff 1}\off =\off c$\nnsp,\oss
the edge-loop relation\dss $L\trf(\trf l_{\dff e_{\dff 1}}\trf)$\dss
is\dss equivalent\dss to\dss \vspace{2pt}
\[
\quad
g_{\trf 1}\dff \cdot\dff g_{\trf 1}
\off =\off 
c
\qff.
\]

\vspace{-10pt}
Let\dss us\dss consider\sss $L\dff(\trf l\trf)$\nnsp.\oss
Since $s_{\dff i}\qff \in\qff \mathcal{D}^{\fff \sim}$
acts\sss on $D$ in\dss the same way\sss as\sss the element\sss of\sss $\mathcal{D}$ denoted\dss by\sss
$s_{\dff i}$\sss before,\oss
the first\dss part\sss of\trs the calculations\dss is\dss the same
and\dss the relation\dss $L\dff(\trf l\trf)$\dss is\vspace{3pt}
\[
\quad
g_{\trf 1}\dff \cdot\dff 
g_{\trf 3}\dff \cdot\dff 
g_{\trf 2}\dff \cdot\dff 
g_{\trf 1}\dff \cdot\dff 
g_{\trf 3}
\off =\off 
s_{\trf 1}\dff \cdot\dff 
s_{\trf 3}\dff \cdot\dff 
s_{\trf 2}\dff \cdot\dff 
s_{\trf 1}\dff \cdot\dff 
s_{\trf 3}
\qff,
\]

\vspace{-9pt}
where\dss $g_{\dff 2}\dff,\pff g_{\dff 3}$\dss are defined\dss by\dss the
formulas\qss (\ref{g-dodecahedron})\qss as before 
and\dss hence\vspace{3pt}
\[
\quad
g_{\trf 1}\dff \cdot\dff 
g_{\trf 3}\dff \cdot\dff 
g_{\trf 2}\dff \cdot\dff 
g_{\trf 1}\dff \cdot\dff 
g_{\trf 3}
\]

\vspace{-33pt}
\[
\quad
=\off
g_{\trf 1}\dff \cdot\dff 
(\trf h\dff \cdot\dff g_{\dff 1}\dff \cdot\dff h^{\dff -\dff 1} \trf)\dff \cdot\dff 
(\trf h^{\dff 2}\dff \cdot\dff g_{\dff 1}\dff \cdot\dff h^{\dff -\dff 2} \trf)\dff \cdot\dff 
g_{\trf 1}\dff \cdot\dff 
(\trf h\dff \cdot\dff g_{\dff 1}\dff \cdot\dff h^{\dff -\dff 1} \trf)
\]

\vspace{-33pt}
\[
\quad
=\off
g_{\trf 1}\dff \cdot\dff 
h\dff \cdot\dff 
g_{\dff 1}\dff \cdot\dff 
h\dff \cdot\dff 
g_{\dff 1}\dff \cdot\dff h^{\dff -\dff 3}\dff \cdot\dff
h\dff \cdot\dff 
g_{\trf 1}\dff \cdot\dff 
h\dff \cdot\dff 
g_{\dff 1}\dff \cdot\dff 
h^{\dff -\dff 1}
\]

\vspace{-6pt}
Since\trs 
$h^{\dff -\dff 3}\off =\off h^{\dff 3}\off =\off c$\dss
and\dss
$c\dff \cdot\dff g_{\dff 1}\off =\off g_{\dff 1}\dff \cdot\dff c$\qss
by\trs $E\dff(\dff e_{\dff 1}\fff,\pff c\trf)$\nnsp,\oss
the\sss last\sss expression\dss is\dss equal\dss to\vspace{6pt}
\[
\quad
g_{\trf 1}\dff \cdot\dff 
h\dff \cdot\dff 
g_{\dff 1}\dff \cdot\dff 
h\dff \cdot\dff 
g_{\dff 1}\dff \cdot\dff 
h\dff \cdot\dff 
g_{\trf 1}\dff \cdot\dff 
h\dff \cdot\dff 
g_{\dff 1}\dff \cdot\dff
h^{\dff -\dff 3}\dff \cdot\dff 
h^{\dff -\dff 1}
\off =\off
(\trf g_{\dff 1}\dff \cdot\dff h\trf)^{\dff 5}\dff \cdot\qff h
\]

\vspace{-6pt}
and\dss hence\qss
$g_{\trf 1}\dff \cdot\dff 
g_{\trf 3}\dff \cdot\dff 
g_{\trf 2}\dff \cdot\dff 
g_{\trf 1}\dff \cdot\dff 
g_{\trf 3}
\off =\off 
(\trf g_{\dff 1}\dff \cdot\dff h\trf)^{\dff 5}\dff \cdot\qff h$\nnsp.\qff\oss
Similarly\halfff,\oss\vspace{4pt}
\[
\quad
s_{\trf 1}\dff \cdot\dff 
s_{\trf 3}\dff \cdot\dff 
s_{\trf 2}\dff \cdot\dff 
s_{\trf 1}\dff \cdot\dff 
s_{\trf 3}
\off =\off 
(\trf s_{\dff 1}\dff \cdot\dff h\trf)^{\dff 5}\dff \cdot\qff h
\qff.
\]

\vspace{-8pt}
It\dss follows\dss that\dss the relation\sss $L\dff(\trf l\trf)$\dss
is\dss equivalent\dss to\vspace{3pt}\vspace{-0.5pt}
\[
\quad
(\trf g_{\dff 1}\dff \cdot\dff h\trf)^{\dff 5}
\off =\off
(\trf s_{\dff 1}\dff \cdot\dff h\trf)^{\dff 5}
\qff.
\]

\vspace{-9pt}
By\qss Lemma\qss \ref{center},\oss which we will\dss prove in a moment\halfff,\pss
$(\trf s_{\dff 1}\dff \cdot\dff h\trf)^{\dff 5}
\off =\dff\off
c$\nnsp.\oss
Since\dss $c\off =\off h^{\dff 3}$\dnsp,\oss
it\dss follows\dss that\dss the\sss group\dss $G^{\dff \sim}$\dss
admits a\sss presentation\dss with\dss two generators\dss 
$g\off =\off g_{\dff 1}$\dss and\dss $r$\dss
and\dss relations\vspace{3pt}
\[
\quad
g^{\dff 2}
\off =\off 
r^{\dff 3}
\off =\off
(\trf g\fff r\trf)^{\dff 5}
\quad\
\mbox{and}\quad\
r^{\dff 6}\off =\dff\off 1 
\qff.
\]

\vspace{-9.25pt}
A conjugation\dss by\dss $r$\dss turns\dss
$r^{\dff 3}
\off =\off
(\trf g\dff r\dff)^{\dff 5}$\dss
into\dss
$r^{\dff 3}
\off =\off
(\trf r\dff g\dff)^{\dff 5}$\dnsp,\qff\oss
leading\dss to\sss the presentation\vspace{3pt}
\[
\quad
g^{\dff 2}
\off =\off 
r^{\dff 3}
\off =\off
(\dff r\fff g\trf)^{\dff 5}
\quad\
\mbox{and}\quad\
r^{\dff 6}\off =\dff\off 1 
\qff.
\]

\vspace{-9.25pt}
A minor\sss modification\dss turns\sss the\sss latter\dss relations in\dss
the classical\dss form of\qss Coxeter\qss \cite{c1}.\oss
Let\trs\vspace{3pt}
\[
\quad
z
\off =\off 
g^{\dff 2}
\off =\off 
r^{\dff 3}
\qff,\quad\
s
\off =\off 
r^{\dff -\dff 1}
\qff,\quad\
\mbox{and}\quad\
t
\off =\off
r\fff g
\qff.
\]

\vspace{-9pt}
In\dss terms of\dss generators\dss $s\fff,\pff t\fff,\pff z$\dss 
the above presentation\dss takes\sss the form\vspace{3pt}
\[
\quad
s^{\dff 3}
\off =\off 
t^{\dff 5}
\off =\off 
(\dff s\dff t\trf)^{\dff 2}
\off =\off
z\quad\
\mbox{and}\quad\
z^{\dff 2}
\off =\off
1
\qff.
\]

\vspace{-9pt}
Coxeter\qss \cite{c1}\qss proved\dss that\dss the relations\qss 
$s^{\dff 3}
\off =\off 
t^{\dff 5}
\off =\off 
(\dff s\dff t\trf)^{\dff 2}
\off =\off
z$\qss
imply\trs 
$z^{\dff 2}\off =\off 1$\dss and\dss hence\sss the relation\dss
$z^{\dff 2}\off =\off 1$\dss is\dss superfluous.\oss
We will\dss return\dss to\sss this remarkable result\sss of\qss Coxeter\dss  
in\dss the
next\sss section.\oss
The relations\dss
$s^{\dff 3}
\off =\off 
t^{\dff 5}
\off =\off 
(\dff s\dff t\trf)^{\dff 2}$\dss 
are\qss Coxeter's\qss relations\dss for\dss $\mathcal{D}^{\dff \sim}$\dnsp.\oss

\mypar{Lemma.}{center}
$(\trf s_{\dff 1}\dff \cdot\dff h\trf)^{\dff 5}
\off =\dff\off
c$\nnsp.\oss\vspace{-0.75pt}

\proof
Let\sss $z$\sss be\sss the center of\trs the face\dss $v\dff w_{\dff 1}\dff a\trf b\dff w_{\dff 3}$\dss
of\trs the dodecahedron\dss $D$\nnsp,\oss and\dss let\sss $y$\sss be\sss the midpoint\sss of\trs
the edge\dss $v\dff w_{\dff 1}$\nnsp.\oss
Then\dss $y\dff v\dff z$\dss is\dss a\sss right-angled\dss triangle.\oss
See\sss the picture.\oss\vspace{12pt}

\begin{figure}[h!]
\hspace*{10.5em}
\includegraphics[width=0.4375\textwidth]{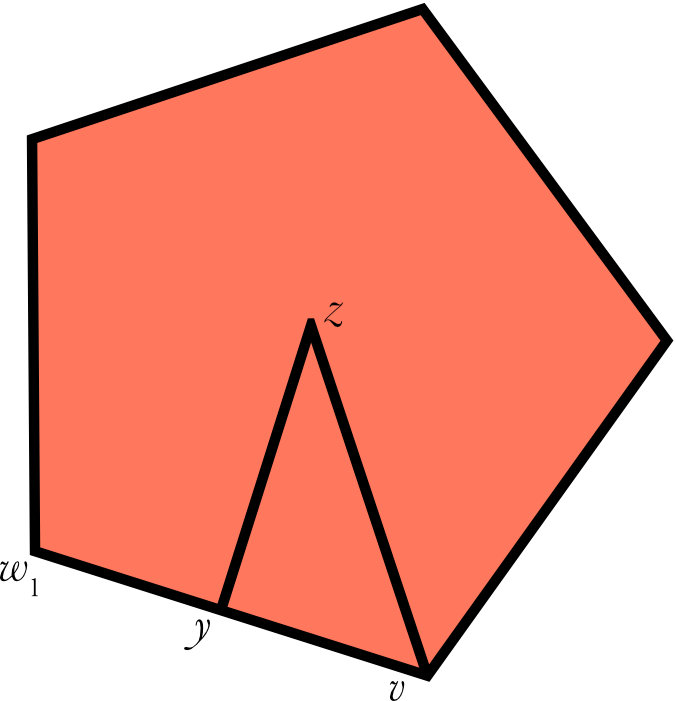}
\end{figure}\vspace{12pt}

Recall\dss that\dss the dodecahedron\dss $D$\dss
has\sss the origin\dss $0\qff \in\qff \rrr^{\dff 3}$\dss as\sss its center\halfff.\oss
Let\dss us\sss consider\dss the radial\dss projection of\trs $D$\dss to\sss
the unit\sss sphere $\mathbb{S}^{\dff 3}$ in\dss $\rrr^{\dff 3}$\dnsp,\oss
and\dss let\trs $y'\dff v'\dff z'$\dss be\sss the radial\dss projection of\trs
the\sss triangle\dss $y\dff v\dff z$\nnsp.\oss
Let\trs $P_{\fff y\fff v}$\nsp,\pss $P_{\fff v\fff z}$\nsp,\pss 
and\dss $P_{\fff z\fff y}$\dss be\sss the planes in\dss $\rrr^{\dff 3}$\dss
and containing\sss $0$\sss and\dss the segments\dss
$y\dff v$\nnsp,\pss $v\dff z$\nsp,\oss and\dss $z\dff y$\dss respectively\halfff.\oss
Let\dss 
$\sigma_{\fff y\fff v}\dff,\off \sigma_{\fff v\fff z}\dff,\off \sigma_{\fff z\fff y}$\dss
be\sss the reflections
in\dss the planes\dss 
$P_{\fff y\fff v}\dff,\off P_{\fff v\fff z}\dff,\off P_{\fff z\fff y}$\dss 
respectively\halfff.\oss
The product\sss of\trs two such\dss reflections\dss 
is\dss equal\dss to\sss the rotation about\dss their\dss line of\dss
intersection\sss by\dss the angle equal\dss to\sss twice\sss the spherical\dss
angle\sss between\dss the corresponding sides 
of\trs $y'\dff v'\dff z'$\dnsp.\oss
The symmetry\sss of\trs the dodecahedron\dss implies\sss that\dss
the products\dss 
$\sigma_{\fff z\fff v}\dff \sigma_{\fff v\fff y}$\nsp,\qss
$\sigma_{\fff v\fff y}\dff \sigma_{\fff y\fff z}$\nsp,\oss and\dss
$\sigma_{\fff y\fff z}\dff \sigma_{\fff z\fff v}$
are clockwise rotations\sss by\dss the angles\sss
$2\dff \pi/3$\nnsp,\dss $\pi$\nnsp,\qss and\sss $2\dff \pi/5$
respectively\halfff.\pss
In\dss turn,\pss this implies\sss that\dss these products are,\qss
respectively\halfff,\qss
the images in $SO\dff (\dff 3\dff)$ of\trs the elements\dss
$h^{\dff -\dff 1}$\dnsp,\qss $s_{\dff 1}$\dss and\sss 
$f\dff \in\qff \spin(\dff 1\dff)$\dnsp,\oss
where $f$ is\dss the clockwise rotation about\dss the\sss line\dss $0\dff z$\dss
by\dss the angle\dss $2\dff \pi/5$\nnsp.\oss
But\vspace{3pt}
\[
\quad
\sigma_{\fff v\fff z}^{\qff 2}
\off =\off
\sigma_{\fff v\fff z}^{\qff 2}
\off =\off
\sigma_{\fff z\fff y}^{\qff 2}
\off =\off
1\quad\
\]

\vspace{-9.5pt}
and\dss hence\vspace{2.5pt}
\[
\quad
(\trf \sigma_{\fff z\fff v}\dff \sigma_{\fff v\fff y} \trf)\dff \cdot\dff 
(\trf \sigma_{\fff v\fff y}\dff \sigma_{\fff y\fff z} \trf)\dff \cdot\dff 
(\trf \sigma_{\fff y\fff z}\dff \sigma_{\fff z\fff v} \trf)\qff
\off =\qff\off
1
\qff.
\]

\vspace{-9pt}
It\dss follows\dss that\dss the image of\trs the product\qss
$h^{\dff -\dff 1}\dff \cdot\dff s_{\dff 1}\dff \cdot\dff f$\pss
in\dss $SO\dff(\dff 3\dff)$\dss is\dss equal\dss to\sss $1$\sss
and\dss hence\sss the product\qss
$h^{\dff -\dff 1}\dff \cdot\dff s_{\dff 1}\dff \cdot\dff f$\pss itself\qss
is\dss equal\dss either\dss to\sss $1$\sss or\trs to\sss $c$\nnsp.\oss

Suppose\sss that\trs
$h^{\dff -\dff 1}\dff \cdot\dff s_{\dff 1}\dff \cdot\dff f\off =\off c$\nnsp.\oss
Recall\dss that\trs
$s_{\dff i}^{\dff -\dff 1}\off =\off s_{\dff i}\dff \cdot\dff c$\dss
and\sss $c$\sss belongs\sss to\sss the center of\trs $\spin(\dff 1\dff)$\nnsp.\oss
It\dss follows\dss that\trs
$(\dff  h^{\dff -\dff 1}\dff \cdot\dff s_{\dff 1}\trf)^{\dff -\dff 1}
\off =\off
s_{\dff 1}\dff \cdot\dff h\dff \cdot\dff c$\dss
and\dss hence\dss
$h^{\dff -\dff 1}\dff \cdot\dff s_{\dff 1}\dff \cdot\dff f\off =\off c$\dss
implies\sss that\vspace{3.5pt}
\[
\quad
f
\off =\off
s_{\dff 1}\dff \cdot\dff h\dff \cdot\dff c
\dff \cdot\dff
c
\off =\off
s_{\dff 1}\dff \cdot\dff h
\qff.
\]

\vspace{-8.5pt}
But\halfff,\pss
clearly\halfff,\pss 
$f^{\qff 5}\qff \in\qff \spin(\dff 1\dff)$\dss 
is\dss a\sss rotation\sss by\dss the angle\dss $2\dff \pi$\dss
and\dss hence\dss $f^{\qff 5}\off =\off c$\dss and\vspace{3.5pt}
\[
\quad
\left(\trf s_{\dff 1}\dff \cdot\dff h\trf\right)^{\dff 5}
\off =\dff\off
f^{\qff 5}
\off =\off
c
\qff.
\]
 
\vspace{-8.5pt}
This proves\sss the\sss lemma\sss modulo\sss the equality\trs
$h^{\dff -\dff 1}\dff \cdot\dff s_{\dff 1}\dff \cdot\dff f\off =\off c$\nnsp.\oss

It\dss remains\sss to prove\sss that\trs
$s_{\dff 1}\dff \cdot\dff h^{\dff -\dff 1}\dff \cdot\dff f\off =\off c$\nnsp.\oss
In\dss fact\halfff,\oss a more general\sss statement\dss is\dss true.\oss
Let\trs $p\fff q\dff r$\dss be a spherical\trs triangle in\dss the
unit\sss sphere $\mathbb{S}^{\dff 3}$\dnsp,\oss
and\dss let\dss $\alpha\fff,\pff \beta\fff,\pff \gamma$\dss
be\sss the interior spherical\sss angles of\trs $p\fff q\dff r$ at\dss the vertices\dss
$p\fff,\pff q\fff,\pff r$\dss
respectively\halfff.\oss
Let\dss us\sss assume\sss that\dss the vertices\dss $p\fff,\pff q\fff,\pff r$\dss
follow\sss in\dss the clockwise order\halfff,\oss
and\trs let\trs 
$s_{\dff p}\dff,\off s_{\dff q}\dff,\off s_{\dff r}
\qff \in\qff \spin(\dff 1\dff)$\dss
be\sss the clockwise rotations about\dss the\sss lines\dss $0\fff p\fff,\pff 0\fff q\fff,\pff 0\fff r$\dss
by\dss the angles\dss
$2\dff \alpha\fff,\pff 2\trf \beta\fff,\pff 2\dff \gamma$\dss
respectively\halfff.\oss
Then\dss\vspace{3pt}
\[
\quad 
s_{\dff p}\dff \cdot\dff s_{\dff q}\dff \cdot\trf s_{\dff r}
\off =\off
c
\qff.
\]

\vspace{-9pt}
This\dss is\dss a special\sss case of\dss a\sss result\sss of\qss Milnor\qss \cite{m2}.\oss
See\qss \cite{m2},\oss the proof\dss of\qss Lemma\qss 3.1.\oss
For\dss the convenience of\trs the reader\dss we reproduce beautiful\qss
Milnor's\qss proof\halfff.\oss
Let\trs
$\sigma_{\fff p\dff q}\dff,\off \sigma_{\fff q\fff r}\dff,\off \sigma_{\fff r\fff p}$\dss
be\sss the reflections\qss 
({\fff}in\dss the sense of\trs the spherical\dss geometry\fff)\qss
in\dss the sides\dss $p\fff q$\nnsp,\qss $q\dff r$\nnsp,\pss and\dss $r\fff p$\dss
of\trs the\sss triangle\dss $p\fff q\dff r$\dss
respectively\halfff.\oss
Then\dss the homomorphism\qss
$\spin(\dff 1\dff)\qff \ttoo\qff SO\dff(\dff 3\dff)$\qss
takes\sss the rotations\dss $s_{\dff p}$\nsp,\qss $s_{\dff q}$\nsp,\oss
and\dss $s_{\dff r}$\qss
to\qss
$\sigma_{\fff r\fff p}\dff \sigma_{\fff p\dff q}$\nsp,\pss
$\sigma_{\fff p\dff q}\dff \sigma_{\fff q\dff r}$\nsp,\oss 
and\dss
$\sigma_{\fff q\fff r}\dff \sigma_{\fff r\fff p}$\dss
respectively\halfff.\pss
As above\vspace{3pt}\vspace{0.375pt}
\[
\quad
\sigma_{\fff p\dff q}^{\qff 2}
\off =\off
\sigma_{\fff q\fff r}^{\qff 2}
\off =\off
\sigma_{\fff r\fff p}^{\qff 2}
\off =\off
1\quad\
\]

\vspace{-9pt}
and\dss hence\vspace{3pt}
\[
\quad
(\trf \sigma_{\fff r\fff p}\dff \sigma_{\fff p\dff q} \trf)\dff \cdot\dff 
(\trf \sigma_{\fff p\dff q}\dff \sigma_{\fff q\dff r} \trf)\dff \cdot\dff 
(\trf \sigma_{\fff q\fff r}\dff \sigma_{\fff r\fff p} \trf)\qff
\off =\qff\off
1
\qff.
\]

\vspace{-9pt}\vspace{0.375pt}
It\dss follows\dss that\dss the image of\trs the product\qss
$s_{\dff p}\dff \cdot\dff s_{\dff q}\dff \cdot\trf s_{\dff r}$\pss
in\dss $SO\dff(\dff 3\dff)$\dss is\dss equal\dss to\sss $1$\sss
and\dss hence\sss the product\qss
$s_{\dff p}\dff \cdot\dff s_{\dff q}\dff \cdot\trf s_{\dff r}$\pss itself\qss
is\dss equal\dss either\dss to\sss $1$\sss or\trs to\sss $c$\nnsp.\oss
Now\dss let\dss us\dss deform\dss the\sss triangle\dss $p\fff q\dff r$
in such a way\dss that\sss all\dss its vertices\sss will\dss tend\dss to\sss
the same point\trs $x\qff \in\qff \mathbb{S}^{\dff 3}$\dnsp.\oss
By\dss the continuity\halfff,\oss during such deformation\dss
the product\trs
$s_{\dff p}\dff \cdot\dff s_{\dff q}\dff \cdot\dff s_{\dff r}$\dss
remains\sss the same.\oss
When all\dss three vertices\dss $p\fff,\pff q\fff,\pff r$\dss
are close\sss to\sss $x$\nnsp,\oss
the\sss triangle\dss is\dss nearly\sss euclidean and\dss hence its
sum of\dss angles\dss is\dss nearly\sss equal\dss to\sss $\pi$\nnsp.\oss
In\dss the\sss limit\dss the
rotations\dss
$s_{\dff p}\dff,\off s_{\dff q}\dff,\off s_{\dff r}
\qff \in\qff \spin(\dff 1\dff)$\dss
turn\sss into\sss the rotations about\dss the\sss line\sss $0\fff x$\sss
by\dss the angles with\dss the sum equal\dss to\sss twice\sss the sum of\dss
angles of\dss a\sss euclidean\dss triangle,\oss i.e.\qss to\dss $2\dff \pi$\nnsp.\oss
Therefore,\oss in\dss the\sss limit\dss the product\trs 
$s_{\dff p}\dff \cdot\dff s_{\dff q}\dff \cdot\dff s_{\dff r}$\dss
is\dss equal\dss to a rotation\dss by\dss the angle\dss $2\dff \pi$\nnsp,\oss
i.e.\qss to\sss $c$\nnsp.\oss
It\dss follows\dss that\trs
$s_{\dff p}\dff \cdot\dff s_{\dff q}\dff \cdot\dff s_{\dff r}
\off =\off
c$\dss
during\dss the whole deformation,\oss 
and\dss hence for\dss the original\dss triangle also.\oss  \eproof

\mysection{Coxeter's\qss implication}{coxeter-implication}

\myuppar{Coxeter's\dss implication\sss in\dss terms of\qss $g\fff,\pff r\fff,\pff z$\nnsp.}
The goal\sss of\trs section\dss is\dss to show\sss how\dss the methods of\qss
Sections\qss \ref{one}\dss --\dss \ref{examples}\qss can\dss be adapted\dss to
prove\sss the remarkable implication\sss due\sss to\dss Coxeter\dss
and\dss mentioned\dss before\dss Lemma\qss \ref{center}.\oss
Coxeter's\dss original\dss proof\dss is\dss quite different\halfff.\oss
See\dss Appendix\qss \ref{coxeter-proof}.\oss
To\sss begin with,\oss we restate\sss this implication\sss in\dss terms of\dss
generators\dss $g\fff,\pff r\fff,\pff z$\dss from\dss Section\qss \ref{examples}.\oss
In\dss these\sss terms it\dss takes\sss the following\dss form\fff:\oss
if\pss $g\fff,\pff r\fff,\pff z$\dss are elements of\dss a\sss group\dss $G$\dss 
and\qss\vspace{3pt}
\begin{equation}
\label{coxeter-relations}
\quad
r^{\dff -\dff 3}
\off =\off 
(\dff r\dff g\trf)^{\dff 5}
\off =\off 
g^{\dff 2}
\off =\off
z
\qff,
\end{equation}

\vspace{-9pt}
then\qss $z^{\dff 2}\off =\off 1$\dnsp.\oss 
Since\sss
$r^{\dff -\dff 3}
\off =\off
g^{\dff 2}
\off =\off
z$\nnsp,\pss the elements\sss
$g$ and $r$ commute with $z$
and\dss hence\sss $z$\sss belongs\sss
to\sss the center\sss of\trs $G$\nnsp.\oss
Therefore\sss the relation\dss
$(\dff r\dff g\trf)^{\dff 5}
\off =\off
z$\dss
is\dss equivalent\dss to\dss
$(\dff g\dff r\trf)^{\dff 5}
\off =\off
z$\nnsp.

\myuppar{The universal\sss example for\dss Coxeter's\dss implication.}
In order\sss to prove\sss this implication,\oss
it\dss is\dss sufficient\dss to consider\sss the universal\sss example,\oss
namely\halfff,\oss the group\dss
$\mathcal{G}$\dss defined\dss by\dss generators\dss
$g\fff,\pff r\fff,\pff z$\trs 
and\dss relations\qss (\ref{coxeter-relations}).\oss 
There\dss is\dss a unique homomorphism\dss
$\varphi\dff \colon\dff
\mathcal{G}\qff \ttoo\qff \mathcal{D}$\dss 
such\dss that\vspace{3pt}\vspace{-0.75pt}
\[
\quad
\varphi\qff \colon\qff
g\off \longmapsto\off s_{\dff 1}\qff,\quad\
r\off \longmapsto\off h\qff,\quad\
z\off \longmapsto\off 1
\qff.
\]

\vspace{-9pt}\vspace{-0.75pt}
Since\sss the group\dss $\mathcal{D}$\dss can\sss be obtained\dss
from\dss $\mathcal{G}$\dss by\sss adding\dss the relation\dss $z\off =\off 1$\nnsp,\oss
the kernel\dss $\ker\dff \varphi$\dss of\dss $\varphi$\dss is\dss
the subgroup\dss normally\dss generated\dss by\sss $z$\nnsp.\oss
Since\sss $z$\sss belongs\sss to\sss the center of\trs $\mathcal{G}$\nnsp,\oss
the cyclic subgroup\sss $Z$\sss generated\dss by\sss $z$\sss
is\dss contained\sss in\dss the center\sss and\dss hence\dss is\dss
equal\dss to\sss the subgroup\dss normally\dss generated\dss by\sss $z$\nnsp.\oss 
Therefore\sss $Z$\sss is\dss
equal\dss to\sss the kernel\dss $\ker\dff \varphi$\dss of\dss $\varphi$\nnsp.\oss

The homomorphism\sss $\varphi$\sss defines an action of\trs $\mathcal{G}$\dss 
on\dss the dodecahedron\dss $D$\dss and on a related\sss polyhedron\sss $T$\dss
to be defined\sss in a\sss moment\halfff.\oss
The\sss latter action\dss is\dss going\dss to be our main\dss tool.\oss

\myuppar{Modifying\dss the dodecahedron.}
Let\sss $x_{\dff 1}$\sss be a point\sss on\dss the edge\dss 
$e_{\dff 1}\off =\off v\dff w_{\trf 1}$\dss of\trs 
$D$\dss near\sss the vertex\dss $v$\nnsp.\oss
It\dss matters only\dss that\sss $x_{\dff 1}$\sss is\dss 
closer\sss to $v$\sss than\dss the midpoint\sss of\dss $e_{\dff 1}$\nsp,\oss
but\dss it\dss is\dss convenient\dss to\sss think\dss that\sss 
$x_{\dff 1}$\sss is\dss fairly\sss close\sss to $v$\nnsp.\oss
Let\dss us\dss consider\dss the orbit\trs $\mathcal{D}\fff x_{\dff 1}$\dss of\dss $x_{\dff 1}$\sss 
under\dss the action of\trs $\mathcal{D}$\nnsp.\oss
If\trs $u\qff \in\qff \mathcal{D}$\dss fixes $x_{\dff 1}$\nsp,\oss
then\dss 
$u\dff(\dff e_{\dff 1}\trf)\off =\off e_{\dff 1}$\nnsp,\oss
and\sss since $x_{\dff 1}$\sss is\dss not\dss the midpoint\sss of\dss $e_{\dff 1}$\nsp,\oss
it\dss follows\sss that\trs $u\off =\off 1$\nnsp.\oss 
Hence\sss the stabilizer of\dss $x_{\dff 1}$\sss in\dss $\mathcal{D}$\dss 
is\dss trivial\sss and\dss
for every\dss two points\dss $y\fff,\off y'\qff \in\qff \mathcal{D}\fff x_{\dff 1}$\dss
of\trs $T$\trs there\dss is\dss a\sss unique element\trs $t\qff \in\qff \mathcal{D}$\dss
such\dss that\trs $t\trf(\trf y\trf)\off =\off y'$\dnsp.\oss 

Let\trs $T$\dss be\sss the convex\sss hull\sss of\trs this orbit\halfff.\oss
It\dss is\dss a\sss polyhedron\dss having\dss
elements of\trs the orbit\trs $\mathcal{D}\fff x_{\dff 1}$\dss as its vertices
and\sss can\dss be obtained\dss from\dss $D$\dss by\sss
cutting\sss off\dss a small\trs triangular\sss pyramid\sss at\sss each vertex.\oss
The pyramid\sss at\sss $v$\sss has\dss 
$v\fff,\qff\off x_{\dff 1}\fff,\qff\off h\dff(\dff x_{\dff 1}\trf)$\nnsp,\pss 
and\sss $h^{\dff 2}\dff(\dff x_{\dff 1}\trf)$
as its\sss vertices;\oss
the other pyramids are\sss the images of\trs this one 
under\dss the rotations of\trs $D$\nnsp.\oss 
Let\dss $Y$\trs be\sss the graph defined\dss by\dss the vertices and edges of\trs $T$\dnsp.\oss
It\dss can\dss be obtained\dss from\dss the graph\dss $X$\dss 
defined\dss by\dss the vertices and edges of\trs $D$\dss
by\sss replacing\sss each vertex\sss by a\qss \emph{small\dss triangle}\qss
as on\dss the following\sss picture.\vspace{12pt}

\begin{figure}[h!]
\hspace*{8em}
\includegraphics[width=0.6\textwidth]{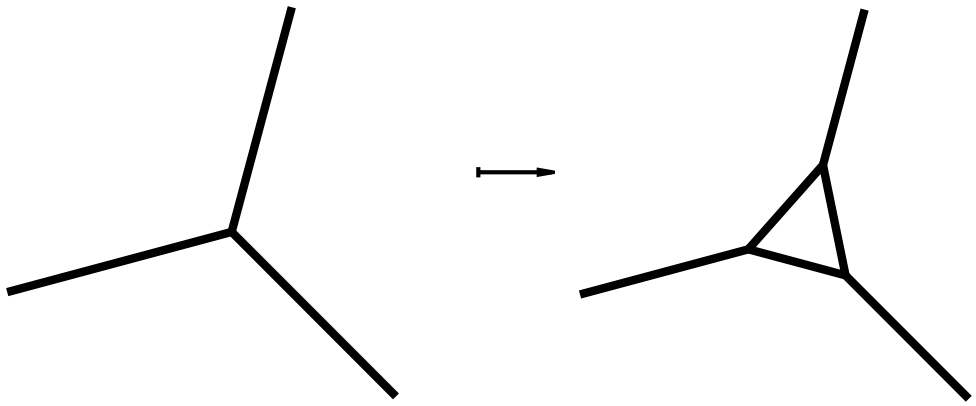}
\end{figure}\vspace{12pt}\vspace{-2.5pt}

The polyhedron\dss $T$\dss and\dss the  graph\dss $Y$\dss have\sss two\sss types of\dss edges.\oss
The edges of\dss the first\dss type correspond\dss to\sss the edges of\trs $X$\dss
and connect\dss two small\dss triangles.\oss
We will\sss call\trs them\dss the\qss \emph{pentagon-edges}.\oss
The edges of\trs the second\dss type are\sss the sides of\trs small\trs triangles
corresponding\dss to\sss the vertices of\trs $X$\nnsp,\oss
i.e.\qss to\sss the vertices of\trs the dodecahedron\dss $D$\nnsp.\oss
We will\sss call\trs them\dss the\qss \emph{triangle-edges}.\oss
Clearly\halfff,\oss pentagon-edges admit\dss inversions,\oss
and\dss triangle-edges do not\qss
(an\sss isometry\sss of\trs $D$\dss interchanging\dss the endpoints of\dss
a\sss triangle-edge\sss has\sss to be orientation-reversing\fff).\oss

\myuppar{Orientations.}
An\qss \emph{oriented edge}\qss of\dss a\sss graph\dss is\dss an edge\sss
together\dss with an\qss \emph{orientation}\qss of\trs this edge,\oss
i.e.\qss a designation of\dss one of\trs its endpoints as\sss its\qss \emph{origin}\qss
and\dss the other as\sss its\qss \emph{target}.\oss
Let\dss us\dss orient\dss the surface of\trs the polyhedron\dss $T$\dss
({\fff}i.e.\qss all\dss faces of\trs $T$\nsp)\qss
in such a way\dss that\dss the positive direction\dss is\dss
the clockwise one when one\sss looks at\dss $T$\sss from outside.\oss
This orientation defines an orientation of\dss each\dss triangular\dss face of\trs $T$\dss
and\dss hence an orientation of\dss every\dss triangle-edge of\trs $Y$\nnsp.\oss
We will\sss call\dss these orientations of\trs triangle-edge\sss the\qss
\emph{clockwise}\qss ones,\oss and\dss the opposite orientations\sss the\qss
\emph{counter-clockwise}\qss ones.\oss

\myuppar{The canonical\dss $\mathcal{D}$\dnsp-scaffolding.}
We need a version of\trs the notion of\dss a regular scaffolding\sss from\dss
Section\qss \ref{presentations}.\oss
Since\sss the vertices of\trs the graph\dss $Y$\dss are elements 
of\trs the orbit\dss $\mathcal{D}\fff x_{\dff 1}$\nsp,\oss
for every oriented edge $\varepsilon$ of\trs $Y$\dss
there\dss is\dss a\sss
unique element\trs $t_{\trf \varepsilon}\qff \in\qff \mathcal{D}$\trs
taking\dss the origin of\trs $\varepsilon$\dss to\sss 
the\sss target of\trs $\varepsilon$\nnsp.\oss
The map\dss
$\varepsilon\off \longmapsto\off t_{\trf \varepsilon}$\dss 
is\dss the\qss \emph{canonical\dss $\mathcal{D}$\dnsp-scaffolding}\qss of\trs $Y$\dnsp.\oss

If\trs an oriented edge $\varepsilon$\sss results\sss from orienting\sss
a\sss pentagon-edge $e$\nnsp,\oss
then\dss 
$t_{\trf \varepsilon}$\dss 
is\dss the unique element\sss of\dss $\mathcal{D}$\dss interchanging\dss 
the endpoints of\dss $e$\nnsp.\oss
In\dss particular\halfff,\oss if\dss $e$\sss corresponds\sss to\sss
the edge\sss $e_{\dff 1}$\sss of\trs $X$\nnsp,\oss
then\dss
$t_{\trf \varepsilon}\off =\off s_{\dff 1}$\nsp.\oss
Clearly\halfff,\oss if\dss $\varepsilon$\sss is\dss an oriented edge,\pss
$u\qff \in\qff \mathcal{D}$\nnsp,\oss
and\dss $\delta\off =\off u\dff(\dff \varepsilon\dff)$\nnsp,\oss
then\dss
$t_{\trf \delta}\off =\off u\dff t_{\trf \varepsilon}\dff u^{\dff -\dff 1}$\dnsp.\oss
We see\sss that\dss the map\dss
$\varepsilon\off \longmapsto\off t_{\trf \varepsilon}$\dss
satisfies analogues of\dss conditions\qss ({\fff}i{\fff})\qss and\qss ({\fff}iii{\fff})\qss
from\dss the definition of\dss regular scaffoldings\qss
(see\dss Section\qss \ref{presentations}).\oss

For each\sss vertex\sss $w$\sss of\dss $X$\dss there\dss is\dss a\sss unique
element\sss $h_{\dff w}\qff \in\qff \mathcal{D}$\sss leaving\sss $w$\sss
invariant\sss and\dss rotating\dss $D$\sss 
coun\-ter-clock\-wise by\dss the angle\dss $2\dff \pi/3$\nnsp.\oss
In\dss particular\halfff,\pss $h_{\dff v}\off =\off h$\nnsp.\oss
Let\sss $\varepsilon$\sss be a\sss triangle-edge in\dss the small\trs
triangle corresponding\dss to $w$\nnsp.\oss
Clearly\halfff,\qss
if\trs the orientation of\sss $\varepsilon$ is\dss the coun\-ter-clockwise one,\oss
then\dss
$t_{\trf \varepsilon}\off =\off h_{\dff w}$\nsp,\oss
and\trs if\qss it\trs is\dss the clockwise one,\oss
then\dss
$t_{\trf \varepsilon}\off =\off h_{\dff w}^{\dff -\dff 1}$\dnsp.\oss
This\sss implies an analogue of\trs the condition\qss ({\fff}ii{\fff})\qss
from\dss the definition of\dss regular scaffoldings,\oss
but\dss we are not\dss going\dss to use it\sss
and\dss leave\sss the\sss task of\dss stating\dss it\dss to\sss the reader\halfff.\oss

\myuppar{The canonical\dss $\mathcal{G}$\dnsp-scaffolding.}
We would\dss like\sss to\sss lift\sss elements $t_{\trf \varepsilon}$\sss to elements
$\tau_{\dff \varepsilon}\qff \in\qff \mathcal{G}$ 
in\sss a canonical\dss way\halfff.\oss
Let $e$ be an edge of\sss $X$\nnsp.\oss
Then$e\off =\off \gamma\trf(\dff e_{\dff 1}\trf)$
for some
$\gamma\qff \in\qff \mathcal{G}$\nnsp.\oss
Let\vspace{3pt}\vspace{-0.2pt}
\[
\quad
g_{\dff e}
\off =\off 
\gamma\qff \cdot\qff 
g\qff \cdot\qff 
\gamma^{\dff -\dff 1}
\qff.
\]

\vspace{-9pt}\vspace{-0.2pt}
Let\dss us\sss check\dss that\dss $g_{\dff e}$\sss does not\sss depend
on\dss the choice of\dss $\gamma$\nnsp.\oss
If\trs
$e
\off =\off
\gamma'\trf(\dff e_{\dff 1}\trf)
\off =\off 
\gamma\trf(\dff e_{\dff 1}\trf)$\nnsp,\oss
then\dss
$\varphi\dff(\trf \gamma'\trf)\dff \cdot\dff \varphi\dff(\trf \gamma\trf)$\dss
leaves\sss $e_{\dff 1}$\sss invariant\sss and\dss hence either\dss\vspace{4.5pt}
\[
\quad
\varphi\dff(\trf \gamma'\trf)
\off =\dff\off 
\varphi\dff(\trf \gamma\trf)\dff,
\hspace*{1.4em}\mbox{or}\hspace*{1.5em}
\varphi\dff(\trf \gamma'\trf)
\off =\dff\off 
\varphi\dff(\trf \gamma\trf)\dff \cdot\dff s_{\dff 1}
\off =\off
\varphi\dff(\trf \gamma\dff \cdot\dff g\trf)
\qff.
\]

\vspace{-7.5pt}
If\trs 
$\varphi\dff(\trf \gamma'\trf)\off =\off \varphi\dff(\trf \gamma\trf)$\nnsp,\oss
then\dss
$\gamma'\off =\off \gamma\dff \cdot\dff \delta$\dss
for some\dss $\delta\qff \in\qff \ker\dff \varphi$\nnsp,\oss
and since\dss  
$\ker\dff \varphi$\dss 
is\dss contained\sss in\dss the center of\trs $\mathcal{G}$\nnsp,\oss
replacing\dss $\gamma$\dss by\dss $\gamma'$\dss does not\sss change\sss $g_{\dff e}$\nsp.\oss
Clearly\halfff,\oss replacing\dss $\gamma$\dss by\dss $\gamma\dff \cdot\dff g$\dss
also does not\sss change\sss $g_{\dff e}$\nsp.\oss
Therefore replacing\dss $\gamma$\dss by\dss $\gamma'$\dss 
does not\sss change\sss $g_{\dff e}$\sss  in\sss both cases.\oss

Next\halfff,\oss let\dss $w$\dss be a vertex of\trs $X$\nnsp.\oss
Then\dss $w\off =\off \gamma\trf(\dff v\trf)$\dss
for some
$\gamma\qff \in\qff \mathcal{G}$\nnsp.\oss
Let\vspace{3pt}\vspace{-0.2pt}
\[
\quad
r_{\fff w}
\off =\off 
\gamma\qff \cdot\qff 
r\qff \cdot\qff 
\gamma^{\dff -\dff 1}
\qff.
\]

\vspace{-9pt}\vspace{-0.2pt}
Let\dss us\sss check\dss that\dss $r_{\fff w}$\sss does not\sss depend
on\dss the choice of\dss $\gamma$\nnsp.\oss
If\trs
$w
\off =\off
\gamma'\trf(\dff v\trf)
\off =\off 
\gamma\trf(\dff v\trf)$\nnsp,\oss
then\dss
$\varphi\dff(\trf \gamma'\trf)\dff \cdot\dff \varphi\dff(\trf \gamma\trf)$\dss
leaves\sss $v$\sss invariant\sss and\dss hence either\dss
$\varphi\dff(\trf \gamma'\trf)
\off =\dff\off 
\varphi\dff(\trf \gamma\trf)$\nnsp,\oss or\vspace{4.5pt}
\[
\quad
\varphi\dff(\trf \gamma'\trf)
\off =\dff\off 
\varphi\dff(\trf \gamma\trf)\dff \cdot\dff h
\off =\off
\varphi\dff(\trf \gamma\dff \cdot\dff r\trf)\dff,
\hspace*{1.4em}\mbox{or}\hspace*{1.5em}
\varphi\dff(\trf \gamma'\trf)
\off =\dff\off 
\varphi\dff(\trf \gamma\trf)\dff \cdot\dff h^{\dff 2}
\off =\dff\off
\varphi\dff(\trf \gamma\dff \cdot\dff r^{\dff 2}\qff)
\qff.
\]

\vspace{-7.5pt}
If\trs 
$\varphi\dff(\trf \gamma'\trf)\off =\off \varphi\dff(\trf \gamma\trf)$\nnsp,\oss
then\dss
replacing\dss $\gamma$\dss by\dss $\gamma'$\dss does not\sss change\sss $r_{\fff w}$\sss
by\dss the same reason as above.\oss
Similarly\halfff,\oss replacing\dss $\gamma$\dss by\dss $\gamma\dff \cdot\dff r$\dss
or\dss $\gamma\dff \cdot\dff r^{\dff 2}$\dss
also does not\sss change\sss $r_{\fff w}$\nsp.\oss
Therefore in every\sss case replacing\dss $\gamma$\dss by\dss $\gamma'$\dss 
does not\sss change\sss $r_{\fff w}$\nsp.\oss

Now\sss we are ready\dss to define\sss elements $\tau_{\dff \varepsilon}$\nsp.\oss
Suppose\sss first\dss that\sss
$\varepsilon$\sss results\sss from orienting\sss
a\sss pentagon-edge $e$\nnsp.\oss
Let\sss $d$\dss be\sss the corresponding\sss edge of\trs $X$\dss
and\dss let\dss us\sss set\trs
$\tau_{\dff \varepsilon}\off =\off g_{\dff d}$\nsp.\oss
Suppose\sss now\dss that\sss
$\varepsilon$\sss results\sss from orienting\sss
a\dss triangle-edge $e$\nnsp.\oss
Let\sss $w$\sss be\sss the vertex of\trs $X$\dss
corresponding\dss to\sss the small\trs triangle having\sss $e$\sss as a side.\oss
In\dss this case we set\trs
$\tau_{\dff \varepsilon}
\off =\off
r_{\fff w}$\dss
if\qss the\sss orientation\sss of\dss $\varepsilon$\sss is\dss the counter-clockwise one,\oss
and\trs
$\tau_{\dff \varepsilon}
\off =\off
r_{\fff w}^{\dff -\dff 1}$\dss
if\qss the\sss orientation\sss of\dss $\varepsilon$\sss is\dss the clockwise one.\oss
It\dss follows\dss immediately\dss from\dss the definitions\sss that\trs
$\varphi\dff(\dff \tau_{\dff \varepsilon}\trf)
\off =\off
t_{\trf \varepsilon}$\dss
for every\sss oriented edge\sss $\varepsilon$\nnsp.\oss
The map\dss
$\varepsilon\off \longmapsto\off \tau_{\dff \varepsilon}$\dss 
is\dss the\qss \emph{canonical\dss $\mathcal{G}$\dnsp-scaffolding}\qss of\trs $Y$\dnsp.\oss

\myuppar{Elements\sss associated\dss with\dss paths\sss in $Y$\dnsp.}
Let\dss $p$\trs be a\sss path\sss in\sss $Y$\dnsp,\oss 
i.e.\qss a\sss sequence\vspace{3pt}
\[
\quad 
p
\off =\off
\{\trf v_{\dff i} \qff\}_{\trf 0\qff \leq\qff i\qff \leq\qff n}
\]

\vspace{-9pt}
of\dss vertices\sss $v_{\dff i}$\sss of\trs $Y$\dss
such\dss that\sss 
$v_{\dff i\dff -\dff 1}$\sss is\dss connected with\sss 
$v_{\dff i}$\sss by\sss an edge for\qss 
$1\qff \leq\qff i\qff \leq\qff n$\nnsp.\oss 
Let\sss $\varepsilon_{\dff i}$\dss be\sss the edge connecting\dss
$v_{\dff i\dff -\dff 1}$ with $v_{\dff i}$\dss
and oriented\dss in such a way\dss that\sss $v_{\dff i\dff -\dff 1}$\sss is\trs
its\dss origin,\oss
where\qss
$1\qff \leq\qff i\qff \leq\qff n$\nnsp.\oss
Let\qss 
$t_{\trf i}\off =\off t_{\qff \varepsilon_{\dff i}}$\qss
and\qss
$\tau_{\trf i}\off =\off \tau_{\qff \varepsilon_{\dff i}}$\nsp,\oss
and\qss let\vspace{3pt}
\[
\quad
\Pi_{\dff \mathcal{D}}\dff(\dff p\trf)
\off =\off
t_{\trf n}\dff \cdot\dff
t_{\trf n\dff -\dff 1}\dff \cdot\dff
\ldots\dff \cdot\dff
t_{\trf 2}\dff \cdot\dff
t_{\trf 1}
\hspace*{1.4em}\mbox{and}\hspace*{1.5em}
\Pi_{\dff \mathcal{G}}\dff(\dff p\trf)
\off =\off
\tau_{\trf n}\dff \cdot\dff
\tau_{\trf n\dff -\dff 1}\dff \cdot\dff
\ldots\dff \cdot\dff
\tau_{\trf 2}\dff \cdot\dff
\tau_{\trf 1}
\qff.
\]

\vspace{-9pt}

Obviously\halfff,\oss the homomorphism\dss
$\varphi\dff \colon\dff
\mathcal{G}\qff \ttoo\qff \mathcal{D}$\dss
takes\dss
$\Pi_{\dff \mathcal{D}}\dff(\dff p\trf)$\dss
to\dss
$\Pi_{\dff \mathcal{G}}\dff(\dff p\trf)$\nnsp.\oss
Clearly\halfff,\pss\vspace{3.5pt}
\[
\quad
t_{\trf n}\dff \cdot\dff
t_{\trf n\dff -\dff 1}\dff \cdot\dff
\ldots\dff \cdot\dff
t_{\trf 2}\dff \cdot\dff
t_{\trf 1}\trf(\trf v_{\trf 0}\trf)
\off =\off
v_{\trf n}
\]

\vspace{-8.5pt}
and\dss hence\dss $\Pi_{\dff \mathcal{D}}\dff(\dff p\trf)$\dss
is\dss the unique element\sss of\trs $\mathcal{D}$\dss
taking\sss $v_{\dff 0}$\sss to\sss $v_{\dff n}$\nsp.\oss
In\dss particular\halfff,\oss
if\trs $p$\sss is\dss a\sss loop,\oss i.e.\pss if\qss
$v_{\dff n}\off =\off v_{\dff 0}$\nsp,\oss
then\dss
$\Pi_{\dff \mathcal{D}}\dff(\dff p\trf)
\off =\off
1$\nnsp.\oss
Since\dss $\ker\dff \varphi$\dss is\dss generated\dss by\sss $z$\nnsp,\oss
in\dss this case\vspace{3pt}
\[
\quad
\Pi_{\dff \mathcal{G}}\dff(\dff p\trf)
\off =\off
z^{\dff k}
\qff,
\]

\vspace{-9pt}
for some integer\dss $k$\nnsp.\oss 
Suppose\sss now\dss that\qss
$p'
\off =\off
\{\trf w_{\dff i} \qff\}_{\trf 0\qff \leq\qff i\qff \leq\qff m}$\qss
is\dss a\sss path starting\sss at\dss the endpoint\sss of\dss $p$\nnsp,\oss
i.e.\qss such\dss that\dss
$w_{\trf 0}\off =\off v_{\dff n}$\nsp.\oss
Let\trs $p' p$\dss be\sss the composition of\dss $p$\sss and\sss $p'$\nnsp,\oss
i.e.\qss the path\vspace{3pt}
\[
\quad
v_{\trf 0}\dff,\off
v_{\trf 1}\dff,\off
\ldots\dff,\off
v_{\fff n}\dff,\off
w_{\trf 1}\dff,\off
\ldots\dff,\off
w_{\fff m}
\qff.
\]

\vspace{-9pt}
Then,\oss obviously\halfff,\oss\vspace{3pt}
\[
\quad
\Pi_{\dff \mathcal{D}}\dff(\dff p' p\trf)
\off =\dff\off
\Pi_{\dff \mathcal{D}}\dff(\dff p'\trf)
\qff \cdot\pff 
\Pi_{\dff \mathcal{D}}\dff(\dff p\trf)
\hspace*{1.4em}\mbox{and}\hspace*{1.5em}
\Pi_{\dff \mathcal{G}}\dff(\dff p' p\trf)
\off =\dff\off
\Pi_{\dff \mathcal{G}}\dff(\dff p'\trf)
\qff \cdot\pff 
\Pi_{\dff \mathcal{G}}\dff(\dff p\trf)
\qff.
\]

\vspace{-7.5pt}
Next\halfff,\oss let\dss $p^{\dff -\dff 1}$\dss be\sss the\sss inverse path\sss of\dss $p$\nnsp,\oss
i.e.\qss the sequence\dss
$v_{\fff n}\dff,\off 
v_{\fff n\dff -\dff 1}\dff,\off 
\ldots\dff,\off 
v_{\trf 1}\dff,\off
v_{\trf 0}$\nsp.\oss
Then\vspace{4.5pt}
\[
\quad
\Pi_{\dff \mathcal{D}}\dff\left(\dff p^{\dff -\dff 1}\trf\right)
\qff \cdot\pff 
\Pi_{\dff \mathcal{D}}\dff\bigl(\dff p\trf\bigr)
\off =\dff\off
\Pi_{\dff \mathcal{D}}\dff\left(\dff p^{\dff -\dff 1} p\trf\right)
\off =\off
1
\qff
\]

\vspace{-7.5pt}
because\dss $p^{\dff -\dff 1} p$\dss is\dss a\sss loop.\oss
The situation\sss is\dss more interesting\dss for\dss $\Pi_{\dff \mathcal{G}}$\nsp.\oss

\mypar{Lemma.}{g-cancellations}
\emph{In\dss the above situation,\oss}\vspace{4.5pt}
\[
\quad
\Pi_{\dff \mathcal{G}}\dff(\dff p^{\dff -\dff 1}\trf)
\qff \cdot\pff 
\Pi_{\dff \mathcal{G}}\dff(\dff p\trf)
\off =\dff\off
\Pi_{\dff \mathcal{G}}\dff(\dff p^{\dff -\dff 1} p\trf)
\off =\off
z^{\dff k}
\qff,
\]

\vspace{-8.5pt}
\emph{where\dss $k$\dss is\dss the number of\qss pentagon-edges
among\dss the edges\dss
$\varepsilon_{\dff 1}\dff,\off
\varepsilon_{\dff 2}\dff,\off
\ldots\dff,\off
\varepsilon_{\dff n}$\nsp.\oss}\vspace{1.125pt}

\proof
For each\dss $i\off =\off 1\fff,\pff 2\fff,\pff \ldots\fff,\pff n$\qss
let\trs $\varepsilon_{\dff i}^{\dff -\dff 1}$\dss be\sss the edge connecting\dss
$v_{\dff i}$ with $v_{\dff i\dff -\dff 1}$
and\dss oriented\dss in such a way\dss that\sss $v_{\dff i}$\sss is\trs
its\dss origin,\oss
and\dss let\qss
$\overline{\tau}_{\dff i}
\off =\off
\tau_{\qff \varepsilon_{\dff i}^{\dff -\dff 1}}$.\oss
Then\vspace{3.5pt}
\[
\quad
\Pi_{\dff \mathcal{G}}\dff(\dff p^{\dff -\dff 1} p\trf)
\off =\dff\off
\overline{\tau}_{\trf 1}\dff \cdot\qff
\ldots\qff \cdot\qff
\overline{\tau}_{\trf n\dff -\dff 1}\dff \cdot\qff
\overline{\tau}_{\dff n}\off \cdot\off
\tau_{\trf n}\dff \cdot\qff
\tau_{\trf n\dff -\dff 1}\dff \cdot\qff
\ldots\qff \cdot\qff
\tau_{\trf 1}
\qff.
\]

\vspace{-8.5pt}
If\trs $\varepsilon_{\dff i}$\dss is\dss a\dss triangle-edge,\oss
then\qss
$\overline{\tau}_{\dff i}
\qff \cdot\qff
\tau_{\dff i}
\off =\off
1$\nnsp,\oss
as\dss it\dss immediately\dss follows\sss from\dss the definitions.\vspace{1.125pt}

In contrast\halfff,\oss
if\trs $\varepsilon_{\dff i}$\dss is\dss a\dss pentagon-edge,\oss
then\qss 
$\overline{\tau}_{\dff i}
\off =\off
\tau_{\dff i}
\off =\off
\gamma\qff \cdot\qff 
g\qff \cdot\qff 
\gamma^{\dff -\dff 1}$\qss
for some\qss $\gamma\qff \in\qff \mathcal{G}$\qss
and\vspace{3pt}
\[
\quad
\overline{\tau}_{\dff i}
\qff \cdot\qff
\tau_{\dff i}
\off =\off
\gamma\qff \cdot\qff 
g\qff \cdot\qff 
\gamma^{\dff -\dff 1}
\qff \cdot\qff
\gamma\qff \cdot\qff 
g\qff \cdot\qff 
\gamma^{\dff -\dff 1}
\]

\vspace{-34pt}
\[
\quad
\phantom{\overline{\tau}_{\dff i}
\qff \cdot\qff
\tau_{\dff i}
\off }
=\off
\gamma\qff \cdot\qff 
g\qff \cdot\qff 
g\qff \cdot\qff 
\gamma^{\dff -\dff 1}
\]

\vspace{-34pt}
\[
\quad
\phantom{\overline{\tau}_{\dff i}
\qff \cdot\qff
\tau_{\dff i}
\off }
=\off
\qff \cdot\qff 
z\qff \cdot\qff 
\gamma^{\dff -\dff 1}
\off =\off
z\qff \cdot\qff
\gamma\qff \cdot\qff  
\gamma^{\dff -\dff 1}
\off =\dff\off
z
\qff,
\]

\vspace{-7.5pt}
where we used\dss the fact\dss that\dss $z$\dss 
belongs\sss to\sss the center of\trs $\mathcal{G}$\nnsp.\oss

Now\dss we can consecutively\qss ``cancel''\qss the products\dss
$\overline{\tau}_{\dff i}
\qff \cdot\qff
\tau_{\dff i}$\dss
in\dss the above expression\dss for\dss the element\trs
$\Pi_{\dff \mathcal{G}}\dff(\dff p^{\dff -\dff 1} p\trf)$\nnsp,\oss
starting\dss with\dss
$\overline{\tau}_{\dff n}
\qff \cdot\qff
\tau_{\dff n}$.\oss
If\trs $\varepsilon_{\dff i}$\dss is\dss a\dss triangle-edge,\oss
then\qss
$\overline{\tau}_{\dff i}
\qff \cdot\qff
\tau_{\dff i}
\off =\off
1$\qss
and\dss
$\overline{\tau}_{\dff i}
\qff \cdot\qff
\tau_{\dff i}$\dss
actually\sss cancels.\oss
If\trs $\varepsilon_{\dff i}$\dss is\dss a\dss pentagon-edge,\oss
then\qss
$\overline{\tau}_{\dff i}
\qff \cdot\qff
\tau_{\dff i}
\off =\off
z$\qss
and\qss
$\overline{\tau}_{\dff i}
\qff \cdot\qff
\tau_{\dff i}$\qss
should\dss be replaced\dss by\dss $z$\dss
and\dss then\dss $z$\dss moved\dss to,\pss say\halfff,\oss
the\sss left\sss of\dss the whole product\halfff.\oss
The\sss latter\dss is\dss possible since\dss $z$\dss 
belongs\sss to\sss the center\halfff.\oss
In\dss the end\sss we will\dss be\sss left\dss with several\dss
factors equal\dss to\dss $z$\nnsp,\oss
one factor\dss for each\dss pentagon-edge.\oss
The\sss lemma\sss follows.\oss  \eproof\vspace{1pt}

\mypar{Lemma.}{rotating-loops}
\emph{Suppose\sss that\qss
$l
\off =\off
\{\trf v_{\dff i} \qff\}_{\trf 0\qff \leq\qff i\qff \leq\qff n}$\qss
is\dss a\dss loop,\oss
and\dss let\sss $l'$\sss be a\sss loop of\qss the form}\vspace{3pt}
\[
\quad
v_{\dff i}\dff,\off
v_{\dff i\dff +\dff 1}\dff,\off
\ldots\dff,\off
v_{\dff n}\dff,\off
v_{\trf 1}\dff,\off
v_{\trf 2}\dff,\off
\ldots\dff,\off
v_{\dff i\dff -\dff 1}\dff,\off
v_{\dff i}
\off,
\]

\vspace{-9pt}
\emph{where\qss 
$0\qff \leq\qff i\qff \leq\qff n\qff -\qff 1$\nnsp.\oss
Then\pss
$\Pi_{\dff \mathcal{G}}\dff(\trf l'\qff)
\off =\off
\Pi_{\dff \mathcal{G}}\dff(\trf l\trf)$\nnsp.\oss}\vspace{1pt}

\proof
Let\trs $\gamma\qff \in\qff \mathcal{G}$\dss be such\dss that\trs
$\gamma\trf(\dff v_{\trf 0}\trf)\off =\dss\off v_{\dff i}$\nsp.\oss
Then\vspace{3pt}
\[
\quad
\Pi_{\dff \mathcal{G}}\dff(\trf l'\qff)
\off =\off
\gamma\qff \cdot\qff 
\Pi_{\dff \mathcal{G}}\dff(\trf l\trf)\qff \cdot\qff 
\gamma^{\dff -\dff 1}
\qff.
\]

\vspace{-9pt}
Since\sss $l$\sss is\dss a\sss loop,\pss
$\Pi_{\dff \mathcal{G}}\dff(\dff l\trf)
\off =\off
z^{\dff k}$\dss
for some\sss integer\sss $k$\sss
and\dss hence\dss
$\Pi_{\dff \mathcal{G}}\dff(\dff l\trf)$\dss
belongs\sss to\sss the center\sss of\trs $\mathcal{G}$\nnsp.\oss
It\dss follows\dss that\trs
$\Pi_{\dff \mathcal{G}}\dff(\trf l'\qff)
\off =\off
\Pi_{\dff \mathcal{G}}\dff(\trf l\trf)$\nnsp.\oss  \eproof\vspace{1pt}

\mypar{Lemma.}{boundaries-of-faces}
\emph{Let\trs $F$\sss be a face of\pss $T$\dnsp,\oss
and\dss let\sss $l$\dss be\sss a\sss loop\dss following\dss
the vertices of\pss $F$\sss in\dss the clockwise\dss direction.\oss
Then\pss
$\Pi_{\dff \mathcal{G}}\dff(\trf l\qff)
\off =\dff\off
z$\nnsp.\oss}

\proof
Let\trs $F$\dss be\sss a\sss triangular\dss face of\trs $T$\dnsp.\oss
Then\dss $F$\dss results from cutting\sss off\dss a\sss triangular\sss
pyramid at\sss some vertex\sss $w$\sss of\trs the dodecahedron\dss $D$\dss
and\dss
$\tau_{\dff \varepsilon}\off =\off r_{\fff w}^{\dff -\dff 1}$\dss
for every\sss clockwise oriented edge\sss $\varepsilon$\sss
connecting\dss two vertices of\trs $F$\nnsp.\oss
It\dss follows\dss that\vspace{3pt}
\[
\quad
\Pi_{\dff \mathcal{G}}\dff(\trf l\qff)
\off =\dff\off
\gamma\qff \cdot\qff 
r^{\dff -\dff 3}\qff \cdot\qff 
\gamma^{\dff -\dff 1}
\]

\vspace{-9pt}
for some\qss $\gamma\qff \in\qff \mathcal{G}$\nnsp.\oss
But\trs 
$r^{\dff -\dff 3}\off =\dff\off z$\dss
and\dss hence\trs
$\Pi_{\dff \mathcal{G}}\dff(\trf l\qff)
\off =\dff\off
z$\dss
in\dss this case.\oss

Let\dss us\sss consider\sss now\dss the\sss ten-sided\dss faces of\trs $T$\dss
contained\sss in\dss pentagonal\trs faces of\trs $D$\nnsp.\oss
Let\trs $F\fff,\off F'$\dss
be\sss two such\dss faces,\oss
and\dss let\dss $l\fff,\pff l'$\qss 
be\sss loops following\dss clockwise\sss
the vertices of\trs $F\fff,\off F'$\dss
respectively\halfff.\oss
Then\dss
$F'\off =\off \gamma\trf(\trf F\trf)$\dss
for some\dss 
$\gamma\qff \in\qff \mathcal{G}$\nnsp.\oss
Since\sss
both\dss loops\sss $l'$\dss and\trs $\gamma\trf(\trf l\trf)$\dss
follow\dss the vertices of\trs $F'$\dss clockwise,\oss
they\sss are related as\sss the loops\ $p'$\sss and\dss $p$\dss
in\dss Lemma\qss \ref{rotating-loops}\qss
and\dss hence\trs
$\Pi_{\dff \mathcal{G}}\dff(\trf l'\qff)
\off =\off
\Pi_{\dff \mathcal{G}}\dff(\trf \gamma\trf(\trf l\trf)\trf)$\nnsp.\oss
At\dss the same\sss time,\oss\vspace{3pt}
\[
\quad
\Pi_{\dff \mathcal{G}}\dff\left(\trf \gamma\trf(\trf l\trf)\qff\right)
\off =\off
\gamma\qff \cdot\qff 
\Pi_{\dff \mathcal{G}}\dff(\trf l\qff)\qff \cdot\qff 
\gamma^{\dff -\dff 1}
\qff.
\]

\vspace{-9pt}
Since\dss $\Pi_{\dff \mathcal{G}}\dff(\trf l\trf)$\dss 
belongs\sss to\sss the center\sss of\trs $\mathcal{G}$\nnsp,\oss
it\dss follows\dss that\trs
$\Pi_{\dff \mathcal{G}}\dff(\trf l'\qff)
\off =\off
\Pi_{\dff \mathcal{G}}\dff(\trf \gamma\trf(\trf l\trf)\trf)
\off =\off 
\Pi_{\dff \mathcal{G}}\dff(\trf l\trf)$\nnsp.\oss
Therefore,\oss it\dss is\dss sufficient\dss to consider only\sss
one\sss such\dss face.\oss

Let\trs $F$\dss be\sss the\sss ten-sided\dss face contained\sss in\dss the face\dss
$v\dff w_{\dff 1}\dff a\trf b\dff w_{\dff 3}$\dss
of\trs $D$\dss
(see\dss Section\qss \ref{examples}).\oss
It\dss  is\dss convenient\dss to rename\sss the vertices\qss
$v\fff,\off w_{\dff 1}\dff,\off a\fff,\off b\fff,\off w_{\dff 3}$\qss
as\qss
$v_{\trf 1}\dff,\off v_{\trf 2}\dff,\off v_{\trf 3}\dff,\off v_{\trf 4}\fff,\off v_{\trf 5}$\qss
respectively\halfff.\oss
Recall\dss that\sss $x_{\dff 1}$\dss is\dss a\sss point\sss on\dss the edge\dss
$e_{\dff 1}\off =\off v\dff w_{\dff 1}$\dss near\dss $v$\nnsp.\oss
Hence\sss $x_{\dff 1}$\dss is\dss a vertex of\trs $F$\dnsp.\oss
Let\vspace{3pt}
\[
\quad
x_{\dff 1}\dff,\off
y_{\dff 1}\dff,\off
x_{\dff 2}\dff,\off
y_{\dff 2}\dff,\off
\ldots\dff,\off
x_{\dff 5}\dff,\off
y_{\dff 5}
\]

\vspace{-9pt}
be\sss the vertices of\trs $F$\dss listed\dss in\dss the clockwise order
along\dss the boundary\sss of\trs $F$\nnsp.\oss
See\sss the picture.\vspace{12pt}

\begin{figure}[h!]
\hspace*{10.5em}
\includegraphics[width=0.4757\textwidth]{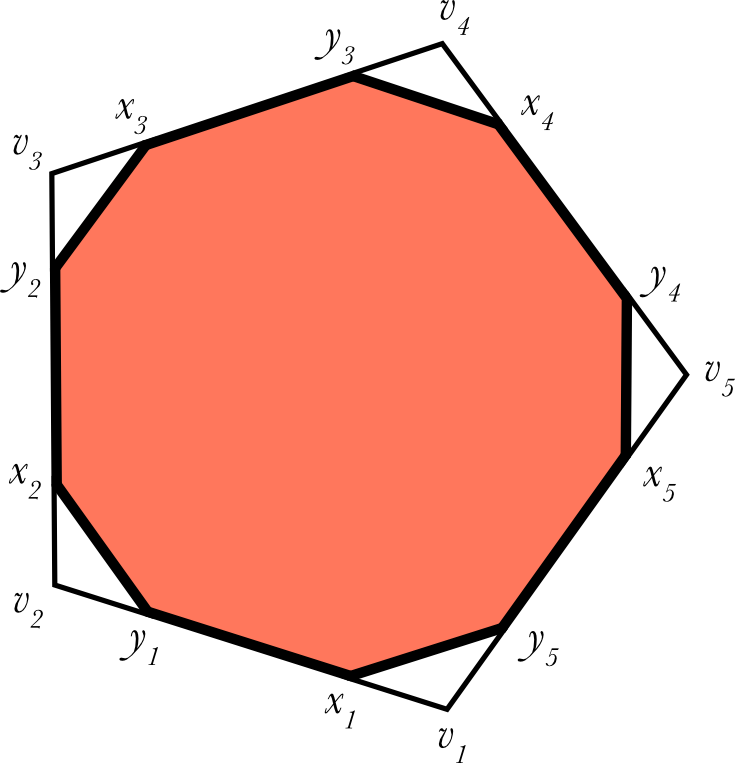}
\end{figure}\vspace{12pt}

Let\trs $c_{\dff i}$\sss is\dss the edge\dss 
$x_{\dff i}\qff y_{\dff i}$\qss
of\trs $Y$\dnsp.\qff\oss
Let\qss
$g_{\dff i}\off =\off g_{\trf c_{\trf i}}$\qss
and\qss
$r_{\dff i}\off =\off r_{\dff v_{\dff i}}$.\oss
Then\vspace{3pt}
\[
\quad
y_{\dff i}
\off =\dff\off 
g_{\dff i}\dff(\dff x_{\dff i}\trf)
\hspace*{1.4em}\mbox{and}\hspace*{1.5em}
x_{\dff i\dff +\dff 1}
\off =\dff\off 
g_{\dff i}\dff(\trf y_{\dff i}\trf)
\]

\vspace{-9pt}
for every\trs $i\off =\off 1\fff,\pff 2\fff,\pff \ldots\fff,\pff 5$\nnsp,\oss
where\sss the subscript\trs $i\qff +\qff 1$\dss is\dss interpreted\dss modulo\dss $5$\nnsp.\oss
Clearly\halfff,\pss $g_{\dff 1}\off =\off g$\dss and\dss $r_{\dff 1}\off =\off r$\nnsp.\oss
Without\sss any\dss loss of\dss generality\dss we can assume\sss that\dss 
the\sss loop\sss $l$\dss begins\sss at\sss $y_{\dff 5}$\sss
and\dss hence\sss $l$\dss is\dss the\sss loop\dss
$y_{\dff 5}\fff,\pff
x_{\dff 1}\fff,\pff
y_{\dff 1}\fff,\pff
x_{\dff 2}\fff,\pff
\ldots\dff,\pff
x_{\dff 5}\fff,\pff
y_{\dff 5}$\nsp.\oss
Then\vspace{3pt}
\[
\quad
\Pi_{\dff \mathcal{G}}\dff(\trf l\qff)
\off =\off
g_{\dff 5}\qff \cdot\qff r_{\dff 5}
\qff \cdot\pff
\ldots
\qff \cdot\pff
g_{\dff 2}\qff \cdot\qff r_{\dff 2}
\qff \cdot\pff
g_{\dff 1}\qff \cdot\qff r_{\dff 1}
\off.
\]

\vspace{-9pt}
The action of\trs the element\trs $g\dff \cdot\dff r\qff \in\qff \mathcal{G}$\dss on\dss $D$\dss
is\dss the same as\sss the action of\trs the element\trs 
$g\dff \cdot\dff h\qff \in\qff \mathcal{D}$\dss from\dss
Section\qss \ref{examples}.\oss
Hence,\oss after\dss renaming\dss the vertices,\oss
the equalities\qss (\ref{rotating-a-face})\qss imply\dss that\vspace{3pt}
\[
\quad
g\dff \cdot\dff r\qff(\dff v_{\dff i}\trf)
\off =\off
v_{\dff i\dff +\dff 1}
\hspace*{1.4em}\mbox{and}\hspace*{1.5em}
g\dff \cdot\dff r\qff(\dff c_{\dff i}\trf)
\off =\off
c_{\dff i\dff +\dff 1}
\]

\vspace{-9pt}
for every\trs $i\off =\off 1\fff,\pff 2\fff,\pff \ldots\fff,\pff 5$\nnsp.\oss
It\dss follows\dss that\vspace{3pt}
\[
\quad
r_{\dff i\dff +\dff 1}
\off =\off
(\trf g\dff \cdot\dff r\trf)^{\dff i}\qff \cdot\qff
r\qff \cdot\qff
(\trf g\dff \cdot\dff r\trf)^{\dff -\dff i}
\hspace*{1.4em}\mbox{and}\hspace*{1.5em}
g_{\dff i\dff +\dff 1}
\off =\off
(\trf g\dff \cdot\dff r\trf)^{\dff i}\qff \cdot\qff
g\qff \cdot\qff
(\trf g\dff \cdot\dff r\trf)^{\dff -\dff i}
\]

\vspace{-9pt}
for every\trs $i\off =\off 0\fff,\pff 1\fff,\pff \ldots\fff,\pff 4$\nnsp.\oss
In\dss turn,\oss this\sss implies\sss that\vspace{3pt}
\[
\quad
g_{\dff i\dff +\dff 1}\dff \cdot\dff r_{\dff i\dff +\dff 1}
\off =\off
(\trf g\dff \cdot\dff r\trf)^{\dff i}\qff \cdot\qff
(\trf g\dff \cdot\dff r\trf)\qff \cdot\qff
(\trf g\dff \cdot\dff r\trf)^{\dff -\dff i}
\off =\off
g\dff \cdot\dff r
\qff
\]

\vspace{-9pt}
for every\trs $i\off =\off 0\fff,\pff 1\fff,\pff \ldots\fff,\pff 4$\dss
and\dss hence\qss \vspace{3pt}
\[
\quad
\Pi_{\dff \mathcal{G}}\dff(\trf l\qff)
\off =\off
(\trf g\dff \cdot\dff r\trf)^{\dff 5}
\off =\off
z 
\qff.
\]

\vspace{-9pt}
This completes\sss the proof\dss of\trs the\sss lemma.\oss  \eproof

\myuppar{Cyclic reorderings and a\sss multiplication of\trs loops.}
Let\dss us\sss say\dss that\dss a\sss loop\dss $l'$\trs is\dss a\qss
\emph{cyclic\sss reordering}\qss of\dss a\dss loop\dss $l$\qss
if\trs the\sss loops\trs $l\fff,\pff l'$\dss are related\sss as\sss in\dss
Lemma\qss \ref{rotating-loops}.\oss
By\qss Lemma\qss \ref{rotating-loops}\qss a\sss cyclic reordering\sss of\dss
a\sss loop\dss $l$\dss does not\dss change\trs
$\Pi_{\dff \mathcal{G}}\dff(\trf l\trf)$\nnsp.\oss

Suppose\sss that\trs $l\fff,\pff l'$\trs are\sss two\sss loops
and\dss $p$\sss is\dss a subpath of\dss a cyclic reordering\sss of\trs $l$\dss 
such\dss that\sss $p^{\dff -\dff 1}$\sss 
is\dss a subpath of\dss a cyclic reordering\sss of\trs $l'$\nnsp.\oss
After\dss replacing\dss $l$\dss and\trs $l'$\dss by\dss their cyclic reorderings,\oss
if\dss necessary\halfff,\oss
we can assume\sss that\trs
$l\off =\off m\fff p$\dss and\dss $l'\off =\off p^{\dff -\dff 1} m'$\dss
for some paths\dss $m\fff,\pff m'$\nnsp.\oss 
This allows\sss to define a new\dss loop\dss
$l\dff \circ_{\dff p}\dff l'
\off =\off
m\dff m'$\nnsp.\oss
Lemma\qss \ref{g-cancellations}\qss implies\sss that\vspace{3pt}\vspace{1.25pt}
\[
\quad
\Pi_{\dff \mathcal{G}}\dff(\trf l'\trf)
\qff \cdot\pff 
\Pi_{\dff \mathcal{G}}\dff(\trf l\trf)
\off =\off
\Pi_{\dff \mathcal{G}}\dff(\trf m\trf)
\qff \cdot\pff 
\Pi_{\dff \mathcal{G}}\dff(\trf p\trf)
\qff \cdot\pff 
\Pi_{\dff \mathcal{G}}\dff(\trf p^{\dff -\dff 1}\trf)
\qff \cdot\pff
\Pi_{\dff \mathcal{G}}\dff(\trf m'\trf) 
\]

\vspace{-36pt}
\[
\quad
\phantom{\Pi_{\dff \mathcal{G}}\dff(\trf l'\trf)
\qff \cdot\pff 
\Pi_{\dff \mathcal{G}}\dff(\trf l\trf)
\off }
=\off
\Pi_{\dff \mathcal{G}}\dff(\trf m\trf)
\qff \cdot\pff 
z^{\dff k}
\qff \cdot\pff
\Pi_{\dff \mathcal{G}}\dff(\trf m'\trf) 
\]

\vspace{-36pt}
\[
\quad
\phantom{\Pi_{\dff \mathcal{G}}\dff(\trf l'\trf)
\qff \cdot\pff 
\Pi_{\dff \mathcal{G}}\dff(\trf l\trf)
\off }
=\off
z^{\dff k}
\qff \cdot\pff
\Pi_{\dff \mathcal{G}}\dff(\trf m\trf)
\qff \cdot\pff 
\Pi_{\dff \mathcal{G}}\dff(\trf m'\trf) 
\off =\off
z^{\dff k}
\qff \cdot\pff
\Pi_{\dff \mathcal{G}}\dff(\trf l\dff \circ_{\dff p}\dff l'\trf)
\qff,
\]

\vspace{-9pt}\vspace{1.25pt}
where\dss $k$\dss is\dss the number of\qss pentagon-edges
among\dss the edges of\trs $p$\nnsp.\oss
It\dss may\dss happen\dss that\trs
$l'\off =\off l^{\dff -\dff 1}$\dss and\dss $p\off =\off l$\nnsp.\oss
In\dss this case\trs 
$l\dff \circ_{\dff p}\dff l'$\qss
is\dss a\sss loop of\trs length $0$ and\trs
$\Pi_{\dff \mathcal{G}}\dff(\trf l\dff \circ_{\dff p}\dff l'\trf)\off =\off 1$\nnsp.\oss

\myuppar{A\sss proof\dss of\qss Coxeter's\dss implication.}
The idea\dss is\dss to multiply\dss identities of\qss
Lemma\qss \ref{boundaries-of-faces}\qss corresponding\dss to all\dss faces of\trs $T$\dss
and\dss then use\sss the cancellation\dss process of\qss
Lemma\qss \ref{g-cancellations}.\oss
The main\dss task\dss is\dss to arrange\sss the multiplication\sss in\sss
such a\sss way\dss that\trs Lemma\qss \ref{g-cancellations}\qss applies.\oss

Let\trs $B$\dss be\sss the boundary\sss of\trs $T$\nnsp.\oss
Suppose\sss that\dss $U\qff \subset\qff B$\dss 
is\dss the union of\trs several\dss faces of\trs $T$\dss
and\dss is\dss homeomorphic\sss to a disc.\oss
Then\dss $U$\dss is\dss bounded\sss in\dss $B$\dss by\sss a\sss polygonal\sss circle\dss $L$\nnsp.\oss
The clockwise orientation of\trs $U$\dss defines an orientation of\trs $L$\nnsp,\oss
which we will\dss call\trs the\dss \emph{$U$\nsp\dnsp-clockwise}\qss orientation.\oss
Choosing\sss a vertex of\trs $L$\dss as\sss the initial\sss vertex and\dss following\dss
the circle\dss $L$\dss clockwise\sss leads\sss to a\sss loop in\dss $Y$\dnsp.\oss
Choosing\sss another\sss vertex as\sss the initial\dss vertex\dss leads\sss to a\sss 
cyclic reordering\sss of\trs this\sss loop.\oss
So,\oss up\sss to a cyclic reordering\dss this\sss loop\dss is\dss well-defined.\oss
We will\sss denote\sss it\dss by\sss $l_{\trf U}$\nnsp.\oss

Suppose now\dss that\trs $U'\qff \subset\qff B$\dss 
is\dss another union of\trs several\dss faces of\dss $T$\sss
and\dss is\dss also homeomorphic\sss to a disc.\oss
Suppose\sss that\halfff,\oss moreover\halfff,\oss
the intersection\dss $U\dff \cap\dff U'$\dss 
is\dss a\sss polygonal\dss path\dss $P$\dss in\dss $B$\nnsp,\oss
i.e.\qss is\dss the union of\dss several\sss geometric edges of\dss $T$\sss
which\dss is\dss homeomorphic either\sss to a segment\sss or\dss to a circle.\oss
If\dss $P$\sss is\dss homeomorphic\sss to a circle,\oss
then\dss $U\dff \cap\dff U'\off =\off L$\dss
and\dss $U\dff \cup\dff U'\off =\off B$\nnsp.\oss
Otherwise\dss $U\dff \cup\dff U'$\dss is\dss homeomorphic\sss to a disc.\oss
Let\dss $p$\dss be\sss the path\sss in\dss $Y$\dss obtained\dss by\dss following\dss
the vertices of\dss $P$\sss in\dss the $U$\nsp\dnsp-clockwise direction.\oss
Then\dss $p^{\dff -\dff 1}$\dss is\dss the path obtained\dss by\dss following\dss
the vertices of\dss $P$\sss in\dss the $U'$\nsp\dnsp-clockwise direction
and\vspace{1.5pt}
\[
\quad
l_{\trf U\qff \cup\qff U'}
\off =\off
l_{\trf U}\qff \circ_{\dff p}\qff l_{\trf U'}
\off.
\]

\vspace{-10.5pt}
It\dss follows\dss that\vspace{1.5pt}
\[
\quad
\Pi_{\dff \mathcal{G}}\dff\bigl(\trf l_{\trf U}\trf\bigr)
\qff \cdot\qff
\Pi_{\dff \mathcal{G}}\dff\bigl(\trf l_{\trf U'}\trf\bigr)
\off =\off
z^{\dff k}
\qff \cdot\qff
\Pi_{\dff \mathcal{G}}\dff\left(\trf l_{\trf U}\qff \circ_{\dff p}\qff l_{\trf U'}\trf\right)
\off =\off
z^{\dff k}
\qff \cdot\qff
\Pi_{\dff \mathcal{G}}\dff\bigl(\trf l_{\trf U\qff \cup\qff U'}\trf\bigr)
\off.
\]

\vspace{-10.5pt}
where\dss $k$\dss is\dss the number of\qss pentagon-edges
among\dss the edges of\trs $p$\nnsp.\oss
If\trs $U\dff \cup\dff U'\off =\off B$\nnsp,\oss
then\dss $l_{\trf U\qff \cup\qff U'}$\dss is\dss a\sss loop of\trs length\sss $0$\sss
and\dss 
$\Pi_{\dff \mathcal{G}}\dff(\trf l_{\trf U\qff \cup\qff U'}\trf)
\off =\off 
1$\nnsp.\oss

The polyhedron\dss $T$\dss has $20$\sss triangular\sss faces 
and\sss $12$\sss ten-sides faces,\oss
corresponding\dss to\sss the vertices 
and\dss faces of\trs the dodecahedron\dss $D$\nnsp.\oss
It\dss is\dss easy\dss to see\sss that\dss the faces\dss
$F_{\dff 1}\dff,\off
F_{\trf 2}\dff,\off
\ldots\dff,\off
F_{\trf 32}$\dss
can\dss be numbered\sss in such a\sss way\dss that\qss 
$U_{\dff i}
\off =\off
F_{\dff 1}\qff \cup\qff \ldots\qff \cup\qff F_{\dff i}$\qss
is\dss homeomorphic\sss to a disc and\dss
$U_{\dff i}\dff \cap\dff F_{\dff i\dff +\dff 1}$\dss
is\dss homeomorphic\sss to a segment\dss 
for every\dss $i\qff \leq\qff 31$\nnsp,\oss
the intersection\dss
$U_{\trf 31}\dff \cap\qff F_{\trf 32}$\dss is\dss equal\dss to\sss the boundary\sss
of\trs $F_{\trf 32}$\nsp,\oss
and\dss $U_{\trf 32}\off =\off B$\nnsp.\oss
Moreover\halfff,\oss every\sss pentagon-edge\dss is\dss con\-tained\sss in one of\trs
the intersections\dss $U_{\dff i}\dff \cap\dff F_{\dff i\dff +\dff 1}$\nsp.\oss
An\dss induction\dss shows\sss that\vspace{3.5pt}
\[
\quad
\Pi_{\dff \mathcal{G}}\dff\bigl(\trf l_{\qff F_{\trf 1}}\trf\bigr)
\qff \cdot\qff
\Pi_{\dff \mathcal{G}}\dff\bigl(\trf l_{\qff F_{\trf 2}}\trf\bigr)
\qff \cdot\qff
\ldots
\qff \cdot\qff
\Pi_{\dff \mathcal{G}}\dff\bigl(\trf l_{\qff F_{\trf i}}\trf\bigr)
\off =\dff\off
z^{\dff k_{\trf i}}
\qff \cdot\pff
\Pi_{\dff \mathcal{G}}\dff\bigl(\trf l_{\qff U_{\trf i}}\trf\bigr)
\qff,
\]

\vspace{-8.5pt}
where\dss $k_{\trf i}$\dss is\dss the number of\qss pentagon-edges contained\dss in\dss
the interior of\trs $U_{\dff i}$\nsp.\oss
By\dss Lemma\qss \ref{boundaries-of-faces}\qss every\dss factor on\dss the left\dss
is\dss equal\dss to\sss $z$\nnsp.\oss
Since\dss 
$\Pi_{\dff \mathcal{G}}\dff\left(\trf l_{\qff U_{\trf 32}}\trf\right)\off =\off 1$\dss
and\dss there are $30$\sss pentagon-edges,\oss
for\dss $i\off =\off 32$\dss the\sss last\sss displayed equality\dss
turns into\dss 
$z^{\dff 32}\off =\off z^{\dff 30}$\dnsp.\oss
Therefore\dss
$z^{\dff 2}\off =\off 1$\nnsp.\oss
This completes\sss the proof\dss of\qss Coxeter's\dss implication.\oss
Note\sss that\dss this proof\dss explains why\sss the exponent\sss is\sss $2$
in\dss the identity\trs
$z^{\dff 2}\off =\off 1$\nnsp:\oss
it\dss is\dss the\dss Euler\dss characteristic 
of\trs the $2$\dnsp-sphere.\oss  \eproof

\newpage
\mysection{Actions\qss with\qss several\qss orbits\qss of\qss vertices}{several}

\myuppar{Generators.}
\emph{We will\sss assume\sss that\sss $X$\sss is\dss connected.}\oss
An\qss \emph{oriented\sss edge}\dss $e$\sss of\trs $X$\sss
is\dss an edge of\dss $X$\dss together with an\qss \emph{orientation},\oss
i.e.\qss an\qss \emph{ordered\dss pair}\qss $e\off =\off (\dff x\fff,\dff y\trf)$\dss
of\dss vertices\dss $x\fff,\qff y$\dss connected\sss by\sss an edge.\oss
The vertices $x$ and $y$ are called\dss the\qss \emph{origin}\qss 
and\dss \dss the\qss \emph{target}\qss of\dss $e$\sss 
and are denoted\dss by\dss
$\aaa\dff(\dff e\dff)$\sss and\dss $\ttt\trf(\dff e\trf)$\sss
respectively\halfff.\oss 
For an oriented edge\sss $e$\sss of\trs $X$\dss we denote\sss by\dss $\overline{e}$\dss
the same edge with\dss the orientation\sss reversed,\oss
so\sss that\sss
$\aaa\dff(\qff \overline{e}\qff)\off =\off \ttt\dff(\dff e\trf)$\dss
and\dss
$\ttt\dff(\qff \overline{e}\qff)\off =\off \aaa\dff(\dff e\trf)$\nnsp.\oss

Let\sss us\sss fix a set\sss $V$ of\dss representatives of\trs
the orbits of\dss the action of\sss $G$ on\dss the set\sss $X_{\dff 0}$ of\dss
vertices of\sss $X$\nnsp.\oss 
Let\sss $E$ be\sss the set\sss of\dss oriented edges $e$ of\sss $X$ such\dss that\sss
$\aaa\dff(\dff e\trf)\qff \in\qff V$\dnsp.\oss
For every\sss $v\qff \in\qff V$\dss let\dss $E_{\dff v}$\dss 
be\sss the set\sss of\dss oriented edges $e$ of\sss $X$ 
such\dss that\sss $\aaa\dff(\dff e\trf)\off =\off v$\nnsp.\oss
Then\vspace{0.75pt}
\[
\quad
E
\off =\off
\bigcup\nolimits_{\qff v\qff \in\qff V}\pff E_{\dff v}
\off.
\]

\vspace{-12pt}\vspace{0.75pt}
For every $e\qff \in\qff E$ the orbit\sss of\sss
$\ttt\dff(\dff e\trf)$ intersects $V$ at\sss a\sss single vertex which we will\sss
denote by $v\dff(\dff e\trf)$\nnsp.\oss
Let\dss us\dss choose\sss for\sss every\dss $e\qff \in\qff E$\dss
an element\sss $s_{\dff e}$\sss such\dss that\sss 
$s_{\dff e}\dff(\trf v\dff(\dff e\trf)\trf)
\pff =\off
\ttt\dff(\dff e\trf)$\nnsp.\oss
Let\sss
$\mathcal{S}
\off =\off 
\{\qff s_{\dff e}\qff |\qff e\qff \in\qff E \qff\}$\nnsp.\oss
For each $e\qff \in\qff E$\dss 
let\sss $g_{\dff e}$ be an abstract symbol\sss 
corresponding\dss to\sss $e$\nnsp.\oss
Let\sss $\mathcal{F}$\sss be\sss the free\sss group  having\sss  
$\mathcal{G}
\off =\off 
\{\qff g_{\dff e}\qff |\qff e\qff \in\qff E \qff\}$
as\sss the set\sss of\dss generators.\oss
Let\vspace{1.5pt}
\[
\quad
H
\off =\off
*_{\fff v\qff \in\qff V}\pff G_{\dff v}
\]

\vspace{-12pt}\vspace{1.5pt}
be\sss the free product\sss of\trs the groups\sss $G_{\dff v}$\nsp,\oss and\sss 
$\psi\dff \colon\dff
\mathcal{F}
\trf \ast\trf 
H 
\qff \ttoo\qff 
G$\sss
be\sss the unique homomorphism\qss 
equal\dss to\sss the inclusion 
$G_{\dff v}\qff \ttoo\qff G$ on every\sss  
$G_{\dff v}$\sss 
and\dss such\dss that\sss 
$\psi\dff(\trf g_{\dff e}\trf)
\off =\off
s_{\dff e}$
for every\dss $e\qff \in\qff E$\nnsp.\oss

\mypar{Lemma.}{transitivity-m} 
\emph{For\sss every\sss vertex $w$ of\pss $X$\dss
there\sss exists\sss   
$g\qff \in\qff \psi\dff(\trf \mathcal{F}\trf)$ 
such\dss that\dss $g\dff(\dff w\trf)\qff \in\qff V$\dnsp.\oss}

\proof
Since\sss $X$\sss is\dss connected,\oss
for every $w\qff \in\qff X_{\dff 0}$\sss
there exists a sequence\dss
$v_{\trf 0}\dff,\off 
v_{\dff 1}\dff,\off
\ldots\dff,\off
v_{\dff k}$\dss
of\dss vertices such\sss that\sss
$w\off =\off v_{\dff k}$\nsp,\dss 
$v\off =\off v_{\trf 0}$\nsp,\oss
and $v_{\fff i}$\sss is\dss connected\sss with $v_{\fff i\dff +\dff 1}$\sss
by\sss an edge of\sss $X$\sss
for each
$i\off =\off 0\fff,\pff 1\fff,\qff \ldots\fff,\qff k\qff -\qff 1$\nnsp.\oss
Arguing\sss by\dss induction,\oss we can assume\sss
that\dss there exists\dss
$g\qff \in\qff \psi\dff(\trf \mathcal{F}\trf)$ 
such\dss that\dss 
$g\trf(\dff v_{\dff k\dff -\dff 1}\trf)\qff \in\qff V$\dnsp.\oss
If\dss also\dss 
$g\trf(\dff v_{\dff k}\trf)\qff \in\qff V$\dnsp,\oss
then we are done.\oss
In any\sss case,\pss
$g\dff(\dff v_{\dff k\dff -\dff 1}\trf)$\sss
is\dss connected\dss with\sss
$g\dff(\dff v_{\dff k}\trf)$\sss
by\sss an edge of\trs $X$\nnsp.\oss
Let\dss $e$\sss be\sss this edge oriented\sss in such a way\sss that\sss
$\aaa\dff(\dff e\trf)
\off =\off
g\dff(\dff v_{\dff k\dff -\dff 1}\trf)$\dss
and\dss
$\ttt\dff(\dff e\trf)
\off =\off
g\dff(\dff v_{\dff k}\trf)$\nnsp.\oss
Then\dss $e\qff \in\qff E_{\dff u}$\nsp,\oss
where\dss
$u\off =\off g\dff(\dff v_{\dff k\dff -\dff 1}\trf)$\nnsp,\oss
and\dss hence\vspace{3pt}
\[
\quad
s_{\dff e}^{\dff -\dff 1}\fff g\qff(\dff w\trf)
\off =\off
s_{\dff e}^{\dff -\dff 1}\fff g\qff(\dff v_{\dff k}\qff)
\off =\off
s_{\dff e}^{\dff -\dff 1}\dff\left(\qff g\qff(\dff v_{\dff k}\trf) \qff\right)
\off =\off
s_{\dff e}^{\dff -\dff 1}\dff\left(\trf \ttt\dff(\dff e\trf)\trf\right)
\off =\off
v\dff(\dff e\trf)
\pff.
\]

\vspace{-9pt}
Clearly\halfff,\pss 
$g\qff \in\qff \psi\dff(\trf \mathcal{F}\trf)$
implies\sss that\trs
$s_{\dff e}^{\dff -\dff 1}\fff g
\qff \in\qff 
\psi\dff(\trf \mathcal{F}\trf)$\nnsp.\oss 
Since,\oss on\dss the other\dss hand,\pss
$v\dff(\dff e\trf)
\qff \in\qff
V$\dnsp,\oss
this completes\sss the induction step.\oss  \eproof

\mypar{Corollary\halfff.}{generators-m} 
\emph{The\sss group\dss $G$\dss is\dss generated\dss by\trs 
$\psi\dff(\trf F_{\dff E}\trf)$\dss and\dss $H$\nnsp.\oss 
Moreover\halfff,\oss 
$G\off =\off \psi\dff(\trf \mathcal{F}\trf)\dff \cdot\dff H$\nnsp.\oss 
In\dss particular\halfff,\oss
the homomorphism\qss 
$\psi\dff \colon\dff 
\mathcal{F}
\trf \ast\dff 
H
\qff \ttoo\qff 
G$\qss is\dss surjective.\oss}

\proof
Let\sss $g\qff \in\qff G$
and 
$v\qff \in\qff V$\dnsp.\oss
Lemma\qss \ref{transitivity-m}\qss implies\sss that
$f\fff g\trf(\dff v\trf)\qff \in\qff V$\sss
for some\sss 
$f\qff \in\qff \psi\dff(\trf \mathcal{F}\trf)$\nnsp.\oss
Since
$f\fff g\trf(\dff v\trf)$
belongs\sss to\sss the $G$\dnsp-orbit\sss of\sss $v$
and $V$ intersects every $G$\dnsp-orbit\sss 
only once,\oss it\dss follows\dss that\sss
$f\fff g\trf(\dff v\trf)\off =\off v$\nnsp,\oss
i.e.\qss that
$f\fff g\qff \in\qff G_{\dff v}$\nsp.\oss
Hence
$g\qff \in\qff f^{\dff -\dff 1}\dff \cdot\dff G_{\dff v}$
and
$g\qff \in\qff \psi\dff(\trf \mathcal{F}\trf)\dff \cdot\dff H$\nnsp.\oss  \eproof

\mypar{Lemma.}{rotated-ends-m}
\emph{Suppose\sss that\sss 
$v\qff \in\qff V$\dnsp,\qss
$e\qff \in\qff E_{\dff v}$\nsp,\qss
$h\qff \in\qff G_{\dff v}$\nsp.\oss
Let\sss $c\off =\off h\trf(\dff e\dff)$\sss
and\sss $w\off =\off v\dff(\dff e\trf)$\nnsp.\oss
Then\dss 
$v\dff(\dff c\trf)
\off =\off
v\dff(\dff e\trf)
\off =\off
w$\sss
and\trs
$s_{\dff c}^{\dff -\dff 1}\dff \cdot\dff h\dff \cdot\dff s_{\dff e}
\off \in\pff G_{\dff w}$\nsp.\oss}

\proof
The assumption\qss
$h\qff \in\qff G_{\dff v}$\qss
implies\sss that\qss
$\ttt\dff(\dff c\trf)\off =\off h\trf(\trf \ttt\dff(\dff e\trf)\trf)$\nnsp.\oss
It\dss follows\dss that\vspace{4.5pt}
\[
\quad 
v\dff(\dff c\trf)
\off =\off
s_{\dff c}^{\dff -\dff 1}\dff(\trf \ttt\dff(\dff c\trf)\trf)
\off =\off
s_{\dff c}^{\dff -\dff 1}\fff \cdot\dff h\trf(\trf \ttt\dff(\dff e\trf)\trf)
\off =\off
s_{\dff c}^{\dff -\dff 1}\fff \cdot\dff h\dff \cdot\dff s_{\dff e}\qff(\trf v\dff(\dff e\trf)\trf)
\qff.
\]

\vspace{-12pt}\vspace{4.5pt}
Hence 
\nsp$v\dff(\dff c\trf)$
belongs\sss to\sss the orbit\sss of $v\dff(\dff e\trf)$\nnsp.\oss
Since
\nsp$v\dff(\dff c\trf)\fff,\pff v\dff(\dff e\trf)\qff \in\qff V$
and\sss orbits\dss in\-ter\-sect $V$ only\sss once,\qss
$v\dff(\dff c\trf)
\off =\off
v\dff(\dff e\trf)$
and\dss hence\dss
$s_{\dff c}^{\dff -\dff 1}\dff \cdot\dff h\dff \cdot\dff s_{\dff e}
\off \in\pff G_{\dff w}$\nsp,\oss
where  
$w\off =\off v\dff(\dff e\trf)\off =\off v\dff(\dff c\trf)$\nnsp.\oss  \eproof

\myuppar{The edge relations.}
Suppose\sss that\trs
$e\qff \in\qff E_{\dff v}$\dss and\dss $h\qff \in\qff G_{\dff v}$\dss
for some\sss vertex\dss $v\qff \in\qff V$\dnsp,\oss
and\dss let\dss 
$w
\off =\off
v\dff(\dff e\trf)$\dss 
and\dss 
$c\off =\off h\trf(\dff e\dff)$\nnsp.\oss
Lemma\qss \ref{rotated-ends-m}\qss implies\sss 
that\qss\vspace{4.5pt}
\[
\quad
k\trf(\dff e\fff,\pff h\trf)
\off =\off
s_{\dff c}^{\dff -\dff 1}\dff \cdot\dff h\dff \cdot\dff s_{\dff e}
\off \in\pff G_{\dff w}
\qff.
\]

\vspace{-12pt}\vspace{4.5pt}
The\qss \emph{edge relation}\qss $E\dff(\dff e\fff,\pff h\trf)$\dss is\dss
the relation\sss 
$g_{\dff c}^{\dff -\dff 1}\dff \cdot\dff h\dff \cdot\dff g_{\dff e}
\off =\off
k\trf(\dff e\fff,\pff h\trf)$\nnsp.\oss
Here\sss $h$\sss is\dss considered as an element\sss of\trs 
the free factor\sss $G_{\dff v}$\sss of\trs $H$\dss
and\dss $k\trf(\dff e\fff,\pff h\trf)$\dss
as an element\sss of\trs the free factor\sss $G_{\dff w}$\nsp.\oss

\mypar{Lemma.}{tracing-a-path-m}
\emph{Suppose\sss that\sss $p$\sss is\dss a\sss path\dss
in\dss the\sss graph\sss $X$\nnsp,\oss 
i.e.\qss that\sss 
$p$ is\dss a\sss sequence of\dss vertices\qss 
$w_{\dff 0}\dff,\off w_{\dff 1}\dff,\off\ldots\dff,\off w_{\dff n}$\dss
of\oss $X$\sss such\dss that\dss
$c_{\dff i}
\off =\off
(\dff w_{\dff i\dff -\dff 1}\dff,\qff w_{\dff i} \trf)$\qss
is\dss an\sss oriented\dss edge of\trs $X$\dss 
for\sss every\dss
$i\qff \geq\qff 1$\nnsp.\oss
If\qss $w_{\trf 0}\qff \in\qff V$\dnsp,\oss
then\dss there\dss is\dss a\sss unique sequence\dss
$e_{\dff 1}\dff,\off e_{\dff 2}\dff,\off\ldots\dff,\off e_{\dff n}
\qff \in\qff
E$\dss
such\dss that}\vspace{4.5pt}
\begin{equation}
\label{path-edge}
\quad
s_{\dff 1}\dff \cdot\dff s_{\dff 2}\dff \cdot\dff
\ldots\dff \cdot\dff 
s_{\dff i\dff -\dff 1}\trf
(\dff e_{\dff i}\trf)
\off =\off
c_{\trf i}
\quad
\mbox{\emph{and}}\quad
\end{equation}

\vspace{-34.5pt}
\begin{equation}
\label{path-step-m}
\quad
s_{\dff 1}\dff \cdot\dff s_{\dff 2}\dff \cdot\dff
\ldots\dff \cdot\dff 
s_{\dff i}\trf
(\trf v\dff(\dff e_{\dff i}\trf)\trf)
\off =\off
w_{\dff i}
\end{equation}

\vspace{-12pt}\vspace{4.5pt}
\emph{for\sss every\dss
$i\off =\off 1\fff,\pff 2\fff,\pff \ldots\fff,\pff n$\nnsp,\oss
where\dss
$s_{\dff i}\off =\off s_{\dff e_{\dff i}}$.\oss
Also,\qss
$\aaa(\dff e_{\dff i}\trf)
\off =\off
v\dff(\dff e_{\dff i\dff -\dff 1}\trf)$\sss 
for\sss $i\qff \geq\qff 2$\nnsp.\oss}

\proof
Let\sss $e_{\dff 1}\off =\off c_{\dff 1}$\nsp.\oss
Suppose\sss that\dss the edges\dss
$e_{\dff 1}\dff,\off e_{\dff 2}\dff,\off\ldots\dff,\off e_{\dff k\dff -\dff 1}
\qff \in\qff E$\dss
are already\sss determined
and\qss (\ref{path-step-m})\qss holds\sss for\dss
$i\qff =\qff k\qff -\qff 1$\nnsp.\oss
Then\qss (\ref{path-edge})\qss
with\sss $i\off =\off k$\sss uniquely determines\sss $e_{\dff k}$\nsp.\oss
Moreover\halfff,\qss
$s_{\dff 1}\dff \cdot\dff s_{\dff 2}\dff \cdot\dff
\ldots\dff \cdot\dff 
s_{\dff k\dff -\dff 1}\trf
(\trf \aaa(\dff e_{\dff k}\trf)\trf)
\off =\off
w_{\dff k\dff -\dff 1}$\nsp.\oss
At\sss the same\sss time\qss (\ref{path-step-m})\qss
for\dss $i\qff =\qff k\qff -\qff 1$\sss
implies\sss that\sss
$s_{\dff 1}\dff \cdot\dff s_{\dff 2}\dff \cdot\dff
\ldots\dff \cdot\dff 
s_{\dff k\dff -\dff 1}\trf
(\trf v\dff(\dff e_{\dff k\dff -\dff 1}\trf)\trf)
\off =\off
w_{\dff k\dff -\dff 1}$\nsp.\oss
It\sss follows\sss that\sss
$\aaa(\dff e_{\dff k}\trf)
\off =\off
v\dff(\dff e_{\dff k\dff -\dff 1}\trf)
\qff \in\qff
V$
and\dss hence\sss $e_{\dff k}$\sss belongs\sss to\sss $E$\sss 
and\dss  
$s_{\dff k}\off =\off s_{\dff e_{\dff k}}$\dss is\dss defined.\oss
Therefore\vspace{3pt}
\[
\quad
s_{\dff 1}\dff \cdot\dff s_{\dff 2}\dff \cdot\dff
\ldots\dff \cdot\dff 
s_{\dff k\dff -\dff 1}\dff \cdot\dff
s_{\dff k}\trf
(\trf v\dff(\dff e_{\dff k}\trf)\trf)
\off =\off
s_{\dff 1}\dff \cdot\dff s_{\dff 2}\dff \cdot\dff
\ldots\dff \cdot\dff 
s_{\dff k\dff -\dff 1}\trf
\bigl(\trf  \ttt(\dff e_{\dff k}\trf)\trf\bigr)
\off =\off
w_{\dff k}
\]

\vspace{-12pt}\vspace{3pt}
where\sss the\sss last\sss equality\sss holds by\sss the definition of $e_{\dff k}$\nsp.\oss
We see\sss that\qss (\ref{path-step-m})\qss holds\sss for\dss 
$i\off =\off k$\nnsp.\oss
An\dss induction\dss by\sss $k$\sss completes\sss the proof\halfff.\oss  \eproof

\myuppar{The\sss pseudo-loop\sss relations.}
A path\sss in\dss $X$\dss 
is\dss said\dss to be a\dss \emph{pseudo-loop}\pss if\trs it\dss begins and ends in\dss $V$\dnsp.\oss
Suppose\sss that\dss the path\sss $p$\sss in\dss Lemma\qss \ref{tracing-a-path-m}\fff\qss
is\dss a pseudo-loop,\oss
i.e.\qss that\trs $w_{\dff 0}\dff,\off w_{\dff n}\qff \in\qff V$\dnsp.\oss
Since every\sss orbit\dss intersects\sss $V$\sss only\sss once,\oss
(\ref{path-step-m})\qss with\dss $i\off =\off n$\dss
implies\sss that 
$v\dff(\dff e_{\dff n}\trf)\off =\off w_{\dff n}$
and\vspace{3pt}
\[
\quad
s_{\dff 1}\dff \cdot\dff s_{\dff 2}\dff \cdot\dff
\ldots\dff \cdot\dff 
s_{\dff n}\trf
(\dff w_{\dff n}\trf)
\off =\off
w_{\dff n}
\qff,
\]

\vspace{-12pt}\vspace{3pt}
i.e.\qss
$s_{\dff 1}\dff \cdot\dff s_{\dff 2}\dff \cdot\dff
\ldots\dff \cdot\dff 
s_{\dff n}
\off \in\off
G_{\dff w_{\dff n}}$\nsp.\oss
Let\trs
$g_{\dff i}\off =\off g_{\dff e_{\dff i}}$.\oss
The\qss \emph{pseudo-loop\sss relation}\qss $L\trf(\trf p\trf)$\dss is\dss
the relation\vspace{1.5pt}
\[
\quad
g_{\dff 1}\dff \cdot\qff g_{\dff 2}
\dff \cdot\qff 
\ldots
\qff \cdot\qff 
g_{\dff n}
\off =\off
s_{\dff 1}\dff \cdot\dff s_{\dff 2}
\dff \cdot\dff 
\ldots
\dff \cdot\dff 
s_{\dff n}
\qff.
\]

\vspace{-12pt}\vspace{1.5pt}
Here\qss  
$g_{\dff 1}\dff \cdot\qff g_{\dff 2}
\dff \cdot\qff 
\ldots
\qff \cdot\qff 
g_{\dff n}
\qff \in\qff
\mathcal{F}$\qss 
and\qss 
$s_{\dff 1}\dff \cdot\dff s_{\dff 2}
\dff \cdot\dff 
\ldots
\dff \cdot\dff 
s_{\dff n}$\qss
is\dss an element\sss of\trs the free factor\sss $G_{\dff w_{\dff n}}$\nsp.\oss
When\sss $p$\dss is\dss actually\sss a\sss loop,\oss
i.e.\qss when\dss $w_{\dff n}\off =\off w_{\dff 0}$\nsp,\oss
we\sss will\sss call\dss the relation $L\trf(\trf p\trf)$ a\qss \emph{loop-relation}.\oss

\myuppar{The edge-loop relations.}
As in\dss Section\qss \ref{one},\oss
for $v\qff \in\qff V$\dss and $e\qff \in\qff E_{\dff v}$\sss 
we denote by\sss $l_{\dff e}$\sss
the\sss loop\dss
$v\fff,\pff t\trf(\dff e\trf)\fff,\pff v$\sss
and\sss call\sss such\dss loops\qss \emph{edge-loops}.\oss
An edge-loop\sss $l_{\dff e}$\sss is\dss a\sss loop and\dss 
we will\sss call\trs the corresponding\sss loop\sss relation\sss $L\trf(\trf l_{\dff e}\trf)$
an\qss \emph{edge-loop}\qss relation.\oss
It\dss has\sss the form\sss 
$g_{\dff e}\dff \cdot\dff g_{\dff a}
\off =\off
s_{\dff e}\dff \cdot\dff s_{\dff a}$\nsp,\oss
where 
$a
\off =\off 
s_{\dff e}^{\dff -\dff 1}\dff(\qff \overline{e} \qff)$\nnsp.\oss
In\dss par\-tic\-u\-lar\halfff,\oss
the orbit\sss of\trs $\ttt\dff(\dff a\trf)$\dss
contains $v$ and\dss hence\dss $v\dff(\dff a\trf)\off =\off v$\nnsp.\oss

\myuppar{Tautological\dss relations.}
For each $v\qff \in\qff V$\trs let\trs $T_{v}$\dss be\sss the set\sss of\dss oriented
edges\dss $e\qff \in\qff E_{\dff v}$\dss such\dss that\trs
$s_{\dff e}\qff \in\qff G_{\dff v}$\nsp.\oss
For\dss $e\qff \in\qff T_{\fff v}$\dss
the\qss \emph{tautological\dss relation}\qss $T\dff(\dff e\trf)$\dss 
is\dss the\sss relation\dss
$g_{\dff e}\off =\off s_{\dff e}$\nsp,\oss
where
$s_{\dff e}$ 
is\dss considered as an element\sss of\trs the
free factor\dss $G_{\dff v}$\nsp.\oss
Let\trs $T\off =\off \bigcup_{\qff v}\fff T_v$\nsp.\oss

\myuppar{Introducing\dss relations\dss in\dss $\mathcal{F}\dff \ast\dff H$\nnsp.}
Let\dss $\mathcal{L}$\dss be
a collection\sss of\trs pseudo-loops,\oss
and\dss let\trs $\mathcal{T}$\dss be a subset\sss of\trs $T$\dnsp.\oss
Let us impose on\dss $\mathcal{F}\dff \ast\trf H$\dss all\dss 
edge relations,\oss 
the pseudo-loop relations\sss $L\trf(\trf p\trf)$\sss such\dss that\sss either\sss
$p\qff \in\qff \mathcal{L}$\dss or\dss
$p$\dss is\dss an edge-loop,\oss
and\dss the\sss tautological\dss relations\dss $T\dff(\dff e\trf)$\dss 
for\dss $e\qff \in\qff \mathcal{T}$\dnsp,\oss
and\dss let\sss $\mathbb{G}$\sss be the resulting\sss group.\oss
Let\sss
$\rho\dff \colon\dff
\mathcal{F}\dff \ast\trf H\qff \ttoo\qff \mathbb{G}$\sss
be\sss the quotient\sss map.\oss
Since all\sss these relations hold\sss in $G$\nnsp,\qss
$\psi$\sss induces a homomorphism\sss
$\varphi\dff \colon\dff
\mathbb{G}\qff \ttoo\qff G$\nnsp.\oss
Clearly,\qss $\psi\off =\off \varphi\dff \circ\dff \rho$\nnsp.\oss
For $v\qff \in\qff V$\sss let\sss
$\mathbb{G}_{\dff v}
\off =\off
\rho\trf(\trf G_{\dff v}\trf)$\nnsp.\oss
Since 
$\psi$ is\dss equal\dss 
to\sss the inclusion $G_{\dff v}\qff \ttoo\qff G$ 
on $G_{\dff v}$\nsp,\oss
the maps 
$G_{\dff v}\qff \ttoo\qff \mathbb{G}_{\dff v}$ 
and
$\mathbb{G}_{\dff v}\qff \ttoo\qff G_{\dff v}$ 
induced\dss by\sss $\rho$\sss and\sss $\varphi$\sss 
respectively are isomorphisms.\oss
For every $e\qff \in\qff E$\sss let\sss 
$\mathfrak{g}_{\dff e}
\off=\off 
\rho\dff(\trf g_{\dff e}\dff)
\qff \in\qff 
\mathbb{G}$\nnsp.\oss
Then\dss 
$\varphi\dff(\trf \mathfrak{g}_{\dff e}\dff)
\off =\off 
\psi\dff(\trf g_{\dff e}\dff)\off =\off s_{\dff e}$\nsp.\oss
Let\sss
$\mathfrak{G}
\off =\off
\rho\trf(\trf \mathcal{G}\trf)$  
and\dss 
$\mathbb{H}
\off =\off 
\rho\dff(\trf H\trf)\qff \subset\qff \mathbb{G}$\nnsp.\oss
Then $\rho$ and $\varphi$\sss induce bijections\sss
$\mathcal{G}\qff \ttoo\qff \mathfrak{G}$
and\sss
$\mathfrak{G}\qff \ttoo\qff \mathcal{S}$
respectively.\oss

\mypar{Lemma.}{exchange-m} 
\emph{Let\dss $v\qff \in\qff V$\dnsp,\sss $e\qff \in\qff E$\nnsp,\qss
and\sss $w\off =\off v\dff(\dff e\trf)$\nnsp.\oss
Then\qss
$\mathbb{G}_{\dff v}\dff \cdot\qff \mathfrak{g}_{\dff e}
\off \subset\off 
\mathfrak{G}\dff \cdot\qff \mathbb{G}_{\dff w}$\nsp.}

\proof
Let\sss $t\qff \in\qff G_{\dff v}$ and\sss $d\off =\off t\trf(\dff e\dff)$\nnsp.\oss 
Since\sss the relation\trs $E\dff(\dff e\fff,\pff t\trf)$\sss 
holds in\sss $\mathbb{G}$\nnsp,\oss\vspace{2pt}
\[
\quad
\rho\trf\left(\trf g_{\dff d}^{\dff -\dff 1}\dff \cdot\dff t\dff \cdot\dff g_{\dff e}\qff\right)
\off =\off
\rho\trf\bigl(\trf k\trf(\dff e\fff,\pff t\trf)\trf\bigr)
\qff.
\]

\vspace{-12pt}\vspace{2pt}
Lemma\qss \ref{rotated-ends-m}\qss implies\sss that\sss
$k\trf(\dff e\fff,\pff h\trf)\qff \in\qff G_{\dff w}$
and\dss hence\sss
$\rho\dff (\trf k\trf(\dff e\fff,\pff t\trf)\trf)\qff \in\qff \mathbb{G}_{\dff w}$\nsp.\oss
It\sss follows\sss that\vspace{3pt}
\[
\quad
\mathfrak{g}_{\dff d}^{\dff -\dff 1} 
\cdot\dff 
\rho\dff (\dff t\trf)
\dff \cdot\dff 
\mathfrak{g}_{\dff e}
\pff \in\pff
\mathbb{G}_{\dff w}
\qff
\]

\vspace{-12pt}\vspace{3pt}
and\dss hence\sss
$\rho\dff (\dff t\trf)
\dff \cdot\dff 
\mathfrak{g}_{\dff e}
\off \in\off
\mathbb{G}_{\dff E}\dff \cdot\qff \mathbb{G}_{\dff w}$\nsp.\oss
The\sss lemma\sss follows.\oss  \eproof

\myuppar{Kozsul\trs models.}
The\qss \emph{Kozsul\dss model}\pss $\mathbb{X}$ of\dss $X$\dss 
has as\sss the set\sss of\dss vertices $\mathbb{X}_{\dff 0}$
the disjoint\sss union\vspace{3pt}
\[
\quad
\mathbb{X}_{\dff 0}
\off =\off
\coprod\nolimits_{\dff v\qff \in\qff V}\qff \mathbb{G}\left/\dff \mathbb{G}_{\dff v}\right.
\off
\]

\vspace{-12pt}\vspace{3pt}
of\dss sets 
$\mathbb{G}/\fff \mathbb{G}_{\dff v}$\sss
with $v\qff \in\qff V$\dnsp.\oss 
The action of $\mathbb{G}$ on sets $\mathbb{G}/\mathbb{G}_{\dff v}$
define an action of $\mathbb{G}$\sss on $\mathbb{X}_{\dff 0}$\nsp.\oss
For $v\qff \in\qff V$ let 
$v^{\fff *}$ be\sss the coset 
$\mathbb{G}_{\dff v}\qff \in\qff \mathbb{G}/ \mathbb{G}_{\dff v}$\sss 
thought\sss as a vertex
of\dss the future graph $\mathbb{X}$\nnsp.\oss
Let\dss 
$f\dff \colon\dff \mathbb{X}_{\dff 0}\qff \ttoo\qff X_{\dff 0}$ 
be\sss the map defined\dss by\sss
$f\dff(\trf \gamma\trf(\trf v^{\fff *}\trf)\trf)
\off =\off
\varphi\dff(\trf \gamma\trf)\dff (\dff v\trf)$\sss
for every $v\qff \in\qff V,\off \gamma\qff \in\qff \mathbb{G}$\nnsp.\oss
Since 
$\varphi\trf(\trf \mathbb{G}_{\dff v}\trf)\off =\off G_{\dff v}$
fixes $v$\nnsp,\oss
the map $f$ is\dss well\sss defined.\oss
As in\dss Section\qss \ref{models},\oss 
the map\sss $f$ is\dss \emph{$\mathbb{G}$\dnsp-equivariant}\qss 
with respect\dss to\sss the natural\sss action of $\mathbb{G}$ on 
$\mathbb{X}_{\dff 0}$ 
and\sss the action on 
$X_{\trf 0}$ via $\varphi$\nnsp. 

Let\sss us define for every 
$x\qff \in\qff \mathbb{X}_{\dff 0}$\sss 
a set\sss 
$N\trf(\dff x\trf)$\dss 
of\qss \emph{neighbors}\qss of $x$\nnsp,\oss 
the vertices\sss to be connected\sss with $x$ by an edge of\sss $\mathbb{X}$
and\sss then check\sss that\sss the resulting\sss relation of\dss being a\qss
\emph{neighbor}\pss is\dss symmetric.\oss 
Every\dss vertex of\trs $\mathbb{X}$\sss has\sss the form
$\gamma\trf(\dff v^{\dff *}\trf)$\nnsp,\oss 
where\sss $\gamma\qff \in\qff \mathbb{G}$
and $v\qff \in\qff V$\dnsp.\oss
For $e\qff \in\qff E_{\dff v}$\sss let\sss us\sss
think about 
$v^{\fff *},\qff\off v\dff(\dff e\trf)^{\fff *}$ 
and\sss 
$\mathfrak{g}_{\dff e}$\dss
as\sss the\sss lifts\dss of\sss 
$v\fff,\pff v\dff(\dff e\trf)$ 
and\sss $s_{\dff e}$
respectively\halfff,\oss
and set\vspace{3pt}
\[
\quad
N\trf(\trf \gamma\trf(\dff v^{\dff *}\trf)\trf) 
\off =\dff\off
\{\pff \gamma\dff \cdot\dff \mathfrak{g}_{\dff e}\trf
(\trf v\dff(\dff e\trf)^{\dff *}\trf)
\pff |\pff
e\qff \in\qff E_{\dff v}
\pff\}
\qff.
\]

\vspace{-12pt}\vspace{3pt}
\mypar{Lemma.}{correctness-of-N-m}
\emph{The definition\sss of\pss $N\trf(\dff x\trf)$\dss is\dss correct\halfff,\oss
i.e.\qss $N\trf(\dff x\trf)$\dss does not\sss depend on\dss the 
choice of\qss $\gamma\fff,\pff v$\qss such\dss that\qss
$x\off =\off \gamma\trf(\trf v^{\fff *}\trf)$\nnsp.\oss}

\proof
Let\dss us\sss prove first\dss that\sss $v$\sss is\dss uniquely\sss determined\dss
by\sss $x$\nnsp.\oss
If\trs
$\gamma\trf (\dff  v^{\fff *}\trf)
\off =\off
\beta\trf (\dff  u^{\fff *}\trf)$\nnsp,\oss
then\vspace{3pt}
\[
\quad
\varphi\dff(\trf \gamma\trf)\dff (\dff v\trf)
\off =\off
f\dff(\trf \gamma\trf(\trf v^{\fff *}\trf)\trf)
\off =\off
f\dff(\qff \beta\trf(\trf u^{\fff *}\trf)\trf)
\off =\off
\varphi\dff(\trf \beta\trf)\dff (\dff u\trf)
\]

\vspace{-12pt}\vspace{3pt}
and\dss hence\sss $v\fff,\pff u$\sss belong\dss
to\sss the same $G$\dnsp-orbit\halfff.\oss
Since\dss $v\fff,\pff u\qff \in\qff V$\dnsp,\oss
this implies\sss that\trs $u\off =\off v$\nnsp.\oss
This proves\sss that\sss $v$\sss is\dss indeed\sss uniquely\sss determined\dss
by\sss $x$\nnsp.\oss
The rest\sss of\dss the proof\dss is\dss similar\sss to\sss the proof\dss
of\trs Lemma\qss \ref{correctness-of-N},\oss
with\dss Lemma\qss \ref{exchange-m}\qss
playing\sss the role of\trs Lemma\qss \ref{exchange}.\oss  \eproof

\mypar{Lemma.}{symmetry-m} 
\emph{Let\qss $x\fff,\pff y\qff \in\qff \mathbb{X}_{\trf 0}$\nsp.\oss 
If\pss $y\qff \in\qff N\trf(\trf x\trf)$\nnsp,\oss
then\qss $x\qff \in\qff N\trf(\trf y\trf)$\nnsp.\oss}

\proof
The proof\dss follows\sss the same route as\sss the proof\dss of\trs
Lemma\qss \ref{symmetry}.\oss
Let\trs us choose\dss $v\qff \in\qff V$\sss and\sss $\gamma\qff \in\qff \mathbb{G}$\dss
such\dss that\sss
$x\off =\off \gamma\trf(\dff v^{\fff *}\trf)$\nnsp.\oss 
If\sss $y\qff \in\qff N\trf(\dff x\trf)$\nnsp,\oss
then\sss \vspace{3pt}
\[
\quad
y
\off =\off 
\gamma
\qff \cdot\qff 
\mathfrak{g}_{\dff e}\qff(\trf v\dff(\dff e\trf)^{\fff *}\trf)
\]

\vspace{-12pt}\vspace{3pt}
for some oriented edge\sss $e\qff \in\qff E_{\dff v}$\nsp.\oss
The edge-loop\sss relation\sss $L\trf(\trf l_{\dff e}\trf)$
has\sss the form\sss
$g_{\dff e}\dff \cdot\dff g_{\dff a}
\off =\off
s_{\dff e}\dff \cdot\dff s_{\dff a}$\sss
with\sss $a\qff \in\qff E_{\dff v\dff(\dff e\trf)}$\dss and\dss
$s_{\dff e}\dff \cdot\dff s_{\dff a}\qff \in\qff G_{\dff v}$\nsp.\oss
By\sss applying\sss $\rho$\sss 
we see\sss that\dss 
$\mathfrak{g}_{\dff e}\dff \cdot\dff \mathfrak{g}_{\dff a}
\off =\off
\rho\trf(\dff s_{\dff e}\dff \cdot\dff s_{\dff a}\trf)
\pff \in\off
\mathbb{G}_{\dff v}$\nsp.\oss
Since\dss $\mathbb{G}_{\dff v}$\dss fixes\sss $v^{\dff *}$\sss 
under\dss the action of\trs $\mathbb{G}$\sss on\dss $\mathbb{X}_{\dff 0}$\nnsp,\oss 
this implies\sss that\vspace{3pt}
\[
\quad
x
\off =\off 
\gamma\trf(\dff v^{\fff *}\trf)
\off =\off 
\gamma\dff \cdot\trf
\rho\trf\left(\dff s_{\dff e}\dff \cdot\dff s_{\dff a}\trf\right)\trf(\dff v^{\fff *}\trf)
\off =\off
\gamma\trf \cdot\trf \mathfrak{g}_{\dff e}\dff \cdot\dff \mathfrak{g}_{\dff a}\trf(\dff v^{\fff *}\trf)
\off =\off
(\trf \gamma\trf \cdot\trf \mathfrak{g}_{\dff e} \trf)\dff \cdot\dff \mathfrak{g}_{\dff a}\dff
(\trf v\dff(\dff a\trf)^{\dff *}\trf)
\qff,
\]

\vspace{-12pt}\vspace{3pt}
where at\sss the\sss last\sss step we used\sss the fact\dss that\sss
$v\dff(\dff a\trf)\off =\off v$\nnsp.\oss
Since
$\dis
y
\off =\off 
\gamma\trf \cdot\trf \mathfrak{g}_{\dff e} \qff
(\qff v\dff(\dff e\trf)^{\fff *}\qff)$
and
$a\qff \in\qff E_{\dff v\dff(\dff e\trf)}$\nsp,\oss
this means\sss that\trs
$x\qff \in\qff N\trf(\trf y\trf)$\nnsp.\oss
The lemma follows.\oss  \eproof

\myuppar{The\sss graph\sss $\mathbb{X}$ and\sss the map\sss 
$f\dff \colon\dff \mathbb{X}\qff \ttoo\qff X$\nnsp.}
By\qss Lemmas\qss \ref{correctness-of-N-m}\qss and\qss \ref{symmetry-m}\qss
the relation\dss $y\qff \in\qff N\trf(\trf x\trf)$\dss
is\dss correctly\sss defined and symmetric.\oss
Therefore we can define\sss the graph\sss $\mathbb{X}$\nnsp,\oss
the action of\dss $\mathbb{G}$ on $\mathbb{X}$ and\sss the $\mathbb{G}$\dnsp-equivariant\sss map
$f\dff \colon\dff \mathbb{X}\qff \ttoo\qff X$
exactly as in\dss Section\qss \ref{models}.\oss

\mypar{Lemma.}{local-isomorphism-m} 
\emph{The map\qss 
$f\dff \colon\dff 
\mathbb{X}\qff \ttoo\qff X$\qss 
is\dss a\sss local\dss isomorphism of\trs graphs,\oss 
i.e.\qss for every\dss vertex\dss
$x$ of\pss $\mathbb{X}$\dss the map\dss $f$\dss maps\sss the set\dss of\qss edges of\pss 
$\mathbb{X}$\dss having\dss $x$\sss as\sss an\sss endpoint\dss bijectively\sss 
onto\sss the set\sss of\qss edges of\pss $X$\dss having\dss
$f\dff (\dff x\trf)$ as\sss an\sss endpoint\halfff.\oss}

\proof
Let\sss $v\qff \in\qff V$\dnsp.\oss
By\sss the definitions,\oss
if\sss $e\qff \in\qff E_{\dff v}$\nsp,\oss
then\vspace{3pt}
\[
\quad
f\trf(\trf \mathfrak{g}_{\dff e}\trf
(\dff v\dff(\dff e\trf)^{\fff *}\trf)\qff)
\off =\off
\varphi\trf(\trf \mathfrak{g}_{\dff e}\trf)\trf
(\dff v\dff(\dff e\trf)\trf)
\off =\off
s_{\dff e}\dff(\dff v\dff(\dff e\trf)\trf)
\off =\off
\ttt(\dff e\trf)
\qff.
\]

\vspace{-12pt}\vspace{3pt}
Since an edge  $e\qff \in\qff E_{\dff v}$\sss
is\dss uniquely determined\dss by $v$ and\sss $\ttt(\dff e\trf)$\nnsp,\oss
this implies\sss that\sss $f$ induces a bijection between\sss the neighbors
of\sss $v^{\fff *}$ and\sss the neighbors of\sss $v$\nnsp.\oss
This proves\sss the\sss lemma for $x\qff \in\qff V$\dnsp.\oss
Now\sss the equivariance of\sss $f$\sss implies\sss the general\sss case.\oss  \eproof

\mypar{Lemma.}{lifted-loops-m}
\emph{Let\sss $p$ be\sss a\sss path\dss
$w_{\dff 0}\dff,\off w_{\dff 1}\dff,\off\ldots\dff,\off w_{\dff n}$\dss
in\qss $X$\nnsp.\oss 
Suppose\sss that\sss $p$\sss is\qss a\dss pseudo-loop,\oss
i.e.\qss that\trs $w_{\dff 0}\dff,\pff w_{\dff n}\qff \in\qff V$\dnsp.\oss
If\pss $p\qff \in\qff \mathcal{L}$\nnsp,\oss 
then\dss the\dss lift\dss of\trs $p$\sss to\qss $\mathbb{X}$\trs 
starting\dss at\qss $w_{\trf 0}^{\fff *}$\dss 
ends\dss at\qss $w_{\fff n}^{\fff *}$\nsp.\oss}

\proof 
Let\dss us\sss consider\dss the sequence of\dss oriented edges\dss
$e_{\dff 1}\dff,\off e_{\dff 2}\dff,\off\ldots\dff,\off e_{\dff n}
\qff \in\qff
E$\dss
determined\dss by\dss $p$\dss
as\sss in\dss Lemma\qss \ref{tracing-a-path-m}.\oss
For
$i\off =\off 1\fff,\qff  \ldots\fff,\qff n$\dss
let\sss\vspace{3pt}
\[
\quad
s_{\dff i}
\off =\off 
s_{\dff e_{\dff i}}
\qff,\quad
g_{\dff i}
\off =\off 
g_{\dff e_{\dff i}}
\qff,\quad
\mathfrak{g}_{\dff i}
\off =\off 
\mathfrak{g}_{\dff e_{\dff i}}
\qff,\quad
z_{\trf i}
\off =\off
\mathfrak{g}_{\dff 1}
\dff \cdot\dff 
\mathfrak{g}_{\dff 2}
\dff \cdot\dff 
\ldots
\dff \cdot\dff 
\mathfrak{g}_{\dff i}\trf
\left(\trf v\dff(\dff e_{\dff i}\trf)^{\fff *}\trf\right)
\qff.
\]

\vspace{-12pt}\vspace{3pt}
Without\sss defining $e_{\dff 0}$
let\sss us\sss set\sss 
$v\dff(\dff e_{\dff 0}\trf) 
\off =\off w_{\dff 0}$\nsp,\dss 
$z_{\dff 0}
\off =\off 
w_{\dff 0}^{\fff *}$\nsp.\oss
Then\trs Lemma\qss \ref{tracing-a-path-m}\qss implies\sss that\sss
$\aaa(\dff e_{\dff i}\trf)
\off =\off
v\dff(\dff e_{\dff i\dff -\dff 1}\trf)$\sss 
for every\sss $i\qff \geq\qff 1$\sss
and\dss hence\vspace{3pt}
\[
\quad
\mathfrak{g}_{\dff i}\trf(\dff v\trf(\dff e_{\dff i}\trf)^{\fff *}\trf)
\off =\off
\mathfrak{g}_{\dff e_{\dff i}}\trf(\dff v\trf(\dff e_{\dff i}\trf)^{\fff *}\trf)
\off \in\off
N\trf(\trf v\trf(\dff e_{\dff i\dff -\dff 1}\trf)^{\fff *}\trf)
\]

\vspace{-12pt}\vspace{3pt}
and\dss hence\sss
$z_{\dff i}\qff \in\qff N\trf(\trf z_{\dff i\dff -\dff 1}\trf)$\sss
for every\trs $i\qff \geq\qff 1$\nnsp.\oss
Therefore\sss
$z_{\qff 0}\dff,\off z_{\trf 1}\dff,\off\ldots\dff,\off z_{\trf n}$\sss
is\dss a path\sss in $\mathbb{X}$\nnsp,\oss
which we will\sss denote by $q$\nnsp.\oss
Since\sss 
$f\dff(\dff z_{\trf i}\trf)
\off =\off
s_{\dff 1}\dff \cdot\dff s_{\dff 2}\dff \cdot\dff
\ldots\dff \cdot\dff 
s_{\dff i}\qff
(\trf v\dff(\dff e_{\dff i}\trf)\trf)$\sss
for every $i$\nnsp,\oss
the equality\qss (\ref{path-step-m})\qss implies\sss that\sss $q$\sss
is\dss the unique\sss lift\sss of\dss $p$\nnsp.\oss
It\dss remains\sss to find\dss where\sss  
$q$\sss ends.\oss
Since $w_{\fff n}$ belongs\sss to\sss the orbit\sss of\dss $v\dff(\dff e_{\dff n}\trf)$
and\sss $w_{\fff n}\qff \in\qff V$\sss because $p$\sss is\dss a pseudo-loop,\qss
$w_{\fff n}
\off =\off 
v\dff(\dff e_{\dff n}\trf)$ 
and\dss hence\sss $q$ ends at\sss
$\mathfrak{g}_{\dff 1}
\dff \cdot\dff 
\mathfrak{g}_{\dff 2}
\dff \cdot\dff 
\ldots
\dff \cdot\dff 
\mathfrak{g}_{\dff n}\trf
(\trf w_{\fff n}^{\fff *}\trf)$\nnsp.\oss
Recall\dss that\dss $L\dff(\trf p\trf)$\dss is\dss the relation\sss 
$g_{\dff 1}
\dff \cdot\dff 
g_{\dff 2}
\dff \cdot\dff 
\ldots
\dff \cdot\dff 
g_{\dff n}
\off =\off
s_{\dff 1}\dff \cdot\dff s_{\dff 2}
\dff \cdot\dff 
\ldots
\dff \cdot\dff 
s_{\dff n}$
and\dss that\sss 
$s_{\dff 1}\dff \cdot\dff s_{\dff 2}\dff \cdot\dff
\ldots\dff \cdot\dff 
s_{\dff n}
\off \in\off
G_{\dff w_{\dff n}}$\nsp.\oss
It\sss follows\sss that\vspace{3pt}
\[
\quad
\mathfrak{g}_{\dff 1}
\dff \cdot\dff 
\mathfrak{g}_{\dff 2}
\dff \cdot\dff 
\ldots
\dff \cdot\dff 
\mathfrak{g}_{\dff n}
\off =\off
\rho\dff(\trf
s_{\dff 1}\dff \cdot\dff s_{\dff 2}
\dff \cdot\dff 
\ldots
\dff \cdot\dff 
s_{\dff n}
\trf)
\off \in\off
\rho\dff(\trf G_{\dff w_{\dff n}}\trf)
\off =\off
\mathbb{G}_{\dff w_{\dff n}}
\pff.
\]

\vspace{-12pt}\vspace{3pt}
Since $\mathbb{G}_{\dff w_{\dff n}}$\sss is\dss the stabilizer of 
$w_{\dff n}^{\fff *}$\dss in $\mathbb{X}$\sss by\dss the construction of\sss $\mathbb{X}$\nnsp,\oss
this implies\sss that\sss the end\sss vertex\sss 
$\mathfrak{g}_{\dff 1}
\dff \cdot\dff 
\mathfrak{g}_{\dff 2}
\dff \cdot\dff 
\ldots
\dff \cdot\dff 
\mathfrak{g}_{\dff n}\trf
(\trf w_{\fff n}^{\fff *}\trf)$
of\sss $q$\sss is\dss equal\dss to $w_{\fff n}^{\fff *}$\nsp.\oss  \eproof

\mypar{Lemma.}{tree-of-representatives}
\emph{There exists a subtree\sss $A$ of\pss $X$\dss
such\dss that\dss the set\sss of\dss vertices of\qss $A$\dss
is\dss a\sss set\sss of\trs representatives\dss of\trs orbits.\oss}

\proof
Let\dss $A$\dss be maximal\sss among\sss
the subtrees of\trs $X$\dss such\dss that\sss 
all\dss their vertices belong\sss to different\sss orbits.\oss
Suppose\sss that\dss there\dss is\dss an orbit\sss of\dss vertices not\dss intersecting\sss $A$\nnsp.\oss
Since\dss $X$\dss is\dss connected,\oss there\dss is\dss a\sss path\dss
$w_{\dff 0}\dff,\off w_{\dff 1}\dff,\off\ldots\dff,\off w_{\fff n}$\dss
in\dss $X$\dss such\dss that\sss $w_{\dff 0}$\sss is\dss a\sss vertex of\dss $A$\sss
and\dss the orbit\sss of\dss $w_{\fff n}$\sss does not\dss intersect\sss $A$\nnsp.\oss
Let\sss $w_{\dff i}$\sss be\sss the first\sss vertex along\dss this path such\dss that\dss
its orbit\sss does not\dss intersect\sss $A$\nnsp.\oss
Then\sss 
$g\dff(\dff w_{\dff i\dff -\dff 1}\trf)$ is\dss a\sss vertex of\dss $A$
for some $g\qff \in\qff G$\nnsp.\oss
The vertex $g\dff(\dff w_{\dff i}\trf)$ is\dss not\sss a\sss vertex of\dss $A$\sss
and\dss is\dss connected\sss with 
$g\dff(\dff w_{\dff i\dff -\dff 1}\trf)$
by an edge\qss
({\fff}because\sss $w_{\dff i}$\sss is\dss connected\sss 
with\sss $w_{\dff i\dff -\dff 1}$\sss by\sss an edge).\oss
Adding\dss to\sss $A$\sss this edge\sss together\sss with\dss the vertex
$g\dff(\dff w_{\dff i}\trf)$
results in a\sss tree properly\sss containing\sss $A$\sss and such\dss that\sss
all\dss its\sss vertices belong\dss to different\sss orbits.\oss
The contradiction with\sss the choice of\sss $A$ completes\sss the proof.\oss  \eproof

\myuppar{Additional\dss assumptions.}
Let\sss us\sss fix a subtree $A$ such as in\dss Lemma\qss \ref{tree-of-representatives}\qss
and\dss take as $V$\sss the set\sss of\dss vertices of $A$\nnsp.\oss
Let $v\qff \in\qff V$ and $e$ be an oriented edge of\sss $A$
such\dss that 
$e\qff \in\qff E_{\dff v}$\nsp.\oss
Then $\ttt\dff(\dff e\trf)\qff \in\qff V$ and\dss hence
$\ttt\dff(\dff e\trf)\off =\off v\dff(\dff e\trf)$\nnsp.\oss
Therefore $s_{\dff e}\off =\off 1$\sss is\dss a\sss legitimate choice,\oss
and\dss we will\sss assume\sss that  $s_{\dff e}\off =\off 1$
for every such $e$\nnsp.\oss
Further\halfff,\oss we will\sss assume\sss that $\mathcal{T}$ contains all\sss
oriented edges of\sss $A$\nnsp.\oss 
Then our relations\sss include\sss
$g_{\dff e}\off =\off 1$\sss for every\sss oriented edge $e$ of\dss $A$\nnsp.\oss
Therefore\vspace{3pt}
\[
\quad
\ttt\dff(\dff e\trf)^{\fff *}
\off =\off 
v\dff(\dff e\trf)^{\fff *}
\off =\off
\mathfrak{g}_{\dff e}\trf\left(\dff v\dff(\dff e\trf)^{\fff *}\trf\right)
\]

\vspace{-12pt}\vspace{3pt}
and\dss hence $\ttt\dff(\dff e\trf)^{\fff *}$
is\dss a neighbor\sss of\sss 
$v^{\fff *}\off =\off \aaa\dff(\dff e\trf)^{\fff *}$ in $\mathbb{X}$\nnsp.\oss
Therefore\sss if\sss $v\fff,\pff w\qff \in\qff V$
are connected\sss by\sss an edge of\sss $A$\nnsp,\oss
then
$v^{\fff *}\dnsp,\off w^{\fff *}
\qff \in\qff
\mathbb{X}_{\dff 0}$
are connected\dss by a\qss (unique)\qss edge of\sss $\mathbb{X}$\nnsp.\oss
Hence\sss the vertices $v^{\fff *}$ with $v\qff \in\qff V$\sss
together\sss with\sss these edges form a subtree $\mathbb{A}$ of\sss $\mathbb{X}$
such\dss that\sss $f$\sss induces an\sss isomorphism\dss
$\mathbb{A}\qff \ttoo\qff A$\nnsp.\oss
The subtree\dss $\mathbb{A}$\dss is\dss a\sss canonical\dss lift\sss of\sss
$A$\sss to\dss $\mathbb{X}$\nnsp.\oss

\mypar{Lemma.}{connectedness-m} 
\emph{The graph\trs $\mathbb{X}$\dss is\dss connected.\oss}

\proof
To begin\dss with,\oss for every $e\qff \in\qff E$\sss
there exits an edge of\sss $\mathbb{X}$ connecting
a\sss vertex of\sss $\mathbb{A}$ with a\sss vertex of\sss
$g_{\dff e}\dff(\trf \mathbb{A}\trf)$\nnsp.\oss
Indeed,\oss if\sss $e\qff \in\qff E_{\dff v}$\nsp,\oss
then\dss the vertex $v^{\fff *}$ of\sss $\mathbb{A}$\sss
is\dss connected\sss with\dss the vertex\dss
$\mathfrak{g}_{\dff e}\trf\left(\dff v\dff(\dff e\trf)^{\fff *}\trf\right)$
of\trs $g_{\dff e}\dff(\trf \mathbb{A}\trf)$\nnsp.\oss
By\sss applying $g_{\dff e}^{\dff -\dff 1}$
we see\sss that\dss there exists also an edge connecting
a\sss vertex\sss of\sss $\mathbb{A}$ with a\sss vertex of\sss
$g_{\dff e}^{\dff -\dff 1}\dff(\trf \mathbb{A}\trf)$\nnsp.\oss
In addition,\oss if\sss $v\qff \in\qff V$ and 
$h\qff \in\qff \mathbb{G}_{\dff v}$\nnsp,\oss
then\sss $v^{\fff *}$\sss is\dss a common\sss vertex of\sss 
$\mathbb{A}$ and $h\dff(\trf \mathbb{A}\trf)$\nnsp.\oss
Every\sss vertex $x$ of\sss $\mathbb{X}$ has\sss the form
$x\off =\off \gamma\qff(\dff v^{\fff *}\trf)$ for\sss some $v\qff \in\qff V$
and $\gamma\qff \in\qff \mathbb{G}$\nnsp.\oss
The element\sss $\gamma$ can\dss be represented 
as a product\sss of\dss the form\vspace{3pt}
\[
\quad
\gamma
\off =\off
\gamma_{\dff 1}\dff \cdot\qff
\gamma_{\dff 2}\dff \cdot\qff
\ldots\qff \cdot\qff
\gamma_{\dff m}
\pff,
\]

\vspace{-12pt}\vspace{3pt}
where each $\gamma_{\dff i}$\sss is\dss either equal\dss to
$\mathfrak{g}_{\dff e}$ or\sss $\mathfrak{g}_{\dff e}^{\dff -\dff 1}$
for some $e\qff \in\qff E$\nnsp,\oss
or\sss belongs\sss to\sss the free factor $\mathbb{G}_{\dff w}$ 
for some $w\qff \in\qff V$\dnsp.\oss
Let\sss us consider\sss subtrees\vspace{3pt}
\[
\quad
\mathbb{A}_{\dff i}
\off =\off
\gamma_{\dff 1}\dff \cdot\qff
\gamma_{\dff 2}\dff \cdot\qff
\ldots\qff \cdot\qff
\gamma_{\dff i}\qff
(\trf \mathbb{A}\trf)
\]

\vspace{-12pt}\vspace{3pt}
of\sss $\mathbb{X}$\nnsp,\oss 
where $i\off =\off 0\fff,\qff 1\fff,\qff \ldots\fff,\qff m$\nnsp.\oss
Then $x$\sss is\dss a\sss vertex of\sss $\mathbb{A}_{\dff m}$
and,\oss by\sss the usual\sss convention,\qss
$\mathbb{A}_{\dff 0}\off =\off \mathbb{A}$\nnsp.\oss
As we\sss just\sss saw,\oss for every $i$
either a vertex of\sss 
$\gamma_{\dff i}\qff
(\trf \mathbb{A}\trf)$
is\dss connected\dss with a vertex of\sss $\mathbb{A}$
by\sss an edge,\oss or\sss
$\gamma_{\dff i}\qff
(\trf \mathbb{A}\trf)$
and $\mathbb{A}$ have a common vertex.\oss
It\dss follows\sss that\dss for every $i$
either a vertex of\sss 
$\mathbb{A}_{\dff i}$\dss
is\dss connected\dss with a vertex of\sss $\mathbb{A}_{\dff i\dff -\dff 1}$\dss
by\sss an edge,\oss or
$\mathbb{A}_{\dff i}$
and $\mathbb{A}_{\dff i\dff -\dff 1}$ have a common vertex.\oss
Clearly\halfff,\oss the union of\trs trees\sss $\mathbb{A}_{\dff i}$\dss
and\dss these connecting\sss edges\dss is\dss connected\sss
and contains $x$ and $\mathbb{A}$\nnsp.\oss
It\sss follows\dss that\sss $\mathbb{X}$\sss is\dss connected.\oss  \eproof

\myuppar{Remark.}
The proof\dss of\trs Lemma\qss \ref{connectedness-m}\qss
is\dss modelled on an argument\dss from\sss the\dss Bass--Serre\trs theory
of\dss groups acting on\sss trees.\oss
Cf.\qss Serre\qss \cite{s},\oss the proof\dss of\trs Theorem\qss 12\qss in\dss Chapter\qss I.

\myuppar{From\dss pseudo-loops\dss to\dss loops.}
Let $q$ be a pseudo-loop,\oss i.e.\qss a path 
$w_{\dff 0}\dff,\qff w_{\dff 1}\dff,\qff \ldots\dff,\qff w_{\dff m}$
in\dss $X$\dss
such\dss that\sss
$w_{\dff 0}\dff,\qff w_{\dff m}\qff \in\qff V$\dnsp.\oss
Since $A$\sss is\dss connected,\oss 
there\dss is\dss a\sss path\sss 
$r$ in $A$ connecting $w_{\dff m}$\sss 
with $w_{\dff 0}$\nsp.\oss
By\sss following\dss $r$ first\sss
and\dss then\sss following\sss $q$\sss
we\sss get\sss a\sss loop $p$ starting\sss and ending at\sss $w_{\dff m}$\nsp.\oss
We will\sss say\sss that\sss such a\sss loop $p$\sss is\dss
a\qss \emph{closure}\qss of\sss $q$\nnsp.\oss
As we will\sss see in a moment,\qss 
$L\trf(\trf q\trf)$\dss is\dss equivalent\dss to $L\trf(\trf p\trf)$
and\dss hence replacing\dss the relation\dss $L\trf(\trf q\trf)$\dss
by\dss $L\trf(\trf p\trf)$\dss
does not\sss affect\sss $\mathbb{G}$\dss and\sss $\mathbb{X}$\nnsp.\oss
This partially\sss motivates choosing\sss for every\sss pseudo-loop
$q\qff \in\qff \mathcal{L}$ a\sss closure $p$\nnsp.\oss
Of\dss course,\oss if\sss $q$\sss is\dss already\sss a\sss loop,\oss
one can\dss take $p\off =\off q$\nnsp.\oss
Let\sss $\mathcal{L}_{\dff c}$ 
be\sss the set\sss of\dss these closures.\oss

\mypar{Lemma.}{closures}
\emph{$L\trf(\trf q\trf)$\sss is\dss equivalent\dss to\dss $L\trf(\trf p\trf)$\nnsp.\oss}

\proof
Let\sss
$e_{\dff 1}\dff,\qff e_{\dff 2}\dff,\qff \ldots\dff,\qff e_{\dff n}
\qff \in\qff E$\nnsp,\oss
where $n\off =\off k\qff +\qff m$\nnsp,\oss
be\sss the sequence of\dss edges related\dss to $p$
as in\dss Lemma\qss \ref{tracing-a-path-m}.\oss
Let\sss
$s_{\dff i}\off =\off s_{\dff e_{\dff i}}$
as\sss in\dss Lemma\qss \ref{tracing-a-path-m},\oss
and\sss let\sss
$g_{\dff i}\off =\off g_{\dff e_{\dff i}}$.\oss
Let\sss
$v_{\dff 0}\dff,\qff v_{\dff 1}\dff,\qff \ldots\dff,\qff v_{\dff k}$
be\sss the vertices of\dss $r$\nnsp.\oss
Since  
$v_{\dff 0}\dff,\qff v_{\dff 1}\dff,\qff \ldots\dff,\qff v_{\dff k}$
are vertices of\sss $A$ and\sss hence\sss belong\sss to\sss $V$\dnsp,\oss
the edges
$e_{\dff 1}\dff,\qff e_{\dff 2}\dff,\qff \ldots\dff,\qff e_{\dff k}$
are\sss the edges of\sss $r$\dnsp.\oss
In\sss more details,\qss 
$\aaa(\dff e_{\dff i}\trf)\off =\off v_{\dff i\dff -\dff 1}$
and\dss
$\ttt(\dff e_{\dff i}\trf)\off =\off v_{\dff i}$
for\dss $i\qff \leq\qff k$\nnsp.\oss
It\dss follows\sss that\sss $s_{\dff i}\off =\off 1$
and $g_{\dff i}\off =\off 1$
is\dss a\sss tautological\sss relation\dss
for $i\qff \leq\qff k$\nnsp.\oss
This implies\sss that\dss 
$L\trf(\trf p\trf)$\dss is\dss equivalent\dss to\sss 
$g_{\dff k\dff +\dff 1}
\dff \cdot\qff 
\ldots
\qff \cdot\qff 
g_{\dff n}
\off =\off
s_{\dff k\dff +\dff 1}
\dff \cdot\qff 
\ldots
\qff \cdot\qff 
s_{\dff n}$\nsp.\oss
The\sss latter\sss relation\dss is\dss nothing\sss else
but\dss the pseudo-loop relation $L\trf(\trf q\trf)$\nnsp.\oss
Therefore\sss the pseudo-loop relation\sss $L\trf(\trf q\trf)$\sss
is\dss equivalent\dss to\sss the\sss loop\sss relation\sss 
$L\trf(\trf p\trf)$\nnsp.\oss  \eproof

\mypar{Theorem.}{simply-connected-m}
\emph{Let\pss $\num{X}^{\dff +}$ be\sss the result\sss of\qss glueing\dss 
$2$\dnsp-cells\sss to\sss the geometric realization\dss $\num{X}$\dss of\trs $X$\dss 
along\dss the geometric realization of\trs all\trs loops\sss of\qss
the form\dss $g\dff(\trf p\trf)$\dss with\dss 
$g\qff \in\qff G$\dss and\dss 
$p\qff \in\qff \mathcal{L}_{\dff c}$\nsp.\oss 
If\pss $\num{X}^{\dff +}$\dss is\dss simply-connected,\oss 
then every\trs lift\trs to\dss $\mathbb{X}$ of\trs a closed\dss path\dss in\dss $X$\dss
is\dss closed.\oss}

\proof 
Let\sss us\sss fix a vertex $v$ of\sss $A$\nnsp.\oss
As in\sss the proof\dss of\trs Theorem\qss \ref{simply-connected},\oss
it\dss is\dss sufficient\sss to prove\sss that\sss the lift\sss
of\dss every\sss path of\trs the form 
$r\dff \cdot\dff g\dff(\trf p\trf)\dff \cdot\dff r^{\dff -\dff 1}$
starting at\sss $v^{\dff *}$ is\dss closed,\oss
where $p\qff \in\qff \mathcal{L}_{\dff c}$\nsp,\oss 
$g\qff \in\qff G$\nnsp,\pss 
and $r$\sss is\dss a path connecting $v$ with 
$g\dff(\dff w\trf)$\nnsp,\oss 
where $w\qff \in\qff V$\dss 
is\dss the starting\sss vertex of\sss $p$\nnsp.\oss
Let\sss $r^{\sim}$\sss be the lift\sss of\dss $r$\sss starting\sss at\dss $v^{\dff *}$\dss 
and\dss $\beta\dff(\trf u^{\dff *}\trf)$ be its endpoint\halfff.\oss
Then\dss\vspace{3pt}
\[
\quad 
f\dff\left(\trf \beta\dff(\trf u^{\dff *}\trf)\trf\right)
\off =\dff\off
\varphi\dff(\trf \beta\trf)\dff(\dff u\trf)
\]

\vspace{-12pt}\vspace{3pt}
is\dss the endpoint\sss of\dss $r$\nnsp,\oss
i.e.\qss is\dss equal\dss to\dss $g\dff(\dff w\trf)$\nnsp.\oss 
In other\dss terms,\qss 
$\varphi\dff(\trf \beta\trf)\trf (\dff u\trf)
\off =\off g\dff(\dff w\trf)$\dss 
and\dss hence\sss the vertices\sss $u$ and $w$\sss belong\sss to\sss the same orbit.\oss 
Since $w\fff,\pff u\qff \in\qff V$\dnsp,\oss
this implies\sss that\trs $u\off =\off w$
and\dss hence\dss
$g^{\dff -\dff 1}\dff \cdot\qff \varphi\dff(\trf \beta\trf)
\qff \in\qff
G_{\dff w}$\nsp.\oss
The rest\sss of\dss the proof\dss is\dss similar\sss to\sss the\sss arguments
in\sss the proof\dss of\trs Theorem\qss \ref{simply-connected},\oss
but\sss with\sss $\mathbb{G}_{\dff v}$ replaced\dss by\sss $\mathbb{G}_{\dff w}$\nsp.\oss  \eproof

\mypar{Corollary\halfff.}{isomorphisms-m}
\emph{Under\dss
the assumptions of\trs the\dss theorem\dss the map\dss 
$f\dff \colon\dff \mathbb{X}\qff \ttoo\qff X$\dss 
is\dss an\dss isomorphism of\qss graphs and\dss the\sss homomorphism\dss
$\varphi\dff \colon\dff \mathbb{G}\qff \ttoo\qff G$\dss
is\dss an\dss isomorphism\sss of\qss groups.\oss}

\proof
Arguing as in\sss the proof\dss of\trs Corollary\qss \ref{isomorphism-complexes},\oss
one can see\sss that\sss $f\dff \colon\dff \mathbb{X}\qff \ttoo\qff X$\sss
is\dss an\sss isomorphism.\oss
It\sss follows\sss that\sss for every $v\qff \in\qff V$\sss
the map
$\mathbb{G}/ \mathbb{G}_{\dff v}
\qff \ttoo\qff
G/ G_{\dff v}$
induced\dss by\sss $\varphi$\sss is\dss a\sss bijection.\oss
But\sss $\varphi$\sss induces\sss isomorphism
$\mathbb{G}_{\dff v}
\qff \ttoo\qff
G_{\dff v}$\sss
and\dss hence\dss $\varphi$\dss is\dss an\sss isomorphism.\oss  \eproof

\myuppar{Implications between edge and edge-loop relations.}
Let\sss $e\qff \in\qff E_{\dff v}$ and $w\off =\off v\dff(\dff e\trf)$\nnsp.\oss
As in\dss Section\qss \ref{implications},\oss the relation\sss
$E\dff(\dff e\fff,\pff t\trf)$ holds\sss if\dss and\dss only\trs if\dss
there\dss is\dss an element\sss $h\qff \in\qff G_{\dff w}$\dss
such\dss that\trs
$t\dff \cdot\dff g_{\dff e}
\off =\off
g_{\dff t\dff(\dff e\trf)}\dff \cdot\dff 
h$\nnsp.\oss
Using\sss this observation and\sss replacing $G_{\dff v}$\sss by $G_{\dff w}$ in\sss the proofs,\oss
we see\sss that\trs Lemmas\qss \ref{edge-product},\oss 
\ref{edge-inverse}\qss and\qss \ref{edge-triple-product}\qss
hold\sss in\sss the present\sss context\sss without\sss any changes in\sss the statements.\oss
In\trs Lemmas\qss \ref{no-inversion-edges},\oss \ref{rotating-edge-loops},\oss 
and\qss \ref{inverting-edge-loops}\qss one needs say\sss that\sss $v$\sss
is\dss an arbitrary\sss element\sss of\sss $V$\dnsp.\oss
There\dss is\dss no need\sss for other changes\sss in\sss the statements,\oss 
and\sss the proofs are arranged\sss in such a way\sss that\sss no changes are needed.\oss
Hence\sss the results of\trs Section\qss \ref{implications}\qss 
hold\sss with\sss trivial\dss modifications for actions 
with several\sss orbits of\dss vertices.\oss

\myuppar{Inversions.}
Following\trs Section\qss \ref{presentations},\oss we will\sss now define an analogue 
of\dss the involution $\iota$\nnsp.\oss
Let\sss $\mathbf{E}$ be\sss the set\sss of\dss orbits of\dss the action of\sss $G$
on\sss the set\sss of\dss oriented edges of\sss $X$\nnsp.\oss
Since $V$\sss is\dss a set\sss of\dss representatives of\dss the action
of\sss $G$ on vertices,\oss one can\sss identify\sss $\mathbf{E}$\sss 
with\sss the disjoint\sss union
of\dss the sets\sss 
$E_{\dff v}\left/\dff G_{\dff v}\right.$\nsp.\oss
The map\dss
$e\qff \longmapsto\qff \overline{e}$\dss 
induces an involution\sss $\iota$\sss on\dss $\mathbf{E}$\nnsp.\oss
The $G$\dnsp-orbit\sss of\dss  
$e$\sss is\dss a\sss fixed\dss point\sss
of\trs the involution\sss $\iota$\sss if\dss and\dss only\trs if\dss the oriented edges
$e$\sss and\dss $\overline{e}$\dss belong\dss to\sss the same $G$\dnsp-orbit\halfff,\oss
i.e.\qss if\trs and\dss only\trs if\qss $\overline{e}\off =\off g\dff(\dff e\trf)$\dss
for some\dss $g\qff \in\qff G$\nnsp.\oss
Such an element\sss $g$\sss is\dss said\sss to be an\qss
\emph{inversion}\qss of\trs $e$\nnsp,\oss
and an oriented edge\dss $e$\dss is\dss said\dss to\qss 
\emph{admit\sss an\dss inversion}\pss
if\trs there exists an\sss inversion of\dss $e$\nnsp.\oss

\myuppar{Scaffoldings.}
In\sss the present\sss context\sss a\pss \emph{scaffolding}\pss 
consists\sss of\dss the disjoint\sss union\sss\vspace{3pt}
\[
\quad
E_{\dff 0}
\off =\off
\coprod\nolimits_{\dff v\qff \in\qff V}\qff E_{\dff v}^{\dff 0}
\off,
\]

\vspace{-12pt}\vspace{3pt}
where each $E_{\dff v}^{\dff 0}$\sss is\dss a family of\dss representatives of\sss
$G_{\dff v}$\dnsp-orbits\sss in $E_{\dff v}$\nsp,\oss
a\dss family\sss of\dss sets\sss $\mathcal{T}_{\fff e}$\sss  
of\trs representatives of\dss cosets in\dss 
$G_{\dff v}\fff\left/\fff G_{\dff e}\right.$\dss for every 
$v\qff \in\qff V$ and\sss 
$e\qff \in\qff E_{\dff v}^{\dff 0}$\nsp,\oss
and\sss a\sss family\dss $s_{\dff e}$\nsp,\dss $e\qff \in\qff E$\dss of\dss elements of\trs $G$\dss
as in\dss Section\qss \ref{several},\oss
i.e.\qss such\dss that\trs 
$s_{\dff e}\dff(\trf v\dff(\dff e\trf)\trf)
\pff =\off
\ttt\dff(\dff e\trf)$\sss
for every $v\qff \in\qff V$ and 
$e\qff \in\qff E_{\dff v}$\nsp.\oss
A\sss scaffolding\dss is\dss said\dss to be\qss \emph{regular}\pss if\trs the following\sss
four conditions hold.\oss\vspace{-5pt}
\begin{itemize}

\item[({\fff}i{\fff})]\quad\
If\dss $e\qff \in\qff E_{\dff 0}$\sss 
admits an\sss inversion,\oss then\sss $s_{\dff e}$\sss is\dss an\sss inversion of\dss $e$\nnsp.\oss

\item[({\fff}ii{\fff})]\quad\
If\dss $e\qff \in\qff E_{\dff 0}$\sss does not\sss 
admits an\sss inversion,\oss then\dss
$a
\off =\off 
s_{\dff e}^{\dff -\dff 1}\dff(\qff \overline{e}\qff)
\pff \in\pff
E_{\dff 0}$\dss
and\dss $s_{\dff a}\off =\off s_{\dff e}^{\dff -\dff 1}$\nnsp.

\item[({\fff}iii{\fff})]\quad\
If\qss $e\qff \in\qff E_{\dff 0}$\nsp,\pss $u\qff \in\qff \mathcal{T}_{\dff e}$\nsp,\oss
and\dss $d\off =\off u\dff(\dff e\trf)$\nnsp,\oss
then\dss 
$s_{\dff d}
\off =\off
u\dff s_{\dff e}\dff u^{\dff -\dff 1}$\nnsp.\oss

\item[({\fff}iv{\fff})]\quad\
The set\sss $E_{\dff 0}$ contains all\sss oriented edges of\sss $A$\nnsp.\oss

\end{itemize}

\vspace{-5pt}
Regular\dss scaffoldings do exists.\oss
The proof\trs is\dss similar\sss to\sss the proof\dss of\trs Lemma\qss \ref{coherent-frames}.\oss

The inclusion $E_{\dff 0}\dff \ttoo\dff E$ induces
a bijection $E_{\dff 0}\dff \ttoo\dff \mathbf{E}$\nnsp.\oss
Hence $\iota$ induces an involution 
$\iota_{\dff 0}$ on\sss $E_{\dff 0}$\nsp.\oss
Let $E_{\dff 1}$ be a set\sss of\dss representatives of\dss orbits of\sss
$\iota_{\dff 0}$\nsp.\oss
If $e\dff \in\dff E_{\dff 0}$ admits an\sss in\-ver\-sion,\oss
then $\iota_{\dff 0}\dff(\dff e\trf)\off =\off e$
and $e\qff \in\qff E_{\dff 1}$\nsp.\oss
Otherwise only\sss one of\trs the edges $e\fff,\off \iota_{\dff 0}\dff(\dff e\trf)$
belongs\sss to $E_{\dff 1}$\nsp.\oss

\mypar{Theorem.}{second-simplification-m}
\emph{Suppose\sss that\dss we are working\trs with a regular\trs scaffolding and\dss
the assumptions of\qss Theorem\qss \ref{simply-connected-m}\qss hold\dss
for a collection\sss $\mathcal{L}_{\dff c}$\sss of\dss closures of\qss pseudo-loops
form a set\sss $\mathcal{L}$\nnsp.\oss
For every\dss $e\qff \in\qff E_{\dff 1}$\dss let\trs $\mathcal{G}_{\dff e}$\dss
be a set\sss of\qss generators of\pss $G_{\dff e}$\nsp.\oss
Then\dss the group\dss $G$\dss can\dss be obtained\dss from\dss $H$\dss
by\sss adding a\sss generator\dss $g_{\dff e}$\dss for every\dss
$e\qff \in\qff E_{\dff 1}$\dss
and\dss the following\dss relations.\oss}\vspace{-5pt}
\begin{itemize}

\item[(${\nsp}E$\nsp)\qff]\quad\
\emph{The edge relations\pss 
$E\dff(\dff e\fff,\qff t\trf)$\qss for\qss $e\qff \in\qff E_{\dff 1}$\dss
and\dss $t\qff \in\qff \mathcal{G}_{\dff e}$\nsp.\oss} 

\item[$(\trf EL\trf)$]\quad\
\emph{The\sss edge-loop relations\pss $L\trf(\trf l_{\dff e}\trf)$\qss 
for edges\dss $e\qff \in\qff E_{\dff 1}$ admitting an inversion.\oss}

\item[(${\nsp}L$\nsp)\qff\halfff]\quad\
\emph{The\qss loop relations\pss
$L\trf(\trf l\qff)$\qss for\qss $l\qff \in\qff \mathcal{L}$\nnsp,\oss
rewritten\sss in\dss
terms of\qss $g_{\dff e}$\dss with\dss
$e\qff \in\qff E_{\dff 1}$\nsp.\oss}

\item[(${\nsp}T$\nsp)\qff]\quad\
\emph{The\sss tautological\dss relation\dss
$g_{\dff e}\off =\off 1$\dss for\sss each edge $e$ of\trs $A$
belonging\sss to $E_{\dff 1}$\nsp.\oss}

\end{itemize}

\vspace{-5pt}
\proof
The proof\dss is\dss similar\sss to\sss the proof\dss of\trs
Theorem\qss \ref{second-simplification}.\oss 
It\dss is\dss based on\dss Corollary\qss \ref{isomorphisms-m}\qss
and\dss on\sss turning\sss relations into definitions\sss
as in\dss Section\qss \ref{presentations}.\oss
We\sss leave details\sss to\sss the reader\halfff.\oss  \eproof

\myuppar{Remark.}
If\trs $X$\dss is\dss a\sss tree and no edge of\trs $X$\sss admits an inversion,\oss
then\dss there are no\sss edge-loop and\dss loop relations
and our presentation\dss is\dss equivalent\dss to\sss the\qss Bass-Serre structure\sss
theorem\qss \cite{s}.\oss
See\dss Serre\qss \cite{s},\oss Section\qss 5.4.\oss

\newpage
\myappend{Coxeter's\qss proof\qss of\pss his\qss implication}{coxeter-proof}

\myapar{Theorem.}{coxeter}
\emph{Let\dss $G$\sss be a\sss group.\qff\oss 
If\pss
$s\fff,\pff t\fff,\pff z\qff \in\qff G$\qss 
and\qss
$s^{\dff 3}
\off =\off 
t^{\dff 5}
\off =\off 
(\dff s\dff t\trf)^{\dff 2}
\off =\off
z$\nnsp,\oss
then\qss $z^{\dff 2}\off =\off 1$\dnsp.}

\prooftitle{Coxeter's\dss proof\oss \textup{\cite{c1}}}
In\dss his book\qss \cite{c2}\qss Coxeter\dss precedes a somewhat\sss 
condensed\sss version of\trs his\sss original\dss proof\pss \cite{c1}\qss by\dss the following\sss
suggestions\sss to\sss the reader\halfff.\oss\vspace{-9pt}

\begin{quoting}
The reader\sss may\dss like\sss to pause\sss here,\oss
and\dss test\dss his own skill\dss before reading\sss on\dff!\oss
The recommended\dss procedure\dss is\dss to express\dss $z$\dss in\sss
various ways,\oss until\sss one of\trs the expressions\dss is\dss
recognized as\sss being equal\dss to its own\sss inverse.\oss
\end{quoting}

\vspace{-9pt}
The present\sss author admits\sss that\dss he did\dss not\dss follow\dss
these suggestions.\oss
Here\dss is\dss the proof\trs from\qss \cite{c1}.\oss
The first\dss part\dss works\sss for\dss the relations\dss
$s^{\dff 3}
\off\dff =\off 
t^{\dff n}
\off =\off 
(\dff s\dff t\trf)^{\dff 2}
\off =\dff\off
z$\qss
with an arbitrary\dss integer\dss $n$\nnsp.\oss
To begin\dss with,\oss
the relations\dss 
$s^{\dff 3}
\off =\off
(\dff s\dff t\trf)^{\dff 2}$\dss
and\dss
$t^{\dff n}
\off =\off
(\dff s\dff t\trf)^{\dff 2}$\dss
imply\dss that\vspace{3pt}\vspace{0.125pt}
\[
\quad
s^{\dff 2}
\off =\off
t\dff s\dff t\quad\
\mbox{and}\quad\
t^{\dff n\dff -\dff 1}
\off =\off
s\dff t\dff s
\]

\vspace{-9pt}\vspace{0.125pt}
and\dss hence\dss
$t\off =\off s^{\dff 2}\qff t^{\dff -\dff 1}\dff s^{\dff -\dff 1}$\qss
and\qss
$s\off =\off t^{\dff -\dff 1}\dff s^{\dff -\dff 1}\dff t^{\dff n\dff -\dff 1}$\dnsp.\qff\oss
It\dss follows\dss that\vspace{5pt}
\[
\quad
s^{\dff 3}
\off =\off
t^{\dff n}
\off =\off
\left(\trf
s^{\dff 2}\qff t^{\dff -\dff 1}\dff s^{\dff -\dff 1}
\qff\right)^{\fff n}
\]

\vspace{-34.5pt}
\[
\quad
\phantom{s^{\dff 3}
\off =\off
t^{\dff n}
\off }
=\off
s^{\dff 2}\qff t^{\dff -\dff 1}\dff s^{\dff -\dff 1}\dff \cdot\dff
s^{\dff 2}\qff t^{\dff -\dff 1}\dff s^{\dff -\dff 1}\dff \cdot\qff
\ldots\qff \cdot\dff
s^{\dff 2}\qff t^{\dff -\dff 1}\dff s^{\dff -\dff 1}
\]

\vspace{-35pt}
\[
\quad
\phantom{s^{\dff 3}
\off =\off
t^{\dff n}
\off }
=\off
s\dff \cdot\dff s\qff t^{\dff -\dff 1}\dff \cdot\dff
s\dff t^{\dff -\dff 1}\dff \cdot\qff
\ldots\qff \cdot\dff
s\dff t^{\dff -\dff 1}\dff \cdot\dff s^{\dff -\dff 1}
\]

\vspace{-34.5pt}
\[
\quad
\phantom{s^{\dff 3}
\off =\off
t^{\dff n}
\off }
=\off
s\dff \cdot\dff
\left(\trf
s\dff t^{\dff -\dff 1}
\qff\right)^{\fff n}\dff \cdot\dff
s^{\dff -\dff 1}
\]

\vspace{-7pt}
and\dss hence\dss 
$s^{\dff 3}
\off =\off
\left(\trf
s\dff t^{\dff -\dff 1}
\qff\right)^{\fff n}$\dnsp.\qff\oss
Similarly\halfff,\oss\vspace{5pt}
\[
\quad
t^{\dff n}
\off =\off
s^{\dff 3}
\off =\off
\left(\trf
s\dff t^{\dff -\dff 1}
\qff\right)^{\fff n}
\off =\off
\left(\trf
t^{\dff -\dff 1}\dff s^{\dff -\dff 1}\dff t^{\dff n\dff -\dff 1}\dff \cdot\dff t^{\dff -\dff 1}
\qff\right)^{\fff n}
\]

\vspace{-33.5pt}
\[
\quad
\phantom{t^{\dff n}
\off =\off
s^{\dff 3}
\off }
=\off
\left(\trf
t^{\dff -\dff 1}\dff s^{\dff -\dff 1}\dff t^{\dff n\dff -\dff 2}
\qff\right)^{\fff n}
\]

\vspace{-34.5pt}
\[
\quad
\phantom{t^{\dff n}
\off =\off
s^{\dff 3}
\off }
=\off
t^{\dff -\dff 1}\dff s^{\dff -\dff 1}\dff t^{\dff n\dff -\dff 2}\dff \cdot\dff
t^{\dff -\dff 1}\dff s^{\dff -\dff 1}\dff t^{\dff n\dff -\dff 2}\dff \cdot\qff
\ldots\qff \cdot\dff
t^{\dff -\dff 1}\dff s^{\dff -\dff 1}\dff t^{\dff n\dff -\dff 2}
\]

\vspace{-35pt}
\[
\quad
\phantom{t^{\dff n}
\off =\off
s^{\dff 3}
\off }
=\off
t^{\dff -\dff 1}\dff \cdot\dff 
s^{\dff -\dff 1}\dff t^{\dff n\dff -\dff 3}\dff \cdot\dff
s^{\dff -\dff 1}\dff t^{\dff n\dff -\dff 3}\dff \cdot\qff
\ldots\qff \cdot\dff
s^{\dff -\dff 1}\dff t^{\dff n\dff -\dff 3}\dff \cdot\dff t
\]

\vspace{-34.5pt}
\[
\quad
\phantom{t^{\dff n}
\off =\off
s^{\dff 3}
\off }
=\off
t^{\dff -\dff 1}\dff \cdot\dff
\left(\trf
s^{\dff -\dff 1}\dff t^{\dff n\dff -\dff 3}
\qff\right)\dff \cdot\dff t
\]

\vspace{-7.5pt}
and\dss hence\dss 
$t^{\dff n}
\off =\off
\left(\trf
s^{\dff -\dff 1}\dff t^{\dff n\dff -\dff 3}
\qff\right)^{\fff n}$\dnsp.\qff\oss
Since\dss
$t\off =\off s^{\dff 2}\qff t^{\dff -\dff 1}\dff s^{\dff -\dff 1}$\dnsp,\qff\oss
for\qss $n\off =\off 5$\qss this implies\sss that\vspace{3.25pt}
\[
\quad
s^{\dff 3}
\off =\off
t^{\trf 5}
\off =\off
\left(\trf
s^{\dff -\dff 1}\dff t^{\trf 2}
\qff\right)^{\fff 5}
\]

\vspace{-33.5pt}
\[
\quad
\phantom{s^{\dff 3}
\off =\off
t^{\dff 5}
\off }
=\off
\left(\trf
s^{\dff -\dff 1}\dff \cdot\qff
\left(\trf
s^{\dff 2}\qff t^{\dff -\dff 1}\dff s^{\dff -\dff 1}
\qff\right)^{\dff 2}
\qff\right)^{\dff 5}
\]

\vspace{-33.5pt}
\[
\quad
\phantom{s^{\dff 3}
\off =\off
t^{\dff 5}
\off }
=\off
\left(\trf
s^{\dff -\dff 1}\qff \cdot\pff
s^{\dff 2}\qff t^{\dff -\dff 1}\dff s^{\dff -\dff 1}\qff \cdot\pff
s^{\dff 2}\qff t^{\dff -\dff 1}\dff s^{\dff -\dff 1}
\qff\right)^{\dff 5}
\]

\vspace{-33.5pt}
\[
\quad
\phantom{s^{\dff 3}
\off =\off
t^{\dff 5}
\off }
=\off
\left(\trf
s\qff t^{\dff -\dff 1}\dff 
s\qff t^{\dff -\dff 1}\dff s^{\dff -\dff 1}
\qff\right)^{\fff 5}
\]

\vspace{-33.5pt}
\[
\quad
\phantom{s^{\dff 3}
\off =\off
t^{\dff 5}
\off }
=\off
s\qff t^{\dff -\dff 1}\dff 
s\qff t^{\dff -\dff 1}\dff s^{\dff -\dff 1}\qff \cdot\pff
s\qff t^{\dff -\dff 1}\dff 
s\qff t^{\dff -\dff 1}\dff s^{\dff -\dff 1}\qff \cdot\pff
\ldots\qff \cdot\pff
s\qff t^{\dff -\dff 1}\dff 
s\qff t^{\dff -\dff 1}\dff s^{\dff -\dff 1}
\]

\vspace{-33.5pt}
\[
\quad
\phantom{s^{\dff 3}
\off =\off
t^{\dff 5}
\off }
=\off
s\qff t\qff \cdot\qff 
t^{\dff -\dff 2}\dff s\qff \cdot\qff
t^{\dff -\dff 2}\dff s\qff \cdot\pff
\ldots\pff \cdot\qff
t^{\dff -\dff 2}\dff s\qff \cdot\qff
t^{\dff -\dff 1}\dff s^{\dff -\dff 1} 
\]

\vspace{-33.5pt}
\[
\quad
\phantom{s^{\dff 3}
\off =\off
t^{\dff 5}
\off }
=\off
s\qff t\dff \cdot\dff 
\left(\trf
t^{\dff -\dff 2}\dff s
\qff\right)^{\dff 5}
t^{\dff -\dff 1}\dff s^{\dff -\dff 1} 
\qff.
\]

\vspace{-8.75pt}
It\dss follows\dss that\qss
$t^{\dff 5}
\off=\off
s^{\dff 3}
\off =\off
t\dff \cdot\dff 
\left(\trf
t^{\dff -\dff 2}\dff s
\qff\right)^{\dff 5}
\dff \cdot\qff
t^{\dff -\dff 1}$\qss
and\trs hence\qss\vspace{1.5pt}
\[
\quad
t^{\dff 5}
\off=\off
\left(\trf
t^{\dff -\dff 2}\dff s
\qff\right)^{\dff 5}
\qff.
\]

\vspace{-10.5pt}
Since also\qss
$t^{\dff 5}
\off =\off
\left(\trf
s^{\dff -\dff 1}\dff t^{\trf 2}
\qff\right)^{\dff 5}$\dnsp,\qff\oss
it\dss follows\dss that\oss
$z^{\trf 2}
\off =\off
t^{\dff 5}\dff \cdot\dff t^{\dff 5}
\dff\off =\off
1$\nnsp.\oss  \eproof

\myappend{Coxeter's\qss implication\qss and\qss universal\qss 
central\qss extensions}{central-extensions}

\myuppar{Universal\sss central\sss extensions.}
Let\dss $G$\dss be a\sss group.\oss
A\qss \emph{central\sss extension}\qss of\trs $G$\dss
is\dss a\sss group\dss $E$\dss together\dss with a\sss surjective\sss homomorphism\dss
$\varphi\dff \colon\dff E\qff \ttoo\qff G$\dss
such\dss that\dss the kernel\sss of\dss $\varphi$\dss is\dss contained\dss in\dss
the center of\trs $G$\nnsp.\oss
A central\sss extension\dss
$\upsilon\dff \colon\dff U\qff \ttoo\qff G$\dss of\trs $G$\dss
is\dss said\dss to be\qss \emph{universal}\pss if\trs for every\sss
central\sss extension\dss
$\varphi\dff \colon\dff E\qff \ttoo\qff G$\dss
there\dss is\dss a\sss one and only\sss one homomorphism\dss
$\eta\dff \colon\dff U\qff \ttoo\qff E$\dss
such\dss that\trs
$\upsilon\off =\off \varphi\dff \circ\dff \eta$\nnsp.\oss
If\dss a universal\sss central\sss extension of\trs $G$\dss exists,\oss
then\dss it\dss is\dss unique up\sss to\sss isomorphism over\dss $G$\nnsp.\oss 
A central\sss extension\dss
$\varphi\dff \colon\dff E\qff \ttoo\qff G$\dss
\emph{splits}\pss if\trs it\sss admits a\qss \emph{section},\oss
i.e.\qss a\sss homomorphism\dss
$s\dff \colon\dff G\qff \ttoo\qff E$\dss
such\dss that\trs $\varphi\dff \circ\dff s\off =\off \id_{\dff G}$\nsp.\oss
See\dss Milnor\qss \cite{m1},\oss Section\qss 5,\pss for a self-contained
exposition of\trs basic facts about\sss universal\sss central\sss extension.\oss

\myapar{Theorem.}{universal-ce}
\emph{A central\sss extension\dss
$\upsilon\dff \colon\dff U\qff \ttoo\qff G$\dss
is\dss universal\qss if\trs and\dss only\trs if\qss
the group\dss $U$\dss is\dss equal\dss to\sss its commutant\trs
$[\dff U\fff,\pff U\dff]$\dss and\sss every\sss central\sss extension of\pss $U$\sss splits.\oss}

\proof
See\qss \cite{m1},\oss Theorem\qss 5.3.\oss  \eproof

\myapar{Lemma.}{universal-perfect}
\emph{As\dss in\qss Section\qss \ref{coxeter-implication},\oss
let\qss $\mathcal{G}$\dss be\sss the group defined\dss by\dss generators\sss
$g\fff,\qff r$ and\dss relations\qss
$g^{\dff 2}
\off =\off 
r^{\dff -\dff 3}
\off =\off
(\trf g\fff r\trf)^{\dff 5}$\dnsp.\oss
Then\dss 
$\mathcal{G}
\off =\off 
[\dff \mathcal{G}\fff,\pff \mathcal{G}\dff]$\nnsp.\oss}

\proof
A presentation of\trs the quotient\dss group\dss 
$\mathcal{G}\fff/[\dff \mathcal{G}\fff,\pff \mathcal{G}\dff]$\dss 
can\sss be obtained\dss by\sss adding\dss
to\sss the relations of\trs $\mathcal{G}$\sss
the commutativity\dss relation\dss $g\fff r\off =\off r g$\nnsp.\oss
These relations\dss imply\dss that\trs\vspace{4.5pt}
\[
\quad
(\trf g\fff r\trf)^{\dff 6}
\off =\off
g^{\dff 6}\dff r^{\dff 6}
\off =\off
r^{\dff -\dff 9}\dff r^{\dff 6}
\off =\off
r^{\dff -\dff 3}
\off =\off
(\trf g\fff r\trf)^{\dff 5}
\]

\vspace{-7.5pt}
and\dss hence $g\fff r\off =\off 1$\nnsp.\oss
It\dss follows\dss that\trs $g^{\dff 3}\off =\off r^{\dff -\dff 3}$\dnsp.\oss
Together\dss with\dss
$g^{\dff 2}\off =\off r^{\dff -\dff 3}$\dss
this implies\sss that\trs $g\off =\off 1$\dss and\dss hence\dss $r\off =\off 1$\nnsp.\oss
Therefore\dss 
$\mathcal{G}\fff/[\dff \mathcal{G}\fff,\pff \mathcal{G}\dff]\off =\off 1$\dss
and\dss hence\dss
$\mathcal{G}
\off =\off 
[\dff \mathcal{G}\fff,\pff \mathcal{G}\dff]$\nnsp.\oss  \eproof

\myapar{Lemma.}{universal-splits}
\emph{Every\sss central\sss extension of\qss the\dss group\dss $\mathcal{G}$\sss splits.\oss}

\proof
It\dss is\dss sufficient\dss to prove\sss that\dss 
for every\sss central\sss extension\dss
$E\qff \ttoo\qff \mathcal{G}$\dss of\trs the\sss group\dss $\mathcal{G}$\dss
the elements\dss
$g\fff,\pff r\qff \in\qff \mathcal{G}$\dss
can\sss be\sss lifted\dss to elements\dss
$\gamma\fff,\pff \rho\qff \in\qff E$\dss
such\dss that\trs
$\gamma^{\trf 2}
\off =\off 
\rho^{\fff -\dff 3}
\off =\off
(\trf \gamma\dff \rho\trf)^{\dff 5}$\dnsp.\oss
Indeed,\oss in\dss this case\dss
$g\off \longmapsto\off \gamma$\nnsp,\pss
$r\off \longmapsto\off \rho$\dss
extends\sss to a section\dss
$\mathcal{G}\qff \ttoo\qff E$\dss
of\qss $E\qff \ttoo\qff \mathcal{G}$\nnsp.\oss
Let\dss us\dss start\dss with\sss arbitrary\dss lifts\dss
$\gamma\fff,\pff \rho\qff \in\qff E$\dss
of\trs the elements\dss
$g\fff,\pff r\qff \in\qff \mathcal{G}$\dss respectively\halfff.\oss
Then\vspace{4.5pt}
\[
\quad
\gamma^{\trf 2}
\off =\off 
a\qff \rho^{\dff -\dff 3}
\off =\off
b\trf (\trf \gamma\dff \rho\trf)^{\dff 5}
\]

\vspace{-7.5pt}
for some elements\dss $a\fff,\pff b$\dss of\trs the\sss kernel\sss
of\trs $E\qff \ttoo\qff \mathcal{G}$\nnsp.\oss
It\dss is\dss sufficient\dss to find\sss elements\dss $x\fff,\pff y$\dss
of\trs this kernel\sss such\dss that\dss the\qss ``corrected''\qss elements\dss
$x\dff g\fff,\pff y\dff r$\dss satisfy\dss the required\dss
relations,\oss i.e.\qss\vspace{4.5pt}
\[
\quad
(\dff x\trf \gamma\trf)^{\trf 2}
\off =\off 
(\trf y\trf \rho \trf)^{\dff -\dff 3}
\off =\off
(\trf x\trf \gamma\dff \cdot\dff y\trf \rho\trf)^{\dff 5}
\]

\vspace{-7.5pt}
Since\dss $x\fff,\pff y$\dss belong\dss to\sss the center\halfff,\oss
these relations are equivalent\dss to\vspace{4.5pt}
\[
\quad
x^{\trf 2}\trf \gamma^{\trf 2}
\off =\off 
y^{\dff -\dff 3}\trf \rho ^{\dff -\dff 3}
\off =\off
x^{\dff 5}\trf y^{\dff 5}\dff (\trf \gamma\trf \rho\trf)^{\dff 5}
\]

\vspace{-10.5pt}
and\dss hence\sss to\vspace{1.5pt}
\[
\quad
\gamma^{\trf 2}
\off =\off 
x^{\trf -\dff 2}\trf y^{\dff -\dff 3}\trf \rho ^{\dff -\dff 3}
\off =\off
x^{\dff 3}\trf y^{\dff 5}\trf (\trf \gamma\trf \rho\trf)^{\dff 5}
\qff.
\]

\vspace{-9pt}
Therefore,\oss it\dss is\dss sufficient\dss to find\sss elements\dss $x\fff,\pff y$\dss
of\trs the\sss kernel\sss such\dss that\vspace{4.5pt}
\[
\quad
x^{\trf -\dff 2}\trf y^{\dff -\dff 3}
\off =\off
a\quad\
\mbox{and}\quad\
x^{\dff 3}\trf y^{\dff 5}
\off =\off 
b
\qff.
\]

\vspace{-9pt}
This\dss is\dss a system of\trs two\sss linear equations in an abelian\dss group\qss
({\fff}the\sss kernel\sss of\qss $E\qff \ttoo\qff \mathcal{G}$\nsp).\oss
The corresponding\sss determinant\dss is\dss
$(\dff -\qff 2\dff)\dff \cdot\dff 5\qff -\qff (\dff -\qff 3\dff)\dff \cdot\dff 5
\off =\off
-\qff 1$\dss 
and\dss hence\sss this system\dss has a solution.\oss
In any\sss case,\oss a direct\sss check shows\sss that\trs
$x\off =\off a^{\dff -\dff 5}\trf b^{\dff -\dff 3}$\dnsp,\pss
$y\off =\off a^{\dff 3}\trf b^{\dff 2}$\dss
is\dss a solution.\oss  \eproof

\myuppar{Two universal\sss central\sss extensions of\trs $\mathcal{D}$\nnsp.}
Lemmas\qss \ref{universal-perfect}\qss and\qss \ref{universal-splits}\qss 
together with\trs Theorem\qss \ref{universal-ce}\qss
immediately\dss imply\dss that\dss the homomorphism\sss
$\varphi\dff \colon\dff
\mathcal{G}\qff \ttoo\qff \mathcal{D}$
from\qss Section\qss \ref{coxeter-implication}\qss is\dss a\sss
universal\sss central\sss extension of\qss $\mathcal{D}$\nnsp.\oss
By\dss the remarks at\dss the beginning of\trs Section\qss \ref{coxeter-implication}\qss
the\sss kernel\sss of\dss $\varphi$\dss is\dss a\sss cyclic\sss group
generated\dss by\dss the element\trs
$z
\off =\off
g^{\dff 2}
\off =\off 
r^{\dff -\dff 3}
\off =\off 
(\dff r\dff g\trf)^{\dff 5}$\dnsp.\oss

On\dss the other\dss hand,\oss the canonical\dss homomorphism\dss
$\spin(\dff 1\dff)\qff \ttoo\qff SO\dff(\dff 3\dff)$\dss
is\dss a central\sss extension of\dss $SO\dff(\dff 3\dff)$\nnsp,\pss
and\dss its\sss kernel\dss is\dss a\sss cyclic\sss group of\dss 
order\sss $2$\nnsp.\oss
Therefore\sss the canonical\dss homomorphism\dss
$\psi\dff \colon\dff
\mathcal{D}^{\fff \sim}\qff \ttoo\qff \mathcal{D}$\dss
from\dss the binary\dss icosahedral\dss group\dss $\mathcal{D}^{\fff \sim}$\dss
to\dss $\mathcal{D}$\dss is\dss
also a central\sss extension\dss with\dss the same cyclic\sss group of\dss
order\sss $2$\sss as its\sss kernel.\oss
The construction of\dss presentations of\trs $\mathcal{D}^{\fff \sim}$\dss
in\dss Section\qss \ref{examples}\qss leads\sss to a homomorphism\dss
$\eta\dff \colon\dff
\mathcal{G}\qff \ttoo\qff \mathcal{D}^{\fff \sim}$\dss
such\dss that\trs
$\varphi\off =\off \psi\dff \circ\dff \eta$\nnsp,\oss
i.e\qss the following\sss diagram\dss is\dss commutative.\oss\vspace{4.25pt}
\[
\quad
\hspace*{3em}
\begin{tikzcd}[column sep=normal, row sep=huge]
\mathcal{G}\off 
\arrow[rr, "{\displaystyle \eta}"] 
\arrow[dr, "{\displaystyle \varphi}"']  
& & \off \mathcal{D}^{\fff \sim} 
\arrow[dl, "{\displaystyle \psi}"] \\
& \mathcal{D}  & \\
\end{tikzcd}
\]

\vspace{-48pt}
The construction also shows\sss that\trs $\eta\dff(\dff z\trf)\off =\off c$\nnsp,\oss
where\sss $c$\sss is\dss the rotation\sss by\dss the angle\sss $2\dff \pi$\nnsp,\oss
the only\dss non-trivial\sss element\sss of\trs the\sss kernel\sss of\sss
$\psi$\nnsp.\oss
Actually\halfff,\pss $z^{\dff 2}\off =\off 1$\dss by\trs Coxeter's\trs implication,\oss
and\dss hence\sss the commutativity\sss of\trs the above diagram\dss implies\sss
that\sss $\eta$\sss is\dss an\sss isomorphism.\oss
Therefore\dss
$\psi\dff \colon\dff
\mathcal{D}^{\fff \sim}\qff \ttoo\qff \mathcal{D}$\dss
is\dss a\sss universal\sss central\sss extension of\trs $\mathcal{D}$\nnsp.\oss

Conversely\halfff,\oss suppose\sss that\dss it\dss is\dss known\dss that\trs
$\psi\dff \colon\dff
\mathcal{D}^{\fff \sim}\qff \ttoo\qff \mathcal{D}$\dss
is\dss a\sss universal\sss central\sss extension.\oss
Since universal\sss central\sss extensions of\trs $\mathcal{D}$ are unique
up\sss to a unique isomorphism over\dss $\mathcal{D}$\nnsp,\oss
this implies\sss that\dss $\eta$\sss is\dss an\sss isomorphism.\oss
Since\dss $c^{\dff 2}\off =\off 1$\dss and\dss
$\eta\dff(\dff z\trf)\off =\off c$\nnsp,\oss
this,\oss in\dss turn,\oss implies\sss that\trs $z^{\dff 2}\off =\off 1$\dss
and\dss hence proves\dss Coxeter's\dss implication.\oss

The fact\dss that\trs
$\psi\dff \colon\dff
\mathcal{D}^{\fff \sim}\qff \ttoo\qff \mathcal{D}$\dss
is\dss a\sss universal\sss central\sss extension 
admits a\sss proof\dss not\dss requiring\dss to\sss know
any\dss presentation\sss of\trs $\mathcal{D}$\dss or\dss
$\mathcal{D}^{\fff \sim}$\dnsp.\oss 
Instead of\dss presentations\sss it\dss uses some basic ideas about\sss
cohomology\sss of\dss groups and\dss the following\dss two\sss theorems.\oss

\myapar{Theorem.}{uce-exists}
\emph{A\sss group\dss $G$\dss admits a\sss universal\sss central\sss extension\dss
if\qss and\dss only\trs if\qss it\dss is\dss equal\dss to its commutant\qss
$[\dff G\fff,\pff G\dff]$\nnsp.\oss}

\proof
See\qss \cite{m1},\oss Theorem\qss 5.7.\oss  \eproof

\myapar{Theorem.}{uce-and-h2}
\emph{If\qss a\sss group\dss $G$\dss admits a\sss universal\sss central\sss extension,\oss
then\dss the\sss kernel\sss of\qss a universal\sss central\sss extension\dss is\dss
isomorphic\dss to\dss $H_{\dff 2}\dff(\trf G\fff,\qff \zzz\trf)$\nnsp.\oss}

\proof
See\qss \cite{m1},\oss Corollary\qss 5.8.\oss  \eproof

\myapar{Theorem.}{uce-binary-icosahedreal}
\emph{The homomorphism\dss
$\psi\dff \colon\dff
\mathcal{D}^{\fff \sim}\dff \ttoo\qff \mathcal{D}$\dss
is\dss a\sss universal\sss central\sss extension of\trs $\mathcal{D}$\nnsp.}

\proof
By\dss the construction,\oss
the homomorphism\dss
$\eta\dff \colon\dff
\mathcal{G}\qff \ttoo\qff \mathcal{D}^{\fff \sim}$\dss
is\dss surjective.\oss
Together\dss with\dss Lemma\qss \ref{universal-perfect}\qss this implies\sss that\dss
the group\dss $\mathcal{D}^{\fff \sim}$\dss is\dss equal\dss to\sss its\sss commutant\halfff.\oss
By\trs Theorem\qss \ref{uce-exists}\qss this\sss implies\sss that\dss there exists
a universal\sss central\sss extension\dss
$\upsilon\dff \colon\dff U\dff \ttoo\qff \mathcal{D}^{\fff \sim}$\nsp\dnsp.\oss 
By\trs Theorem\qss \ref{universal-ce}\qss every\sss 
central\sss extension of\trs $U$\dss splits.\oss
If\dss $\upsilon$\dss is\dss an\sss isomorphism,\oss
then every\sss 
central\sss extension of\trs $\mathcal{D}^{\fff \sim}$\dss splits.\oss
Since\sss the group\dss $\mathcal{D}^{\fff \sim}$\dss is\dss equal\dss to\sss its\sss
commutant\halfff,\oss in\dss this case\dss Theorem\qss \ref{universal-ce}\qss
implies\sss that\dss
$\psi\dff \colon\dff
\mathcal{D}^{\fff \sim}\dff \ttoo\qff \mathcal{D}$\dss
is\dss a\sss universal\sss central\sss extension of\trs $\mathcal{D}$\nnsp.\oss

It\dss remains\sss to prove\sss that\trs
$\upsilon\dff \colon\dff U\dff \ttoo\qff \mathcal{D}^{\fff \sim}$
is\dss an\sss isomorphism.\oss
By\trs Theorem\qss \ref{uce-and-h2}\qss
the\sss kernel\sss of\dss $\upsilon$\sss is\dss isomorphic\sss to\dss
$H_{\dff 2}\dff(\trf \mathcal{D}^{\fff \sim},\qff \zzz\trf)$\nnsp.\oss
Therefore it\dss is\dss sufficient\dss to prove\sss that\trf
$H_{\dff 2}\dff(\trf \mathcal{D}^{\fff \sim},\qff \zzz\trf)
\off =\off
0$\nnsp.\oss

Recall\dss that\dss $\mathcal{D}^{\fff \sim}$\dss is\dss a\sss subgroup of\trs
$\spin(\dff 1\dff)$\dss and\sss consider\dss the quotient\trs
$Q
\off =\off
\spin(\dff 1\dff)/\dff\mathcal{D}^{\fff \sim}$\nsp\dnsp.\oss 
The canonical\dss map\dss
$\spin(\dff 1\dff)\qff \ttoo\qff Q$\dss is\dss a covering\sss space.\oss
Since\trs $\spin(\dff 1\dff)$\dss is\dss
homeomorphic\sss to\sss the sphere\sss $S^{\dff 3}$\dnsp,\oss
this implies\sss that\trs
$\pi_{\dff 1}\dff(\dff Q\dff)\off =\off \mathcal{D}^{\fff \sim}$\dss
and\dss
$\pi_{\dff 2}\dff(\dff Q\dff)\off =\off 0$\nnsp.\oss
In\dss turn,\oss this implies\sss that\dss by\dss glueing\dss to\sss $Q$\sss
cells of\dss dimension\qss $\geq\qff 4$\dss one can construct\sss
a space\sss $Q^{\dff +}$\sss such\dss that\trs\vspace{3.5pt}
\[
\quad
\pi_{\dff 1}\dff(\dff Q^{\dff +}\dff)
\off =\off 
\pi_{\dff 1}\dff(\dff Q\dff)
\off =\off 
\mathcal{D}^{\fff \sim}\quad\
\mbox{and}\quad\
\pi_{\dff i}\dff(\dff Q^{\dff +}\dff)
\off =\off
0\quad\
\mbox{for}\quad\
i\qff \geq\qff 2
\qff.
\]

\vspace{-8.5pt}
By\sss one of\trs the definitions,\pss
$H_{\dff n}\dff(\trf \mathcal{D}^{\fff \sim},\qff \zzz\trf)
\off =\off
H_{\dff n}\dff(\dff Q^{\dff +},\qff \zzz\trf)$\dss
for such\sss $Q^{\dff +}$\dss and every\sss $n$\nnsp.\oss
Since\sss glueing\sss cells of\dss dimension\qss $\geq\qff 4$\dss
does not\sss affect\dss the homology\dss groups 
in dimensions\qss $\leq\qff 2$\nnsp,\oss\vspace{3pt}
\[
\quad
H_{\dff 1}\dff(\dff Q\fff,\qff \zzz\trf)
\off =\off
H_{\dff 1}\dff(\dff Q^{\dff +},\qff \zzz\trf)
\off =\off
H_{\dff 1}\dff(\trf \mathcal{D}^{\fff \sim},\qff \zzz\trf)\quad\
\mbox{and}
\]

\vspace{-36pt}
\[
\quad
H_{\dff 2}\dff(\dff Q\fff,\qff \zzz\trf)
\off =\off
H_{\dff 2}\dff(\dff Q^{\dff +},\qff \zzz\trf)
\off =\off
H_{\dff 2}\dff(\trf \mathcal{D}^{\fff \sim},\qff \zzz\trf)
\qff.
\]

\vspace{-9pt}
By\sss another definition,\pss
$H_{\dff 1}\dff(\dff G\fff,\qff \zzz\trf)
\off =\off
G\fff/[\dff G\fff,\pff G\dff]$\dss
for every\dss group\dss $G$\nnsp.\oss
Since\sss the group\dss $\mathcal{D}^{\fff \sim}$\dss 
is\dss equal\dss to\sss its\sss commutant\halfff,\oss
it\dss follows\sss that\trs
$H_{\dff 1}\dff(\dff Q\fff,\qff \zzz\trf)
\off =\off
0$\nnsp.\oss
By\trs Poincar\'{e}\dss duality\sss applied\dss to\sss
the $3$\dnsp-manifold\sss $Q$\nnsp,\oss
this implies\sss that\trs
$H_{\dff 2}\dff(\dff Q\fff,\qff \zzz\trf)
\off =\off
0$\dss
and\dss hence\dss
$H_{\dff 2}\dff(\trf \mathcal{D}^{\fff \sim},\qff \zzz\trf)
\off =\off
0$\nnsp.\oss  \eproof

\myuppar{Three proofs\sss of\trs Coxeter's\dss implication.}
Of\dss course,\oss Coxeter's\dss original\dss proof\halfff,\oss
presented\dss in\dss Appendix\qss \ref{coxeter-proof},\oss
is\dss the most\sss elementary\sss one,\oss
but\dss it\dss hardly\dss explains\sss why\dss the result\dss is\dss true.\oss
The proof\trs in\dss Section\qss \ref{coxeter-implication}\qss
is\dss still\dss fairly\sss elementary\fff.\oss
With\dss the exception of\dss a reference\sss to\sss the elementary\dss
theory of\dss $CW$\dnsp-complexes in\dss the proof\dss of\qss
Theorem\qss \ref{simply-connected},\oss
it\dss uses\sss geometrical\sss and\dss topological\dss ideas\sss
only\dss to arrange calculations with generators and\sss relations in
a\sss transparent\sss way\halfff.\oss
In contrast\halfff,\oss 
the proof\dss based on\dss the\sss theory\sss of\dss central\sss extensions
and\trs Theorem\qss \ref{uce-binary-icosahedreal}\qss
relies on\sss fairly\sss sophisticated\dss algebraic and\dss topological\dss tools.\oss
Somewhat\dss mysteriously\halfff,\oss
the proofs presented\sss in\dss Section\qss \ref{coxeter-implication}\qss
and\sss in\dss this appendix\sss suggest\sss different\dss
topological\dss reasons for\sss the exponent\dss in\dss Coxeter's\dss implication\dss
being\sss equal\dss to\sss $2$\nnsp.\oss 
In\dss the first\dss proof\trs the number\sss $2$ appears as\dss 
Euler\dss characteristic of\trs $S^{\dff 2}$\dnsp,\oss
while in\dss the second one as\sss the order of\trs the\sss fundamental\dss group\dss
$\pi_{\dff 1}\dff(\trf SO\dff(\dff 3\dff)\trf)$\nnsp,\oss i.e.\qss
the number of\dss sheets\sss of\trs the covering\dss
$\spin(\dff 1\dff)\qff \ttoo\qff SO\dff(\dff 3\dff)$\nnsp.\oss

\myappend{Cayley\qss diagrams\qss and\qss scaffoldings}{cayley-diagrams}

\myuppar{Cayley\dss diagrams.}
Let\sss $G$\sss be a\sss group\sss together\dss with 
a finite set\sss $S$\sss generating\sss $G$\nnsp.\oss
The\qss \emph{Cayley\dss diagram}\pss of\trs $G\fff,\pff S$\dss is\dss a\sss directed\sss
graph\dss with edges\sss labeled\dss by\sss elements of\dss $S$\nnsp.\oss
It\dss is\dss defined as follows.\oss
The vertices are elements of\trs $G$\nnsp.\oss 
For every\dss pair\dss $g\fff,\pff s$\nnsp,\oss
where\dss $g\qff \in\qff G$\dss and\dss $s\qff \in\qff S$\nnsp,\oss
there\dss is\dss an edge connecting\dss $g$\sss with\sss $g\fff s$\nnsp,\oss
directed\dss from\sss $g$\sss to\sss $g\fff s$\nnsp,\oss
and\dss labeled\dss by\sss $s$\nnsp.\oss
There are no other edges.\oss
If\trs $s^{\dff 2}\off =\off 1$\dss for some\dss $s\qff \in\qff S$\nnsp,\oss
then\sss for every\sss vertex\sss $g\qff \in\qff G$\dss the\sss two\sss directed edges,\oss
the one going\dss from\sss $g$\sss to\sss $g\fff s$\nnsp,\oss
and\dss the one going\dss from\sss $g\fff s$\dss to\dss
$g\fff s\fff s\off =\off g$\nnsp,\oss
are usually\sss replaced\dss by\sss single undirected edge connecting\dss
$g$\sss with\sss $g\fff s$\sss and\dss labeled\dss by\sss $s$\nnsp.\oss

\myuppar{Cayley\dss diagrams\sss for symmetric sets of\dss generators.}
A subset\sss $S\qff \subset\qff G$\dss
is\dss called\qss \emph{symmetric}\pss if\trs $S$\dss is\dss invariant\dss under\dss the
involution\dss $g\qff \longmapsto\qff g^{\dff -\dff 1}$\dnsp.\oss
Suppose\sss that\sss $S$\sss is\dss a symmetric set\sss of\dss generators of\trs $G$\dss
and consider\dss the\dss Cayley\dss diagram of\trs $G\fff,\pff S$\nnsp.\oss
In\dss this case,\oss if\trs $g\qff \in\qff G$\dss is\dss connected\dss to\dss
$h\qff \in\qff G$\dss by an edge directed\dss from\sss $g$\sss to\sss $h$\sss
and\dss labeled\dss by\sss $s\qff \in\qff S$\nnsp,\oss
i.e.\qss if\trs $h\off =\off g\fff s$\nnsp,\oss
then $h$ is\dss connected\sss to $g$\sss
by an edge directed\dss from $h$\sss to\sss $g$\sss
and\dss labeled\dss by\sss $s^{\dff -\dff 1}$\dnsp.\oss
If\qss $s^{\dff 2}\off =\off 1$\nnsp,\oss
then\sss both\sss these edges
can\sss be\sss replaced\dss by\sss an undirected edge\sss 
labeled\dss by\sss $s$\nnsp.\oss

In\dss the case of\dss a symmetric set\sss $S$\sss of\dss generators
one can view\sss the\dss Cayley\dss diagram somewhat\sss differently\halfff.\oss
Recall\dss that\sss an\trs \emph{orientation}\qss of\dss an edge in an undirected\sss graph\dss
is\dss a designation of\dss one of\trs its endpoints as\sss its\qss \emph{origin}\qss
and\dss the other as\sss its\qss \emph{target}\qss
(see\dss Sections\qss \ref{coxeter-implication},\dss \ref{several}).\oss
Let\dss $C$\dss be\sss the undirected\sss graph\sss having\dss $G$\dss as\sss
the set\sss of\dss vertices and\dss having an edge connecting\trs
$g\fff,\pff h\qff \in\qff G$\trs if\trs and\dss only\trs if\qss
$g\dff h^{\dff -\dff 1}\qff \in\qff S$\nnsp.\oss
Then\sss for every oriented edge $\varepsilon$ of\trs $C$\dss
there\dss is\dss a\sss
unique element\trs $c_{\trf \varepsilon}\qff \in\qff S$\trs
such\dss that\trs $t\off =\off o\dff c_{\trf \varepsilon}$\nsp,\oss
where\sss $o$\sss it\dss the origin and\sss $t$\sss is\dss the\sss target of\trs this edge.\oss
The graph $C$ together\sss with\sss the map
$\varepsilon\qff \longmapsto\qff c_{\trf \varepsilon}$
is\dss simply\sss another\sss form of\trs the\dss Cayley\dss diagram of\sss $G\fff,\pff S$\nnsp.\qss
At\dss the same\sss time\sss this map\sss
looks very\sss similar\sss to\sss the scaffoldings from\dss Section\qss \ref{coxeter-implication}.

\myuppar{Cayley\dss diagrams and\dss scaffoldings.}
Suppose\sss that\dss $G$\sss acts\qss (from\dss the\sss left\halfff)\qss
on a set\dss $Z$\nnsp.\oss
Suppose\sss that\trs $z\qff \in\qff Z$\dss is\dss such\dss that\dss the map\dss
$g\off \longmapsto\off g\dff(\dff z\trf)$\sss is\dss bijective.\oss
Then one can\sss use\sss this map\sss to\sss identify\dss $G$\dss 
with\dss the orbit\dss $G\dff z$\nnsp.\oss
This identification\dss turns\dss $C$\dss
into a\sss graph\dss $Y$\dss having\dss the orbit\dss $G\dff z$\dss as\sss
the set\sss of\dss vertices.\oss
Let\dss us\dss look at\dss what\dss happens with\sss the map\dss
$\varepsilon\off \longmapsto\off c_{\trf \varepsilon}$\dss
under\dss this identification.\oss
Let\dss $\varepsilon$\dss be\sss the oriented edge of\trs $C$\dss
with\dss the origin\sss $g$\sss and\dss the\sss target\sss $g\fff s$\nnsp.\oss
Then\dss $c_{\trf \varepsilon}\off =\off s$\nnsp.\oss
The identification of\dss $G$\sss with\sss $G\dff z$\sss turns\dss $\varepsilon$\dss
into\sss the oriented edge of\trs $Y$\nnsp,\oss 
still\sss denoted\dss by\dss $\varepsilon$\nnsp,\oss 
having\sss $g\dff(\dff z\trf)$\sss as\sss its origin\dss
and\dss
$g\fff s\dff(\dff z\trf)
\off =\off
g\fff s\dff g^{\dff -\dff 1}\dff \cdot\dff g\dff(\dff z\trf)$\dss
as its\sss target.\oss
Clearly\halfff,\vspace{3pt}
\[
\quad
t_{\qff \varepsilon}
\off =\off
g\fff s\dff g^{\dff -\dff 1}
\off =\off
g\dff c_{\trf \varepsilon}\trf g^{\dff -\dff 1}
\]

\vspace{-12pt}\vspace{3pt}
is\dss the unique element\sss of\sss $G$\sss
taking\ $g\dff(\dff z\trf)$\sss to $g\fff s\dff(\dff z\trf)$\nnsp.\oss
If\sss $s^{\dff 2}\off =\off 1$\nnsp,\oss
then $t_{\qff \varepsilon}$ does not\sss depend on\dss the orientation
of\sss $\varepsilon$\nnsp,\oss
as a\sss trivial\sss calculation shows.\oss
We see\sss that\trs the identification of\sss $G$ with $G\dff z$\sss
turns\sss
$\varepsilon\pff \longmapsto\pff c_{\trf \varepsilon}$\sss
into an analogue\sss
$\varepsilon\pff \longmapsto\pff t_{\qff \varepsilon}$\sss 
of\sss regular scaffoldings\sss from\dss Section\qss \ref{presentations}.

In\dss the case of\trs the action of\trs the group\dss $\mathcal{D}$\dss
on\dss the dodecahedron\dss $D$\dss one can\dss take as\sss $z$\sss
the point\sss $x_{\dff 1}$\sss from\dss Section\qss \ref{coxeter-implication}\qss
and\dss the set\trs 
$\{\qff s_{\dff 1}\dff,\off h\dff,\off h^{\dff -\dff 1} \qff\}$\dss
as\dss $S$\nnsp.\oss
Then\dss $Y$\dss turns out\dss to\sss be\sss the graph denoted\dss by\dss $Y$\dss
in\dss Section\qss \ref{coxeter-implication},\oss and\dss 
$\varepsilon\off \longmapsto\off t_{\qff \varepsilon}$\dss 
to be\sss the canonical\dss $\mathcal{D}$\dnsp-scaffolding of\trs $Y$\dnsp.\oss

\begin{flushright}
October\qss 29,\oss 2023
 
https\halfff:/\!/\hspace*{-0.06em}nikolaivivanov.com

E-mail\halfff:\oss nikolai.v.ivanov{\fff}@{\dff}icloud.com,\oss ivanov{\fff}@{\dff}msu.edu

Department\sss of\qss Mathematics,\oss Michigan\sss State\sss University
\end{flushright}
}

\end{document}